\documentclass[11pt,twoside]{article}
\newcommand{\ifims}[2]{#1} 
\newcommand{\ifAMS}[2]{#1}   
\newcommand{\ifau}[4]{#1}  
\newcommand{\ifbook}[2]{#1}   
\newcommand{\ifsupnorm}[2]{#2} 
\newcommand{\ifLaplace}[2]{#1} 
\newcommand{\ifNL}[2]{#2}  
\newcommand{\ifapp}[2]{#1}  
\newcommand{\ifAnya}[2]{#1}  
\newcommand{\ifadap}[2]{#2}  
\newcommand{\iffourG}[2]{#1}  
\newcommand{\ifsqnorm}[2]{#2}  
\newcommand{\ifMLE}[2]{#1}  

\usepackage{tocloft}

\advance\cftsecnumwidth 0.3em\relax
\advance\cftsubsecnumwidth 0.3em\relax

\usepackage{ifthen}
\usepackage[ampersand]{easylist}
\usepackage{amsmath,amssymb,amsthm}
\usepackage[bb=stixtwo]{mathalpha}
\usepackage{mathtools}
\usepackage{natbib}
\usepackage{epsfig,graphicx}
\usepackage{comment}
\usepackage{color}
\usepackage{srcltx}
\usepackage[mathscr]{eucal}
\usepackage[math]{easyeqn}
\usepackage{etoolbox}
\usepackage{nccmath}
\usepackage{hyperref}
\hypersetup{
            colorlinks,
            linkcolor=hookersgreen,
            linktoc=blue,
            citecolor=ultramarine,
            urlcolor=black,
            filecolor=black
            }
\usepackage{nameref}
\usepackage{multirow}

\usepackage{accents}
\usepackage{mathrsfs}
\usepackage{bbm}
\usepackage{algorithm}
\usepackage{algpseudocode}

\ifims{
\textheight=22cm
\textwidth=14.8cm
\topmargin=0pt
\oddsidemargin=1.0cm
\evensidemargin=1.0cm
\linespread{1.3}

}{ 
}

\numberwithin{equation}{section}
\numberwithin{figure}{section}
\newcounter{example}[section]
\numberwithin{example}{section}
\newcounter{remark}[section]
\numberwithin{remark}{section}
\newtheorem{theorem}{Theorem}[section]
\newtheorem{proposition}[theorem]{Proposition}
\newtheorem{lemma}[theorem]{Lemma}
\newtheorem{corollary}[theorem]{Corollary}

\newtheorem{exmp}[example]{Example}
\newtheorem{rmrk}[remark]{Remark}
\newenvironment{example}{\begin{exmp}\rm}{\end{exmp}}
\newenvironment{remark}{\begin{rmrk}\rm}{\end{rmrk}}

\ifbook{
    \newcommand{\Chapter}[1]{\section{#1}}
    \newcommand{\Section}[1]{\subsection{#1}}
    \newcommand{\Subsection}[1]{\subsubsection{#1}}
    \def\Chname{Section }
    \def\chname{section }
    
  }
  {
    \newcommand{\Chapter}[1]{\chapter{#1}}
    \newcommand{\Section}[1]{\section{#1}}
    \newcommand{\Subsection}[1]{\subsection{#1}}
    \def\Chname{Chapter}
    
  }

\renewcommand{\(}{$\,}
\renewcommand{\)}{\,$}

\def\nquad{\hspace{-1cm}}
\def\eqdef{\stackrel{\operatorname{def}}{=}}
\def\eqd{\stackrel{\operatorname{d}}{=}}
\def\tod{\stackrel{d}{\longrightarrow}}

\DeclareMathAlphabet{\mathbbmsl}{U}{bbm}{bx}{sl}

\DeclareMathSymbol{\Alpha}{\mathalpha}{operators}{"41}
\DeclareMathSymbol{\Beta}{\mathalpha}{operators}{"42}
\DeclareMathSymbol{\Epsilon}{\mathalpha}{operators}{"45}
\DeclareMathSymbol{\Zeta}{\mathalpha}{operators}{"5A}
\DeclareMathSymbol{\Eta}{\mathalpha}{operators}{"48}
\DeclareMathSymbol{\Iota}{\mathalpha}{operators}{"49}
\DeclareMathSymbol{\Kappa}{\mathalpha}{operators}{"4B}
\DeclareMathSymbol{\Mu}{\mathalpha}{operators}{"4D}
\DeclareMathSymbol{\Nu}{\mathalpha}{operators}{"4E}
\DeclareMathSymbol{\Omicron}{\mathalpha}{operators}{"4F}
\DeclareMathSymbol{\Rho}{\mathalpha}{operators}{"50}
\DeclareMathSymbol{\Tau}{\mathalpha}{operators}{"54}
\DeclareMathSymbol{\Chi}{\mathalpha}{operators}{"58}
\DeclareMathSymbol{\omicron}{\mathord}{letters}{"6F}

\newcommand{\cc}[1]{\mathscr{#1}}
\newcommand{\bb}[1]{\boldsymbol{#1}}

\DeclareFontFamily{U}{mathx}{\hyphenchar\font45}
\DeclareFontShape{U}{mathx}{m}{n}{
<5><6><7><8><9><10>
<10.95><12><14.4><17.28><20.74><24.88>
mathx10
}{}
\DeclareSymbolFont{mathx}{U}{mathx}{m}{n}
\DeclareFontSubstitution{U}{mathx}{m}{n}
\DeclareMathAccent{\widebar}{0}{mathx}{"73}

\renewcommand{\bar}[1]{\widebar{#1}}
\renewcommand{\hat}[1]{\widehat{#1}}
\renewcommand{\tilde}[1]{\widetilde{#1}}

\makeatletter
\def\mathcenterto#1#2{\mathclap{\phantom{#1}\mathclap{#2}}\phantom{#1}}
\let\old@widetilde\widetilde
\def\widetildeto#1#2{\mathcenterto{#2}{\old@widetilde{\mathcenterto{#1}{#2\,}}}}
\let\old@widehat\widehat
\def\widehatto#1#2{\mathcenterto{#2}{\old@widehat{\mathcenterto{#1}{#2\,}}}}
\makeatother

\newcommand{\thankstitle}[1]{\ifthenelse{\equal{#1}{}}{}{\thanks{#1}}}
\newcommand{\thanksau}[1]{\ifthenelse{\equal{#1}{}}{}{\thanks{#1}}}

\newcommand{\aua}[6]
{\def\authora{#1}
\def\runauthora{#2}
\def\addressa{#3}
\def\emaila{#4}
\def\affiliationa{#5}
\def\thanksa{#6}}

\def\theauthors{
\ifau{ 
  \author{
    \authora
    \thanksau{\thanksa}
    \\[5.pt]
    \addressa \\
    \texttt{ \emaila}
  }
}
{  
  \author{
    \authora
    \thanksau{\thanksa}
    \\[5.pt]
    \addressa \\
    \texttt{ \emaila}
    \and
    \authorb
    \thanksau{\thanksb}
    \\[5.pt]
    \addressb \\
    \texttt{ \emailb}
  }
}
{   
  \author{
    \authora
    \thanksau{\thanksa}
    \\[5.pt]
    \addressa \\
    \texttt{ \emaila}
    \and
    \authorb
    \thanksau{\thanksb}
    \\[5.pt]
    \addressb \\
    \texttt{ \emailb}
    \and
    \authorc
    \thanksau{\thanksc}
    \\[5.pt]
    \addressc \\
    \texttt{ \emailc}
  }
} {   
  \author{
    \authora
    \thanksau{\thanksa}
    \\[5.pt]
    \addressa \\
    \texttt{ \emaila}
    \and
    \authorb
    \thanksau{\thanksb}
    \\[5.pt]
    \addressb \\
    \texttt{ \emailb}
    \and
    \authorc
    \thanksau{\thanksc}
    \\[5.pt]
    \addressc \\
    \texttt{ \emailc}
    \and
    \authord
    \thanksau{\thanksd}
    \\[5.pt]
    \addressd \\
    \texttt{ \emaild}
  }
}
}

\renewcommand{\Gamma}{\varGamma}
\renewcommand{\Pi}{\varPi}
\renewcommand{\Sigma}{\varSigma}
\renewcommand{\Delta}{\varDelta}
\renewcommand{\Lambda}{\varLambda}
\renewcommand{\Psi}{\varPsi}
\renewcommand{\Phi}{\varPhi}
\renewcommand{\Theta}{\varTheta}
\renewcommand{\Omega}{\varOmega}
\renewcommand{\Xi}{\varXi}
\renewcommand{\Upsilon}{\varUpsilon}

\def\argmax{\operatornamewithlimits{argmax}}
\def\argmin{\operatornamewithlimits{argmin}}

\def\av{\bb{a}}

\def\cv{\bb{c}}

\def\ev{\bb{e}}

\def\uv{\bb{u}}
\def\vv{\bb{v}}
\def\wv{\bb{w}}
\def\xv{\bb{x}}

\def\zv{\bb{z}}

\def\Av{\bb{A}}
\def\Bv{\bb{B}}

\def\Sv{\bb{S}}

\def\Uv{\bb{U}}

\def\Xv{\bb{X}}
\def\Yv{\bb{Y}}
\def\Zv{\bb{Z}}

\def\epsv{\bb{\varepsilon}}

\def\gammav{\bb{\gamma}}

\def\omegav{\bb{\omega}}

\def\xiv{\bb{\xi}}

\def\Psiv{\bb{\Psi}}

\def\sumi{\sum_{i=1}^{n}}

\usepackage{color}

\definecolor{blue(pigment)}{rgb}{0.2, 0.2, 0.6}
\definecolor{ultramarine}{rgb}{0.07, 0.04, 0.56}
\definecolor{darkspringgreen}{rgb}{0.09, 0.45, 0.27}
\definecolor{hookersgreen}{rgb}{0.0, 0.44, 0.0}
\definecolor{hgreen}{rgb}{0.09, 0.46, 0.2}
\definecolor{plum(traditional)}{rgb}{0.56, 0.27, 0.52}
\definecolor{purple(html/css)}{rgb}{0.5, 0.0, 0.5}
\definecolor{magenta(dye)}{rgb}{0.79, 0.08, 0.48}

\def\grad{\sigma}
\def\accu{\mu}

\def\jJ{J}

\def\CGPz{\CONSTi_{\cgp}}

\def\cgp{w}

\def\hspm{\hspace{1pt}}
\def\AvGP{\Avm_{\!\GP}}

\def\AFN{\mathbbmsl{U}}
\def\Avm{\bb{M}}

\def\deta{\tau}

\def\nEO{N}

\def\Xiv{\bb{\Xi}}
\def\AFN{\mathbb{Z}}
\def\BFN{\mathbb{B}}
\def\afv{\bb{s}}
\def\afn{\mathbb{z}}
\def\dmax{\kappa}

\def\PREC{\Ups}

\def\bvn{\bb{\mu}}
\def\svn{\bb{\phi}}

\def\GaussD{\gaussv_{\DPTG}}

\def\AFN{\mathbb{U}}

\def\Prec{\upsv}
\def\Precs{\Prec^{*}}
\def\Err{\mathcal{E}}

\def\VBH{S}
\def\VBHt{\tilde{\VBH}}

\def\leave{c}
\def\IIm{I^{\leave}}

\def\TG{\Gamma}
\def\TGD{\mathbbm{J}}

\def\trT{M}
\def\trTv{\bb{\trT}}
\def\trTvt{\tilde{\trTv}}

\def\Pmuvp{\Phi}

\def\muHa{\muH_{\alpH}}
\def\muHc{\muH_{c}}

\def\perm{\pi}

\def\KH{\KS}

\def\alpH{\alpha}
\def\alpHm{\rho}

\def\Gauss{\mathcal{G}}
\def\Gauss{\zeta}
\def\Gaussv{\bb{\Gauss}}

\def\Gaussvb{\bar{\Gaussv}}

\def\tensco{\tau}
\def\Tens{\mathcal{T}}
\def\Tenst{\tilde{\Tens}}
\def\TensG{\mathbb{T}}

\def\TensU{\Tens}

\def\charf{\mathbb{f}}

\def\supH{\lambda}
\def\rexH{\epsilon}

\def\Egs{\E_{\gaussv}}
\def\Pgs{\P_{\gaussv}}

\def\EUV{\E_{\, \UVL} \,}
\def\EX{\mathcal{E}}

\def\tarp{\theta}
\def\tarpv{\bb{\tarp}}
\def\tarps{\tarp^{*}}

\def\deta{\kappa}
\def\detam{\deta_{0}}

\def\rsmall{\varrho}
\def\Psimean{\bar{\Psiv}}

\def\hmax{\mathsf{c}}

\def\hL{h}

\def\smpa{\smp_{0}}

\def\lgd{f}

\def\PfL{\P_{\lgd}}

\def\dagg{\prime}

\def\amax{\nu}
\def\Cond{\,\, \Big| \, }

\def\Matr{\mathfrak{M}}

\def\weights{\weight^{*}}

\def\vH{\vA}

\def\Eta{\mathcal{H}}

\def\nbin{N}

\def\ic{i'}

\def\HVB{\mathbbmsl{V}} 

\def\HL{\mathbb{m}}

\def\HLH{\mathbb{h}}

\def\smp{s}

\def\dltw{\delta}

\def\dltwb{\omega}

\def\dltwu{\tau}
\def\dltwd{\dltw^{\dagg}}
\def\dltwbd{\dltwb^{\dagg}}

\def\dltwaa{\CONSTi}

\def\dltwhat{\hat{\dltw}}

\def\dlt{\delta}

\def\dltt{\tilde{\dlt}}

\def\vv{\bb{v}}

\def\II{\mathcal{I}}

\def\DFL{\mathbb{D}}

\def\DFN{\DVL}

\def\Projc{\Proj^{\perp}}

\def\R{\mathbbmsl{R}}
\def\E{\mathbbmsl{E}}

\def\P{\mathbbmsl{P}}

\def\kappa{\varkappa}

\def\Bernoulli{\operatorname{Bernoulli}}

\def\diag{\operatorname{diag}}

\def\Fr{\operatorname{Fr}}

\def\ND{\mathcal{N}}

\def\Var{\operatorname{Var}}

\def\T{\top}
\def\tr{\operatorname{tr}}

\def\cond{\, \big| \,}

\def\nsize{{n}}

\def\sumi{\sum_{i=1}^{\nsize}}

\def\ex{\mathrm{e}}

\def\Id{\mathbbmsl{I}}
\def\Ind{\operatorname{1}\hspace{-4.3pt}\operatorname{I}}

\def\alp{\alpha}

\def\avn{\av}

\def\bias{\mathsf{b}}

\def\biasD{\bias_{\scalebox{0.666}{$\DPN$}}}

\def\B{\cc{B}}
\def\BB{\mathbbmsl{B}}

\def\BBHt{\tilde{\BBH}}
\def\BHG{\mathcal{B}}

\def\BBB{\cc{B}}

\def\BBH{B}

\def\Bv{\bb{\BB}}

\def\bp{b}

\def\cdens{\phi}

\def\CA{\mathcal{A}}

\def\CGP{w}
\def\CGPa{\CGP_{0}}
\def\CGPs{\CGP_{*}}

\def\CONST{\mathtt{C} \hspace{0.1em}}
\def\CONSTi{\mathtt{C}}

\def\CONSTIF{\CONSTi_{\IF}}

\def\CONSTgmb{\CONSTi_{\cdens}}
\def\CONSTchar{\CONSTi_{\charf}}

\def\const{\mathsf{c}}

\def\constg{\kappa}

\def\DP{D}

\def\DPN{D}

\def\DPH{\DP}

\def\DVL{\mathbb{D}}

\def\dimA{\mathbb{p}}

\def\dimVG{\dimA_{\GP}}

\def\dimL{\dimA}

\def\dimH{\dimA}

\def\dimp{p}

\def\dimG{\dimA_{\GP}}

\def\dimQ{\dimA_{\QP}}

\def\dimq{q}

\def\dimm{m}

\def\dimD{\dimA_{\scalebox{0.666}{${\DPN}$}}}

\def\dimm{M}

\def\dens{f}

\def\Eta{\cc{H}}

\def\err{\diamondsuit}

\def\errS{\mathcal{E}}
\def\errSt{\tilde{\errS}}
\def\errSv{\errS}

\def\eps{\epsilon}			
\def\eps{\varepsilon}

\def\fs{f}

\def\fn{g}

\def\fG{f_{\GP}}

\def\fGu{h}

\def\gaussv{\bb{\gauss}}
\def\gauss{\gamma}

\def\gaussF{\gaussv_{\IFl}}

\def\Gauss{\zeta}
\def\Gaussv{\bb{\Gauss}}
\def\GaussG{\mathcal{G}}

\def\gm{\mathtt{g}}

\def\gmb{\gm}

\def\gmn{\gm}

\def\gp{g}

\def\GP{G}

\def\GPKS{\mathcal{K}}

\def\GPT{\mathcal{G}}

\def\hKS{\hat{\KS}}

\def\IF{\mathbbmsl{F}}
\def\IFGP{\IF_{\GP}}
\def\IFN{\IFL}

\def\IFL{{\mathbbmss{F}}}

\def\IFL{\mathbb{F}}

\def\imi{\mathtt{i}}

\def\jc{j'}

\def\KS{A}

\def\Kappa{\cc{K}}

\def\kc{k'}

\def\LT{L}
\def\LGP{\LT_{\GP}}

\def\LL{\cc{L}}

\def\loss{\wp}

\def\ldens{\ell}

\def\muH{\mu}

\def\norms{\circ} 
\def\normG{\alpha}

\def\nuH{\nu}

\def\omegav{\bb{\phi}}

\def\pent{\operatorname{pen}}

\def\Pone{P}
\def\Pdom{\mu_{0}}
\def\PDOM{\bb{\mu}_{0}}

\def\Psimean{\bar{\Psi}}
\def\Proj{\Pi}

\def\priord{\pi}

\def\QP{Q}

\def\rhot{t}

\def\rhou{b}

\def\rhogmn{\varrho}

\def\rhoH{\rho}

\def\riskt{\cc{R}}

\def\rr{\mathtt{r}}

\def\rrGP{\rr_{\GP}}

\def\rrB{\rr_{\scalebox{0.666}{${\BBH}$}}}
\def\rrD{\rr_{\scalebox{0.666}{${\DPN}$}}}

\def\rrn{\rr}

\def\rrL{\rr}

\def\rrTG{\rr}

\def\Sigmat{\tilde{\Sigma}}
\def\Sigmah{\hat{\Sigma}}

\def\Spsi{S}

\def\thetas{\theta^{*}}

\def\Tau{T}

\def\uvd{\uv^{\circ}}

\def\uvdd{\uv}

\def\ups{\upsilon}
\def\upsv{\bb{\ups}}

\def\upsvd{\upsv^{\circ}}
\def\upsvs{\upsv^{*}}

\def\upsvr{\breve{\upsv}}

\def\upsvn{\upsvd}

\def\ups{\upsilon}
\def\upsv{\bb{\ups}}

\def\upsd{\ups^{\circ}}

\def\upss{\ups^{*}}

\def\UV{\mathcal{U}}

\def\UVL{\UV}

\def\UVz{\UV}

\def\Ups{\varUpsilon}
\def\Upsd{\Ups^{\circ}}

\def\vA{\mathtt{v}}

\def\VP{V}

\def\weight{w}

\def\wv{\bb{w}}

\def\ww{w}

\def\WV{\mathcal{W}}

\def\xivGP{\xiv_{\GP}}

\def\xvs{\xv^{*}}

\def\xx{\mathtt{x}}
\def\xxc{\xx_{c}}

\def\xxs{\xx_{\norms}}

\def\xxe{\xx_{\ex}}

\def\XX{\cc{X}}

\def\yy{\mathtt{y}}

\def\zq{z}
\def\zqc{\zq_{c}}

\def\zqs{\bar{\zq}}

\def\zqe{\mathsf{z}}

\def\zz{\mathfrak{z}}

\def\zzH{\zz}
\def\thetitle{Finite sample expansions and risk bounds in high-dimensional SLS models}
\def\theruntitle{Finite sample expansions and risk bounds in high-dimensional SLS models}

\def\theabstract{
This note extends the results of classical parametric statistics like Fisher and Wilks theorem to modern setups
with a high or infinite parameter dimension, limited sample size, and possible model misspecification. 
We consider a special class of stochastically linear smooth (SLS) models satisfying three major conditions: the stochastic component
of the log-likelihood is linear in the model parameter and the expected 
negative log-likelihood is a smooth and convex function.
For a (quasi) penalized maximum likelihood estimators (pMLE), 
we establish three types of results: 
(1) concentration in a small vicinity of the ``truth'';
(2) Fisher and Wilks expansions;
(3) risk bounds.
In all results, the remainder is given explicitly and can be evaluated in terms of the effective sample size
and effective parameter dimension which allows us to identify the so-called \emph{critical parameter dimension}.
The results are also dimension and coordinate-free. 
The obtained finite sample expansions are of special interest because they can be used not only for obtaining the risk bounds
but also for inference, studying the asymptotic distribution, analysis of resampling procedures, etc.
The main tool for all these expansions is the so-called ``basic lemma'' about linearly perturbed optimization. 
Despite their generality, all the presented bounds are nearly sharp and the classical asymptotic results can be obtained
as simple corollaries. 
Our results indicate that the use of advanced fourth-order expansions 
allows to relax the critical dimension condition \( \dimA^{3} \ll n \) from \cite{SpLaplace2022} to 
\( \dimA^{3/2} \ll n \).
Examples for classical models like logistic regression, log-density and precision matrix estimation 
illustrate the applicability of general results.
We also indicate how the standard rate results from nonparametric statistics can be derived from the obtained risk bounds.
}

\def\kwdp{62F10,62E17}
\def\kwds{62J12}

\def\thekeywords{penalized MLE, Fisher and Wilks expansions, risk bounds}

\def\thankstitle{}

\aua
{Vladimir Spokoiny}
{Spokoiny, V. }
{
    Weierstrass Institute and HU Berlin,  
    \\
    Mohrenstr. 39, 10117 Berlin, Germany
    }
{spokoiny@wias-berlin.de}
{Weierstrass-Institute and Humboldt University Berlin}
{
  Financial support by the German Research Foundation (DFG) through the Collaborative Research Center 1294 ``Data assimilation'' is gratefully acknowledged.
}

\begin{document}
\thispagestyle{empty}
{
\title{\thetitle}
\theauthors

\maketitle
\begin{abstract}
{\footnotesize \theabstract}
\end{abstract}

\ifAMS
    {\par\noindent\emph{AMS 2010 Subject Classification:} Primary \kwdp. Secondary \kwds}
    {\par\noindent\emph{JEL codes}: \kwdp}

} 

\tableofcontents

\Chapter{Introduction}
\label{Sgenintr}


%

This paper presents some general results describing 
the properties 
of the \emph{penalized Maximum Likelihood Estimator} (pMLE).
Our starting point is a parametric model assumption about the distribution \( \P \) of the data \( \Yv \).
This distribution is described by a parameter \( \upsv \) 
from a parametric set \( \Ups \).
The quality of the data fit by a parameter \( \upsv \) is measured by
a fidelity (empirical risk) function \( L(\Yv,\upsv) \).
Its population counterpart (expectation w.r.t. the true data distribution \( \P \))
defines the risk function.
Our leading example is given by the negative log-likelihood function 
\( L(\upsv) = - \log \frac{d\P_{\upsv}}{d\PDOM}(\Yv) \) under 
the parametric assumption that 
\( \P \) belongs a given parametric family \( (\P_{\upsv} \, , \upsv \in \Ups) \) 
dominated by a sigma-finite measure \( \PDOM \).
This assumption is usually an idealization of reality
and the true distribution \( \P \) is not an element of \( (\P_{\upsv}) \).
However, a parametric assumption, even being wrong, may appear to be very useful, because it yields 
the method of estimation.
Namely, the (quasi) MLE \( \tilde{\upsv} \) is defined by minimizing 
the fidelity (empirical risk) function \( L(\upsv) = L(\Yv,\upsv) \)
over the parameter set \( \Ups \):
\begin{EQA}
	\tilde{\upsv}
	&=&
	\argmin_{\upsv \in \Ups} L(\upsv) .
\label{tuauv}
\end{EQA}
For a penalty function \( \pent_{\GP}(\upsv) \) on \( \Ups \), the \emph{penalized MLE} 
\( \tilde{\upsv}_{\GP} \) is defined by minimizing the penalized fidelity 
\( \LGP(\upsv) = L(\upsv) + \pent_{\GP}(\upsv) \):
\begin{EQA}
	\tilde{\upsv}_{\GP}
	&=&
	\argmin_{\upsv \in \Ups} \LGP(\upsv)
	=
	\argmin_{\upsv \in \Ups} \bigl\{ L(\upsv) + \pent_{\GP}(\upsv) \bigr\} .
\label{tuGauULGuauUp}
\end{EQA}
A typical example of \( \pent_{\GP} \) is a quadratic penalty:
\begin{EQA}
	\pent_{\GP}(\upsv)
	&=&
	\frac{1}{2} \| \GP \upsv \|^{2}
\label{vfgiubhygtfdehwjkie}
\end{EQA}
for a symmetric \( \dimp \times \dimp \) positive definite matrix \( \GP \in \Matr_{\dimp} \).
A prominent ridge regression also known as Tikhonov regularization
corresponds to \( \GP^{2} = \lambda \Id_{\dimp} \).

\Section{Classical parametric theory}

The classical Fisher parametric theory assumes that 
\( \Ups \) is a subset of a finite-dimensional Euclidean space \( \R^{\dimp} \),
the underlying data distribution \( \P \) indeed belongs to the considered parametric family 
\( (\P_{\upsv}) \), that is, \( \Yv \sim \P = \P_{\upsvs} \) for some \( \upsvs \in \Ups \).
In addition, some regularity of the family \( (\P_{\upsv}) \), or, equivalently, 
of the log-likelihood function \( - L(\upsv) \) is assumed. 
This, in particular, enables us to 
apply the third order Taylor expansion of \( L(\upsv) \) around the point of maximum \( \tilde{\upsv} \)
and to obtain Fisher expansions
\begin{EQA}[rcccl]
	\tilde{\upsv} - \upsvs
	& \approx &
	- \IF^{-1} \nabla L(\upsvs)
	\, ,
	\qquad
	- 2 L(\tilde{\upsv}) + 2 L(\upsvs)
	& \approx &
	\| \IF^{-1/2} \nabla L(\upsvs) \|^{2}
	.
\label{tumusaxFm1nLus}
\end{EQA}
Here  \( \IF = \IF(\upsvs) \) is the total Fisher information at \( \upsvs \) 
defined as the Hessian of the expected negative log-likelihood function \( \E L(\upsv) \):
\begin{EQA}
	\IF(\upsv)
	&=&
	\nabla^{2} \E L(\upsv) .
\label{IFuvmn2ELuv}
\end{EQA}
Under standard parametric assumptions, \( \IF(\upsv) \) is symmetric positive definite, \( \IF(\upsv) \in \Matr_{\dimp} \).
Moreover, if the data \( \Yv \) is generated as a sample of independent random variables 
\( Y_{1},\ldots,Y_{n} \), then the log-likelihood has an additive structure:
\( L(\upsv) = \sumi \ell(Y_{i},\upsv) \).
This allows to establish asymptotic standard normality of the standardized score 
\( \xiv \eqdef - \IF^{-1/2} \nabla L(\upsvs) \) and hence, to state Fisher and Wilks Theorems:
as \( n \to \infty \)
\begin{EQ}[rclcl]
	\IF^{1/2} \bigl( \tilde{\upsv} - \upsvs \bigr)
	& \approx &
	\xiv
	& \tod &
	\gaussv,
	\\
	- 2 L(\tilde{\upsv}) + 2 L(\upsvs)
	& \approx &
	\| \xiv \|^{2}
	& \tod &
	\| \gaussv \|^{2} \sim \chi^{2}_{\dimp} \, ,
\label{2Lm2Lusg2sc2p}
\end{EQ}
where \( \gaussv \) is a standard Gaussian vector in \( \R^{\dimp} \) and 
\( \chi^{2}_{\dimp} \) is a chi-squared distribution with \( \dimp \) degrees of freedom.
These results are fundamental and build the basis for most statistical applications 
like analysis of variance, canonical and correlation analysis, 
uncertainty quantification and hypothesis testing etc.
We refer to \cite{vdV00} and \cite{LeCa2006} for a comprehensive discussion and a historical overview of the related results
including the general LAN theory by L. Le Cam.
\Section{Extensions}
Modern statistical problems require to extend the classical results in several directions.

\paragraph{Model misspecification and bias}
Very often, the underlying data generating measure \( \P \) 
is not an element of the family \( (\P_{\upsv} \, , \upsv \in \Ups) \).
This means that the used log-likelihood function \( - L(\upsv) \) is not necessarily a true log-likelihood.
In particular, the condition \( \E \exp\{  - L(\upsv) \} = 1 \) does not hold.
This enables us to incorporate e.g. minimum contrast estimation, method of moments, etc.
The target of estimation \( \upsvs \) has to be redefined 
as the maximizer of the expected log-likelihood:
\begin{EQA}
	\upsvs
	& \eqdef &
	\argmin_{\upsv \in \Ups} \E L(\upsv) 
\label{pd47787tyhrujeirfjk}
\end{EQA}
leading to some modelling bias as the distance between \( \P \) and \( \P_{\upsvs} \).
This also concerns the use of a penalty leading to some \emph{penalization bias}.
When operating with the penalized loss \( \LGP(\upsv) \), the target of estimation 
becomes
\begin{EQA}
	\upsvs_{\GP}
	& \eqdef &
	\argmin_{\upsv \in \Ups} \E \LGP(\upsv) ,
\label{upssaruELuaiU}
\end{EQA}
which might be significantly different from \( \upsvs \).
This requires to carefully evaluate the penalization bias \( \upsvs_{\GP} - \upsvs \).

\paragraph{Finite samples, general likelihood, effective sample size}
Another important issue is a possibility of relaxing the assumption of i.i.d. or independent observations
which ensures an additive structure of the function \( L(\upsv) \).
We consider a general likelihood function, its structure does not need to be specified. 
In particular, we do not assume independent or progressively dependent observations and additive structure 
of the log-likelihood. 
We can even proceed with just one observation.
However, for stating our results about accuracy of estimation, we need
a notion of \emph{effective sample size} \( n \).
This is given via the so-called Fisher information matrix.
Everywhere we use the notation 
\begin{EQA}
	\IF(\upsv) 
	&=& 
	\nabla^{2} \E L(\upsv) ,
	\qquad
	\IF_{\GP}(\upsv) 
	= 
	\nabla^{2} \E \LGP(\upsv) 
	\, .
\label{e456yhgfr567uhgf}
\end{EQA}
We also write
\(
	\IF
	=
	\IF(\upsvs_{\GP}) \),
	\(
	\IF_{\GP} 
	= 
	\IF_{\GP}(\upsvs_{\GP})
\).
If the \( Y_{i} \)'s are i.i.d. then \( \IF(\upsv) \) is proportional to \( n \).
Therefore, the value 
\( n = \| \IF^{-1} \|^{-1} \) serves as a ``sample size''.

\paragraph{Effective parameter dimension and critical dimension}
One more important issue is the parameter dimension \( \dimp \).
The classical theory assumes \( \dimp \) fixed and \( n \) large.
We aim at relaxing both conditions by allowing a large/huge/infinite parameter dimension
and a small or moderate \( n \).
It appears that all the results below rely on the so-called 
\emph{effective dimension} \( \dimVG \) defined as
\begin{EQA}
	\dimVG
	& \eqdef &
	\tr \bigl\{ \IFGP^{-1} \Var(\nabla L(\upsvs_{\GP})) \bigr\} .
\label{nswe45678lkjmnbvc}
\end{EQA}
This quantity replaces the original dimension \( \dimp \) and it can be small or moderate even for \( \dimp \)
infinite. 
\cite{OsBa2021} used similar definitions in context of M-estimation.
One of the main intentions of our study is to understand the range of applicability of 
the mentioned results in terms of the effective parameter dimension \( \dimVG \) and the effective sample size \( n \).
It appears that most of the results ahead about concentration of the pMLE \( \tilde{\upsv}_{\GP} \) apply under the condition 
\( \dimVG \ll n \) which 
replaces the classical signal-to-noise relation: the effective number of parameters to be estimated is smaller in order
than the effective sample size.
More advanced results like Fisher and Wilks expansions and sharp risk bounds for a low dimensional sub-vector of \( \upsv \) 
may require stronger conditions \( \dimVG^{2} \ll n \) or \( \dimVG^{3/2} \ll n \). 

An essential feature of our results is their dimension-free and coordinate-free form.
The true parametric dimension \( \dimp \) can be very large, it does not appear in the error terms.
Also, we do not use any spectral decomposition or sequence space structure, in particular, we do not require that 
the Fisher information matrix \( \IF \) and the penalty matrix \( \GP^{2} \) are diagonal or can be jointly diagonalized.

\Section{Main steps of study}
Now we briefly describe our setup and the main focus of our analysis.                                                                                                                                                    
Below 
we limit ourselves to a special class of \emph{stochastically linear smooth} (SLS) statistical models.  
The major feature of such models is that the stochastic component 
\( \zeta(\upsv) = L(\upsv) - \E L(\upsv) \) of the log-likelihood \( L(\upsv) \) is linear in parameter \( \upsv \).
We also assume that the expected fidelity \( \E L(\upsv) \) is a convex and smooth function of the parameter \( \upsv \).
This class includes popular Generalized Linear Models but it is much larger. 
In particular, by extending the parameter space, one can consider 
many nonlinear models including nonlinear regression or 
nonlinear inverse problems as a special case of SLS; see \cite{Sp2019NLI}.
The assumption of stochastic linearity helps to avoid heavy tools of empirical process theory 
which is typically used in the analysis of pMLE \( \tilde{\upsv}_{\GP} \);
see e.g. \cite{BiMa1998}, \cite{vdV00}, \cite{geer2000}, \cite{Kosorok}, \cite{LeCa2006}, \cite{nickl_2015} among many others.
We only need some accurate deviation bounds for quadratic forms of the errors; see Section~\ref{Sdevboundgen} in the appendix.
Our aim is to establish possibly sharp and accurate results under realistic assumptions 
on a SLS model.
The study includes several steps.

\paragraph{Concentration of the pMLE}
The first step of our analysis is to establish a concentration result 
for the pMLE \( \tilde{\upsv}_{\GP} \) defined by minimization of \( \LGP(\upsv) \).
If the expected fidelity \( \E \LGP(\upsv) \) is strictly convex and smooth in \( \upsv \) then 
\( \tilde{\upsv}_{\GP} \) well concentrates 
in a small elliptic vicinity \( \CA_{\GP} \) of the ``target'' \( \upsvs_{\GP} \) from \eqref{upssaruELuaiU}:
\begin{EQA}
	\P\bigl( \| \IF_{\GP}^{1/2} (\tilde{\upsv}_{\GP} - \upsvs_{\GP}) \| > 3\rrGP/2 \bigr)
	& \leq &
	3 \ex^{-\xx} ,
\label{uyhhjkigvbghjugfg}
\end{EQA}
where \( \rrGP^{2} \approx \dimVG \).
Similar concentration bounds can be found in \cite{OsBa2021} for the M-estimator.
The result becomes sensible provided that \( \dimVG \ll n \) with \( \nsize^{-1} \asymp \| \IF_{\GP}^{-1} \| \).
In the classical parametric theory, such results about concentration of pMLE involve some advanced tools 
from the empirical process theory.
The use condition \nameref{Eref} about linearity of the stochastic component \( \zeta(\upsv) = L(\upsv) - \E L(\upsv) \)
allows to reduce the analysis to deviation bounds of the quadratic form \( \| \IF_{\GP}^{-1/2} \nabla \zeta \|^{2} \);
cf. condition \nameref{EU2ref}.
Section~\ref{Sdevboundgen} presents several results in this direction under different assumptions on the stochastic gradient 
\( \nabla \zeta \).

\paragraph{3S Fisher and Wilks expansions}
Having established the concentration of \( \tilde{\upsv}_{\GP} \in \CA_{\GP} \), 
we can restrict the analysis to this vicinity and 
 use the Taylor expansion of the penalized fidelity function 
\( \LGP(\upsv) \).
This helps to derive rather precise approximations for
\( \tilde{\upsv}_{\GP} - \upsvs_{\GP} \) and \( \LGP(\tilde{\upsv}_{\GP}) - \LGP(\upsvs_{\GP}) \):
\begin{EQ}[rcl]
    \bigl\| \IF_{\GP}^{1/2} \bigl( \tilde{\upsv}_{\GP} - \upsvs_{\GP} \bigr) + \xivGP \bigr\|
    & \leq &
    \frac{3\dltwu_{3}}{4} \, \bigl\| \xivGP \bigr\|^{2} ,
    \\
    \biggl| 
    \LGP(\tilde{\upsv}_{\GP}) - \LGP(\upsvs_{\GP}) 
    + \frac{1}{2} \bigl\| \xivGP \bigr\|^{2}
    \biggr|
    & \leq &
    \dltwu_{3} \, \bigl\| \xivGP \bigr\|^{3} \, ,
\label{3d3Af12DGi}
\end{EQ}
where 
\begin{EQA}
	\xivGP 
	& \eqdef & 
	- \IF_{\GP}^{-1/2} \nabla \LGP(\upsvs_{\GP}) 
	=
	- \IF_{\GP}^{-1/2} \nabla \zeta ,
\label{f78jjkonbh984e3edr4}
\end{EQA}
and \( \nabla \zeta = \nabla \zeta(\upsv) \) does not depend on \( \upsv \) due to linearity of
\( \zeta(\upsv) = \LGP(\upsv) - \E \LGP(\upsv) \).
The accuracy of approximation is controlled by the value \( \dltwu_{3} \)
which describes the accuracy of the third-order Taylor expansion of the function \( \fG(\upsv) = \E \LGP(\upsv) \)
in terms of the third directional derivatives of \( \fG \).
In typical examples \( \dltwu_{3} \asymp \sqrt{1/n} \). 
The presented results require \( \dltwu_{3}^{2} \, \dimVG \ll 1 \) which again leads to the condition 
\( \dimVG \ll n \).
The first result in \eqref{3d3Af12DGi} about the pMLE \( \tilde{\upsv}_{\GP} \) will be referred to as 
\emph{the Fisher expansion}, while the second one about \( \LGP(\tilde{\upsv}_{\GP}) \) is called 
\emph{the Wilks expansion}.
These two expansions provide a finite sample analog of the asymptotic statements \eqref{2Lm2Lusg2sc2p}
and are informative even in the classical parametric situation.
In fact, under standard assumptions, the normalized score vector \( \xivGP \) is asymptotically normal 
\( \ND(0,\Sigma_{\GP}) \) with \( \Sigma_{\GP} = \IF_{\GP}^{-1/2} \VP^{2} \IF_{\GP}^{-1/2} \in \Matr_{\dimp} \) 
and \( \VP^{2} = \Var\bigl( \nabla L(\upsv) \bigr) \in \Matr_{\dimp} \).
Stochastic linearity implies that the matrix \( \VP^{2} \) does not depend on the point \( \upsv \).
If the model is correctly specified, then \( \Sigma_{\GP} \) approaches the identity as \( n \to \infty \),
and we obtain the classical results \eqref{2Lm2Lusg2sc2p}.
Note that the use of stochastic linearity allows us to obtain much more accurate bounds than 
in \cite{SP2011} or \cite{SP2013_rough}. 
The derived finite sample expansions for the loss \( \tilde{\upsv}_{\GP} - \upsvs_{\GP} \) and the excess loss
\( \LGP(\tilde{\upsv}_{\GP}) - \LGP(\upsvs_{\GP}) \) are our main results which enable us 
not only for establishing the risk bounds like in \cite{OsBa2021} but also for making finite sample inference
and studying the asymptotic behavior of the estimator \( \tilde{\upsv}_{\GP} \).

\paragraph{3S risk bounds}
The loss of \( \tilde{\upsv}_{\GP} \) can be naturally expanded as
\begin{EQA}
	\tilde{\upsv}_{\GP} - \upsvs
	&=&
	\tilde{\upsv}_{\GP} - \upsvs_{\GP} + \upsvs_{\GP} - \upsvs \, .
\label{fughw83ur89hiv45edyf}
\end{EQA}
Due to the Fisher expansion \eqref{3d3Af12DGi}, 
\begin{EQA}
	\tilde{\upsv}_{\GP} - \upsvs_{\GP}
	& \approx &
	- \IF_{\GP}^{-1} \nabla \zeta \, .
\label{jdf56w3yhfi8f8iher4rdk}
\end{EQA}
Similarly we derive an expansion of the bias:
\begin{EQA}
	\upsvs_{\GP} - \upsvs
	& \approx &
	- \IF_{\GP}^{-1} \GP^{2} \upsvs \, .
\label{yd6yhwesduyujt63gty}
\end{EQA}
Putting together these two expansions leads to the so-called bias-variance decomposition of the squared risk:
for any linear mapping \( \QP \colon \R^{\dimp} \to \R^{\dimq} \)
\begin{EQA}
	\E \bigl\| \QP (\tilde{\upsv}_{\GP} - \upsvs) \bigr\|^{2}
	& \approx &
	\riskt_{\QP} \, ,
\label{vfog9ggiog98ryuj3ey}
\end{EQA}
where \( \riskt_{\QP} \) is the squared risk in the approximating linear model:
\begin{EQA}
	\riskt_{\QP}
	& \eqdef &
	\bigl\| \QP \, \IF_{\GP}^{-1} (\nabla \zeta + \GP^{2} \upsvs) \bigr\|^{2}
	=
	\tr \Var(\QP \IF_{\GP}^{-1} \nabla \zeta) + \bigl\| \QP \, \IF_{\GP}^{-1} \GP^{2} \upsvs \bigr\|^{2} \, .
\label{7nev6hdeivtewgvdfivc}
\end{EQA}
Theorem~\ref{TQFiWibias} provides sufficient conditions allowing to state a sharp risk bound:
\begin{EQA}
	(1 - \alp_{\QP})^{2} \riskt_{\QP}
	\leq 
	\E \bigl\| \QP (\tilde{\upsv}_{\GP} - \upsvs) \bigr\|^{2}
	& \leq &
	(1 + \alp_{\QP})^{2} \riskt_{\QP} \, .
\label{8kjeiv7ytew2wegvfuy}
\end{EQA}
Of course, this result is only meaningful if \( \alp_{\QP} \ll 1 \).
It appears that this value strongly depends on the dimension \( \dimq \) of the mapping \( \QP \).
If \( \QP = \Id_{\dimp} \) or \( \QP = \IF_{\GP}^{1/2} \) then \( \dimG \ll n \) is sufficient to ensure \( \alp_{\QP} \ll 1 \).
In the case of a low dimensional target with \( \dimq \asymp 1 \),
the condition \( \alp_{\QP} \ll 1 \) translates into \( \dimG^{2} \ll n \).

\paragraph{4S expansions and risk bounds}
The critical dimension condition \( \dimG^{2} \ll n \) can be very limiting.
Fourth-order smoothness conditions on \( \fG(\upsv) \) allow us to improve the accuracy of expansion \eqref{fughw83ur89hiv45edyf} 
by accounting for the third-order term and thus, relax the critical dimension bound.
Consider the third-order tensor
\( \TensU(\uv) = \frac{1}{6} \langle \nabla^{3} \fs(\upsvs_{\GP}), \uv^{\otimes 3} \rangle \).
Let \( \nabla \TensU(\uv) = \frac{1}{2} \langle \nabla^{3} \fs(\upsvs_{\GP}), \uv^{\otimes 2} \rangle \) be its gradient. 
Define the vectors \( \svn_{\GP} \) and \( \bvn_{\GP} \) by
\begin{EQ}[rcl]
	\svn_{\GP} 
	&=&
	\IF_{\GP}^{-1} \bigl\{ \nabla \zeta + \nabla \TensU(\IF_{\GP}^{-1} \nabla \zeta) \bigr\}
	\, ,
\label{8vfjvr43223efryfuweefgi}
	\\
	\bvn_{\GP}
	&=&
	\IF_{\GP}^{-1} \GP^{2} \upsvs + \IF_{\GP}^{-1} \nabla \Tens(\IF_{\GP}^{-1} \GP^{2} \upsvs) \, .
\end{EQ}
Theorem~\ref{Teff41} states the following bound: 
\begin{EQA}
	\| \QP \, (\tilde{\upsv}_{\GP} - \upsvs + \svn_{\GP} + \bvn_{\GP}) \|
	& \leq &
	\| \QP \IF_{\GP}^{-1/2} \| \, 
	\Bigl( \frac{\dltwu_{4}}{2} + \dltwu_{3}^{2} \Bigr) \,  
	\bigl( \| \IF_{\GP}^{-1/2} \nabla \zeta \|^{3} + \bias_{\GP}^{3} \bigr) ,
\label{0mkvhgjnrwwe3u8gtygi}
\end{EQA}
where \( \bias_{\GP} = \| \IF_{\GP}^{-1/2} \GP^{2} \upsvs \| \) and \( \dltwu_{4} \) controls fourth directional
derivatives of \( \fG \).
Typically \( \dltwu_{4} \asymp n^{-1} \) and \eqref{0mkvhgjnrwwe3u8gtygi}
is an improvement of \eqref{3d3Af12DGi} because the full dimensional error term in the right-hand side of \eqref{0mkvhgjnrwwe3u8gtygi} 
is of order \( \dimG^{3/2}/n \) compared to \( \dimG^{2}/n \) in \eqref{3d3Af12DGi}.
Therefore, the corrections from \eqref{8vfjvr43223efryfuweefgi} improves the critical dimension condition
from \( \dimG^{2} \ll n \) to \( \dimG^{3/2} \ll n \).
An interesting question of using a higher order expansion of \( \fG \) for a further relaxation of the critical dimension
condition is still open because even for 4S case, a closed-form solution of the corresponding 4S approximation problem
is not available.


\paragraph{Tools}
The presented results are based on two kinds of statements.
The results about concentration of the pMLE heavily rely on deviation bounds for quadratic forms of a centered and 
standardized score vector. 
Such results are collected in Section~\ref{Sdevboundgen}.
We separately study the cases of Gaussian errors, sub-Gaussian errors, and sub-exponential errors. 
The other important technical ingredient is the theory of perturbed optimization.
The main result of Theorem~\ref{PFiWigeneric2}
describes the solution of a convex optimization problem after a linear perturbation.
This result only relies on the smoothness and convexity of the objective function.
Section~\ref{Slocalsmooth} presents this and similar results.  

\Section{Examples}
Section~\ref{SgenBounds} presents some general theoretical results which improve and extend the similar results from the earlier paper \cite{SpLaplace2022}.
Section~\ref{SexamplesSLS} illustrates how the general conditions of Section~\ref{SgenBounds} can be checked 
for the classical setups like logistic regression, log-density, and precision matrix estimation.
All the mentioned examples are particular cases of Generalized Linear Models (GLM).
However, the SLS approach goes far beyond the GLM setup.
In particular, the paper \cite{Sp2024} explains how the so-called calming device can be used to bring 
a nonlinear regression problem to the SLS setup. 
Similarly, one can consider models like deep neuronal networks, nonlinear inverse problems, etc.
One more class of examples is given by error-in-operator models.
This class includes random design regression, instrumental regression, functional data analysis, 
diffusion, and McKean-Vlasov models, etc.
The calming trick applies here as well; see \cite{PuNoSp2024} for the case of a high-dimensional random design.  
The other examples include effective dimension reduction, 
Gaussian mixture estimation,
low-rank matrix recovery,
covariance and precision matrix estimation, 
smooth functional estimation, 
among others.
However, a rigorous treatment of each problem requires a separate study with a careful check
of the conditions and specific results and will be done elsewhere.

%


\Chapter{Properties of the MLE \( \tilde{\upsv} \) for SLS models}
\label{SgenBounds}
This \chname collects general results about concentration and expansion of the quasi MLE in the SLS setup
which substantially improve the bounds from \cite{SP2013_rough} and \cite{SpLaplace2022}.
%
We assume to be given a random function \( L(\upsv) \), \( \upsv \in \Ups \subseteq \R^{\dimp} \),
\( \dimp < \infty \).
This function can be viewed as negative log-likelihood or loss/empirical risk.
Consider in parallel two optimization problems defining 
the MLE \( \tilde{\upsv} \) and its population counterpart (the background truth) \( \upsvs \):
\begin{EQA}[rcl]
	\tilde{\upsv} 
	&=& 
	\argmin_{\upsv} L(\upsv),
	\qquad
	\upsvs 
	=
	\argmin_{\upsv} \E L(\upsv),
	\qquad
\label{tuGauLGususGE}
\end{EQA}
Define the Fisher information matrix \( \IF(\upsv) \eqdef \nabla^{2} \E L(\upsv) \) 
and denote \( \IF = \IF(\upsvs) \). 

\Section{Basic conditions}
\label{Scondgeneric}
Now we present our major conditions.
The most important one is about linearity of the stochastic component 
\( \zeta(\upsv) = L(\upsv) - \E L(\upsv) = L(\upsv) - \E L(\upsv) \).

\medskip
\begin{description}
    \item[\label{Eref} \( \bb{(\zeta)} \)]
      \textit{The stochastic component \( \zeta(\upsv) = L(\upsv) - \E L(\upsv) \) of the process \( L(\upsv) \) is linear in 
      \( \upsv \). 
      We denote by \( \nabla \zeta \equiv \nabla \zeta(\upsv) \in \R^{\dimp} \) its gradient
      }.
\end{description}

Below we assume some concentration properties of the stochastic vector \( \nabla \zeta \).
More precisely, we require that \( \nabla \zeta \) obeys the following condition.

\begin{description}
\item[\label{EU2ref}\( \bb{(\nabla \zeta)} \)]
	\textit{There exists \( \VP^{2} \geq \Var(\nabla \zeta) \) such that  
	for all considered  \( \BBH \in \Matr_{\dimp} \) and \( \xx > 0 \)
	}
\begin{EQA}
	\P\bigl( \| \BBH^{1/2} \VP^{-1} \nabla \zeta \| \geq \zq(\BBH,\xx) \bigr)
	& \leq &
	3 \ex^{-\xx} ,
\label{2emxGPm12nz122}
	\\
	\zq^{2}(\BBH,\xx)
	& \eqdef &
	\tr \BBH + 2 \sqrt{\xx \, \tr \BBH^{2}} + 2 \xx \| \BBH \| \,  .
\label{34rtyghuioiuyhgvftid}
\end{EQA}
\end{description}

This condition can be effectively checked if the errors in the data exhibit sub-gaussian or sub-exponential behavior; see 
\ifsqnorm{Section~\ref{SdevboundnonGauss}.}{\cite{Sp2023c}, \cite{Sp2023d}.}
The important special case corresponds to \( \BBH = \IF^{-1/2} \VP^{2} \IF^{-1/2} \) 
and \( \xx \approx \log n \) leading to the bound
\begin{EQA}
	\P\bigl( \| \IF^{-1/2} \nabla \zeta \| > \zq(\BBH,\xx) \bigr)
	& \leq &
	3/n .
\label{udyvfeyejff6777dj23}
\end{EQA}
The  value \( \dimA = \tr(\IF^{-1} \VP^{2}) \) can be called the \emph{effective dimension}; see \cite{SP2013_rough}.

We also assume that the loss function \( L(\upsv) \) or, equivalently, its deterministic part 
\( \E L(\upsv) \) is a convex function.

\medskip
\begin{description}
    \item[\label{LLref} \( \bb{(\mathcal{C})} \)]
      \textit{The function \( \E L(\upsv) \) is convex on \( \Ups \) which is open and convex set in \( \R^{\dimp} \).      
      }
\end{description}
\medskip

\ifLaplace{}{
In Section~\ref{Spostconcentr} we consider a stronger condition of semi-concavity of \( \E L(\upsv) \). }
%

Later we will also need some smoothness conditions on the function \( f(\upsv) = \E L(\upsv) \)
within a local vicinity of the point \( \upsvs \).
The notion of locality is given in terms of a metric tensor \( \DPN \in \Matr_{\dimp} \).
In most of the results later on, one can use \( \DPN = \IF^{1/2} \).
In general, we only assume \( \DPN^{2} \leq \dmax^{2} \IF \) for some \( \dmax > 0 \).
%
Introduce the error of the second-order Taylor approximation at a point \( \upsv \) in a direction \( \uv \) by
\begin{EQ}[rcl]
	\dltw_{3}(\upsv,\uv) 
	&=& 
	f(\upsv + \uv) - f(\upsv) - \langle \nabla f(\upsv), \uv \rangle 
	- \frac{1}{2} \langle \nabla^{2} f(\upsv), \uv^{\otimes 2} \rangle , 
	\\
	\dltwd_{3}(\upsv,\uv) 
	&=&
	\langle \nabla f(\upsv + \uv), \uv \rangle - \langle \nabla f(\upsv), \uv \rangle 
	- \langle \nabla^{2} f(\upsv), \uv^{\otimes 2} \rangle \, .
\label{dltw3vufuv12f2g}
\end{EQ}
Second order smoothness means a bound of the form
\begin{EQA}
	|\dltw_{3}(\upsv,\uv)| 
	& \leq & 
	\dltwb \| \DPN \uv \|^{2} \, ,
	\quad
	|\dltwd_{3}(\upsv,\uv)| 
	\leq 
	\dltwbd \| \DPN \uv \|^{2} \, ,
	\qquad
	\| \DPN \uv \| \leq \rr \, ,
\label{7jduusyswjdtcrtebfjvu}
\end{EQA}
for some radius \( \rr \) and small constants \( \dltwb \) and \( \dltwbd \).
These quantities can be effectively bounded under smoothness
conditions \nameref{LL3tref}, \nameref{LLsT3ref}, or \nameref{LLtS3ref} given in Section~\ref{Slocalsmooth}.
For instance, under \nameref{LL3tref}, by Lemma~\ref{LdltwLa3t}, it holds for a small constant \( \dltwu_{3} \)
\begin{EQA}
	\dltwbd
	& \leq &
	\dltwu_{3} \, \rr \, ,
	\qquad
	\dltwb
	\leq 
	\dltwu_{3} \, \rr/ 3 .
\label{swdy6fwqd6qtxcvbdyfdtw}
\end{EQA}
Also under \nameref{LLtS3ref}, the same bounds apply with \( \dltwu_{3} = \hmax_{3} \, n^{-1/2} \); see Lemma~\ref{LdltwLaGP}.

The class of models satisfying the conditions \nameref{Eref}, 
\nameref{EU2ref}, and \nameref{LLref}
with a smooth function \( f(\upsv) = \E L(\upsv) \) will be referred to as \emph{stochastically linear smooth} (SLS). 
This class includes linear regression, generalized linear models (GLM), and log-density models; 
see \cite{SpPa2019}, \cite{OsBa2021}\ifapp{ or Section~\ref{SGBvM} later.}{ or \cite{SpLaplace2022}.}
However, this class is much larger.
For instance, nonlinear regression can be adapted to the SLS framework 
by an extension of the parameter space; see 
\ifNL{Section~\ref{Snoninverse}.}{\cite{Sp2024}.}

\Section{Concentration of the MLE \( \tilde{\upsv} \). 2S-expansions}
\label{SgenpMLE}
This section discusses properties of the MLE \( \tilde{\upsv} = \argmin_{\upsv} L(\upsv) \)
under second-order smoothness conditions.
%
Fix \( \xx > 0 \) and define with \( \VP^{2} \) from \nameref{EU2ref} and \( \BBH = \IF^{-1/2} \VP^{2} \IF^{-1/2} \) 
\begin{EQA}[c]
	\UVz
	\eqdef 
	\bigl\{ \uv \colon \| \IF^{1/2} \uv \| \leq \tfrac{4}{3} \rrB \bigr\},
	\qquad
	\rrB \eqdef \zq(\BBH,\xx) ,
\label{rm1ucDGu2r0DGu}
\end{EQA}
where \( \zq(\BBH,\xx) \) is given by \eqref{34rtyghuioiuyhgvftid}.
By \nameref{EU2ref}, on a random set \( \Omega(\xx) \) with 
\( \P(\Omega(\xx)) \geq 1 - 3 \ex^{-\xx} \), it holds \( \| \IF^{-1/2} \nabla \zeta \| \leq \rr \).
Further, for the metric tensor \( \DPN \) from \eqref{7jduusyswjdtcrtebfjvu}, define
\begin{EQA}[c]
	\dltwb
    \eqdef 
    \sup_{\uv \in \UVz}
    \frac{2 |\dltw_{3}(\upsvs,\uv)|}{\| \DPN \uv \|^{2}} \,\, ,
    \qquad
    \dltwbd
    \eqdef 
    \sup_{\uv \in \UVz} \frac{|\dltwd_{3}(\upsvs,\uv)|}{\| \DPN \uv \|^{2}} \,\, . 
\label{dtb3u1DG2d3GP}
\end{EQA}

\begin{proposition}
\label{PconcMLEgenc}
Suppose
\nameref{Eref},
\nameref{EU2ref},
and 
\nameref{LLref}.
Let also \( \DPN^{2} \leq \dmax^{2} \IF \) and \( \dltwbd \, \dmax^{2} < 1/4 \);
see \eqref{dtb3u1DG2d3GP}.
Then on \( \Omega(\xx) \), 
it holds 
\begin{EQA}
	\| \IF^{1/2} (\tilde{\upsv} - \upsvs) \|  
	& \leq &
	\frac{4}{3} \, \rrB \, ,
	\qquad
	\| \DPN (\tilde{\upsv} - \upsvs) \|  
	\leq
	\frac{4 \dmax}{3} \, \rrB
	\, . 
\label{rhDGtuGmusGU0}
\end{EQA}
\end{proposition}

\begin{proof}
Apply Proposition~\ref{Pconcgeneric} to 
\( \fs(\upsv) = \E L(\upsv) \), \( \amax = 3/4 \), and \( \Av = \nabla \zeta \).
\end{proof}

Concentration of \( \tilde{\upsv} \) around \( \upsvs \) 
can be used to establish a version of the Fisher expansion for 
the estimation error \( \tilde{\upsv} - \upsvs \) and the Wilks expansion for the excess 
\( L(\tilde{\upsv}) - L(\upsvs) \).
The result substantially improves the bounds from \cite{OsBa2021} for M-estimators and follows by Proposition~\ref{PFiWigeneric}.

\begin{theorem}
\label{TFiWititG}
Assume the conditions of Proposition~\ref{PconcMLEgenc}.
Then on \( \Omega(\xx) \!\) 
\begin{EQ}[rcl]
    2 L(\tilde{\upsv}) - 2 L(\upsvs) 
    + \bigl\| \IF^{-1/2} \nabla \zeta \bigr\|^{2}
    & \leq &
    \frac{\dltwb}{1 + \dmax^{2} \hspm \dltwb} \bigl\| \DPN \IF^{-1} \nabla \zeta \bigr\|^{2} \, ,
    \\
    2 L(\tilde{\upsv}) - 2 L(\upsvs) 
    + \bigl\| \IF^{-1/2} \nabla \zeta \bigr\|^{2}
    & \geq &
    - \frac{\dltwb}{1 - \dmax^{2} \hspm \dltwb} \bigl\| \DPN \IF^{-1} \nabla \zeta \bigr\|^{2} .
\label{3d3Af12DGttG}
\end{EQ}
Also
\begin{EQ}[rcl]
    \bigl\| \DPN \bigl( \tilde{\upsv} - \upsvs + \IF^{-1} \nabla \zeta \bigr) \bigr\|
    & \leq &
    \frac{2\sqrt{\dltwb}}{1 - \dmax^{2} \hspm \dltwb} \, \bigl\| \DPN \IF^{-1} \nabla \zeta \bigr\| \, ,
    \\
    \bigl\| \DPN \bigl( \tilde{\upsv} - \upsvs \bigr) \bigr\|
    & \leq &
    \frac{1 + 2 \sqrt{\dltwb}}{1 - \dmax^{2} \hspm \dltwb} \, \bigl\| \DPN \IF^{-1} \nabla \zeta \bigr\| \, .
\label{DGttGtsGDGm13rG}
\end{EQ}
\end{theorem}

\Section{Expansions and risk bounds under third-order smoothness}
\label{SFiWiexs3}

The results of Theorem~\ref{TFiWititG} can be refined under third-order smoothness conditions.
%
%
Namely, Proposition~\ref{Pconcgeneric2} yields the following Wilks expansion for the MLE \( \tilde{\upsv} \).

\begin{theorem}
\label{TFiWititG2}
Assume 
\nameref{Eref},
\nameref{EU2ref},
and 
\nameref{LLref}.
Let also \nameref{LL3tref} hold at \( \upsvs \)  
with a metric tensor \( \DPN \) and values \( \rr \) and \( \dltwu_{3} \) satisfying 
\begin{EQA}
	\DPN^{2} \leq \dmax^{2} \, \IF , 
	\quad \rrn \geq \frac{4\dmax}{3} \, \rrB ,
	\quad 
	\dltwu_{3} \, \dmax^{3} \, \rrB
	& < &
	\frac{1}{4} \, ,
\label{yxdhewndu7jwnjjuuMLE}
\end{EQA}
for \( \rrB \) from \eqref{rm1ucDGu2r0DGu}.
Then on \( \Omega(\xx) \), it holds 
\begin{EQA}
	\| \IF^{1/2} (\tilde{\upsv} - \upsvs) \|  
	& \leq &
	\frac{4}{3} \rrB \, ,
	\qquad
	 \| \DPN (\tilde{\upsv} - \upsvs) \|  
	\leq 
	\frac{4\dmax}{3} \, \rrB \, ,
\label{rhDGtuGmusGU0a2MLE}
\end{EQA}
and
\begin{EQ}[rcl]
    \Bigl| 2 L(\tilde{\upsv}) - 2 L(\upsvs) + \| \IF^{-1/2} \nabla \zeta \|^{2} \Bigr|
    & \leq &
    \frac{\dltwu_{3}}{2} \, \| \DPN \IF^{-1} \nabla \zeta \|^{3} 
    \, .
\label{3d3Af12DGttG2}
\end{EQ}
\end{theorem}

Under \nameref{LLsT3ref}, Proposition~\ref{PFiWigeneric2} yields an advanced Fisher expansion.
Define
\begin{EQA}[c]
	\BBH_{\DPN} = \DPN \IF^{-1} \VP^{2} \IF^{-1} \DPN,
	\\
	\dimD
	\eqdef
	\tr \BBH_{\DPN} \, ,
	\quad
	\rrD \eqdef \zq(\BBH_{\DPN},\xx)
	\leq 
	\sqrt{\tr \BBH_{\DPN}} + \sqrt{2\xx \, \| \BBH_{\DPN} \|} \, ;
\label{y7s7d7d77dfdy7fuegue3j}
\end{EQA}
cf. \eqref{34rtyghuioiuyhgvftid}.
By \nameref{EU2ref}, it holds 
\( \P(\| \DPN \, \IF^{-1} \nabla \zeta \| > \rrD) \leq 3\ex^{-\xx} \).
The result follows by limiting to the set \( \Omega(\xx) \) on which \( \| \DPN \, \IF^{-1} \nabla \zeta \| \leq \rrD \)
and by applying Proposition~\ref{PFiWigeneric2}.

\begin{theorem}
\label{TFiWititG3}
Assume 
\nameref{Eref},
\nameref{EU2ref},
and 
\nameref{LLref}.
Let \nameref{LLsT3ref} hold at \( \upsvs \) with a metric tensor 
\( \DPN \) and values \( \rr \) and \( \dltwu_{3} \) 
satisfying
\begin{EQA}[c]
	\DPN^{2} \leq \dmax^{2} \hspm \IF ,
	\quad
	\rr \geq \frac{3}{2} \, \rrD \, ,
	\quad
	\dltwu_{3} \, \dmax^{2} \hspm \rrD < \frac{4}{9} \, ,
\label{8difiyfc54wrboeMLE}
\end{EQA}
where \( \rrD \) is from \eqref{y7s7d7d77dfdy7fuegue3j}.
With \( \Omega(\xx) = \{ \| \DPN \, \IF^{-1} \nabla \zeta \| \leq \rrD \} \), 
it holds \( \P(\Omega(\xx)) \geq 1 - 3 \ex^{-\xx} \) and on \( \Omega(\xx) \)
\begin{EQA}[rcl]
    \| \DPN^{-1} \IF (\tilde{\upsv} - \upsvs + \IF^{-1} \nabla \zeta) \|
    & \leq &
    \frac{3\dltwu_{3}}{4} \| \DPN \, \IF^{-1} \nabla \zeta \|^{2} 
	\, .
\label{DGttGtsGDGm13rG22}
\end{EQA}
%
\end{theorem}

Expansion \eqref{DGttGtsGDGm13rG22} yields accurate risk bounds.
\begin{theorem}
\label{TQFiWi}
Assume \nameref{Eref},
\nameref{EU2ref},
and 
\nameref{LLref}.
Let \( \fs(\upsv) = \E L(\upsv) \) satisfy \nameref{LLsT3ref} at \( \upsvs \) with some
\( \DPN \), \( \rr \), and \( \dltwu_{3} \).
Let also
\begin{EQA}[c]
	\DPN^{2} \leq \dmax^{2} \, \IF \, ,
	\qquad 
	\rr \geq \frac{3}{2} \rrD \, ,
	\qquad
	\dmax^{2} \dltwu_{3} \, \rrD < \frac{4}{9} \, ;
\label{7deyedhfg5563w6rfygjLR}
\end{EQA}
see \eqref{y7s7d7d77dfdy7fuegue3j}.
For any linear mapping \( \QP \colon \R^{\dimp} \to \R^{\dimq} \), 
it holds on \( \Omega(\xx) \)
\begin{EQA}
	\| \QP (\tilde{\upsv} - \upsvs + \IF^{-1} \nabla \zeta) \|
	& \leq &
	\| \QP \IF^{-1} \DPN \| \, \frac{3\dltwu_{3}}{4} \, \| \DPN \IF^{-1} \nabla \zeta \|^{2}
	\, .
	\qquad
	\quad
\label{g25re9fjfregdndg}
\end{EQA}
Also, introduce 
\begin{EQA}[c]
	\riskt_{\QP} \eqdef \E \{ \| \QP \IF^{-1} \nabla \zeta \|^{2} \Ind_{\Omega(\xx)} \} 
	\leq 
	\dimQ  \, 
\label{7djhed8cjfct534etgdhdyQP}
\end{EQA}
with \( \dimQ \eqdef \E \| \QP \IF^{-1} \nabla \zeta \|^{2} = \tr \Var(\QP \IF^{-1} \nabla \zeta) \).
Then
\begin{EQA}
	\E \bigl\{ \| \QP (\tilde{\upsv} - \upsvs) \| \Ind_{\Omega(\xx)} \bigr\} 
	& \leq &
	\riskt_{\QP}^{1/2} 
	+ \| \QP \IF^{-1} \DPN \| \, \frac{3\dltwu_{3}}{4} \, \dimD \, .
\label{EtuGus11md3GQP}
\end{EQA}
Further, assume \( \E \bigl\{ \| \DPN \IF^{-1} \nabla \zeta \|^{4} \Ind_{\Omega(\xx)} \bigr\} \leq \CONSTi_{4}^{2} \, \dimD^{2} \) 
and define
\begin{EQ}[rcl]
	\alp_{\QP}
	& \eqdef & 
	\frac{\| \QP \IF^{-1} \DPN \| \, (3/4) \dltwu_{3} \, \CONSTi_{4} \, \dimD } {\sqrt{\riskt_{\QP}}}
	\, . 
\label{6dhx6whcuydsds655srew}
\end{EQ}
If \( \alp_{\QP} < 1 \) then
\begin{EQA}
	(1 - \alp_{\QP})^{2} \riskt_{\QP} 
	\leq 
	\E \bigl\{ \| \QP \, (\tilde{\upsv} - \upsvs) \|^{2} \Ind_{\Omega(\xx)} \bigr\}
	& \leq &
	(1 + \alp_{\QP})^{2} \riskt_{\QP} \, .
\label{EQtuGmstrVEQtGQP}
\end{EQA}
\end{theorem}

\Section{Effective and critical dimension in ML estimation}
\label{ScritdimMLE}
This section discusses the important question of the critical parameter dimension 
still ensuring the validity of the presented results.
To be more specific, we only consider the 3S-results of Theorem~\ref{TFiWititG3}.
Also, assume \( \dmax \equiv 1 \).
The important constant \( \dltwu_{3} \) is identified by \nameref{LLtS3ref}: \( \dltwu_{3} = \hmax_{3} /\sqrt{n} \),
where the scaling factor \( n \) means the sample size.
It can be defined as the smallest eigenvalue of the Fisher operator \( \IF \).

First, we discuss the case \( \QP = \DPN = \IF^{1/2} \).
It appears that in this full dimensional situation, all the obtained results apply and are meaningful under the condition \( \dimA \ll n \),
where \( \dimA = \tr(\BBH) \) for \( \BBH = \IF^{-1/2} \VP^{2} \IF^{-1/2} \) is the \emph{effective dimension} of the problem.
Indeed,  \( \rrD^{2} = \rrB^{2} \approx \tr(\BBH) = \dimA \), and
condition \eqref{8difiyfc54wrboeMLE} requires \( \dltwu_{3} \, \rrD \ll 1 \) which can be spelled out as \( \dimA \ll n \).
Expansion \eqref{DGttGtsGDGm13rG22} means
\begin{EQA}
	\| \IF^{1/2} (\tilde{\upsv} - \upsvs) \|
    & \leq &
    \| \IF^{-1/2} \nabla \zeta \| + \frac{3\dltwu_{3}}{4} \| \IF^{-1/2} \nabla \zeta \|^{2} \, ,
\label{Y6DF66FCC6dseeEGHBE}
\end{EQA}
and the second term on the right-hand side of this bound is smaller 
than the first one under the same condition \( \dltwu_{3} \, \rrD \ll 1 \).
Similar observations apply to bound \eqref{EQtuGmstrVEQtGQP} of Theorem~\ref{TQFiWi} 
which is meaningful only if \( \alp_{\QP} \) in \eqref{6dhx6whcuydsds655srew} is small.
As \( \riskt_{\QP} \approx \dimQ = \dimA \), the condition \( \dltwu_{3} \, \rrD \ll 1 \) implies \( \alp_{\QP} \ll 1 \)
and hence, the bound \eqref{EQtuGmstrVEQtGQP} is sharp.
We conclude that the main properties of the MLE \( \tilde{\upsv} \) 
are valid under the condition \( \dimA \ll n \) meaning sufficiently many observations 
per effective number of parameters.

\ifNL{The situation changes drastically if \( \QP \) is not full-dimensional as e.g. in semiparametric estimation,
when \( \QP \) projects onto a low-dimensional target component.
We will see in \Chname \ref{SsemiMLE} that in this case,
\( \alp_{\QP} \ll 1 \) requires \( \dimA^{2} \ll n \).}{}
\iffourG{
An interesting question about a further improvement of the error term in \eqref{g25re9fjfregdndg}
will be discussed in the next section.}{}

\iffourG{
\Section{Bounds under fourth-order smoothness}
\label{SMLE4}
This section explains how the accuracy of the expansions for MLE can be improved
and the critical dimension condition can be relaxed under 
fourth-order smoothness of \( \fs(\upsv) = \E L(\upsv) \).

Consider the third-order tensor
\( \TensU(\uv) = \frac{1}{6} \langle \nabla^{3} \fs(\upsvs), \uv^{\otimes 3} \rangle \) and its gradient
\( \nabla \TensU(\uv) = \frac{1}{2} \langle \nabla^{3} \fs(\upsvs), \uv^{\otimes 2} \rangle \). 
Define a random vector \( \svn \) by
\begin{EQ}[rcl]
	\svn 
	&=&
	\IF^{-1} \nabla \zeta + \IF^{-1} \nabla \TensU(\IF^{-1} \nabla \zeta) 
	\, .
\label{8vfjvr43f8khg54ed54}
\end{EQ}
The next result shows that the use of \( \svn \) in place of \( \IF^{-1} \nabla \zeta \)
allows to improve the accuracy of the Fisher expansion \eqref{DGttGtsGDGm13rG22} and of the Wilks expansion \eqref{3d3Af12DGttG2}.

\begin{theorem}
\label{Teffp4s}
Assume \nameref{Eref}, \nameref{LLref}, and \nameref{EU2ref}.
Let \nameref{LLsT3ref} and \nameref{LLsT4ref} hold at \( \upsvs \) and 
\begin{EQA}[c]
	\DPN^{2} \leq \dmax^{2} \, \IF \, ,
	\;\; 
	\rr \geq \frac{3}{2} \rrD \, ,
	\;\;
	\dmax^{2} \dltwu_{3} \, \rrD < \frac{4}{9} \, ,
	\;\;
	\dmax^{2} \dltwu_{4} \, \rrD^{2} < \frac{1}{3} 
	\, ,
\label{7deyedhfg5563w6rfygjLR4}
\end{EQA}
with \( \rrD \) from \eqref{y7s7d7d77dfdy7fuegue3j}.
Then \( \svn \) from \eqref{8vfjvr43f8khg54ed54} fulfills on \( \Omega(\xx) \)
\begin{EQA}[rcl]
\label{yfjvh9h4f5e53tghfugu44}
	\| \DPN^{-1} \IF ( \tilde{\upsv} - \upsvs + \svn) \|
    & \leq &
    \Bigl( \frac{\dltwu_{4}}{2} + \dmax^{2} \dltwu_{3}^{2} \Bigr) \, \| \DPN \IF^{-1} \nabla \zeta \|^{3} \, ,
	\\
	\| \DPN^{-1} \IF \, (\svn - \IF^{-1} \nabla \zeta) \|
	&=&
	\| \DPN^{-1} \nabla \TensU(\IF^{-1} \nabla \zeta) \|
	\leq 
	\frac{\dltwu_{3}}{2} \, \| \DPN \IF^{-1} \nabla \zeta \|^{2}
	\, ,
	\qquad
\label{iuvchycvf6e64rygh3224}
\end{EQA}
and
\begin{EQA}
	&& \nquad
	\bigl| 
		L(\tilde{\upsv}) - L(\upsvs) 
		+ \frac{1}{2} \| \IF^{-1/2} \nabla \zeta \|^{2} + \TensU(\IF^{-1} \nabla \zeta) 
	\bigr|
	\\
	& \leq &
	\frac{\dltwu_{4} + 4 \dmax^{2} \dltwu_{3}^{2}}{8} \| \DPN \IF^{-1} \nabla \zeta \|^{4} 
    + \frac{\dmax^{2} (\dltwu_{4} + 2 \dmax^{2} \dltwu_{3}^{2})^{2}}{4} \, \| \DPN \IF^{-1} \nabla \zeta \|^{6} \, .
    \qquad
\label{87dfsjqweudsfjht7w3}
\end{EQA}
\end{theorem}

\begin{proof}
See Theorem~\ref{Pconcgeneric4} with \( \Av = \nabla \zeta \),
\( \avn = - \svn \), and \( \IFN = \IF \).
\end{proof}

The obtained expansion yields the bound on the loss and risk of \( \tilde{\upsv} \).
Define 
\begin{EQA}[rcl]
	\riskt_{\QP}
	& \eqdef &
	\E \bigl\{ \| \QP \IF^{-1} \nabla \zeta \|^{2} \Ind_{\Omega(\xx)} \bigr\} \, ,
\label{7dhdrw2bfvu78u78e4ndwQ}
	\\
\label{hdvje39bugebfh8edx}
	\riskt_{\QP,2} 
	& \eqdef & 
	\E \bigl\{ \| \QP \svn \|^{2} \Ind_{\Omega(\xx)} \bigr\} \, .
\end{EQA}

\begin{theorem}
\label{Teff41}
Assume the conditions of Theorem~\ref{Teffp4s} and let
\begin{EQA}[c]
	\E \bigl\{ \| \DPN \IF^{-1} \nabla \zeta \|^{k} \Ind_{\Omega(\xx)} \bigr\} \leq \CONSTi_{k}^{2} \, \dimD^{k/2} ,
	\qquad
	k=3,4,6 
	\, .
\label{6hjdfv8e6hyefyeheew7sk}
\end{EQA}
Then it holds for any linear mapping \( \QP \)
\begin{EQ}[rcl]
	&& \nquad
	\E \bigl\{ \| \QP \, (\tilde{\upsv} - \upsvs) \| \Ind_{\Omega(\xx)} \bigr\}
	\leq 	 
	\E \bigl\{ \| \QP \svn \| \Ind_{\Omega(\xx)} \bigr\}
	+ \| \QP \IF^{-1} \DPN \| \, 
	\Bigl( \frac{\dltwu_{4}}{2} + \dmax^{2} \dltwu_{3}^{2} \Bigr) \, \CONSTi_{3}^{2} \, \dimD^{3/2} \, ,
\label{0mkvhgjnrw3dfwe3u8gtygE1}
	\\
	&& \nquad
	\Bigl| \E \bigl\{ \| \QP \svn \| \Ind_{\Omega(\xx)} \bigr\}
	- \E \bigl\{ \| \QP \IF^{-1} \, \nabla \zeta \| \Ind_{\Omega(\xx)} \bigr\} 
	\Bigr|
	\leq 
	\| \QP \IF^{-1} \DPN \| \, \frac{\dltwu_{3}}{2} \, \dimD  \, .
\label{udtgecthwjdytdehduuc6}
\end{EQ}
With \( \riskt_{\QP,2} \) from \eqref{hdvje39bugebfh8edx}, let
\begin{EQA}
	\alp_{\QP,2}
	& \eqdef &
	\frac{\| \QP \IF^{-1} \DPN \| \, ( \dltwu_{4}/2 + \dmax^{2} \dltwu_{3}^{2}) \, \CONSTi_{6} \, \dimD^{3/2}}
		 {\sqrt{\riskt_{\QP,2}}} 
	< 1 
	\, .
\label{yfhcvhched6chejdrte}
\end{EQA}
Then
\begin{EQA}
	\bigl( 1 - \alp_{\QP,2} \bigr)^{2} \riskt_{\QP,2}
	\leq 
	\E \bigl\{ \| \QP \, (\tilde{\upsv} - \upsvs) \|^{2} \Ind_{\Omega(\xx)} \bigr\}
	& \leq &
	\bigl( 1 + \alp_{\QP,2} \bigr)^{2} \riskt_{\QP,2} \, .
\label{6shx76whnjvyehfbvyfh}
\end{EQA}
If another constant \( \alp_{\QP,1} < 1 \) ensures 
\begin{EQ}[rcl]
	&&
	\| \QP \IF^{-1} \DPN \| \, \frac{\dltwu_{3}}{2} \, \CONSTi_{4} \, \dimD 
	\leq 
	\alp_{\QP,1} \, \sqrt{\riskt_{\QP}} \, 
\label{6dhx6whcuydsds655srew4}
\end{EQ}
with \( \riskt_{\QP} \) from \eqref{7dhdrw2bfvu78u78e4ndwQ} then
\begin{EQA}
	\riskt_{\QP} (1 - \alp_{\QP,1})^{2} 
	\leq 
	\riskt_{\QP,2}
	& \leq &
	\riskt_{\QP} (1 + \alp_{\QP,1})^{2} \, .
\label{EQtuGmstrVEQtGQ2}
\end{EQA}
\end{theorem}

\begin{proof}
Rescaling of \( \DPN \) reduces the proof to \( \dmax = 1 \).
Theorem~\ref{Teffp4s} yields
\begin{EQA}
	\| \QP \, (\tilde{\upsv} - \upsvs + \svn) \|
	& \leq &
	\| \QP \IF^{-1} \DPN \| \, 
	\Bigl( \frac{\dltwu_{4}}{2} + \dltwu_{3}^{2} \Bigr) \,  
	\| \DPN \IF^{-1} \nabla \zeta \|^{3} \, ,
	\qquad
	\qquad
\label{0mkvhgjnrwwe3u8gtyg}
	\\
	\| \QP \{ \svn - \IF^{-1} \nabla \zeta \} \|
	& \leq &
	\frac{\dltwu_{3}}{2} \, \| \QP \IF^{-1} \DPN \| \, \| \DPN \IF^{-1} \nabla \zeta \|^{2} \, .
\label{vuedy766t4e3bfvyt6e}
\end{EQA}
Now \eqref{0mkvhgjnrw3dfwe3u8gtygE1} follows from \eqref{6hjdfv8e6hyefyeheew7sk} with \( k=3 \).
Next, we study the quadratic risk of \( \tilde{\upsv} \).
Define \( \epsv_{\QP} = \QP (\tilde{\upsv} - \upsvs + \svn) \). 
By \eqref{0mkvhgjnrwwe3u8gtyg} 
\begin{EQA}
	\sqrt{\E ( \| \epsv_{\QP} \|^{2} \Ind_{\Omega(\xx)} )}
	& \leq &
	\| \QP \IF^{-1} \DPN \| \, 
	\frac{\dltwu_{4}}{2} \,  
	\sqrt{\E \| \DPN \IF^{-1} \nabla \zeta \|^{6} \Ind_{\Omega(\xx)}} 
	\leq 
	\alp_{\QP,2} \sqrt{\riskt_{\QP,2}} \, ,
\label{d8ew3jfjhyvb6rt6543ejhvu}
\end{EQA}
and \eqref{6shx76whnjvyehfbvyfh} follows. 
Further, denote 
\begin{EQA}[rcccl]
	\loss_{\QP}
	& \eqdef &
	\QP \IF^{-1} \nabla \zeta \, ,
	\qquad
	\delta_{\QP}
	& \eqdef &
	\QP (\IF^{-1} \nabla \zeta + \svn) \, .
\label{7ejvuejfumfmrgvytweof7}
\end{EQA}
By definition, \( \riskt_{\QP} = \E \bigl\{ \| \loss_{\QP} \|^{2} \Ind_{\Omega(\xx)} \bigr\} \), 
\( \riskt_{\QP,2} = \E \bigl\{ \| \loss_{\QP} + \delta_{\QP} \|^{2} \Ind_{\Omega(\xx)} \bigr\} \), and 
\begin{EQA}
	\riskt_{\QP,2} - \riskt_{\QP}
	& = &
	\E \bigl\{ \| \delta_{\QP} \|^{2} \Ind_{\Omega(\xx)} \bigr\}
	+ 2 \E \bigl\{ \langle \loss_{\QP} , \delta_{\QP} \rangle \Ind_{\Omega(\xx)} \bigr\} .
\label{7jde7jevyreteyv8rnbe}
\end{EQA}
Also \eqref{iuvchycvf6e64rygh3224} and \eqref{6dhx6whcuydsds655srew4} imply
\begin{EQA}
	\sqrt{\E \bigl( \| \delta_{\QP} \|^{2} \Ind_{\Omega(\xx)} \bigr) }
	& \leq &
	\| \QP \IF^{-1} \DPN \| \, \frac{\dltwu_{3}}{2} \, 
	\sqrt{\E \| \DPN \IF^{-1} \nabla \zeta \|^{4} \Ind_{\Omega(\xx)} } 
	\\
	& \leq &
	\| \QP \IF^{-1} \DPN \| \, \frac{\dltwu_{3}}{2} \, \CONSTi_{4} \, \dimD 
	\leq 
	\alp_{\QP,1} \, \sqrt{\riskt_{\QP}} \, .
\label{yfhf73hjf9bryrnvbir}
\end{EQA}
This proves \eqref{EQtuGmstrVEQtGQ2}.
\end{proof}

\begin{remark}
As \( \| \DPN \, \IF^{-1} \nabla \zeta \| \leq \rrD \) on \( \Omega(\xx) \), it holds 
\begin{EQA}
	\E \bigl( \| \DPN \IF^{-1} \nabla \zeta \|^{4} \Ind_{\Omega(\xx)} \bigr)
	& \leq &
	\rrD^{2} \, \E \bigl( \| \DPN \IF^{-1} \nabla \zeta \|^{2} \Ind_{\Omega(\xx)} \bigr)
	\leq 
	\rrD^{2} \, \dimD \, .
\label{6hdnweyftegwgfywggvf636}
\end{EQA}
If \( \rrD^{2} \approx \dimD \), then  \( \CONSTi_{4} \approx 1 \) in \eqref{6dhx6whcuydsds655srew}.
\end{remark}

The results of Theorem~\ref{Teff41} enable us to improve the issue of \emph{critical dimension}.
For simplicity, let \( \QP = \DPN = \IF^{1/2} \).
Then the derived bounds are meaningful if 
\begin{EQA}
	( \dltwu_{4} + \dltwu_{3}^{2} ) \, \dimD^{3/2} 
	&=&
	o(1) .
\label{ufcjnciwf7fnrvuehj}
\end{EQA} 
Assuming \( \dltwu_{4} \asymp 1/n \) and \( \dltwu_{3}^{2} \asymp 1/n \),
we obtain the critical dimension condition
\( \dimD^{3/2} \ll n \) which is weaker than \( \dimD^{2} \ll n \).
Condition \eqref{6dhx6whcuydsds655srew4} ensuring equivalence of \( \riskt_{\QP,2} \) and \( \riskt_{\QP} \)
requires \( \dltwu_{3} \, \dimD \ll \riskt_{\QP} \) as in the 3S case.
}{}

\iffourG{
\Section{Fourth-order expansion and bias correction}
\label{SMLE4b}

Expansion \eqref{8vfjvr43f8khg54ed54} and \eqref{yfjvh9h4f5e53tghfugu44}
\begin{EQA}
	\tilde{\upsv} - \upsvs
	& \approx &
	- \svn
	=
	- \IF^{-1} \nabla \zeta - \IF^{-1} \nabla \TensU(\IF^{-1} \nabla \zeta)
\label{dfv8edefye3hejxs2qw}
\end{EQA}
involves the term \( \IF^{-1} \nabla \TensU(\IF^{-1} \nabla \zeta) \) which is quadratic in \( \IF^{-1} \nabla \zeta \).
If the dimension \( \dimp \) is large, one can expect that this quadratic form concentrates around its expectation.
This suggests a third-order correction of the MLE \( \tilde{\upsv} \):
\begin{EQA}
	\hat{\upsv}
	&=&
	\tilde{\upsv} + \IF^{-1} \, \E \, \nabla \TensU(\IF^{-1} \nabla \zeta),
	\\
	\hat{\upsv} - \upsvs
	& \approx &
	- \IF^{-1} \nabla \zeta 
	- \IF^{-1} \bigl\{ \nabla \TensU(\IF^{-1} \nabla \zeta) - \E \, \nabla \TensU(\IF^{-1} \nabla \zeta) \bigr\} .
\label{tsgqhyug8g900o4kqtdje}
\end{EQA}
To avoid technical issues and complicated notation, we do not distinguish between
\( \E \, \nabla \TensU(\IF^{-1} \nabla \zeta) \) and 
\( \E \, \bigl\{ \nabla \TensU(\IF^{-1} \nabla \zeta) \Ind_{\Omega(\xx)} \bigr\} \),
where \( \Omega(\xx) \) is a random set of dominating probability shown in Theorem~\ref{Teff41}.
Define 
\begin{EQ}[rcl]
\label{hdvje39bug53ebfh8edx3}
	\riskt_{\QP,3} 
	& \eqdef & 
	\E \bigl\{ \| \QP (\svn - \E \svn) \|^{2} \Ind_{\Omega(\xx)} \bigr\} \, ;
\end{EQ}
cf. \eqref{hdvje39bugebfh8edx}.
Obviously \( \riskt_{\QP,3} \leq \riskt_{\QP,2} \).
The next result is an immediate extension of Theorem~\ref{Teff41}.

\begin{theorem}
\label{Teff4bc}
Assume the conditions of Theorem~\ref{Teffp4s} and let
\begin{EQA}[c]
	\E \bigl\{ \| \DPN \IF^{-1} \nabla \zeta \|^{k} \Ind_{\Omega(\xx)} \bigr\} \leq \CONSTi_{k}^{2} \, \dimD^{k/2} ,
	\qquad
	k=3,4,6 
	\, .
\label{6hjdfv8e6hyefyeheew7skbc}
\end{EQA}
Then it holds for any linear mapping \( \QP \)
\begin{EQA}[rcl]
	&& \nquad
	\E \bigl\{ \| \QP \, (\hat{\upsv} - \upsvs) \| \Ind_{\Omega(\xx)} \bigr\}
	\\
	& \leq &
	\E \bigl\{ \| \QP (\svn - \E \svn) \| \Ind_{\Omega(\xx)} \bigr\}
	+ \| \QP \IF^{-1} \DPN \| \, 
	\Bigl( \frac{\dltwu_{4}}{2} + \dmax^{2} \dltwu_{3}^{2} \Bigr) \, \CONSTi_{3}^{2} \, \dimD^{3/2} \, ,
\label{0mkvhgjnrw3dfwe3u8gtygE1bc}
\end{EQA}
With \( \riskt_{\QP,3} \) from \eqref{hdvje39bug53ebfh8edx3}, let
\begin{EQA}
	\alp_{\QP,3}
	& \eqdef &
	\frac{\| \QP \IF^{-1} \DPN \| \, ( \dltwu_{4}/2 + \dmax^{2} \dltwu_{3}^{2}) \, \CONSTi_{6} \, \dimD^{3/2}}
		 {\sqrt{\riskt_{\QP,3}}} 
	< 1 
	\, .
\label{yfhcvhched6chejdrtebc}
\end{EQA}
Then
\begin{EQA}
	\bigl( 1 - \alp_{\QP,3} \bigr)^{2} \riskt_{\QP,3}
	\leq 
	\E \bigl\{ \| \QP \, (\hat{\upsv} - \upsvs) \|^{2} \Ind_{\Omega(\xx)} \bigr\}
	& \leq &
	\bigl( 1 + \alp_{\QP,3} \bigr)^{2} \riskt_{\QP,3} \, .
\label{6shx76whnjvyehfbvyfhbc}
\end{EQA}
\end{theorem}
}{}

\Chapter{Penalization, bias-variance decomposition}
\label{SPENbias}
This section explains how the results for SLS models can be extended to the penalized maximum likelihood approach.

\Section{Penalization bias}
\label{Ssmoothbias}

A common approach for improving the performance of MLE is based on regularization or penalization. 
The objective function \( L(\upsv) \) is extended by including a penalty term \( \pent_{\GP}(\upsv) \) 
which is responsible for complexity (roughness) of the parameter \( \upsv \).
A typical example to keep in mind is \( \pent_{\GP}(\upsv) = \| \GP \upsv \|^{2}/2 \)
for a penalization matrix \( \GP^{2} \).
Penalization by \( \pent_{\GP}(\upsv) \) can gradially improve  stability and numerical properties of the estimator,
however, it leads to a change of the ``truth'' \( \upsvs \), and hence, to some bias.
This section describes the bias caused by a smooth penalty.
Define the penalized MLE \( \tilde{\upsv}_{\GP} \)
\begin{EQA}
	\tilde{\upsv}_{\GP}
	& \eqdef &
	\argmin_{\upsv} \LGP(\upsv)
	=
	\argmin_{\upsv} \{ L(\upsv) + \pent_{\GP}(\upsv) \} \, .
\label{vdyedbwdtf6g7h8y6ie}
\end{EQA}
Compared with \eqref{tuGauLGususGE}, consider three optimization problems
\begin{EQA}[rcl]
	\tilde{\upsv}_{\GP} 
	&=& 
	\argmin_{\upsv} \LGP(\upsv),
	\qquad
	\upsvs_{\GP} 
	=
	\argmin_{\upsv} \E \LGP(\upsv),
	\qquad
	\upsvs 
	=
	\argmin_{\upsv} \E L(\upsv) .
	\qquad
\label{tuGauLGususGEgp}
\end{EQA}
Due to 
{Proposition~\ref{PconcMLEgenc}}, 
the penalized MLE  \( \tilde{\upsv}_{\GP} \) estimates rather \( \upsvs_{\GP} \) 
than \( \upsvs \).
This section describes the bias \( \upsvs_{\GP} - \upsvs \) caused by penalization.

Define the penalized Fisher information \( \IF_{\GP}(\upsv) = \nabla^{2} \E \LGP(\upsv) \) and
introduce \( \AvGP(\upsv) \eqdef \nabla \!\pent_{\GP}(\upsv) \).
Set \( \IF_{\GP} = \IF_{\GP}(\upsvs_{\GP}) \),
\begin{EQA}[rcccccl]
	\IF_{\GP} &=& \IF_{\GP}(\upsvs_{\GP}) \, ,
	\qquad
	\AvGP 
	& \eqdef &
	\nabla \!\pent_{\GP}(\upsvs) ,
	\quad
	\biasD
	& \eqdef &
	\| \DPN \IF_{\GP}^{-1} \AvGP \| 
	\, .
\label{ghdrd324ee4ew222}
\end{EQA}
For a quadratic penalty \( \pent_{\GP}(\upsv) = \| \GP \upsv \|^{2}/2 \), this results in
\begin{EQA}[rcccl]
	\AvGP
	&=&
	\GP^{2} \upsvs,
	\qquad
	\biasD
	& \eqdef &
	\| \DPN \IF_{\GP}^{-1} \GP^{2} \upsvs \| 
	\, .
\label{ghdrd324ee4ew222q}
\end{EQA}
Proposition~\ref{Pbiaspen} yields the following result.

\begin{proposition}
\label{Lvarusetb} 
Let \( \fG(\upsv) = \E \LGP(\upsv) \) satisfy \nameref{LLsT3ref} at \( \upsvs_{\GP} \) with some
metric tensor \( \DPN \) and values \( \rr \) and \( \dltwu_{3} \) such that
\begin{EQA}[c]
	\DPN^{2} \leq \dmax^{2} \, \IF_{\GP} \, ,
	\qquad 
	\rr \geq 3 \biasD/2 \, ,
	\qquad
	\dltwu_{3} \, \dmax^{2} \, \biasD < 4/9 ,
\label{7deyedhfg5563w6rfygj}
\end{EQA}
for \( \biasD \) from \eqref{ghdrd324ee4ew222}. 
Then 
\begin{EQA}[rcl]
\label{11ma3eaelDebQ}
	\| \DPN^{-1} \IF_{\GP} (\upsvs_{\GP} - \upsvs + \IF_{\GP}^{-1} \AvGP) \|
	& \leq &
	\frac{3\dltwu_{3}}{4} \, \biasD^{2} 
	\, .
\end{EQA}
The same bounds apply with \( \IF_{\GP}(\upsvs) \) in place of \( \IF_{\GP} = \IF_{\GP}(\upsvs_{\GP}) \).
\end{proposition}

\iffourG{
The accuracy of the expansions for the bias of a pMLE can be improved
under fourth-order smoothness of \( \fG(\upsv) = \E L_{\GP}(\upsv) \).
%
Consider the third-order tensor
\( \TensU(\uv) = \frac{1}{6} \langle \nabla^{3} \fs(\upsvs_{\GP}), \uv^{\otimes 3} \rangle \) and its gradient
\( \nabla \TensU(\uv) = \frac{1}{2} \langle \nabla^{3} \fs(\upsvs_{\GP}), \uv^{\otimes 2} \rangle \). 
Define a vector \( \bvn_{\GP} \) by
\begin{EQ}[rcl]
	\bvn_{\GP}
	&=&
	\IF_{\GP}^{-1} \AvGP + \IF_{\GP}^{-1} \nabla \Tens(\IF_{\GP}^{-1} \AvGP) 
	\, .
\label{8vfjvr43f8khg54ed5G}
\end{EQ}
Theorem~\ref{Pbiaspen} exlains the impact of using \( \bvn_{\GP} \) in place of \( \IF_{\GP}^{-1} \AvGP \).

\begin{theorem}
\label{Teffp4}
Assume \nameref{LLref}.
Let \nameref{LLsT3ref} and \nameref{LLsT4ref} hold at \( \upsvs_{\GP} \) and 
\begin{EQA}[c]
	\DPN^{2} \leq \dmax^{2} \, \IF_{\GP} \, ,
	\;\; 
	\rr \geq \frac{3}{2} \biasD \, ,
	\;\;
	\dmax^{2} \dltwu_{3} \, \biasD < \frac{4}{9} \, ,
	\;\;
	\dmax^{2} \dltwu_{4} \, \biasD^{2} < \frac{1}{3} 
	\, ,
\label{7deyedhfg5563w6rfygjLR4}
\end{EQA}
with \( \biasD = \| \DPN \IF_{\GP}^{-1} \AvGP \| \).
Then 
\begin{EQA}
\label{0mkvhgjnrw3dfwer45}
	\| \DPN^{-1} \IF_{\GP} (\upsvs_{\GP} - \upsvs + \bvn_{\GP}) \|
	& \leq &
	\Bigl( \frac{\dltwu_{4}}{2} + \dmax^{2} \dltwu_{3}^{2} \Bigr) \, \biasD^{3} \, ,
	\qquad
	\\
	\| \DPN^{-1} \IF_{\GP} (\bvn_{\GP} - \IF_{\GP}^{-1} \AvGP) \|
	& \leq &
	\frac{\dltwu_{3}}{2} \biasD^{2} \, .
\label{0mkvhgjnrw3dfwe3u8gtysg}
\end{EQA}
\end{theorem}
}{}

\Section{Loss and risk of the pMLE. Bias-variance decomposition}
\label{Slossrisksb}
Now we combine the previous results about the stochastic term \( \tilde{\upsv}_{\GP} - \upsvs_{\GP} \)
and the bias term \( \upsvs_{\GP} - \upsvs \) to obtain sharp bounds 
on the loss and risk of the pMLE \( \tilde{\upsv}_{\GP} \).

\begin{theorem}
\label{TQFiWibias}
Assume \nameref{Eref},
\nameref{EU2ref},
and 
\nameref{LLref}.
Let \( \fG(\upsv) = \E \LGP(\upsv) \) satisfy \nameref{LLsT3ref} at \( \upsvs_{\GP} \) with some
\( \DPN \), \( \rr \), and \( \dltwu_{3} \).
With \( (\rrD \vee \biasD) \eqdef \max\{ \rrD , \biasD \} \), assume
\begin{EQA}[c]
	\DPN^{2} \leq \dmax^{2} \, \IF_{\GP} \, ,
	\qquad 
	\rr \geq \frac{3}{2} (\rrD \vee \biasD) \, ,
	\qquad
	\dmax^{2} \dltwu_{3} \, (\rrD \vee \biasD) < \frac{4}{9} \, ;
\label{7deyedhfg5563w6rfygjLR}
\end{EQA}
see \eqref{y7s7d7d77dfdy7fuegue3j} and \eqref{ghdrd324ee4ew222}.
For any linear mapping \( \QP \colon \R^{\dimp} \to \R^{\dimq} \), 
it holds on \( \Omega(\xx) \)
\begin{EQA}
	\| \QP (\tilde{\upsv}_{\GP} - \upsvs + \IF_{\GP}^{-1} \nabla \zeta + \IF_{\GP}^{-1} \AvGP) \|
	& \leq &
	\| \QP \IF_{\GP}^{-1} \DPN \| \, \frac{3\dltwu_{3}}{4} \, \bigl( \| \DPN \IF_{\GP}^{-1} \nabla \zeta \|^{2} + \biasD^{2} \bigr) \,
	\, .
	\qquad
	\quad
\label{g25re9fjfregdndgb}
\end{EQA}
Also, introduce \( \dimQ \eqdef \E \| \QP \IF_{\GP}^{-1} \nabla \zeta \|^{2} = \tr \Var(\QP \IF_{\GP}^{-1} \nabla \zeta) \) and
\begin{EQA}[c]
	\riskt_{\QP} \eqdef \E \{ \| \QP \IF_{\GP}^{-1} (\nabla \zeta + \AvGP) \|^{2} \Ind_{\Omega(\xx)} \} 
	\leq 
	\dimQ + \| \QP \IF_{\GP}^{-1} \AvGP \|^{2} \, .
\label{7djhed8cjfct534etgdhdy}
\end{EQA}
Then
\begin{EQA}
	\E \bigl\{ \| \QP (\tilde{\upsv}_{\GP} - \upsvs) \| \Ind_{\Omega(\xx)} \bigr\} 
	& \leq &
	\riskt_{\QP}^{1/2} 
	+ \| \QP \IF_{\GP}^{-1} \DPN \| \, \frac{3\dltwu_{3}}{4} \bigl( \dimD + \biasD^{2} \bigr) \, .
\label{EtuGus11md3GQ}
\end{EQA}
Further, assume \( \E \bigl\{ \| \DPN \IF_{\GP}^{-1} \nabla \zeta \|^{4} \Ind_{\Omega(\xx)} \bigr\} \leq \CONSTi_{4}^{2} \, \dimD^{2} \) 
and define
\begin{EQ}[rcl]
	\alp_{\QP}
	& \eqdef & 
	\frac{\| \QP \IF_{\GP}^{-1} \DPN \| \, (3/4) \dltwu_{3} \, (\CONSTi_{4} \, \dimD + \biasD^{2})} {\sqrt{\riskt_{\QP}}}
	\, . 
\label{6dhx6whcuydsds655srewb}
\end{EQ}
If \( \alp_{\QP} < 1 \) then
\begin{EQA}
	(1 - \alp_{\QP})^{2} \riskt_{\QP} 
	\leq 
	\E \bigl\{ \| \QP \, (\tilde{\upsv}_{\GP} - \upsvs) \|^{2} \Ind_{\Omega(\xx)} \bigr\}
	& \leq &
	(1 + \alp_{\QP})^{2} \riskt_{\QP} \, .
\label{EQtuGmstrVEQtGQ}
\end{EQA}
\end{theorem}

\begin{proof}
It holds by \eqref{DGttGtsGDGm13rG22} and \eqref{11ma3eaelDebQ}
\begin{EQ}[rcl]
	\| \QP (\tilde{\upsv}_{\GP} - \upsvs_{\GP} + \IF_{\GP}^{-1} \nabla \zeta \|
	& \leq &
	\| \QP \IF_{\GP}^{-1} \DPN \| \, \frac{3\dltwu_{3}}{4} \, \| \DPN \IF_{\GP}^{-1} \nabla \zeta \|^{2} \, ,
	\\
	\| \QP (\upsvs_{\GP} - \upsvs + \IF_{\GP}^{-1} \AvGP) \|
	& \leq &
	\| \QP \IF_{\GP}^{-1} \DPN \| \, \frac{3\dltwu_{3}}{4} \, \biasD^{2} \, ,
\label{u87dcudenhedst6cbweq2q}
\end{EQ}
and hence,
\begin{EQA}
	\| \QP (\tilde{\upsv}_{\GP} - \upsvs + \IF_{\GP}^{-1} \nabla \zeta + \IF_{\GP}^{-1} \AvGP) \|
	& \leq &
	\| \QP \IF_{\GP}^{-1} \DPN \| \, \frac{3\dltwu_{3}}{4} \, \bigl( \| \DPN \IF_{\GP}^{-1} \nabla \zeta \|^{2} + \biasD^{2} \bigr)
\label{ydy6fc6fv6e43yte46fhuQ}
\end{EQA}
yielding \eqref{g25re9fjfregdndgb} and \eqref{EtuGus11md3GQ}.
Further, define 
\begin{EQA}
	\epsv_{\GP}
	& \eqdef &
	\QP \bigl\{ \tilde{\upsv}_{\GP} - \upsvs + \IF_{\GP}^{-1} (\nabla \zeta + \AvGP) \bigr\} \, .
\label{dfu8nj3eyvnhedd6ywhwhvfu7}
\end{EQA}
It holds by \eqref{u87dcudenhedst6cbweq2q}
\begin{EQA}
	\E^{1/2} \bigl\{ \| \epsv_{\QP} \|^{2} \Ind_{\Omega(\xx)} \bigr\}
	& \leq &
	\| \QP \IF_{\GP}^{-1} \DPN \| \frac{3\dltwu_{3}}{4} \, 
	\bigl\{ \E^{1/2} \bigl( \| \DPN \IF_{\GP}^{-1} \nabla \zeta \|^{4} \Ind_{\Omega(\xx)} \bigr) + \biasD^{2} \bigr\} 
	\leq 
	\alp_{\QP} \, \riskt_{\QP}^{1/2} \, ,
\label{d7hw3fujgv76444gdfu7eQ}
\end{EQA}
and therefore,
\begin{EQA}
	&& \nquad
	\E^{1/2} \bigl\{ \| \QP \, (\tilde{\upsv}_{\GP} - \upsvs) \|^{2} \Ind_{\Omega(\xx)} \bigr\}
	= 
	\E^{1/2} \bigl\{ \| \QP \IF_{\GP}^{-1} (\nabla \zeta + \AvGP) - \epsv_{\QP} \|^{2} \Ind_{\Omega(\xx)} \bigr\} 
	\\
	& \leq &
	\E^{1/2} \bigl\{ \| \QP \IF_{\GP}^{-1} (\nabla \zeta + \AvGP) \|^{2} \Ind_{\Omega(\xx)} \bigr\}
	+ \E^{1/2} \bigl\{ \| \epsv_{\QP} \|^{2} \Ind_{\Omega(\xx)} \bigr\}
	\leq 
	(1 + \alp_{\QP}) \, \riskt_{\QP}^{1/2} \, .
\label{7dnmswewe8fyhbdfc7er}
\end{EQA}
This yields \eqref{EQtuGmstrVEQtGQ}.
\end{proof}

\begin{remark}
\label{RremainderD}
The condition \( \DPN^{2} \leq \dmax^{2} \IF_{\GP} \) implies \( \| \QP \IF_{\GP}^{-1} \DPN \| \leq \dmax^{2} \| \QP \DPN^{-1} \| \)
which can be used in the remainder for all risk bounds. 
\end{remark}

\begin{remark}
Due to \eqref{EQtuGmstrVEQtGQ} 
\begin{EQA}
	\E \bigl\{ \| \QP \, (\tilde{\upsv}_{\GP} - \upsvs) \|^{2} \Ind_{\Omega(\xx)} \bigr\}
	& = &
	\bigl( \dimQ + \| \QP \IF_{\GP}^{-1} \AvGP \|^{2} \bigr) \, \bigl\{ 1 + o(1) \bigr\}.
	\qquad
\label{EQtuGuvus2tr1o1}
\end{EQA}
This relation is usually referred to as ``bias-variance decomposition''.
Our bound is sharp in the sense that for the special case of linear models, 
\eqref{EQtuGuvus2tr1o1} becomes equality.
Under the so-called ``small bias'' condition \( \| \QP \IF_{\GP}^{-1} \AvGP \|^{2} \ll \dimQ \), 
the impact of the bias induced by penalization is negligible. 
\ifNL{
The relation \( \| \QP \IF_{\GP}^{-1} \AvGP \|^{2} \asymp \dimQ \) is called ``bias-variance trade-off'',
it leads to minimax rate of estimation;
see Section~\ref{Spriorexample}.}{}
\end{remark}

If the constant \( \alp_{\QP} \) from \eqref{6dhx6whcuydsds655srewb} satisfies \( \alp_{\QP} \ll 1 \) then 
by \eqref{EQtuGmstrVEQtGQ}, 
\( \E \bigl\{ \| \QP \, (\tilde{\upsv}_{\GP} - \upsvs) \|^{2} \Ind_{\Omega(\xx)} \bigr\} = (1 + o(1)) \riskt_{\QP} \).
\iffourG{
Now we explain how the accuracy of the expansions for pMLE can be improved
and the critical dimension condition can be relaxed under 
fourth-order smoothness of \( \fG(\upsv) = \E L_{\GP}(\upsv) \).
%
Putting together the results on the stochastic component \( \tilde{\upsv}_{\GP} - \upsvs_{\GP} \) and 
on the bias \( \upsvs_{\GP} - \upsvs \) yields the bound on the loss and risk of the estimator \( \tilde{\upsv}_{\GP} \).
Define 
\begin{EQ}[rcl]
	\riskt_{\QP}
	& \eqdef &
	\E \bigl\{ \| \QP \IF_{\GP}^{-1} (\nabla \zeta + \AvGP) \|^{2} \Ind_{\Omega(\xx)} \bigr\} \, ,
\label{7dhdrw2bfvu78u78e4ndw}
	\\
\label{hdvje39bug53ebfh8edx}
	\riskt_{\QP,2} 
	& \eqdef & 
	\E \bigl\{ \| \QP (\svn_{\GP} + \bvn_{\GP}) \|^{2} \Ind_{\Omega(\xx)} \bigr\} \, .
\end{EQ}

\begin{theorem}
\label{Teff4s}
Assume the conditions of Theorem~\ref{Teffp4s} and Theorem~\ref{Teffp4} and let
\begin{EQA}[c]
	\E \bigl\{ \| \DPN \IF_{\GP}^{-1} \nabla \zeta \|^{k} \Ind_{\Omega(\xx)} \bigr\} \leq \CONSTi_{k}^{2} \, \dimD^{k/2} ,
	\qquad
	k=3,4,6 
	\, .
\label{6hjdfv8e6hyefyeheew7skb}
\end{EQA}
Then it holds for any linear mapping \( \QP \)
\begin{EQA}
	&& \nquad
	\E \bigl\{ \| \QP \, (\tilde{\upsv}_{\GP} - \upsvs) \| \Ind_{\Omega(\xx)} \bigr\}
	\\
	& \leq &	 
	\E \bigl\{ \| \QP (\svn_{\GP} + \bvn_{\GP}) \| \Ind_{\Omega(\xx)} \bigr\}
	+ \| \QP \IF_{\GP}^{-1} \DPN \| \, 
	\Bigl( \frac{\dltwu_{4}}{2} + \dmax^{2} \dltwu_{3}^{2} \Bigr) \, \bigl( \CONSTi_{3}^{2} \, \dimD^{3/2} + \biasD^{3} \bigr) ,
\label{0mkvhgjnrw3dfwe3u8gtygE1b}
	\\
	&& \nquad
	\Bigl| \E \bigl\{ \| \QP (\svn_{\GP} + \bvn_{\GP}) \| \Ind_{\Omega(\xx)} \bigr\}
	- \E \bigl\{ \| \QP \IF_{\GP}^{-1} \, \nabla \zeta + \QP \IF_{\GP}^{-1} \, \AvGP \| \Ind_{\Omega(\xx)} \bigr\} 
	\Bigr|
	\\
	& \leq &
	\| \QP \IF_{\GP}^{-1} \DPN \| \, \frac{\dltwu_{3}}{2} \, \bigl( \dimD + \biasD^{2} \bigr) \, .
\label{udtgecthwjdytdehduuc6}
\end{EQA}
With \( \riskt_{\QP,2} \) from \eqref{7dhdrw2bfvu78u78e4ndw}, let
\begin{EQA}
	\alp_{\QP,2}
	& \eqdef &
	\frac{\| \QP \IF_{\GP}^{-1} \DPN \| \, ( \dltwu_{4}/2 + \dmax^{2} \dltwu_{3}^{2}) \, (\CONSTi_{6} \, \dimD^{3/2} + \biasD^{3})}
		 {\sqrt{\riskt_{\QP,2}}} 
	< 1 
	\, .
\label{yfhcvhched6chejdrte}
\end{EQA}
Then
\begin{EQA}
	\bigl( 1 - \alp_{\QP,2} \bigr)^{2} \riskt_{\QP,2}
	\leq 
	\E \bigl\{ \| \QP \, (\tilde{\upsv}_{\GP} - \upsvs) \|^{2} \Ind_{\Omega(\xx)} \bigr\}
	& \leq &
	\bigl( 1 + \alp_{\QP,2} \bigr)^{2} \riskt_{\QP,2} \, .
\label{6shx76whnjvyehfbvyfhb}
\end{EQA}
If another constant \( \alp_{\QP,1} < 1 \) ensures 
\begin{EQ}[rcl]
	&&
	\| \QP \IF_{\GP}^{-1} \DPN \| \, \frac{\dltwu_{3}}{2} \, \bigl( \CONSTi_{4} \, \dimD + \biasD^{2} \bigr)
	\leq 
	\alp_{\QP,1} \, \sqrt{\riskt_{\QP}} \, 
\label{6dhx6whcuydsds655srew4b}
\end{EQ}
with \( \riskt_{\QP} \) from \eqref{7dhdrw2bfvu78u78e4ndw} then
\begin{EQA}
	\riskt_{\QP} (1 - \alp_{\QP,1})^{2} 
	\leq 
	\riskt_{\QP,2}
	& \leq &
	\riskt_{\QP} (1 + \alp_{\QP,1})^{2} \, .
\label{EQtuGmstrVEQtGQ2b}
\end{EQA}
\end{theorem}

\begin{proof}
Rescaling of \( \DPN \) reduces the proof to \( \dmax = 1 \).
Theorem~\ref{Teffp4} yields
\begin{EQA}
	\| \QP \, (\tilde{\upsv}_{\GP} - \upsvs + \svn_{\GP} + \bvn_{\GP}) \|
	& \leq &
	\| \QP \IF_{\GP}^{-1} \DPN \| \, 
	\Bigl( \frac{\dltwu_{4}}{2} + \dltwu_{3}^{2} \Bigr) \,  
	\bigl( \| \DPN \IF_{\GP}^{-1} \nabla \zeta \|^{3} + \biasD^{3} \bigr) ,
	\qquad
	\qquad
\label{0mkvhgjnrwwe3u8gtygb}
	\\
	\| \QP \{ \svn_{\GP} - \bvn_{\GP} - \IF_{\GP}^{-1} (\nabla \zeta + \AvGP) \} \|
	& \leq &
	\frac{\dltwu_{3}}{2} \, \| \QP \IF_{\GP}^{-1} \DPN \| \, \bigl( \| \DPN \IF_{\GP}^{-1} \nabla \zeta \|^{2} + \biasD^{2} \bigr) \, .
\label{vuedy766t4e3bfvyt6e}
\end{EQA}
Now \eqref{0mkvhgjnrw3dfwe3u8gtygE1b} follows from \eqref{6hjdfv8e6hyefyeheew7skb} with \( k=3 \).
Next, we study the quadratic risk of \( \tilde{\upsv}_{\GP} \).
Define \( \epsv_{\QP} = \QP (\tilde{\upsv}_{\GP} - \upsvs - \svn_{\GP} + \bvn_{\GP}) \). 
By \eqref{0mkvhgjnrwwe3u8gtygb} 
\begin{EQA}
	\sqrt{\E ( \| \epsv_{\GP} \|^{2} \Ind_{\Omega(\xx)} )}
	& \leq &
	\| \QP \IF_{\GP}^{-1} \DPN \| \, 
	\Bigl( \frac{\dltwu_{4}}{2} + \dltwu_{3}^{2} \Bigr) \,  
	\Bigl( \sqrt{\E \| \DPN \IF_{\GP}^{-1} \nabla \zeta \|^{6} \Ind_{\Omega(\xx)}}  + \biasD^{3} \Bigr)
	\leq 
	\alp_{\QP,2} \sqrt{\riskt_{\QP,2}} \, ,
\label{d8ew3jfjhyvb6rt6543ejhvu}
\end{EQA}
and \eqref{6shx76whnjvyehfbvyfhb} follows. 
Further, denote 
\begin{EQA}
	\loss_{\QP}
	& \eqdef &
	\QP \IF_{\GP}^{-1} (\nabla \zeta + \AvGP) ,
	\\
	\delta_{\QP}
	& \eqdef &
	\QP (\IF_{\GP}^{-1} \nabla \zeta - \svn_{\GP}) 
	+ \QP (\IF_{\GP}^{-1} \AvGP - \bvn_{\GP}) .
\label{7ejvuejfumfmrgvytweof7}
\end{EQA}
By definition, \( \riskt_{\QP} = \E \bigl\{ \| \loss_{\QP} \|^{2} \Ind_{\Omega(\xx)} \bigr\} \), 
\( \riskt_{\QP,2} = \E \bigl\{ \| \loss_{\QP} + \delta_{\QP} \|^{2} \Ind_{\Omega(\xx)} \bigr\} \), and 
\begin{EQA}
	\riskt_{\QP,2} - \riskt_{\QP}
	& = &
	\E \bigl\{ \| \delta_{\QP} \|^{2} \Ind_{\Omega(\xx)} \bigr\}
	+ 2 \E \bigl\{ \langle \loss_{\QP} , \delta_{\QP} \rangle \Ind_{\Omega(\xx)} \bigr\} .
\label{7jde7jevyreteyv8rnbe}
\end{EQA}
Also \eqref{6hjdfv8e6hyefyeheew7skb}, \eqref{0mkvhgjnrw3dfwe3u8gtysg}, and \eqref{6dhx6whcuydsds655srew4b} imply
\begin{EQA}
	\sqrt{\E \bigl( \| \delta_{\QP} \|^{2} \Ind_{\Omega(\xx)} \bigr) }
	& \leq &
	\| \QP \IF_{\GP}^{-1} \DPN \| \, \frac{\dltwu_{3}}{2} \, 
	\Bigl( \sqrt{\E \| \DPN \IF_{\GP}^{-1} \nabla \zeta \|^{4} \Ind_{\Omega(\xx)} } + \biasD^{2} \Bigr)
	\\
	& \leq &
	\| \QP \IF_{\GP}^{-1} \DPN \| \, \frac{\dltwu_{3}}{2} \, \bigl( \CONSTi_{4} \, \dimD + \biasD^{2} \bigr)
	\leq 
	\alp_{\QP,1} \, \sqrt{\riskt_{\QP}} \, .
\label{yfhf73hjf9bryrnvbir}
\end{EQA}
This proves \eqref{EQtuGmstrVEQtGQ2b}.
\end{proof}
}{}

\ifsupnorm{\input nanonorm}{}
%
\def\GPp{\GP_{\priord}}
\def\GPb{\GP_{1}}
\def\DPp{\DP_{\priord}}
\def\projp{\text{\Large$\pi$}}
\def\Projc{\mathcal{P}}
\def\dimmc{\dimm'}
\def\dimms{\dimm^{*}}
\def\CGP{w}
\def\CGPw{\tau}
\def\CGPa{\CGP_{0}}
\def\CGPb{\CGP_{1}}
\def\CGPs{\CGPw^{*}}
\def\Fclass{\mathcal{F}}
\def\CONSTGPT{\CONSTi_{\GPT}}

\Section{Bias-variance trade-off and oracle risk bounds}
\label{Spriorexample}
This section presents some examples of choosing \( \GP^{2} \) for 
achieving the ``bias-variance trade-off''
and obtaining rate optimal results.
Let pMLE \( \tilde{\upsv}_{\GP} \), its population counterpart \( \upsvs_{\GP} \),
and the background true parameter \( \upsvs \) be given by \eqref{tuGauLGususGE}.
Theorem~\ref{TQFiWibias} yields the following bound for the risk 
\( \riskt_{\QP} \) of \( \tilde{\upsv}_{\GP} \) given by \eqref{7djhed8cjfct534etgdhdy}:
\begin{EQA}
	\E \| \QP (\tilde{\upsv}_{\GP} - \upsvs) \|^{2}
	& \approx &
	\riskt_{\QP}
	\approx
	\dimA_{\QP} + \| \QP \IF_{\GP}^{-1} \GP^{2} \upsvs \|^{2} \, .
\label{bnjmb3ed6w2jh21fg7fj}
\end{EQA}
This suggest to select the operator \( \GP^{2} \) by forcing the ``bias-variance trade-off'' 
\( \dimA_{\QP} \asymp \| \QP \IF_{\GP}^{-1} \GP^{2} \upsvs \|^{2} \).
Later in this section, we illustrate this relation through popular examples of regularization by projection or by roughness penalty.
For any considered choice of penalization \( \GP^{2} \), we assume the conditions of Propositions~\ref{PconcMLEgenc} and \ref{Lvarusetb}
to be fulfilled. 
To simplify the analysis, we also assume \( \VP^{2} = \DPN^{2} = \IF \) and consider two specific choices of \( \QP \):
prediction/response loss with \( \QP = \IF^{1/2} \) and
estimation loss with \( \QP = \Id \). 
In the latter case, we focus on a direct problem with a bounded condition number of \( \IF \).

\Subsection{Projection estimation}
\label{SprojGP}
Consider the class of projection estimators given by a set of sub-spaces \( \{ \II_{\dimm} \} \) of the parameter space \( \R^{\dimp} \).
For each \( m \), only projection \( \Proj_{\dimm} \upsv \) on the subspace \( \II_{\dimm} \) is considered 
but there is no any additional penalization.
Formally, this corresponds to the diagonal matrix \( \GP_{\dimm}^{2} \) with \( \dimm \) diagonal elements equal to zero,
and the remaining ones equal to infinity.
Later we everywhere use the sub-index \( \dimm \) in place of \( \GP_{\dimm} \).
It appears that \( \IF_{\dimm}^{-1}(\upsv) \, \GP_{\dimm}^{2} = \{ \IF(\upsv) + \GP_{\dimm}^{2} \}^{-1} \GP_{\dimm}^{2} \) 
for any \( \upsv \in \Ups \) is nothing but the orthogonal projector \( \Projc_{\dimm} = \Id - \Proj_{\dimm} \) 
on the orthogonal subspace \( \II_{\dimm}^{c} \).
Also,
\begin{EQA}
	\IF_{\dimm}^{-1} \, \GP_{\dimm}^{2} \upsvs
	&=&
	\Projc_{\dimm} \upsvs
	=
	\upsvs - \Proj_{\dimm} \upsvs \, .
\label{yeufyy3hfyfhwmdygneum}
\end{EQA}
Similarly, \( \IF_{\dimm}^{-1} \IF = \Proj_{\dimm} \) and \( \IF_{\dimm}^{-1} \IF \, \IF_{\dimm}^{-1} = \IF_{\dimm}^{-1} \). 
As \( \DPN^{2} = \VP^{2} = \IF \), this leads to 
\begin{EQA}[c]
	\dimA_{\QP,\dimm} = \tr(\QP \IF_{\dimm}^{-1} \VP^{2} \IF_{\dimm}^{-1} \QP^{\T}) 
	= 
	\tr (\QP \IF_{\dimm}^{-1} \QP^{\T}) 
	\, .
\label{uef8efhjfjfbufb8874uf}
\end{EQA}
In particular,
\begin{EQA}[c]
	\dimA_{\QP,\dimm} 
	= 
	\begin{cases}
		\tr (\Proj_{\dimm}) = \dimm, 	& \QP = \IF^{1/2} , \\
		\tr(\IF_{\dimm}^{-1}), 		& \QP = \Id .
	\end{cases}
\label{uef8efhjfjfbufb8874ufpr}
\end{EQA}
For the corresponding risk \( \riskt_{\QP,\dimm} \) from \eqref{bnjmb3ed6w2jh21fg7fj}, we obtain by Theorem~\ref{TQFiWibias}
\begin{EQA}
	\riskt_{\QP,\dimm}
	& \approx &
	\begin{cases}
		\dimm + \| \IF^{1/2} \Projc_{\dimm} \upsvs \|^{2} & \QP = \IF^{1/2} , \\
		\tr(\IF_{\dimm}^{-1}) + \| \Projc_{\dimm} \upsvs \|^{2}, 		& \QP = \Id .
	\end{cases}
	 \, .
\label{hgfcuygeuhfui3w23}
\end{EQA}
The optimal (or oracle) choice of \( \dimm \) can be given by minimization of the risk \( \riskt_{\QP,\dimm} \):
\begin{EQA}
	\dimms
	& \eqdef &
	\argmin_{\dimm} \riskt_{\QP,\dimm}
	\, .
\label{d7fjfjjvytv435y6ewee}
\end{EQA}
%
A standard way of obtaining the minimax rate of estimation is based on the 
approximation theory for functional spaces.
One assumes that \( \upsvs \) belongs to a special set \( \Fclass \) like a Sobolev or Besov ball, and 
\begin{EQA}
	\| \upsv - \Proj_{\dimm} \upsv \|
	& \leq &
	\rho_{\dimm} \, , 
	\qquad
	\upsv \in \Fclass \, ,
\label{edydjd7u73ef6fhej7hj}
\end{EQA}
where the \( \rho_{\dimm} \)'s are \( \Fclass \)-specific and decrease to zero as \( \II_{\dimm} \) increase.
As an example, consider a ``smooth'' signal \( \upsvs \) from a Sobolev ball \( \BBB(\smpa,\CGPa) \):
\begin{EQA}
	\BBB(\smpa,\CGPa)
	& \eqdef &
	\Bigl\{ \upsv = (\ups_{j}) \in \R^{\dimp} \colon \sum_{j \geq 1} j^{2\smpa} \ups_{j}^{2} \leq \CGPa \Bigr\} 
\label{BBBsttsjj1sj2S}
\end{EQA}
with \( \smpa > 0 \) and \( \CGPa \asymp 1 \).
Then for any \( \upsv \in \BBB(\smpa,\CGPa) \)
\begin{EQA}
	\| \Projc_{\dimm} \upsv \|^{2}
	& \leq &
	\dimm^{-2\smpa} \sum_{j > \dimm} j^{2\smpa} \ups_{j}^{2}
	\leq 
	\rho_{\dimm}^{2}
	=
	\CGPa \, \dimm^{-2\smpa} \, .
\label{ysgduy27828e3h2ewsud}
\end{EQA}
We additionally assume that \( \IF \geq n \Id \), where \( n \) is a scaling parameter meaning the sample size.
For \( \QP = \Id_{\dimp} \), it holds
\begin{EQA}[c]
	\tr(\IF_{\dimm}^{-1})
	\leq 
	\dimm/n \, .
\end{EQA}
Therefore, the \( \upsvs \)-dependent choice \eqref{d7fjfjjvytv435y6ewee} can be replaced by the \( \Fclass \)-specific choice
\begin{EQA}
	\dimms
	&=&
	\argmin_{\dimm} \{ \dimm/n + \rho_{\dimm}^{2} \} \, .
\label{8fifivugvy5435tdjbuuyr}
\end{EQA}
Typically the solution to this problem satisfies the balance relation \( \dimms \asymp n \rho_{\dimms}^{2} \)
leading to the risk \( \riskt_{\QP,\dimms} \asymp \dimms/n \).
For the case of a Sobolev ball, \( \rho_{\dimm}^{2} = \CGPa \dimm^{-2\smpa} \), 
and the trade-off relation reads as \( \dimm/n \asymp \dimm^{-2\smpa} \).
This leads to the standard rule of thumb 
\( \dimms \asymp n^{1/(2 \smpa+1)} \) and \( \riskt_{\QP,\dimms} \asymp n^{-2\smpa/(2\smpa+1)} \).

For the prediction loss with \( \QP = \IF^{1/2} \), the situation is similar as long as 
a \emph{direct} problem with a bounded condition number \( \CONSTIF = \lambda_{\max}(\IF)/\lambda_{\min}(\IF) \) is considered.
Assume that \( n \Id \leq \IF \leq \CONSTIF \, n \Id \), where \( n = \lambda_{\min}(\IF) \). 
As \( \IF \geq n \, \Id_{\dimp} \), we obtain 
\( \| \IF^{1/2} \Projc_{\dimm} \upsvs \|^{2} \geq n \, \rho_{\dimm}^{2} \).
Therefore, the optimal choice of \( \dimm \) can be reduced to minimization of 
\( \dimm + n \, \rho_{\dimm}^{2} \) which coincides with \eqref{8fifivugvy5435tdjbuuyr}.
For the case of a Sobolev ball with \( \rho_{\dimm}^{2} = \CGPa \dimm^{-2\smpa} \),
this yields \( \dimms \asymp n^{1/(2 \smpa+1)} \) and \( \riskt_{\dimms} \asymp n^{1/(2\smpa+1)} \).


\Subsection{Roughness penalty}
In this section we discuss another more general class of penalizing families 
\( \GPT = \{ \GP^{2} \} \) with polynomially growing eigenvalues and
show that the risk of each \( \tilde{\upsv}_{\GP} \) 
can be decomposed and analyzed as in the case projection estimation with a proper choice
of the projection sub-space.
Assume as earlier that \( \VP^{2} = \IF \).
For any \( \QP \) and any \( \GP^{2} \in \GPT \), it holds
\begin{EQA}[c]
	\dimA_{\QP} = \tr(\QP \IF_{\GP}^{-1} \VP^{2} \IF_{\GP}^{-1} \QP^{\T}) 
	= 
	\tr (\QP \IF_{\GP}^{-1} \IF \, \IF_{\GP}^{-1} \QP^{\T}) 
\label{uef8efhjfjfbufb8GP}
\end{EQA}
and
\begin{EQA}[c]
	\dimA_{\QP} 
	= 
	\begin{cases}
		\tr (\IF^{1/2} \IF_{\GP}^{-1} \IF^{1/2})^{2} , 	& \QP = \IF^{1/2} , \\
		\tr(\IF_{\GP}^{-2} \IF), 		& \QP = \Id_{\dimp} \, .
	\end{cases}
\label{uef8efhjfjfbufbefi8GP}
\end{EQA}
Similarly
\begin{EQA}
	\bias_{\QP}
	&=&
	\| \QP \IF_{\GP}^{-1} \GP^{2} \upsvs \|^{2}
	=
	\begin{cases}
		\| \IF^{1/2} \IF_{\GP}^{-1} \GP^{2} \upsvs \|^{2} , 	& \QP = \IF^{1/2} , \\
		\| \IF_{\GP}^{-1} \GP^{2} \upsvs \|^{2}, 		& \QP = \Id_{\dimp} .
	\end{cases}
\label{d7ewjd6jr4g8ejefgghi4}
\end{EQA}
The aim is to describe these quantities and the related risk bounds in terms of the spectral characteristics of the penalty-to-signal matrix
\( \BB \eqdef \IF^{-1/2} \GP^{2} \IF^{-1/2} \).
Later we assume that \( \BB \) fulfills the polynomial growth condition on its spectrum.

\begin{description}
\item[\( \bb{(\operatorname{poly})} \)\label{GPTref}]
	\emph{Let \( \bp_{1}^{2} \leq \ldots \leq \bp_{\dimp}^{2} \) be increasing eigenvalues
	of \( \BB \).
	Then for some fixed constant \( \CONSTi_{\BB} \) and all \( \dimm < \dimp \) 
	}
\begin{EQA}
	\sum_{j = \dimm+1}^{\dimp} \bp_{j}^{-4}
	& \leq &
	\CONSTi_{\BB} \, \dimm \, \bp_{\dimm+1}^{-4} ,
	\qquad
	\sum_{j=1}^{\dimm} \bp_{j}^{4}
	\leq 
	\CONSTi_{\BB} \, \dimm \, \bp_{\dimm}^{4} \, .
\label{dsujjue67he3jv7nmeuej}
\end{EQA}

\end{description}

Condition \eqref{dsujjue67he3jv7nmeuej} assumes that \( \bp_{j}^{2} \) grow at least as \( j^{2\smpa} \) for
\( \smpa > 1/4 \).
The constant \( \CONSTi_{\BB} \) depends on \( \smpa \) only.

\begin{lemma}
\label{LGPgrowth}
Let \( \BB = \IF^{-1/2} \GP^{2} \IF^{-1/2} \) with increasing eigenvalues 
\( \bp_{1}^{2} \leq \ldots \leq \bp_{\dimp}^{2} \) satisfying \nameref{GPTref}.
Then for any \( \dimm \)
\begin{EQA}
	\tr (\IF^{1/2} \IF_{\GP}^{-1} \IF^{1/2})^{2}
	& \leq &
	\Bigl( 1 + \frac{\CONSTi_{\BB}}{\bp_{\dimm+1}^{4}} \Bigr) \, \dimm 
	\, .
\label{r8uerf6yghe5erjkrf7e3j}
\end{EQA}
In particular, if \( \dimm_{\GP} \) is the largest \( m \) such that \( \bp_{\dimm}^{2} \leq 1 \) then
\begin{EQA}
	\tr (\IF^{1/2} \IF_{\GP}^{-1} \IF^{1/2})^{2}
	& \leq &
	(1 + \CONSTi_{\BB}) \dimm_{\GP} 
	\, .
\label{r8uerf6yghe5erjkrf7e3jGP}
\end{EQA}
If \( \IF \geq n \Id_{\dimp} \) then 
\begin{EQA}
	\tr (\IF_{\GP}^{-2} \IF)
	& \leq &
	n^{-1} \tr (\IF^{1/2} \IF_{\GP}^{-1} \IF^{1/2})^{2}
	\leq 
	(1 + \CONSTi_{\BB}) \dimm_{\GP}
	\, .
\label{8d8d8d8h3rf76fyj34ffg}
\end{EQA}
\end{lemma}

\begin{proof}
Without loss of generality, assume that the matrix 
\( \BB = \IF^{-1/2} \GP^{2} \IF^{-1/2} \) is diagonal, that is, 
\( \BB = \diag\{ \bp_{1}^{2}, \ldots, \bp_{\dimp}^{2} \} \).
It holds by \eqref{dsujjue67he3jv7nmeuej} for any \( \dimm \)
\begin{EQA}
	\tr (\IF^{1/2} \IF_{\GP}^{-1} \IF^{1/2})^{2}
	&=&
	\tr (\Id_{\dimp} + \BB)^{-2}
	\leq 
	\sum_{j=1}^{\dimm} \frac{1}{(1 + \bp_{j}^{2})^{2}}
	+ \sum_{j=\dimm+1}^{\dimp} \frac{1}{(1 + \bp_{j}^{2})^{2}}
	\\
	& \leq &
	\dimm + \sum_{j=\dimm+1}^{\dimp} \frac{1}{\bp_{j}^{4}}
	\leq 
	\dimm + \CONSTi_{\BB} \, \dimm \, \bp_{\dimm+1}^{-4}
\label{uf8f73jr7jewlfgu7e4hjfk}
\end{EQA}
and the first bound follows.
Further, the definition of \( \dimm_{\GP} \) yields
\( \bp_{\dimm_{\GP}+1}^{4} \geq	1 \)
which reduces the second bound to the first one.
\end{proof}

Now we evaluate the bias term. 

\begin{lemma}
\label{LbiasGP}
With \( \Projc_{\GP} = \BB(\Id_{\dimp} + \BB)^{-1} \), it holds
\begin{EQA}
	\| \IF_{\GP}^{-1} \GP^{2} \upsvs \|^{2}
	& = &
	\| \IF^{-1/2} \Projc_{\GP} \, \IF^{1/2} \upsvs \|^{2} \, ,
	\\
	\| \IF^{1/2} \IF_{\GP}^{-1} \GP^{2} \upsvs \|^{2}
	& = &
	\| \Projc_{\GP} \, \IF^{1/2} \upsvs \|^{2} 
	\, .
\label{yytxstx6xt6wr5ww543}
\end{EQA}
If \( \IF \) and \( \BB \) or, equivalently, \( \IF \) and \( \GP^{2} \) commute then 
\begin{EQA}
	\| \IF_{\GP}^{-1} \GP^{2} \upsvs \|^{2}
	& = &
	\| \Projc_{\GP} \upsvs \|^{2}
	\, .
\label{8d8d8d8h3uduuyt444}
\end{EQA}
\end{lemma}
\begin{proof}
Statement \eqref{yytxstx6xt6wr5ww543} follows from the identity
\( \IF_{\GP}^{-1} \GP^{2} = \IF^{-1/2} \Projc_{\GP} \, \IF^{1/2} \).
\end{proof}

Now we summarize our finding in the case \( \QP = \IF^{1/2} \) and
\( \QP = \Id_{\dimp} \).

\begin{proposition}
\label{PriskGP}
Let \( \BB = \IF^{-1/2} \GP^{2} \IF^{-1/2} \) satisfy \nameref{GPTref}
and \( \IF \geq n \Id_{\dimp} \).
If \( \dimm_{\GP} \) is the largest \( m \) such that \( \bp_{\dimm}^{2} \leq 1 \) then
\begin{EQA}
	\riskt_{\QP}
	& \leq &
	\begin{cases}
		( 1 + \CONSTi_{\BB}) \dimm_{\GP} + \| \Projc_{\GP} \, \IF^{1/2} \upsvs \|^{2} \, , &  \QP = \IF^{1/2} \, , \\
		n^{-1} ( 1 + \CONSTi_{\BB}) \dimm_{\GP} + \| \IF^{-1/2} \Projc_{\GP} \, \IF^{1/2} \upsvs \|^{2} \, , &  \QP = \Id_{\dimp} \, .
	\end{cases}
\label{ydfjgv7ejed7urkh744u}
\end{EQA}
\end{proposition}

For illustration, let us consider an interesting case when \( \upsvs \) is \( \GPb \)-smooth for \( \GPb^{2} \leq \GP^{2} \), 
that is, \( \GPb \)-smoothness is less restrictive then \( \GP \)-smoothness.

\begin{lemma}
\label{LsmoothGPa}
Let \( \upsvs \) be \( \GPb \)-smooth for some \( \GPb \in \GPT \), that is, 
\( \| \GPb \upsvs \|^{2} \leq \CGPb \).
Let also \( \GP \), \( \IF \), and \( \GPb \) commute and hence, have the same eigenspaces, and 
\( (\gp_{1,j}^{2}) \) be the ordered eigenvalues of \( \GPb^{2} \).
Moreover, let the ratio \( \bp_{j}^{4}/\gp_{1,j}^{2} \) grow with \( j \).
Then
\begin{EQA}
	\| \Projc_{\dimm_{\GP}} \upsvs \|^{2}
	& \leq &
	\frac{1}{\gp_{1,\dimm_{\GP}}^{2}} \, \| \GPb \upsvs \|^{2} .
\label{ydfjfd8uuy655tehdf}
\end{EQA}
\end{lemma}

\begin{proof}
As \( \gp_{1,j}^{2} \) and \( \bp_{j}^{4}/\gp_{1,j}^{2} \) grow with \( j \) and \( \bp_{\dimm} \leq 1 \)
for \( \dimm \leq \dimm_{\GP} \), it holds 
\begin{EQA}
	\| \Projc_{\dimm} \upsv \|^{2}
	&=&
	\sum_{j=1}^{\dimp} \frac{\bp_{j}^{4}}{(1+\bp_{j}^{2})^{2}} \, \ups_{j}^{2}
	\leq 
	\sum_{j=1}^{\dimm} \bp_{j}^{4} \ups_{j}^{2} + \sum_{j=\dimm+1}^{\dimp} \ups_{j}^{2}
	\leq 
	\sum_{j=1}^{\dimm} \frac{\bp_{j}^{4}}{\gp_{1,j}^{2}} \gp_{1,j}^{2} \ups_{j}^{2} 
	+ \sum_{j=\dimm+1}^{\dimp} \frac{\gp_{1,j}^{2}}{\gp_{1,j}^{2}} \ups_{j}^{2}
	\\
	& \leq &
	\frac{\bp_{\dimm}^{4}}{\gp_{1,\dimm}^{2}} \sum_{j=1}^{\dimm} \gp_{1,j}^{2} \ups_{j}^{2}
	+ \frac{1}{\gp_{1,\dimm+1}^{2}} \sum_{j=\dimm+1}^{\dimp} \gp_{1,j}^{2} \ups_{j}^{2}
	\leq 
	\gp_{1,\dimm}^{-2} \| \GPb \upsv \|^{2} \, .
\label{8djfcfjcvjvcht5gr7yhe}
\end{EQA}
This implies the result.
\end{proof}

\medskip
We conclude that a roughness penalty \( \GP^{2} \) satisfying \nameref{GPTref} yields nearly the same risk as the projection estimator 
with a special \( \GP^{2} \)-dependent choice \( \II_{\dimm_{\GP}} \) of the corresponding sub-space.
This reduces the problem of risk minimization to the case of projection estimation considered earlier.

\Subsection{An example}

Consider a particular example when \( \{ \GP^{2} \} \) is a univariate family of penalizing matrices 
\( \GP^{2} \) of the form \( \GP^{2} = \CGPw \GPb^{2} \)
for \( \GPb^{2} = \diag\{ \gp_{1}^{2}, \ldots, \gp_{\dimp}^{2} \} \) 
with \( \gp_{j}^{2} = h(j) \) for a strictly increasing function \( h(\cdot) \).
Later we assume \( \IF = n \Id \).
Each value \( \CGPw \) identifies 
the spectral cut-off value \( \dimm_{\CGPw} \) which solves \( \CGPw \bp_{\dimm}^{2} \approx 1 \) or, equivalently, \( \CGPw \gp_{\dimm}^{2} \approx n \).
This leads to
\begin{EQA}
	\dimm_{\CGPw}
	& \approx &
	h^{-1}(n/\CGPw) .
\label{d7xcyhjkkf7643ddhrusnsdh}
\end{EQA}
Now we study the bias term beginning from the case when 
\( \upsvs \) is \( \GPb \)-smooth:
\( \| \GPb \upsvs \|^{2} \leq \CGPb \).
Then Lemma~\ref{LsmoothGPa} implies
\( \| \Projc_{\dimm_{\CGPw}} \upsvs \|^{2} \leq n^{-1} \| \GP \upsvs \|^{2} \leq n^{-1} \CGPw \| \GPb \upsvs \|^{2} \).
This yields the upper bound for the risk:
\begin{EQA}
	\riskt_{\CGPw}
	& \approx &
	\dimm_{\CGPw}/n + \| \Projc_{\dimm_{\CGPw}} \upsvs \|^{2}
	\leq 
	n^{-1} h^{-1}(n/\CGPw) + n^{-1} \CGPw \, \CGPb \, .
\label{hdycyuccuje47ctdge6}
\end{EQA}
The optimal/oracle choice \( \CGPs \) of \( \CGPw \) is obtained by minimization of this expression 
w.r.t. \( \CGPw \) leading to
\begin{EQA}
	\CGPs
	&=&
	\argmin_{\CGPw} \bigl\{ h^{-1}(n/\CGPw) + \CGPw \, \CGPb \bigr\} \, .
\label{usd8w6fjebcje6cjdebdj}
\end{EQA}
%
For instance, if \( \gp_{j}^{2} = h(j) = j^{2\smpa} \) then \( h^{-1}(j) = j^{1/(2\smpa)} \),
\( \dimm_{\CGPw} \approx (n/\CGPw)^{1/(2\smpa)} \),
and
\begin{EQA}
	\riskt_{\CGPw}
	& \leq &
	\dimm_{\CGPw}/n + \| \Projc_{\dimm_{\CGPw}} \upsvs \|^{2}
	\leq 
	n^{-1} \bigl\{ (n/\CGPw)^{1/(2\smpa)} + \CGPw \, \CGPb \bigr\} \, .
\label{hdycyuccuje47ctdge6}
\end{EQA}
This yields
\begin{EQA}
	\CGPs
	& \asymp &
	n^{1/(2\smpa+1)} \CGPb^{-2\smpa/(1 + 2\smpa)} ,
	\qquad
	\riskt^{*}
	\asymp
	n^{-1} \CGPs \, \CGPb
	\asymp
	n^{-1} (n \, \CGPb)^{1/(2\smpa+1)} \, .
\label{ucg4e6fhcvje48cjvgetdxbd}
\end{EQA}
The case when \( \upsvs \) is not \( \GPb \)-smooth is a bit more involved
because there is no minimax solution over a class of signals \( \upsvs \).
The \( \upsvs \)-dependent choice of \( \dimms \) follows \eqref{d7fjfjjvytv435y6ewee}
and \( \CGPs = n / {\dimms}^{2\smpa} \).

\Section{Inverse problems. Ridge penalty and spectral cut-off}
This section discusses risk bounds for inverse problem setup. 
It also addresses an important question of choosing a proper penalty
to obtain optimal estimation accuracy.
It appears that the answer is different depending on the relation between 
degree of ill-posedness and smoothness of the signal.

First we specify the inverse problem setup.
Its main feature is that the Fisher information operator 
\( \IF = \IF(\upsvs) = \nabla^{2} L(\upsvs) \) is not well-posed,
its conditional number is very large or infinite. 
This makes the choice of the penalization \( \pent_{\GP} \) especially important.
To be more specific, we focus on the estimation problem with \( \QP = \Id_{\dimp} \)
and on a quadratic penalty \( \pent_{\GP}(\upsv) = \| \GP \upsv \|^{2}/2 \) yielding \( \AvGP = \GP^{2} \upsvs \).

We assume that all the conditions of Theorem~\ref{TQFiWibias} are fulfilled.
Under the critical dimension condition, the leading term of the risk is given by \( \riskt \) from
\eqref{7djhed8cjfct534etgdhdy}.
For a smooth operator \( \IF \), a simple ridge penalization does a good job.
In the other cases, one should apply one or another model reduction technique.
We discuss the spectral cut-off method in combination with the approximation spaces setup. 

\Subsection{Ridge penalty and a smooth operator}
An important example of penalty choice is a ridge penalty
\( \GP^{2} = \gp^{2} \Id_{\dimp} \).
It is basis and coordinate free and enforces the ``benign overfitting'' phenomenon
in high-dimensional regression;
see \cite{baLoLu2020,ChMo2022,NoPuSp2024} and references therein.
Later we consider the estimation problem with \( \QP = \Id_{\dimp} \).
Also, we assume \( \Var(\nabla \zeta) \leq \IF \).
This condition can be relaxed in many ways.

Our study reveals an interesting phase transition effect. 
The use of ridge penalty leads to accurate results in the situation 
when the operator \( \IF \) is ``more regular'' than the signal \( \upsvs \).
In this case, ridge penalization enforces nearly the same effect as 
a spectral cut-off method.
Moreover, with a properly chosen ridge parameter \( \gp^{2} \), one can achieve 
the bias-variance trade-off in estimation of the target parameter.
The situation changes dramatically if the operator \( \IF \) is not smoother than the signal. 
This includes the case of a direct problem when the conditional number of the operator 
\( \IF \) is bounded by a fixed constant. 
It is well known that the ridge penalty is not efficient in this case, 
and a model reduction technique should be applied. 

In the rest of this section, we assume a smooth operator \( \IF \).
Our study includes the case with \( \dimp = \infty \).
Introduce the ordered eigenvalues \( \nEO_{1} \geq \nEO_{2} \geq \ldots \geq \nEO_{\dimp} \) of \( \IF \).
Also, define \( \IF_{\GP} = \IF + \gp^{2} \Id_{\dimp} \).
By Theorem~\ref{TQFiWibias},
the squared risk of \( \tilde{\upsv}_{\GP} \) can be approximated by \( \riskt \) with
\begin{EQA}[rcl]
	\riskt
	&=& 
	\tr \Var\bigl( \IF_{\GP}^{-1} \nabla \zeta \bigr)
	+ \| \IF_{\GP}^{-1} \, \GP^{2} \upsvs \|^{2} 
	\, .
\label{jcvyyvb7n99hj3jdfuycv}
\end{EQA}
Introduce an operator \( \BBH \eqdef \bigl( \Id_{\dimp} + \gp^{-2} \IF \bigr)^{-1} \).
It can be viewed as an approximation of the projector
in \( \R^{\dimp} \) on the subspace defined by the inequality
\( \IF \geq \gp^{2} \Id_{\dimp} \).
This subspace is spanned by the eigenvectors \( \ev_{j} \)
corresponding to \( \nEO_{j} \geq \gp^{2} \).
The dimension of this space is given by 
\begin{EQA}[c]
	\jJ_{\gp} \eqdef \max\{ j \colon \nEO_{j} \geq \gp^{2} \} .
\label{dufhiw3jeio26vte3hjb}
\end{EQA}
This relation can be inverted: 
given \( \jJ \), the corresponding ridge factor \( \gp^{2} \) can be given by 
\( \gp^{2} = \nEO_{\jJ} \).

\begin{lemma}
\label{Lvarb}
Let 
\( \nEO_{1} \geq \nEO_{2} \geq \ldots \geq \nEO_{\dimp} \) 
be the ordered eigenvalues of \( \IF \) while \( \ev_{j} \) be the corresponding eigenvectors.
Let also \( \Var(\nabla \zeta) \leq \IF \).
With \( \jJ = \jJ_{\gp} \) from \eqref{dufhiw3jeio26vte3hjb},
\begin{EQA}[rcl]
	\tr\Var (\IF_{\GP}^{-1} \nabla \zeta)
	& \leq &
	\biggl( \sum_{j=1}^{\jJ} \frac{1}{\nEO_{j}} + \gp^{-4} \sum_{j=\jJ+1}^{\dimp} \nEO_{j} \biggr)
	\, ,
	\\
	\gp^{4} \| \IF_{\GP}^{-1} \, \upsvs \|^{2}
	& \leq &
	\bigl\| \BBH \upsvs \bigr\|^{2}
	=
	\sum_{j=1}^{\dimp} \frac{1}{(\gp^{-2} \nEO_{j} + 1)^{2}} \langle \upsvs, \ev_{j} \rangle^{2}
	\, .
\end{EQA}
\end{lemma}
\begin{proof}
By \( \Var(\nabla \zeta) \leq \IF \)
\begin{EQA}[rcl]
	&& \nquad
	\tr \Var (\IF_{\GP}^{-1} \nabla \zeta)
	\leq 
	\tr \bigl\{ (\IF + \gp^{2} \Id_{\dimp})^{-1} \IF (\IF + \gp^{2} \Id_{\dimp})^{-1} \bigr\}
	\\
	&=&
	\sum_{j=1}^{\dimp} \frac{\nEO_{j}}{(\nEO_{j} + \gp^{2})^{2}}
	\leq 
	\biggl( \sum_{j=1}^{\jJ} \nEO_{j}^{-1} + \gp^{-4} \sum_{j=\jJ+1}^{\dimp} \nEO_{j} \biggr)
	\, ,
\end{EQA}
and the first result follows.
The second one is obvious.
\end{proof}

The results can be further simplified if the values \( \nEO_{j} \) decay polynomially:
\( \nEO_{j} \asymp \nEO_{1} \, j^{-2s} \) for \( s > 1/2 \).
Then for any \( \jJ \), we may use
\begin{EQA}[c]
	\sum_{j=1}^{\jJ} \frac{1}{\nEO_{j}} \leq \CONSTi_{1} \, \frac{\jJ}{\nEO_{\jJ}}  \, ,
	\qquad
	\sum_{j=\jJ+1}^{\dimp} \nEO_{j} \leq \CONSTi_{2} \, \jJ \nEO_{\jJ}  \, .
\label{udvcghue3j3vc654rw2yb}
\end{EQA}

\begin{proposition}
\label{Pvarb}
Let \( \GP^{2} = \gp^{2} \Id_{\dimp} \) and \( \Var(\nabla \zeta) \leq \IF \).
Assume \eqref{udvcghue3j3vc654rw2yb}, 
and \( \jJ = \jJ_{\gp} \); see \eqref{dufhiw3jeio26vte3hjb}.
Then
\begin{EQA}[rcl]
	\tr \Var (\IF_{\GP}^{-1} \nabla \zeta)
	& \leq &
	(\CONSTi_{1} + \CONSTi_{2}) \frac{\jJ}{\nEO_{\jJ}}
	\, .
	\qquad
\label{ckjivhiedi93wikvuyb}
\end{EQA}
Further, let \( \upsvs \) satisfy the smoothness condition
\begin{EQA}[c]
	\sum_{j=1}^{\dimp} \cgp_{j}^{2} \langle \upsvs, \ev_{j} \rangle^{2}
	\leq 
	1,
\label{dcyhvc55t8tg764e33rfb}
\end{EQA}
where \( \cgp_{j} \) grow while \( \cgp_{j} \nEO_{j} \) decrease with \( j \).
With \( \jJ = \jJ_{\gp} \) from \eqref{dufhiw3jeio26vte3hjb}, it holds
\begin{EQA}[rcl]
	\bigl\| \bigl( \gp^{-2} \IF + \Id_{\dimp} \bigr)^{-1} \upsvs \bigr\|^{2}
	& \leq &
	\max_{j \leq \dimp} \frac{\cgp_{j}^{-2}}{(\gp^{-2} \nEO_{j} + 1)^{2}}
	\leq 
	\frac{1}{\cgp_{\jJ}^{2}}
	\, 
\label{ckjivhiedi93wikvuy2b}
\end{EQA}
yielding
\begin{EQA}
	\riskt
	& \lesssim &
	\frac{\jJ}{\nEO_{\jJ}} + \frac{1}{\cgp_{\jJ}^{2}} 
	\, .
\label{c8f78vtv4c4tehfivefb}
\end{EQA}
\end{proposition}

\begin{proof}
Statement \eqref{ckjivhiedi93wikvuyb} follows Lemma~\ref{Lvarb} and \eqref{udvcghue3j3vc654rw2yb}.
Further, 
\begin{EQA}[rcl]
	\| \BBH \upsvs \|^{2}
	&=&
	\sum_{j=1}^{\dimp} \frac{1}{(\gp^{-2} \nEO_{j} + 1)^{2}} \langle \upsvs, \ev_{j} \rangle^{2}
	=
	\sum_{j=1}^{\dimp} \frac{\cgp_{j}^{-2}}{(\gp^{-2} \nEO_{j} + 1)^{2}} \, \cgp_{j}^{2} \langle \upsvs, \ev_{j} \rangle^{2}
	\\
	& \leq &
	\max_{j \leq \dimp} \frac{\cgp_{j}^{-2}}{(\gp^{-2} \nEO_{j} + 1)^{2}}
	\sum_{j=1}^{\dimp} \cgp_{j}^{2} \langle \upsvs, \ev_{j} \rangle^{2}
	\leq 
	\biggl\{ \min_{j \leq \dimp} (\gp^{-2} \nEO_{j} \cgp_{j} + \cgp_{j}) \biggr\}^{-2} \, 
	\, .
\end{EQA}
As the values \( \nEO_{j} \cgp_{j} \) decrease and \( \cgp_{j} \) increase with \( j \),
it holds for any \( \jJ \)
\begin{EQA}[c]
	\min_{j \leq \jJ} (\gp^{-2} \nEO_{j} \cgp_{j} + \cgp_{j})
	\geq 
	\gp^{-2} \nEO_{\jJ} \cgp_{\jJ} \, ,
	\quad
	\min_{j > \jJ} (\gp^{-2} \nEO_{j} \cgp_{j} + \cgp_{j})
	\geq 
	\cgp_{\jJ+1}
	\, .
	\qquad
\end{EQA}
Therefore, the choice of \( \jJ \) by \( \gp^{-2} \nEO_{\jJ} \approx 1 \)
yields
\begin{EQA}[c]
	\min_{j \leq \dimp} (\gp^{-2} \nEO_{j} \cgp_{j} + \cgp_{j})
	\geq 
	\cgp_{\jJ}
	\, ,
\end{EQA}
and \eqref{ckjivhiedi93wikvuy2b} follows as well.
\end{proof}

Bound \eqref{ckjivhiedi93wikvuy2b} requires that the values \( \cgp_{j} \) 
grow while \( \nEO_{j} \cgp_{j} \) decrease.
The latter means that the operator \( \IF \) is more regular than the signal \( \upsvs \).
Now, consider the opposite case.
For instance, let the eigenvalues \( \nEO_{j} \) of \( \IF \) decrease slowly with 
\( j \) or remain significantly positive. 
Then inversion of \eqref{dufhiw3jeio26vte3hjb} for small \( \gp^{2} \) can be 
problematic. 
The corresponding cut-off index \( \jJ \) will be too large and cannot be used 
to mimic the bias-variance trade-off.
As a consequence, a ridge penalization is not efficient in such situations.
The following section explains how this situation can be handled using the spectral cut-off method.


\Subsection{Spectral cut-off}
Spectral cut-off is a standard tool in model reduction for inverse problems. 
It assumes that the eigenvector decomposition of \( \IF \) 
with ordered eigenvalues \( \nEO_{1} \geq \nEO_{2} \geq \ldots \geq \nEO_{\jJ} \) and the corresponding eigenvectors
\( \ev_{j} \) are available.
By \( \Proj_{\jJ} \) we denote the canonical projector in the source space \( \R^{\dimp} \) onto the first \( \jJ \) coordinates \( \ev_{j} \).
Also, we assume that smoothness properties of the unknown source signal \( \upsvs \) can be given in the canonical basis 
formed by these eigenvectors \( \ev_{j} \); see \eqref{dcyhvc55t8tg764e33rfb}.
For a fixed cut-off parameter \( \jJ \), consider the truncation penalty \( \GP_{\jJ}^{2} = \diag\{ \GP_{1}^{2},\dots,\GP_{\dimp}^{2} \} \)
with \( \gp_{j}^{2} = 0 \) for \( j \leq \jJ \) and \( \gp_{j}^{2} = \infty \) for \( j > \jJ \).
Effectively, this penalty enforces \( \tilde{\ups}_{\GP,j} = \upss_{\GP,j} = 0 \) for \( j > \jJ \).

\begin{proposition}
\label{LvarEOspect}
Let \( \nEO_{1} \geq \nEO_{2} \geq \ldots \geq \nEO_{\dimp} \) 
be the ordered eigenvalues of \( \IF \) and \( \Var(\nabla \zeta) \leq \IF \).
Let also \( \GP_{\jJ}^{2} \) be the spectral cut-off penalty at \( \jJ \). 
Then
\begin{EQA}[rcccl]
	\tr\Var\bigl\{ (\IF + \GP^{2})^{-1} \nabla \zeta \bigr\}
	& \leq &
	\sum_{j=1}^{\jJ} \frac{1}{\nEO_{j}} 
	\, ,
	\qquad
	\| \IF_{\GP}^{-1} \, \GP^{2} \upsvs \|
	& = &
	\| (\Id_{\dimp} - \Proj_{\jJ}) \upsvs \|
	\, 
\end{EQA}
yielding
\begin{EQA}[c]
	\riskt
	\leq 
	\sum_{j=1}^{\jJ} \frac{1}{\nEO_{j}} + \| (\Id_{\dimp} - \Proj_{\jJ}) \upsvs \|^{2}
	\, .
\end{EQA}
Under \eqref{udvcghue3j3vc654rw2yb} and \eqref{dcyhvc55t8tg764e33rfb}, bound \eqref{c8f78vtv4c4tehfivefb} applies. 
\end{proposition}

In the contrary to ridge penalization, 
the penalty coefficients \( \gp_{j}^{2} \) vanish for \( j \leq \jJ \).
This improves the bound on the penalization bias, the assumption of decay of \( \cgp_{j} \nEO_{j} \) is not required.

A proper choice of the cut-off parameter \( \jJ \) ensures a nearly optimal accuracy of estimation. 
This approach is widely used in linear inverse problems.
However, availability of the SVD for \( \IF \) is a severe limitation.
This issue can be resolved using the \emph{approximation space} approach.

\Subsection{Approximation spaces and truncation penalties}
This setup assumes specific basis in the parameter space \( \R^{\dimp} \)
suitable for describing the smoothness properties of \( \upsvs \) and regularity of \( \IF \) simultaneously.
Without loss of generality, we apply the canonical basis in the source space \( \R^{\dimp} \).
By \( \Proj_{\jJ} \) we denote the canonical projector in the source space \( \R^{\dimp} \) onto the first \( \jJ \) coordinates.

Consider a truncation penalty \( \GP^{2} = \diag\{ \GP_{1}^{2},\dots,\GP_{\dimp}^{2} \} \) with 
 \( \GP_{j}^{2} = 0 \) for \( j \leq \jJ \) and
\( \GP_{j}^{2} = \infty \) for \( j > \jJ \).
Effectively, this penalty enforces \( \tilde{\ups}_{\GP,j} = \upss_{\GP,j} = 0 \) for \( j > \jJ \).

\Subsection{Estimation risk}
A \( \jJ \)-truncation penalty reduces the original problem to the non-penalized MLE 
\begin{EQA}[rcccl]
	\tilde{\upsv}_{\jJ}
	&=&
	\argmax_{\upsv_{\jJ}} \LL_{\jJ}(\upsv_{\jJ})
	\, ,
	\qquad
	\upsvs_{\jJ}
	&=&
	\argmax_{\upsv_{\jJ}} \E \LL_{\jJ}(\upsv_{\jJ})
	\, .
\label{8sdi2ineuiujqdydtw3hychj}
\end{EQA}
The results of Theorem~\ref{TQFiWibias} apply to this setup.
Moreover, many terms entering the formulation of this theorem can be specified in more detail.
First, we discuss the stochastic component. 
Denote 
\begin{EQA}[c]
	\IF_{\jJ} \eqdef \Proj_{\jJ} \IF \Proj_{\jJ} 
	\, . 
\end{EQA}
By \( \Var(\nabla \zeta) \leq \IF \)
\begin{EQA}[rcl]
	\Var (\IF_{\GP}^{-1} \nabla \zeta)
	&=&
	\Var (\IF_{\jJ}^{-1} \nabla \zeta_{\jJ})
	\leq 
	\IF_{\jJ}^{-1} 
	\, .
\label{du7dyfu3ijhrktruncb}
\end{EQA}
Obviously, \( \tr \IF_{\jJ}^{-1} \leq \tr \IF^{-1} \).
Therefore, truncation reduces the effective dimension and improves the stochastic component in the risk. 
It remains to evaluate the bias \( \upsvs_{\jJ} - \upsvs \) caused by this truncation
under the regularity conditions on the operator \( \IF \) and the parameter \( \upsvs \).
%
Regularity of \( \IF \) will be described using the quantities
\begin{EQ}[rcl]
	\nEO_{j}
	& \eqdef &
	\lambda_{j}(\IF_{j})
	\, .
\label{vidf8vccil3l3wf6rbvyb}
\end{EQ}
Here 
\( \lambda_{j}(\BBH) \) means the \( j \)th largest eigenvalue of the matrix \( \BBH \) and
\( \IF_{j} = \Proj_{j} \IF \Proj_{j} \); see \eqref{du7dyfu3ijhrktruncb}.
When considering \( \BBH = \IF_{j} \) as a \( j \)-dimensional 
matrix, the value \( \nEO_{j} \) corresponds to its smallest eigenvalue. 
If the basis vectors \( \ev_{j} \) in \( \R^{\dimp} \) are defined as the ordered eigenvectors of \( \IF \) then
\( \nEO_{j} \) coincide with the \( j \)th eigenvalue of \( \IF \).
For mildly/severely ill-posed problems, these values rapidly decrease with \( j \).

\begin{proposition}
\label{Pappspaceb}
Assume \eqref{vidf8vccil3l3wf6rbvyb}.
Then for any \( \jJ \leq \dimp \) 
\begin{EQA}[rcccl]
\label{fubynbunmn5323hgftwEOa}
	\| \IF_{\jJ}^{-1} \, \GP_{\jJ}^{2} \upsvs \|
	& = &
	\| (\Id_{\dimp} - \Proj_{\jJ}) \upsvs \| 
	\, ,
	\qquad
	\tr(\IF_{\jJ}^{-1})
	& \leq &
	\sum_{j=1}^{\jJ} \frac{1}{\nEO_{j}} 
	\, . 
\end{EQA}
\end{proposition}

\Subsection{Rate of estimation in inverse problems}

This section illustrates the obtained results by showing that 
a ridge penalty for a smooth operator or a truncation penalty with properly
selected parameter \( \jJ \) lead to 
rate optimal accuracy of estimation.
As in the previous sections, we fix an approximation spaces setup.
Regularity of the operator \( \IF \) is described by decreasing sequences 
\( \nEO_{1} \geq \nEO_{2} \geq \ldots \geq \nEO_{\dimp} \) ensuring \eqref{vidf8vccil3l3wf6rbvyb}.
Moreover, we only discuss the case of a mildly ill-posed problem when the values \( \nEO_{j} \) decrease
polynomially; see \eqref{udvcghue3j3vc654rw2yb}.
We also assume a Sobolev smoothness of the signal \( \upsv \) as in \eqref{dcyhvc55t8tg764e33rfb}
with \( \cgp_{j} \) polynomially increasing.
A popular special case is given by \( \nEO_{j} \approx \nEO_{1} \, j^{-2s} \) 
and \( \cgp_{j}^{2} \approx \CGPz \, j^{2\beta} \).
The ridge penalty requires \( s \geq \beta \).
Truncation penalty enables us to relax this condition. 

\begin{proposition}
\label{PrateEiO}
Let the operator \( \IF \) fulfill \eqref{vidf8vccil3l3wf6rbvyb} with 
\( \nEO_{j} \geq \nEO_{1} \, j^{-2s} \) for all \( j > 1 \).
Further, let \( \upsvs \) follow \eqref{dcyhvc55t8tg764e33rfb} with \( \cgp_{j}^{2} = \CGPz \, j^{2\beta} \).
Define \( \jJ \) by \( \jJ \asymp (\nEO_{1}/\CGPz)^{1/(1 + 2\beta + 2s)} \).
Then the risk \( \riskt_{\jJ} \) of the estimator \( \tilde{\upsv}_{\jJ} \) satisfies
\begin{EQA}[c]
	\riskt_{\jJ}
	\lesssim 
	\CGPz^{-\frac{2s+1}{1 + 2\beta + 2s}} \nEO_{1}^{-\frac{2\beta}{1 + 2\beta + 2s}}
	\, .
\end{EQA}
\end{proposition}
\begin{proof}
Define \( \jJ \) by \( \jJ \asymp (\nEO_{1}/\CGPz)^{1/(1 + 2\beta + 2s)} \).
\eqref{vidf8vccil3l3wf6rbvyb} implies
\begin{EQA}[c]
	\tr(\DPN_{\jJ}^{-2})
	\leq 
	\sum_{j=1}^{\jJ} \frac{1}{\nEO_{j}} 
	\leq 
	\frac{1}{\nEO_{1}} \sum_{j=1}^{\jJ} j^{2s}
	\leq 
	\frac{1}{2s+1} \frac{\jJ^{2s+1}}{\nEO_{1}}
	\, .
\end{EQA}
Further, as \( \cgp_{j}^{2} = \CGPz \, j^{2\beta} \) increase with \( j \), it holds
\begin{EQA}[c]
	\| (\Id_{\dimp} - \Proj_{\jJ}) \upsvs \|^{2}
	=
	\sum_{j=\jJ+1}^{\dimp} \langle \upsvs, \ev_{j} \rangle^{2}
	\leq 
	\cgp_{\jJ}^{-2}
	\sum_{j=\jJ+1}^{\dimp} \cgp_{j}^{2} \langle \upsvs, \ev_{j} \rangle^{2}
	\leq 
	\cgp_{\jJ}^{-2}
	= 
	\CGPz^{-1} \, \jJ^{-2\beta}
	\, .
\end{EQA}
This yields
\begin{EQA}[c]
	\riskt
	\lesssim  
	\tr(\DPN_{\jJ}^{-2}) + \cgp_{\jJ}^{-2}
	\leq 
	\nEO_{1}^{-1} \jJ^{2s+1} + \CGPz^{-1} \, \jJ^{-2\beta}
	\lesssim
	\CGPz^{-\frac{2s+1}{1 + 2\beta + 2s}} \nEO_{1}^{-\frac{2\beta}{1 + 2\beta + 2s}}
\end{EQA}
as stated.
\end{proof}

\Chapter{Examples of parametric models}
\label{SexamplesSLS}
This \chname illustrates the general notions on the particular examples including 
logistic regression, log-density estimation, and precision matrix estimation.
We mainly check the general conditions.
This enables us to apply the results of \Chname \ref{SgenBounds}.

\Section{Log-density estimation}
\label{SGBvM}
Suppose we are given a random sample \( X_{1},\ldots,X_{\nsize} \) in \( \R^{d} \).
The density model assumes that all these random variables are independent 
identically distributed from some measure \( \Pone \) with a density 
\( \dens(\xv) \) with respect to a \( \sigma \)-finite measure \( \Pdom \) in \( \R^{d} \).
This density function is the target of estimation.
By definition, the function \( \dens \) is non-negative, measurable, and integrates to one:
\( \int \dens(\xv) \, d\Pdom(\xv) = 1 \).
%
Here and in what follows, the integral \( \int \) without limits means the integral over the whole space \( \R^{d} \).
If \( \dens(\cdot) \) has a smaller support \( \XX \), one can restrict integration to this set.
%
Below we parametrize the model by a linear decomposition of the log-density function.
Let \( \bigl\{ \psi_{j}(\xv), \, j=1,\ldots,\dimp \bigr\} \) with \( \dimp \leq \infty \)
be a collection of functions in \( \R^{d} \) (a dictionary).
%
For each \( \upsv = (\ups_{j}) \in \R^{\dimp} \), define   
\begin{EQA}
    \ldens(\xv,\upsv)
    &\eqdef &
    \ups_{1} \psi_{1}(\xv) + \ldots + \ups_{\dimp} \psi_{\dimp}(\xv) - \cdens(\upsv)
    =
    \bigl\langle \Psiv(\xv), \upsv \bigr\rangle - \cdens(\upsv),
\label{logdenssumj}
\end{EQA}
where \( \Psiv(\xv) \) is a vector with components \( \psi_{j}(\xv) \).
Let \( \cdens(\upsv) \) be given by 
\begin{EQA}
    \cdens(\upsv)
    & \eqdef &
    \log \int \ex^{\langle \Psiv(\xv), \upsv \rangle} \, d\Pdom(\xv) .
\label{gtdelointTPxm0}
\end{EQA}
It is worth stressing that the data point \( \xv \) only enters in the linear term \( \bigl\langle \Psiv(\xv), \upsv \bigr\rangle \)
of the log-likelihood \( \ldens(\xv,\upsv) \).
The function \( \cdens(\upsv) \) is entirely model-driven.
Below we restrict \( \upsv \) to a subset \( \Ups \) in \( \R^{\dimp} \) such that 
\( \cdens(\upsv) \) is well defined and the integral 
\( \int \ex^{\langle \Psiv(\xv), \upsv \rangle} \, d\Pdom(\xv) \) is finite.
Linear log-density modeling assumes  
\begin{EQA}
  	\log \dens(\xv)
  	&=&
  	\ldens(\xv,\upsvs)
	=
	\bigl\langle \Psiv(\xv), \upsvs \bigr\rangle - \cdens(\upsvs)
\label{dnesdxtsl}
\end{EQA}  
for some \( \upsvs \in \Ups \).
A nice feature of such representation is that the function \( \log \dens(\xv) \), in contrary to the density itself, does not need to be non-negative.
Another important benefit of using the log-density is that the stochastic part of the corresponding log-likelihood is \emph{linear} w.r.t. the parameter \( \upsv \).
With \( \Spsi = \sumi \Psiv(X_{i}) \), the negative log-likelihood \( L(\upsv) \) reads as
\begin{EQA}
    L(\upsv) 
    &=& 
    - \sumi \bigl\langle \Psiv(X_{i}), \upsv \bigr\rangle + \nsize \cdens(\upsv) 
    = 
    - \langle \Spsi, \upsv \rangle + \nsize \cdens(\upsv) .
\label{density_likelihood}
\end{EQA}
The truth can be defined as its population counterpart:
\begin{EQA}
    \upsvs
    & = &
    \argmin_{\upsv \in \Ups} \E L(\upsv)
    =
    \argmin_{\upsv \in \Ups} 
      	\bigl\{ - \langle \E \Spsi, \upsv \rangle + \nsize \cdens(\upsv) \bigr\} 
	=
    \argmin_{\upsv \in \Ups} 
      	\bigl\{ - \langle \Psimean, \upsv \rangle + \cdens(\upsv) \bigr\} ,
	\qquad
\label{ttsT1dedef}
\end{EQA}
where \( \Psimean = \nsize^{-1} \E \Spsi \).
This yields the identity
\begin{EQA}
	\nabla \cdens(\upsvs)
	&=&
	\Psimean .
\label{cgsuwywdwdy2whq21ef2w}
\end{EQA}
For a given penalty operator \( \GP^{2} \),
the penalized loss \( \LGP(\upsv) \) reads as
\begin{EQA}
    \LGP(\upsv) 
    &=& 
    L(\upsv) + \frac{1}{2} \| \GP \upsv \|^{2}
    = 
    - \langle \Spsi, \upsv \rangle + \nsize \cdens(\upsv) 
    + \frac{1}{2} \| \GP \upsv \|^{2} .
\label{density_likelihoodG}
\end{EQA}
The penalized MLE 
\( \tilde{\upsv}_{\GP} \) and its population counterpart \( \upsvs_{\GP} \) are defined as
\begin{EQ}[rcccl]
    \tilde{\upsv}_{\GP}
    &=&
    \argmin_{\upsv \in \Ups} L_{\GP}(\upsv) ,
    \qquad
    \upsvs_{\GP}
    &=&
    \argmin_{\upsv \in \Ups} \E L_{\GP}(\upsv) .
\label{tsatTqLGttT122}
\end{EQ}
We are interested in sufficient conditions on the model which enables us to apply the general results of Section~\ref{SgenBounds}
for quantifying the error terms \( \tilde{\upsv}_{\GP} - \upsvs_{\GP} \), \( \upsvs_{\GP} - \upsvs \), and the corresponding risk
\( \E \| \QP (\tilde{\upsv}_{\GP} - \upsvs) \|^{2} \).

\subsection*{Assumptions}
First note that the generalized linear structure of the model automatically yields 
conditions \nameref{LLref} and \nameref{Eref}.
Indeed, convexity of \( \cdens(\cdot) \) implies that 
\( \E L(\upsv) = - \langle \E \Spsi, \upsv \rangle + \nsize \cdens(\upsv) \) is convex. 
Further, for the stochastic component \( \zeta(\upsv) = L(\upsv) - \E L(\upsv) \), it holds 
\begin{EQA}
\label{stochDensDef}
    \nabla \zeta(\upsv)
    &=&
    \nabla \zeta
    =
    \Spsi - \E \Spsi
    = 
    \sumi \bigl[ \Psiv(X_{i}) - \E \, \Psiv(X_{i}) \bigr],
\label{ztLtEltnztSPXi}
\end{EQA}
and \nameref{Eref} follows.
Further, the representation 
\( \E L(\upsv) = - \langle \E \Spsi, \upsv \rangle + \nsize \cdens(\upsv) \) implies
\begin{EQA}
	\IF(\upsv)
	&=&
	\nabla^{2} \E L(\upsv)
	=
	\nabla^{2} L(\upsv)
	=
	\nsize \nabla^{2} \cdens(\upsv) .
\label{IFttn2ELtnn2cdt}
\end{EQA}
To simplify our presentation, we assume that \( X_{1},\ldots,X_{\nsize} \) are indeed i.i.d.
and the density \( \dens(\xv) \) can be represented in the form \eqref{dnesdxtsl} 
for some parameter vector \( \upsvs \).
This can be easily extended to a non i.i.d. case at the cost of more complicated notations.
Then \( \Psimean = \nsize^{-1} \E \Spsi = \E \, \Psiv(X_{1}) \).
Moreover, by \eqref{gtdelointTPxm0},
\( \nabla^{2} \cdens(\upsvs) = \Var \bigl\{ \Psiv(X_{1}) \bigr\} \) and
\begin{EQA}
	\Var(\nabla \zeta)
	&=&
	\nsize \, \nabla^{2} \cdens(\upsvs) 
	=
	\IF(\upsvs) .
\label{3094857ythgfjdkleo}
\end{EQA}
For any \( \upsv \in \Ups \) and \( \rsmall > 0 \), 
consider the elliptic set \( \B_{\rsmall}(\upsv) \subset \R^{\dimp} \) with
\begin{EQA}
	\B_{\rsmall}(\upsv)
	& \eqdef &
	\bigl\{ \uv \in \R^{\dimp} \colon \langle \nabla^{2} \cdens(\upsv), \uv^{\otimes 2} \rangle \leq \rsmall^{2} \bigr\} .
\label{f765re276tf7we6t8erw63}
\end{EQA}

Assume the following conditions.

\begin{description}
\item[\label{upsvsref} \( \bb{(\fs)} \)]
\( X_{1},\ldots,X_{\nsize} \) are i.i.d. from a density \( \fs \) satisfying
\( \log \fs(\xv) = \Psiv(\xv)^{\T} \upsvs - \cdens(\upsvs) \).

\item[\label{Upsref} \( \bb{(\Ups)} \)]
The set \( \Ups \) is open and convex, 
the value \( \cdens(\upsv) \) from \eqref{gtdelointTPxm0} is finite for all \( \upsv \in \Ups \),
\( \upsvs \) from \eqref{ttsT1dedef} 
is an internal point in \( \Ups \) such that 
\( \B_{2\rsmall}(\upsvs) \subset \Ups \) for a fixed \( \rsmall > 0 \).

\item[\label{cdensref} \( \bb{(\cdens)} \)]
For \( \upsv \in \B_{\rsmall}(\upsvs) \) and all \( \uv \) with
\( \langle \nabla^{2} \cdens(\upsv), \uv^{\otimes 2} \rangle \leq 4\rsmall^{2} \),
it holds
\begin{EQA}
	\exp \{ \cdens(\upsv + \uv) - \cdens(\upsv) - \langle \nabla \cdens(\upsv),\uv \rangle \}
	& \leq &
	\CONSTi_{\rsmall} 
	\,  .
\label{vgtcgdcbwcte32vfurdf}
\end{EQA}
\end{description}

Introduce a measure \( P_{\upsv} \) by the relation:
\begin{EQA}
    \frac{d P_{\upsv}}{d\Pdom}(\xv)
    &=&
    \exp\bigl\{ \bigl\langle \Psiv(\xv), \upsv \bigr\rangle - \cdens(\upsv) 
    \bigr\} .
\label{dPdlxeLtSn}
\end{EQA}
Identity \eqref{gtdelointTPxm0} ensures that \( P_{\upsv} \) is a probabilistic measure.
Moreover, under \eqref{dnesdxtsl}, the data generating measure \( \P \) coincides with \( P_{\upsvs}^{\otimes \nsize} \).

\begin{description}
\item[\label{PsiX4ref} \( \bb{(\Psiv_{4})} \)]
	\emph{There are \( \dltwaa_{\Psi,3} \geq 0 \) and \( \dltwaa_{\Psi,4} \geq 3 \) such that 
	for all \( \upsv \in \B_{\rsmall}(\upsvs) \) and \( \zv \in \R^{\dimp} \)
\begin{EQA}
	\bigl| E_{\upsv} \bigl\langle \Psiv(X_{1}) - E_{\upsv} \Psiv(X_{1}),\zv \bigr\rangle^{3} \bigr|
	& \leq &
	\dltwaa_{\Psi,3} \, E_{\upsv}^{3/2} \bigl\langle \Psiv(X_{1}) - E_{\upsv} \Psiv(X_{1}),\zv \bigr\rangle^{2} ,
	\\
	E_{\upsv} \bigl\langle \Psiv(X_{1}) - E_{\upsv} \Psiv(X_{1}),\zv \bigr\rangle^{4}
	& \leq &
	\dltwaa_{\Psi,4} \,  E_{\upsv}^{2} \bigl\langle \Psiv(X_{1}) - E_{\upsv} \Psiv(X_{1}),\zv \bigr\rangle^{2} .
\label{dcief8yfe7eegywddwfeyg}
\end{EQA}
	}
\end{description}

In fact, conditions \nameref{cdensref} and \nameref{PsiX4ref} follow from \nameref{Upsref} and can be considered as a kind of definition 
of important quantities \( \CONSTi_{\rsmall} \), \( \dltwaa_{\Psi,3} \), and \( \dltwaa_{\Psi,4} \) 
which will be used for describing the smoothness properties of \( \cdens(\upsv) \).


\begin{lemma}
\label{LS3dens}
Assume \nameref{upsvsref}, \nameref{Upsref}, \nameref{cdensref}, and \nameref{PsiX4ref}, and 
let \( \rr \leq \rsmall \, \sqrt{\nsize/2} \).
Then, for any \( \upsv \in \B_{\rsmall}(\upsvs) \), 
the function \( f(\upsv) = \E_{\upsvs} L(\upsv) \) satisfies \nameref{LLtS3ref} and \nameref{LLtS4ref}
with \( \hL(\upsv) = \cdens(\upsv) - \langle \nabla \cdens(\upsvs), \upsv \rangle \), 
\( \HL^{2}(\upsv) = \nabla^{2} \cdens(\upsv) \),
and constants \( \hmax_{3} \) and \( \hmax_{4} \) satisfying 
\begin{EQA}
	\hmax_{3}
	&=&
	\dltwaa_{\Psi,3} \, (\dltwaa_{\Psi,4} \, \CONSTi_{\rsmall})^{3/4} \, ,
	\qquad
	\hmax_{4}
	=
	(\dltwaa_{\Psi,4} - 3) \, \dltwaa_{\Psi,4} \, \CONSTi_{\rsmall} \, ,
	\, .
\label{7dgs5rrtgvsdhsdrsdvg}
\end{EQA}
\end{lemma} 

\begin{proof}
Fix \( \upsv \in \B_{\rsmall}(\upsvs) \).
With \( P_{\upsv} \) defined by \eqref{dPdlxeLtSn},
it holds \( E_{\upsv} \Psiv(X_{1}) = \nabla \cdens(\upsv) \) and 
\( \Var_{\upsv}(\Psiv(X_{1})) = \nabla^{2} \cdens(\upsv) \).
Further, if \( \uv \in \B_{2\rsmall}(\upsv) \) and  
\( \upsv + \uv \in \Ups \), then 
\begin{EQA}
	\cdens(\upsv + \uv) 
	&=& 
	\log E_{0} \exp\{ \langle \Psiv(X_{1}),\upsv + \uv \rangle \} 
	=
	\log E_{\upsv} \exp\bigl\{ \bigl\langle \Psiv(X_{1}), \uv \bigr\rangle + \cdens(\upsv) \bigr\} .
\label{dfor4u7tgjhkgtirt}
\end{EQA}
Define \( \epsv = \Psiv(X_{1}) - E_{\upsv} \Psiv(X_{1}) \).
By \( E_{\upsv} \Psiv(X_{1}) = \nabla \cdens(\upsv) \) 
\begin{EQA}
	\log E_{\upsv} \exp ( \langle \epsv,\uv \rangle )
	&=&
	\cdens(\upsv + \uv) - \cdens(\upsv) - \langle E_{\upsv} \Psiv(X_{1}),\uv \rangle
	\\
	&=&	
	\cdens(\upsv + \uv) - \cdens(\upsv) - \langle \nabla \cdens(\upsv),\uv \rangle .
\label{kjhr8dfrikweyyweyerdd}
\end{EQA}
By \nameref{cdensref}
\begin{EQA}
	\cdens(\upsv;\uv)
	& = &
	\cdens(\upsv + \uv) - \cdens(\upsv) - \langle \nabla \cdens(\upsv),\uv \rangle
	\leq 
	\frac{\const_{\CONSTi_{\rsmall}}}{4} \| \HL(\upsv) \uv \|^{2} .
\label{ucfnmdcyeytrdf6edye4t6}
\end{EQA}
Define for \( |t| \leq 1 \) 
\begin{EQA}
	\chi(t)
	& \eqdef &
	\log E_{\upsv} \exp ( t \langle \epsv,\uv \rangle )
	=
	\cdens(\upsv + t \uv) - \cdens(\upsv) - \langle \nabla \cdens(\upsv), t \uv \rangle ,
\label{fg8yw8wf78wq87weqf87we8f}
\end{EQA}
By \nameref{PsiX4ref} with \( \dltwaa_{\Psi,4} \geq 3 \)
\begin{EQA}
	\bigl| \chi^{(3)}(0) \bigr|
	&=&
	\bigl| E_{\upsv} \langle \epsv,\uv \rangle^{3} \bigr|
	\leq 
	\dltwaa_{\Psi,3} \, E_{\upsv}^{3/2} \langle \epsv,\uv \rangle^{2} \, ,
	\\
	\bigl| \chi^{(4)}(0) \bigr|
	&=&
	\bigl| E_{\upsv} \langle \epsv,\uv \rangle^{4} 
	- 3 E_{\upsv}^{2} \langle \epsv,\uv \rangle^{2} \bigr| 
	\leq 
	(\dltwaa_{\Psi,4} - 3) E_{\upsv}^{2} \langle \epsv,\uv \rangle^{2} .
\label{hdudyue2g2w5dgeh3etyd}
\end{EQA}
This implies 
\begin{EQA}
\label{0det53543wyf7rryem}
	| \langle \nabla^{3}\cdens(\upsv), \uv^{\otimes 3} \rangle |
	& \leq &
	\dltwaa_{\Psi,3} \, \langle \nabla^{2} \cdens(\upsv), \uv^{\otimes 2} \rangle^{3/2} \, ,
	\\
	| \langle \nabla^{4}\cdens(\upsv), \uv^{\otimes 4} \rangle |
	& \leq &
	(\dltwaa_{\Psi,4} - 3) \, \langle \nabla^{2} \cdens(\upsv), \uv^{\otimes 2} \rangle^{2} .
\label{0det53543wyf7rryedm}
\end{EQA}
Further, by \eqref{kjhr8dfrikweyyweyerdd} 
\begin{EQA}
	\nabla^{2} \cdens(\upsv + \uvdd)
	&=&
	\nabla^{2} \log E_{\upsv} \ex^{\langle \epsv,\uvdd \rangle} 
	=
	\frac{E_{\upsv} \{ \epsv \epsv^{\T} \ex^{\langle \epsv,\uvdd \rangle} \}}
		 {(E_{\upsv} \, \ex^{\langle \epsv,\uvdd \rangle})^{2}} 
	- \frac{E_{\upsv} \{ \epsv \, \ex^{\langle \epsv,\uvdd \rangle} \} \, 
			E_{\upsv} \{ \epsv \, \ex^{\langle \epsv,\uvdd \rangle} \}^{\T} }
		 {( E_{\upsv} \, \ex^{\langle \epsv,\uvdd \rangle} )^{2}}
\label{cskjoihnbve4r5ujiko}
\end{EQA}
and it follows by \eqref{0det53543wyf7rryedm} and \eqref{vgtcgdcbwcte32vfurdf} in view of 
\( E_{\upsv} \, \ex^{\langle \epsv,\uvdd \rangle} \geq 1 \)
for any \( \zv \in \R^{\dimp} \)
\begin{EQA}
	\bigl\langle \nabla^{2} \cdens(\upsv + \uvdd), \zv^{\otimes 2} \bigr\rangle
	& \leq &
	E_{\upsv} \bigl\{ \langle \epsv,\zv \rangle^{2} \ex^{\langle \epsv,\uvdd \rangle} \bigr\}
	\\
	& \leq &
	E_{\upsv}^{1/2} \langle \epsv,\zv \rangle^{4} \, \, E_{\upsv}^{1/2} \ex^{2 \langle \epsv,\uvdd \rangle} 
	\leq 
	\sqrt{\dltwaa_{\Psi,4} \, \CONSTi_{\rsmall}} \,\, \bigl\langle \nabla^{2} \cdens(\upsv), \zv^{\otimes 2} \bigr\rangle
\label{dpe034985tuykhlgpd}
\end{EQA}
This implies 
\begin{EQA}
	\frac{\langle \nabla^{2} \cdens(\upsv + \uvdd), \zv^{\otimes 2} \rangle}
		 {\langle \nabla^{2} \cdens(\upsv), \zv^{\otimes 2} \rangle} 
	=
	\frac{\| \HL(\upsv + t \uv) \zv \|^{2}}{\| \HL(\upsv) \zv \|^{2}}
	& \leq &
	\sqrt{\dltwaa_{\Psi,4} \, \CONSTi_{\rsmall}} 
	\, .
\label{vwehuf8ew87fuwyfe2}
\end{EQA}
%
Now we are prepared to finalize the check of \nameref{LLtS3ref} and \nameref{LLtS4ref}.
Let \( \upsv \in \B_{\rsmall}(\upsvs) \).
For any \( \uv \) with \( \| \HL(\upsv) \uv \| \leq \rr/\sqrt{\nsize} \leq \rsmall \), by \eqref{0det53543wyf7rryem} and \eqref{vwehuf8ew87fuwyfe2}
\begin{EQA}
    \frac{|\langle \nabla^{3} \cdens(\upsv + t \uv),\zv^{ \otimes 3} \rangle|}{\| \HL(\upsv) \zv \|^{3}}
    & \leq &	
    \frac{\dltwaa_{\Psi,3} \, \| \HL(\upsv + t \uv) \zv \|^{3}}{\| \HL(\upsv) \zv \|^{3}}
    \leq 
    \dltwaa_{\Psi,3} \, (\dltwaa_{\Psi,4} \, \CONSTi_{\rsmall})^{3/4} \, ,
\label{cuhywe8w3endsd72teu2hwd}
\end{EQA}
and \nameref{LLtS3ref} follows with 
\( \hmax_{3} = \dltwaa_{\Psi,3} \, (\dltwaa_{\Psi,4} \, \CONSTi_{\rsmall})^{3/4} \). 
The proof of \nameref{LLtS4ref} is similar. 
%
\end{proof}

Now we check \nameref{EU2ref}.
To be more specific, consider the deviation bound for \( \IF_{\GP}^{-1/2} \nabla \zeta = \Sv - \E \Sv \), where
\( \IF_{\GP} = \IF + \GP^{2} \) for \( \IF = \IF(\upsvs) \) and a penalty operator \( \GP^{2} \).
Define
\begin{EQA}
	\dimG
	& \eqdef &
	\tr(\IF_{\GP}^{-1} \IF) ,
	\qquad
	\rrGP 
	=
	\sqrt{\dimG} + \sqrt{2\xx} \, .
\label{7djxdcwesrtcywe3bdyu7sa}
\end{EQA}

\begin{lemma}
\label{LEUdens}
Assume the conditions of Lemma~\ref{LS3dens}. 
Then \nameref{EU2ref} holds with \( \VP^{2} = 2 \nsize \nabla^{2} \cdens(\upsvs) \) for \( \xx \leq (\rsmall \sqrt{\nsize/2} - \sqrt{\dimG})^{2}/4 \).
\end{lemma}
 
\begin{proof}
I.i.d. structure of \( \Sv = \sum_{i} X_{i} \) and \eqref{3094857ythgfjdkleo} yield 
\( \Var(\Sv) = \nsize \nabla^{2} \cdens(\upsvs) \).
Further, for any \( \uv \in \B_{\rsmall}(\upsvs) \), again by the i.i.d. assumption
and by \eqref{kjhr8dfrikweyyweyerdd}
\begin{EQA}
	\nsize^{-1} \log \E_{\upsvs} \exp\bigl\{ \langle \nabla \zeta ,\uv \rangle \bigr\}
	&=&
	\log \E_{\upsvs} \ex^{\langle \epsv, \uv \rangle} 
	=
	\cdens(\upsvs + \uv) - \cdens(\upsvs) - \langle \nabla \cdens(\upsvs),\uv \rangle .
\label{gw7fere9rh8uyegrh978e9}
\end{EQA}
For \( \rrGP \leq \rsmall \, \nsize^{1/2} \), consider all \( \uv \) with
\( \nsize \langle \nabla^{2} \cdens(\upsvs), \uv^{\otimes 2} \rangle \leq \rrGP^{2} \). 
If \( \hmax_{3} \, \rrGP \leq 3 \nsize^{1/2} \), 
then by \nameref{LLtS3ref} and \eqref{gtcdsftdfvtwdsefhfdvfrvsewseGP} of Lemma~\ref{LdltwLaGP}
\begin{EQA}
	\cdens(\upsvs + \uv) - \cdens(\upsvs) - \langle \nabla \cdens(\upsvs),\uv \rangle
	& \leq &
	\frac{1 + \hmax_{3} \, \rrGP \, \nsize^{-1/2}/3}{2} \, \langle \nabla^{2} \cdens(\upsvs), \uv^{\otimes 2} \rangle
	\leq 
	\langle \nabla^{2} \cdens(\upsvs), \uv^{\otimes 2} \rangle .
\label{d0if897erfe87y9u}
\end{EQA}
This implies \eqref{expgamgm} with \( \VP^{2} = 2 \nsize \nabla^{2} \cdens(\upsvs) \),
\( \gm = \rsmall \sqrt{\nsize/2} \) and thus, the deviation bound
\eqref{PxivbzzBBroB} of Theorem~\ref{Tdevboundgm} implies \nameref{EU2ref} for 
\( \xx \leq \xxc \leq (\rsmall \sqrt{\nsize/2} - \sqrt{\dimG})^{2}/4 \).
\end{proof}

\Section{Histogramm estimation}
Let \( X_{1},\ldots,X_{n} \) be i.i.d. with values in \( \XX \subseteq \R^{d} \).
Consider a partition \( \XX = \cup_{j=1}^{\dimp} E_{j} \) for non-overlapping sets \( E_{j} \) and define
\( \tarps_{j} = \P(X_{i} \in E_{j}) \), \( \psi_{j}(\xv) = \Ind(\xv \in E_{j}) \), and
\begin{EQA}
	S_{j}
	&=&
	\sumi \Ind(X_{i} \in E_{j})
	=
	\sumi \psi_{j}(X_{i}) .
\label{sicuy9isudiowel3eujhe}
\end{EQA}
W.l.o.g. assume that all \( \tarps_{j} \) are positive and \( \tarps_{j} \leq 1/2 \).
Set \( \upss_{j} = \log \tarps_{j} \).
The negative penalized log-likelihood reads 
\begin{EQA}
	\LGP(\upsv)
	&=&
	- \sum_{j=1}^{\dimp} S_{j} \ups_{j} + n \cdens(\upsv) + \frac{1}{2} \| \GP \upsv \|^{2} \, ,
\label{iw0klesd98c7yw6wyshvdid}
\end{EQA}
where
\begin{EQA}
	\cdens(\upsv) 
	&=& 
	\log(\ex^{\ups_{1}} + \ldots + \ex^{\ups_{\dimp}}) .
\label{dicvjcvy67er646yrfhwjh}
\end{EQA}
This model can be viewed as a special case of log-density.
Penalization by \( \| \GP \upsv \|^{2}/2 \) is important in the cases when some of \( \tarps_{j} \) are very close to zero
leading to large negative values of \( \ups_{j} \).
Later we consider a special case of a ridge penalization with \( \GP^{2} = \gp^{2} \Id_{\dimp} \).

\begin{lemma}
\label{Lhistphi}
For any \( \upsv = (\ups_{1} , \ldots , \ups_{\dimp})^{\T} \), 
it holds with \( \tarpv = \ex^{- \cdens(\upsv)} ( \ex^{\ups_{1}} , \ldots , \ex^{\ups_{\dimp}} )^{\T} \)
\begin{EQA}
	\nabla^{2} \cdens(\upsv)
	&=&
	\diag(\tarpv) - \tarpv \tarpv^{\T} \, .
\label{jd7wjwjiv5rrtw76387jeej}
\end{EQA}
\end{lemma}

\begin{proof}
It holds for all \( j \neq m \leq \dimp \)
\begin{EQA}[rcccl]
	\frac{\partial}{\partial \ups_{j}} \cdens(\upsv)
	&=&
	\frac{\ex^{\ups_{j}}}{\ex^{\ups_{1}} + \ldots + \ex^{\ups_{\dimp}}} 
	&=&
	\tarp_{j} 
	\, ,
\label{8y8diewjoiw66cuywjdfjh87}
	\\
	\frac{\partial^{2}}{\partial \ups_{j}^{2}} \cdens(\upsv)
	&=&
	\frac{\ex^{\ups_{j}}}{\ex^{\ups_{1}} + \ldots + \ex^{\ups_{\dimp}}} 
	- \frac{\ex^{2\ups_{j}}}{(\ex^{\ups_{1}} + \ldots + \ex^{\ups_{\dimp}})^{2}} 
	&=&
	\tarp_{j} - \tarp_{j}^{2}
	\, ,
	\\
	\frac{\partial^{2}}{\partial \ups_{j} \, \partial \ups_{m}} \cdens(\upsv)
	&=&
	- \frac{\ex^{\ups_{j}} \, \ex^{\ups_{m}}}{(\ex^{\ups_{1}} + \ldots + \ex^{\ups_{\dimp}})^{2}} 
	&=&
	\tarp_{j} \, \tarp_{m}
	\, ,
\label{yd6yh2jvc8tqgc7djw}
\end{EQA}
and the assertion follows.
\end{proof}

Later we drop the quadratic term \( \tarpv \tarpv^{\T} \) in \( \cdens(\upsv) \) and set
\( \HLH^{2}(\upsv) = \diag(\tarpv) \).

\begin{lemma}
\label{Lcdenshi}
Fix \( \upsv \in \R^{\dimp} \) and \( \tarpv = \ex^{- \cdens(\upsv)} ( \ex^{\ups_{1}} , \ldots , \ex^{\ups_{\dimp}} )^{\T} \).
Then it holds for \( \HLH^{2}(\upsv) = \diag(\tarpv) \) and any \( \uv \) with \( \| \uv \|_{\infty} \leq \const \)
\begin{EQA}
	\ex^{-\const} \, \HLH(\upsv)
	& \leq &
	\HLH(\upsv + \uv)
	\leq 
	\ex^{\const} \, \HLH(\upsv) .
\label{iud78ujkies8wei3f98e3k}
\end{EQA}
Moreover, \( \cdens(\upsv;\uv) = \cdens(\upsv + \uv) - \cdens(\upsv) - \langle \nabla \cdens(\upsv), \uv \rangle \) fulfills
\begin{EQA}
	\cdens(\upsv;\uv)
	& \leq & 
	\frac{\ex^{2\const}}{2} \| \HLH(\upsv) \uv \|^{2} \, .
\label{tdcjhc76yetey4rghc7ye}
\end{EQA}
\end{lemma}

\begin{proof}
Clearly \( \| \uv \|_{\infty} \leq \const \) implies 
\( \ex^{-\const} \ex^{\ups_{j}} \leq \ex^{\ups_{j} + u_{j}} \leq \ex^{\const} \ex^{\ups_{j}} \) and
the first assertion follows by definition of \( \HLH(\upsv) \).
Further, with some \( t \in [0,1] \)
\begin{EQA}
	\cdens(\upsv;\uv)
	&=&
	\cdens(\upsv + \uv) - \cdens(\upsv) - \langle \nabla \cdens(\upsv), \uv \rangle
	=
	\frac{1}{2} \langle \nabla^{2} \cdens(\upsv + t \uv), \uv^{\otimes 2} \rangle
	\leq 
	\frac{\ex^{2\const}}{2} \| \HLH(\upsv) \uv \|^{2}
\label{fciub6yhuii99wendy}
\end{EQA}
and \eqref{tdcjhc76yetey4rghc7ye} follows as well.
\end{proof}

The next step is to bound the higher order derivatives of \( \cdens(\upsv) \).

\begin{lemma}
\label{Lcdens34hi}
For any \( \uv \in \R^{\dimp} \), define a r.v. \( U \) with \( \P_{\upsv}(U = u_{j}) = \tarp_{j} \).
Then
\begin{EQA}[rcccl]
	\bigl\langle \nabla \cdens(\upsv), \uv \bigr\rangle
	&=&
	\sum_{j} u_{j} \tarp_{j} 
	&=&
	\E_{\upsv} U
	\, ,
	\\
	\bigl\langle \nabla^{2} \cdens(\upsv), \uv^{\otimes 2} \bigr\rangle
	&=&
	\sum_{j} u_{j}^{2} \tarp_{j} - \Bigl( \sum_{j} u_{j} \tarp_{j} \Bigr)^{2}
	&=&
	\E_{\upsv} (U - \E_{\upsv} U)^{2}
	\, ,
\label{kdyyudd6ehjd8femfof}
\end{EQA}
and
\begin{EQA}
	\bigl\langle \nabla^{3} \cdens(\upsv), \uv^{\otimes 3} \bigr\rangle
	&=&
	\sum_{j} u_{j}^{3} \tarp_{j}
	- 3 \sum_{j} u_{j} \tarp_{j} \sum_{j} u_{j}^{2} \tarp_{j} 
	+ \Bigl( \sum_{j} u_{j} \tarp_{j} \Bigr)^{3}
	=
	\E_{\upsv} (U - \E_{\upsv} U)^{3} 
	\, .
\label{odu8c87cd78d7df766w3}
\end{EQA}
Moreover, 
\begin{EQA}[rcccl]
	\bigl| \E_{\upsv} U \bigr|
	& \leq &
	\| \HLH(\upsv) \uv \| ,
	\qquad
	\E_{\upsv} (U - \E_{\upsv} U)^{2}
	& \leq &
	\E_{\upsv} U^{2}
	= 
	\| \HLH(\upsv) \uv \|^{2} ,
\label{6djhc78w3jdkl9hgb84rj}
\end{EQA}
and for any \( \gp^{2} \geq 0 \), it holds with \( \HL^{2}(\upsv) = \HLH^{2}(\upsv) + \gp^{2} \Id_{\dimp} \)
\begin{EQA}
	\bigl\langle \nabla^{3} \cdens(\upsv), \uv^{\otimes 3} \bigr\rangle
	& \leq &
	\frac{1}{(\tarp_{\min} + \gp^{2})^{1/2}} \| \HL(\upsv) \uv \|^{3} + 3 \| \HLH(\upsv) \uv \|^{3} \, .
\label{s8d8wdiwjdaiqdeihdksdd}
\end{EQA}
\end{lemma}

\begin{proof}
It holds by \( \tarp_{j} = \ex^{\ups_{j}}/(\ex^{\ups_{1}} + \ldots + \ex^{\ups_{\dimp}}) \)
with \( \bar{u} \eqdef \E_{\upsv} U = \sum_{j} u_{j} \tarp_{j} \)
\begin{EQA}[rcccl]
	\bigl\langle \nabla \cdens(\upsv), \uv \bigr\rangle
	&=&
	\frac{\sum_{j} u_{j} \ex^{\ups_{j}}}{\ex^{\ups_{1}} + \ldots + \ex^{\ups_{\dimp}}} 
	&=&
	\bar{u}
	\, ,
	\\
	\bigl\langle \nabla^{2} \cdens(\upsv), \uv^{\otimes 2} \bigr\rangle
	&=&
	\frac{\sum_{j} u_{j}^{2} \ex^{\ups_{j}}}{\ex^{\ups_{1}} + \ldots + \ex^{\ups_{\dimp}}} 
	- \frac{(\sum_{j} u_{j} \ex^{\ups_{j}})^{2}}{(\ex^{\ups_{1}} + \ldots + \ex^{\ups_{\dimp}})^{2}} 
	&=&
	\E_{\upsv} (U - \bar{u})^{2}
	\, ,
\label{kdyyudd6ehjd8femfof}
\end{EQA}
and \eqref{6djhc78w3jdkl9hgb84rj} follows by \( \E_{\upsv} U^{2} = \sum_{j} u_{j}^{2} \tarp_{j} \) and by definition of \( \HLH(\upsv) \).
Similarly,
\begin{EQA}
	\bigl\langle \nabla^{3} \cdens(\upsv), \uv^{\otimes 2} \bigr\rangle
	&=&
	\frac{\sum_{j} u_{j}^{3} \ex^{\ups_{j}}}{\ex^{\ups_{1}} + \ldots + \ex^{\ups_{\dimp}}} 
	- \frac{3\sum_{j} u_{j} \ex^{\ups_{j}} \,\, \sum_{j} u_{j}^{2} \ex^{\ups_{j}}}{(\ex^{\ups_{1}} + \ldots + \ex^{\ups_{\dimp}})^{2}} 
	+ \frac{2(\sum_{j} u_{j} \ex^{\ups_{j}})^{3}}{(\ex^{\ups_{1}} + \ldots + \ex^{\ups_{\dimp}})^{3}}
	\\
	&=&
	3 \E_{\upsv} U^{3} - 3 \bar{u} \, \E_{\upsv} U^{2} + 2 \bar{u}^{3}
	=
	\E_{\upsv} (U - \bar{u})^{3}
	\, .
\label{kdyyudd6ehjd8fem3333}
\end{EQA}
By \eqref{6djhc78w3jdkl9hgb84rj}
\begin{EQA}
	\bigl| 3 \bar{u} \, \E_{\upsv} U^{2} - 2 \bar{u}^{3}  \bigr|
	& \leq &
	\bigl| \bar{u} \, \E_{\upsv} U^{2} \bigr| + 2 \bigl| \bar{u} \, (\E_{\upsv} U^{2} - \bar{u}^{2}) \bigr|
	\leq 
	3 \bigl| \bar{u} \, \E_{\upsv} U^{2} \bigr|
	\leq 
	3 \| \HLH(\upsv) \uv \|^{3} \, .
\label{7dejwekdfbyy2wnedyg7y7e3h}
\end{EQA}
Further, for any \( \gp^{2} \geq 0 \)
\begin{EQA}
	\Bigl| \sum_{j} u_{j}^{3} \tarp_{j} \Bigr|
	& \leq &
	\Bigl| \sum_{j} u_{j}^{3} (\tarp_{j} + \gp^{2}) \Bigr|
	\leq 
	\frac{1}{(\tarp_{\min} + \gp^{2})^{1/2}} \sum_{j} |u_{j}|^{3} (\tarp_{j} + \gp^{2})^{3/2} 
	\\
	& \leq &
	\frac{1}{(\tarp_{\min} + \gp^{2})^{1/2}} \, \Bigl| \sum_{j} u_{j}^{2} (\tarp_{j} + \gp^{2}) \Bigr|^{3/2} \, ,
\label{uye7fyeihfiwe8wrd8w2lwu}
\end{EQA}
and \eqref{s8d8wdiwjdaiqdeihdksdd} follows.
\end{proof}


\Section{Precision matrix estimation}
Let \( \Xv_{1},\ldots,\Xv_{\nsize} \) be i.i.d. zero mean Gaussian vector in \( \R^{\dimp} \) 
with a covariance \( \Sigma \): \( \Xv_{i} \sim \ND(0,\Sigma) \).
Our goal is to estimate the corresponding precision matrix \( \Prec = \Sigma^{-1} \).
Later we identify the matrix \( \Prec \) with the point in the linear subspace \( \PREC \) of \( \R^{\dimp \times \dimp} \)
composed by symmetric matrices. 
The ML approach leads to the negative log-likelihood
\begin{EQA}
	L(\Prec)
	&=&
	\frac{1}{2} \sumi \langle \hKS_{i}, \Prec \rangle - \frac{\nsize}{2} \log \det(\Prec)
\label{uchdy7c7fj4ib8jelvgewr}
\end{EQA}
with \( \hKS_{i} = \Xv_{i} \Xv_{i}^{\T} \).
Here and later \( \langle A,B \rangle \) means \( \tr(A B) \) for \( A, B \in \PREC \).
The corresponding MLE \( \tilde{\Prec} \) minimizes \( L(\Prec) \):
\begin{EQA}
	\tilde{\Prec}
	&=&
	\argmin_{\Prec \in \PREC} L(\Prec) .
\label{c8jr9fvc6hocejvkborencxy}
\end{EQA}
The target of estimation \( \Precs \) can be defined as its population counterpart:
\begin{EQA}
	\Precs
	&=&
	\argmin_{\Prec \in \PREC} \E L(\Prec) .
\label{c8jr9fvc6hocejvkbotarg}
\end{EQA}
Now introduce a quadratic penalization on \( \Prec \) in the form \( \| \GPKS \Prec \|_{\Fr}^{2}/2 \)
for a linear operator \( \GPKS \) on the space \( \PREC \).
One typical example corresponds to the case with
\begin{EQA}
	\| \GPKS \Prec \|_{\Fr}^{2}
	&=&
	\sum_{m=1}^{\dimp} \| \GPKS_{m} \Prec_{m} \|^{2} ,
	\qquad
	\Prec = (\Prec_{1},\ldots,\Prec_{\dimp}) ,
\label{uev7eun3evbkbdr3fevunwd}
\end{EQA}
for a family of linear mappings \( \GPKS_{1},\ldots,\GPKS_{\dimp} \) in \( \R^{\dimp} \).
The corresponding penalized MLE \( \tilde{\Prec}_{\GPKS} \) is defined by minimizing 
the penalized loss \( L_{\GPKS}(\Prec) = L(\Prec) + \| \GPKS \Prec \|_{\Fr}^{2} \):
\begin{EQA}
	\tilde{\Prec}_{\GPKS}
	&=&
	\argmin_{\Prec} L_{\GPKS}(\Prec)
	=
	\argmin_{\Prec} \biggl\{ L(\Prec) + \frac{1}{2} \| \GPKS \Prec \|_{\Fr}^{2} \biggr\}
	\\
	&=&
	\argmin_{\Prec} \biggl\{ 
		\sumi \langle \hKS_{i}, \Prec \rangle - \nsize \log \det(\Prec) + \| \GPKS \Prec \|_{\Fr}^{2}
	\biggr\} .
\label{gdhwjduvyeebdchyndec}
\end{EQA}
Define also the penalized target \( \Precs_{\GPKS} \) as
\begin{EQA}
	\Precs_{\GPKS}
	&=&
	\argmin_{\Prec} \biggl\{ \E L(\Prec) + \frac{1}{2} \| \GPKS \Prec \|_{\Fr}^{2} \biggr\}
	=
	\argmin_{\Prec} \biggl\{ 
		\nsize \langle \Sigma, \Prec \rangle - \nsize \log \det(\Prec) + \| \GPKS \Prec \|_{\Fr}^{2}
	\biggr\} .
\label{9xjcxy3wjvbhfobogte5tv}
\end{EQA}
We intend to state some sharp bounds on the loss and risk of \( \tilde{\Prec}_{\GPKS} \) by applying
the general results of Section~\ref{SgenBounds}.
Model \eqref{uchdy7c7fj4ib8jelvgewr} is a special case of an exponential family.
Therefore, the basic assumptions \nameref{Eref} and \nameref{LLref} are fulfilled
automatically.
Next, we check the smoothness properties of \( \E L(\Prec) \) in terms of the Gatoux derivatives.

\begin{lemma}
\label{LPrecsmooth}
Let \( \Prec \in \PREC \) be positive definite.
For any \( \zv \in \PREC \) and \( \Uv = \Prec^{-1/2} \, \zv \, \Prec^{-1/2} \), it holds
\begin{EQA}
	\frac{d^{2}}{dt^{2}} \, \E L(\Prec + t \, \zv) \bigg|_{t=0}
	&=&
	\frac{\nsize}{2} \, \tr \Uv^{2}
	=
	\frac{\nsize}{2} \, \tr\{ (\Prec^{-1} \zv)^{2} \} .
\label{hf7vnjbgu437bjh9r5}
\end{EQA}
Similarly
\begin{EQA}
	\frac{d^{3}}{dt^{3}} \E L(\Prec + t \zv) \bigg|_{t=0}
	&=&
	- \nsize \tr\Uv^{3} \, ,
	\qquad
	\frac{d^{4}}{dt^{4}} \E L(\Prec + t \zv) \bigg|_{t=0}
	=
	3 \nsize \tr\Uv^{4} \, .
	\qquad
\label{8vjkvct6ehjb8ejuetb8t4r}
\end{EQA}
\end{lemma}

\begin{proof}
Fix some \( \zv \in \PREC \).
It holds by \eqref{uchdy7c7fj4ib8jelvgewr} 
with \( \Uv = \Prec^{-1/2} \, \zv \, \Prec^{-1/2} \)
\begin{EQA}
	&& \nquad
	\frac{d^{2}}{dt^{2}} \E L(\Prec + t \, \zv) \bigg|_{t=0}
	=
	- \frac{\nsize}{2} \, \frac{d^{2}}{dt^{2}} \, \log \det(\Prec + t \, \zv) \bigg|_{t=0}
	\\
	&=&
	- \frac{\nsize}{2} \, \frac{d^{2}}{dt^{2}} \, \log \det(\Id_{\dimp} + t \, \Uv) \bigg|_{t=0}
	=
	\frac{\nsize}{2} \, \tr \Uv^{2}
	=
	\frac{\nsize}{2} \, \| \Prec^{-1/2} \, \zv \, \Prec^{-1/2} \|_{\Fr}^{2} \, .
\label{ufiw2nfvybyfn3g893mjfu}
\end{EQA}
This formula can easily be checked when \( \Uv \) is diagonal, the general case 
is reduced to this one by an orthogonal transform.
\eqref{8vjkvct6ehjb8ejuetb8t4r} can be checked similarly.
\end{proof}

\noindent
Bounds \eqref{8vjkvct6ehjb8ejuetb8t4r} help to check condition \nameref{LLsT3ref} and \nameref{LLsT4ref}.

\begin{lemma}
\label{LS3S4Prec}
For \( \Prec \in \PREC \) positive definite, define \( \DFN^{2}(\Prec) \) by 
\begin{EQA}
	\| \DFN(\Prec) \zv \|_{\Fr}^{2}
	&=&
	\frac{\nsize}{2} \| \Prec^{-1/2} \, \zv \, \Prec^{-1/2} \|_{\Fr}^{2} \, ,
	\qquad
	\zv \in \PREC \, .
\label{ydjd67544rgduidf78rh}
\end{EQA}
Then \nameref{LLsT3ref} and \nameref{LLsT4ref} are fulfilled with \( \rr < n/2 \) and
\begin{EQA}
	\dltwu_{3}
	&=&
	\sqrt{8} \bigl( 1 - \sqrt{2 \rr^{2}/\nsize} \bigr)^{-3} n^{-1/2} \, ,
	\qquad
	\dltwu_{4}
	=
	12 \bigl( 1 - \sqrt{2 \rr^{2}/\nsize} \bigr)^{-4} n^{-1} \, .
\label{9fjdftdf4e5wgdfgyeu}
\end{EQA}
\end{lemma}

\begin{proof}
Consider \( \uv \in \PREC \) such that \( \| \DFN(\Prec) \uv \|_{\Fr} \leq \rr \).
Fix \( \zv \in \PREC \) and define \( \Uv = (\Prec+\uv)^{-1/2} \zv \, (\Prec+\uv)^{-1/2} \).
Then by \eqref{8vjkvct6ehjb8ejuetb8t4r}
and by \( \bigl| \tr \Uv^{3} \bigr| \leq \bigl( \tr \Uv^{2} \bigr)^{3/2} \), the function \( \fs(\Prec) \) satisfies
\begin{EQA}
	\bigl| \langle \nabla^{3} \fs(\Prec + \uv), \zv^{\otimes 3} \rangle \bigr|
	& \leq &
	\nsize \bigl| \tr \UV^{3} \bigr|
	\leq 
	\nsize \bigl( \tr \Uv^{2} \bigr)^{3/2}
	\leq 
	n \| (\Prec+\uv)^{-1/2} \zv \, (\Prec+\uv)^{-1/2} \|_{\Fr}^{3}
\label{6hjdsy6dhw8df7gwtrdghv}
\end{EQA}
Further, \( \| \DFN(\Prec) \uv \|_{\Fr} \leq \rr \) implies
\begin{EQA}
	\| \Prec^{-1/2} (\Prec+\uv) \Prec^{-1/2} - \Id \|^{2}
	&=&
	\| \Prec^{-1/2} \uv \, \Prec^{-1/2} \|^{2}
	\leq 
	\| \Prec^{-1/2} \uv \, \Prec^{-1/2} \|_{\Fr}^{2}
	\leq 
	2 \rr^{2}/\nsize \, 
\label{ufehedjbv87rjfhxrrs}
\end{EQA}
yielding 
\begin{EQA}
	\| (\Prec+\uv)^{-1/2} \Prec^{1/2} \|^{2}
	&=&
	\| (\Id + \Prec^{-1/2} \uv \, \Prec^{-1/2})^{-1} \|
	\leq 
	\frac{1}{1 - \sqrt{2 \rr^{2}/\nsize}} \, .
\label{jduwjklb6554e34ftds}
\end{EQA}
Therefore, for any \( \zv \in \PREC \)
\begin{EQA}
	\| (\Prec+\uv)^{-1/2} \zv \, (\Prec+\uv)^{-1/2} \|_{\Fr}
	& \leq &
	\| (\Prec+\uv)^{-1/2} \Prec \|^{2} \,\, \| \Prec^{-1/2} \zv \, \Prec^{-1/2} \|_{\Fr}
	\\
	& \leq &
	\frac{1}{1 - \sqrt{2 \rr^{2}/\nsize}} \, \| \Prec^{-1/2} \zv \, \Prec^{-1/2} \|_{\Fr} \, .
\label{hc7hjejifog0ejfytkh3}
\end{EQA}
Therefore,
\begin{EQA}
	\bigl| \langle \nabla^{3} \fs(\Prec + \uv), \zv^{\otimes 3} \rangle \bigr|
	& \leq &
	\nsize \bigl( 1 - \sqrt{2 \rr^{2}/\nsize} \bigr)^{-3} \| \Prec^{-1/2} \zv \, \Prec^{-1/2} \|_{\Fr}^{3}
\label{jw32gwnduvf5hwdtfnryr}
\end{EQA}
and condition \nameref{LLsT3ref} is fulfilled with 
\( \dltwu_{3} = \sqrt{8} \bigl( 1 - \sqrt{2 \rr^{2}/\nsize} \bigr)^{-3} n^{-1/2} \).
Similarly, one can check \nameref{LLsT4ref}.
\end{proof}

By \eqref{hf7vnjbgu437bjh9r5}, the Fisher matrix 
\( \IF(\Prec) \eqdef - \nabla^{2} \E L(\Prec) \) is a linear operator in \( \PREC \) with
\begin{EQA}
	&& \nquad
	\langle \IF(\Prec), \uv^{\otimes 2} \rangle
	=
	\frac{\nsize}{2} \, \tr \Uv^{2}
	=
	\frac{\nsize}{2} \, \tr\{ (\Prec^{-1} \uv)^{2} \} \, ,
	\qquad
	\Uv = \Prec^{-1/2} \, \zv \, \Prec^{-1/2} \, .
\label{ufiw2nfvybyfn3g8IF}
\end{EQA}
The penalized Fisher information operator \( \IF_{\GPKS} = \IF(\Precs) + \GPKS^{2} \) is given by
\begin{EQA}
	\| \IF_{\GPKS}^{1/2} \uv \|_{\Fr}^{2}
	&=&
	\frac{\nsize}{2} \tr\{ (\Sigma \uv)^{2} \} + \| \GPKS \uv \|_{\Fr}^{2} \, ;
	\qquad
	\uv \in \PREC \, .
\label{8jc53lv9cjfyenvjcorvkbkKS}
\end{EQA}

The next condition to be verified is \nameref{EU2ref}.
\begin{lemma}
\label{LEUPrec}
Let \( \IF_{\GPKS} = \IF(\Precs) + \GPKS^{2} \).
Define also
\begin{EQA}
	\BBH_{\GPKS}^{2}
	& \eqdef &
	\IF_{\GPKS}^{-1/2} \Sigma^{2} \IF_{\GPKS}^{-1/2} \, .
\label{d6hwneuv8vflbi83j}
\end{EQA}
If \( \dimH_{\GPKS} < \nsize/8 \), then \nameref{EU2ref} is fulfilled with
\begin{EQA}
	\rr_{\GPKS}(\xx)
	& \eqdef &
	\sqrt{\dimH_{\GPKS}} + \sqrt{2 \xx} \, ,
	\qquad
	\dimH_{\GPKS} 
	\eqdef 
	(\tr \BBH_{\GPKS})^{2} + \tr \BBH_{\GPKS}^{2} \, .
\label{0cyd5c56fyb84fnv73fv}
\end{EQA}
\end{lemma}

\begin{proof}
The stochastic component \( \zeta(\Prec) = L(\Prec) - \E L(\Prec) \) reads
\begin{EQA}
	\zeta(\Prec) 
	&=& 
	L(\Prec) - \E L(\Prec)
	=
	\frac{1}{2} \sumi \langle \Err_{i}, \Prec \rangle
\label{jcucutehvjuvjrencfsdki}
\end{EQA}
with \( \Err_{i} = \hKS_{i} - \E \hKS_{i} = \Xv_{i} \Xv_{i}^{\T} - \Sigma \).
Clearly \( \zeta(\Prec) \) is linear in \( \Prec \) and condition \nameref{Eref} 
is fulfilled. 
Moreover, for any direction \( \uv \) in the parameter space \( \Ups \),
\begin{EQA}
	\langle \nabla \zeta,\uv \rangle
	&=&
	\frac{1}{2} \sumi \langle \Err_{i}, \uv \rangle
	=
	\frac{1}{2} \sumi (\Xv_{i}^{\T} \uv \Xv_{i} - \langle \Sigma,\uv \rangle) 
	=
	\frac{1}{2} \sumi \bigl\{ \gaussv_{i}^{\T} \Uv \gaussv_{i} - \tr (\Uv) \bigr\}
\label{hyhewuc7jffjud7egyvee}
\end{EQA}
with \( \gaussv_{i} = \Sigma^{-1/2} \Xv_{i} \) standard Gaussian and 
\( \Uv = \Sigma^{1/2} \uv \Sigma^{1/2} \).
By Lemma~\ref{Gaussmoments}
\begin{EQA}
	\Var\langle \nabla \zeta,\uv \rangle
	&=&
	\frac{\nsize}{4} \Var\bigl( \gaussv_{1}^{\T} \Uv \gaussv_{1} \bigr)
	=
	\frac{\nsize}{2} \tr(\Uv^{2}) 
	=
	\frac{\nsize}{2} \tr(\Sigma \uv)^{2} .
\label{7xhcu3ob0fjf48jvu3erlvew3y}
\end{EQA}
In particular, for \( \Prec = \Precs = \Sigma^{-1} \), 
the operator \( \IF = \IF(\Precs) \) coincides with \( \Var(\nabla \zeta) \):
\begin{EQA}
	\| \IF^{1/2} \uv \|_{\Fr}^{2}
	&=&
	\frac{\nsize}{2} \tr(\Sigma \uv)^{2} .
\label{8jc53lv9cjfyenvjcorvkbk}
\end{EQA}
This is in agreement with the fact that under the true parametric assumption, it holds
\( - \nabla^{2} \E L(\Precs) = \Var\bigl( \nabla \zeta \bigr) \).
One can easily check for any \( \uv \in \PREC \)
\begin{EQA}
	\Var \langle \IF_{\GPKS}^{-1/2} \nabla \zeta, \uv \rangle
	&=&
	\Var \langle \nabla \zeta, \IF_{\GPKS}^{-1/2} \uv \rangle
	=
	\frac{\nsize}{2} \tr(\uv \IF_{\GPKS}^{-1/2} \Sigma^{2} \IF_{\GPKS}^{-1/2} \uv) 
	=
	\frac{\nsize}{2} \tr(\uv \, \BBH_{\GPKS}^{2} \, \uv)
\label{cdhd7c3h3ycckwxu7cjeen}
\end{EQA}
with \( \BBH_{\GPKS} \) from \eqref{d6hwneuv8vflbi83j}.
Moreover, Theorem~\ref{TFrobGaussA} yields under \( \dimH_{\GPKS} < \nsize/8 \)
\begin{EQA}
	\P\bigl( \| \IF_{\GPKS}^{-1/2} \nabla \zeta \|_{\Fr} \geq \rr_{\GPKS}(\xx) \bigr)
	& \leq &
	3 \ex^{-\xx} ;
\label{ycncdh4fd87hfibi78ehv}
\end{EQA}
see \eqref{0j3wgyrfg76hy5tjhtdeert}, and  \nameref{EU2ref} is verified.
\end{proof}



\def\deta{\delta}

\Section{Anisotropic logistic regression}
\label{Saniligot}
This section considers a popular logistic regression model. 
It is widely used e.g. in binary classification in machine learning for binary classification or in binary response models in econometrics.
The results presented here can be viewed as an extension of \cite{SpPa2019} and \cite{SpLaplace2022}.

Suppose we are given a vector of independent observations/labels \( \Yv = (Y_{1},\ldots, Y_{\nsize})^{\T} \) 
and a set of the corresponding feature vectors \( \Psiv_{i} \in \R^{\dimp} \).
Each binary label \( Y_{i} \) is modelled as a Bernoulli random variable with the parameter 
\( \thetas_{i} = \P(Y_{i} = 1) \).
Logistic regression reduces this model to linear regression for the canonical parameter 
\( \upss_{i} = \log \frac{\thetas_{i}}{1-\thetas_{i}} \) in the form \( \upss_{i} = \langle \Psiv_{i}, \upsvs \rangle \),
where \( \upsv \) is the parameter vector in \( \R^{\dimp} \).
The corresponding negative log-likelihood reads 
\begin{EQA}
	L(\upsv)
	&=&
	- \sumi \Bigl\{ Y_{i} \, \langle \Psiv_{i}, \upsv \rangle 
	- \cdens\bigl( \langle \Psiv_{i}, \upsv \rangle \bigr) \Bigr\}
\label{LuYiPsiGLM}
\end{EQA}
with \( \cdens(\ups) = \log\bigl( 1 + \ex^{\ups} \bigr) \).
This function \( \cdens(\cdot) \) is convex with \( \cdens''(\ups) = \frac{\ex^{\ups}}{(1 + \ex^{\ups})^{2}} \).
We also assume that the \( \Psiv_{i} \)'s are deterministic, otherwise, we condition on the design.
A penalized MLE \( \tilde{\upsv}_{\GP} \) is defined by minimization
of the penalized log-likelihood \( \LGP(\upsv) = L(\upsv) + \| \GP \upsv \|^{2}/2 \)
for the quadratic penalty \( \| \GP \upsv \|^{2}/2 \):
\begin{EQA}
	\tilde{\upsv}_{\GP} 
	&=& 
	\argmin_{\upsv \in \R^{\dimp}} \LGP(\upsv) .
\label{vb8rj4ru7yfueow3e}
\end{EQA}
The truth and the penalized truth are 
defined via the expected log-likelihood
\begin{EQ}[rcccl]
	\upsvs 
	&=& 
	\argmin_{\upsv \in \R^{\dimp}} \E L(\upsv) ,
	\quad
	\upsvs_{\GP} 
	&=& 
	\argmin_{\upsv \in \R^{\dimp}} \E \LGP(\upsv) .
\label{vb8rj4ru7yfueow3eEE}
\end{EQ}
%
The Fisher matrix \( \IF(\upsv) \) at \( \upsv \) is given by
\begin{EQ}[rcccl]
	\IF(\upsv)
	&=&
	\nabla^{2} \E L(\upsv)
	=
	\sumi  \weight_{i}(\upsv) \, \Psiv_{i} \, \Psiv_{i}^{\T} ,
	\quad
	\weight_{i}(\upsv)
	& \eqdef &
	\cdens''\bigl( \langle \Psiv_{i}, \upsv \rangle \bigr) 
	\, .
	\qquad
\label{IFGusiphpPiuvT}
\end{EQ}
We also write 
\begin{EQA}
	\IF_{\GP}(\upsv) &=& \IF(\upsv) + \GP^{2} ,
	\qquad
 	\IF_{\GP} = \IF_{\GP}(\upsvs_{\GP}) = \IF(\upsvs_{\GP}) + \GP^{2} .
\label{tudxc7c7cycy55dfgetcct}
\end{EQA}
Alternatively one can define \( \IF_{\GP} = \IF_{\GP}(\upsvs) \).
Later we need a metric tensor \( \DPN \) which defines a local vicinity of \( \upsvs \).
We assume \( \IF \leq \DPN^{2} \leq \IF_{\GP} \).
If \( \IF \) is well posed, then \( \DPN^{2} = \IF \) is a default choice.

Further, we discuss the stochastic component of \( L(\upsv) \)
\begin{EQA}
	\zeta(\upsv) 
	&=& 
	L(\upsv) - \E L(\upsv) 
	= 
	- \sumi \bigl( Y_{i} - \E Y_{i} \bigr) \langle \Psiv_{i},\upsv \rangle 
\label{7fjvuehd6vgwesdrdgenhy}
\end{EQA}
It is obviously linear in \( \upsv \) with
\begin{EQA}
	\nabla \zeta
	&=&
	- \sumi \bigl( Y_{i} - \E Y_{i} \bigr) \Psiv_{i} \, .
\label{nabzesiYiEYiPi}
\end{EQA}
Further, with \( Y_{i} \sim \Bernoulli(\thetas_{i}) \) 
\begin{EQA}
	\Var(\nabla \zeta)
	&=&
	\sumi \Var(Y_{i}) \, \Psiv_{i} \Psiv_{i}^{\T} 
	=
	\sumi \thetas_{i} (1 - \thetas_{i}) \, \Psiv_{i} \Psiv_{i}^{\T} .
\label{V2siVeiPiT2}
\end{EQA}
If model \eqref{LuYiPsiGLM} is correctly specified, that is, 
\( \thetas_{i} = \ex^{\langle \Psiv_{i},\upsvs \rangle}/(1 + \ex^{\langle \Psiv_{i},\upsvs \rangle}) \), 
then \( \Var(\nabla \zeta) = \IF(\upsvs) \); see \eqref{IFGusiphpPiuvT}.
Instead, we assume \( \Var(\nabla \zeta) \leq \DPN^{2} \).
This can be relaxed to \( \Var(\nabla \zeta) \leq \CONST \DPN^{2} \) for some fixed constant \( \CONST \).
%
%
Now we check the general conditions from \Chname\ref{SgenBounds}.
Convexity of \( \cdens(\cdot) \) yields concavity of \( L(\upsv) \) and thus, \nameref{LLref}.
Condition \nameref{Eref} is granted by \eqref{nabzesiYiEYiPi}.
For checking the other conditions, we need some regularity of the design 
\( \Psiv_{1},\ldots,\Psiv_{\nsize} \).
Let \( \upsvd \) be either \( \upsvs \), or \( \upsvs_{\GP} \).
Define an elliptic vicinity \( \Upsd \) of \( \upsvd \) as
\begin{EQA}
	\Upsd 
	&=& 
	\{ \upsv \colon \| \DPN \, (\upsv - \upsvd) \| \leq \rrL \} 
	\, .
	\qquad
\label{txxctde43d4s4srgwcrdf}
\end{EQA}

\begin{description}
\item[\label{Psivdref} \( (\Psiv^{\circ}) \)]

	\emph{For some \( \detam > 0 \)}
\begin{EQA}[rcl]
\label{diuq3rvyewrhjdsfasd}
	\max_{i \leq \nsize} \| \DPN^{-1} \, \Psiv_{i} \|
	& \leq &
	\detam \, .
\end{EQA}
\end{description}

It appears that  all the results from Section~\ref{SgenBounds}
provided that \( \detam \) is sufficiently small.
We also assume for some \( \deta \geq 0 \)
\begin{EQA}
	\sumi \langle \Psiv_{i}, \zv \rangle^{4} \, \weight_{i}(\upsv) 
	& \leq &
	\deta^{2} \, \| \DPN \zv \|^{4} \, ,
	\qquad
	\zv \in \R^{\dimp} \, .
\label{k23CPpHkuf2}
\end{EQA}
Later we show that this condition follows from \nameref{Psivdref} for \( \deta^{2} \leq \sqrt{\ex} \, \detam^{2} \).
However, this is a conservative upper bound, condition \eqref{k23CPpHkuf2} may apply with much smaller values of \( \deta \).
As the final accuracy bound depends on \( \deta \) rather then on \( \detam \), we keep \eqref{k23CPpHkuf2}
as a separate condition.

The penalty \( \| \GP \upsv \|^{2}/2 \) in the pMLE \( \tilde{\upsv}_{\GP} \) results in some bias.
It can be measured by
\begin{EQA}[rcl]
	\bias_{\DPN}
	& \eqdef &
	\| \DPN \IF_{\GP}^{-1} \GP^{2} \upsvs \| 
	\, .
\label{ghdrd324ee4ewlll}
\end{EQA}

\begin{proposition}
\label{PGLMall}
Consider the logistic regression model \eqref{LuYiPsiGLM} and let \( \Var(\nabla \zeta) \leq \DPN^{2} \); 
see \eqref{IFGusiphpPiuvT} and \eqref{V2siVeiPiT2}.
With \( \BBH = 2 (\DPN \, \IF_{\GP}^{-1} \DPN)^{2} \), define 
\( \rr_{\DPN} = \zq(\BBH,\xx) \leq \sqrt{\dimA_{\DPN}} + \sqrt{2 \xx \, \lambda_{\DPN}} \) for 
\( \dimA_{\DPN} = \tr (\BBH) \), \( \lambda_{\DPN} = \| \BBH \| \).
Assume \nameref{Psivdref} with \( \detam \) satisfying 
\begin{EQA}[c]
	(\rr_{\DPN} \vee \bias_{\DPN}) \, \detam \leq 1/3, 
	\qquad 
	1 / (2\detam) \geq \sqrt{2\xx} + \sqrt{\dimA_{\DPN}/\lambda_{\DPN}} 
	\, .
\label{t56dfkjd6ffg5t54453}
\end{EQA}
Then \eqref{k23CPpHkuf2} holds for \( \deta^{2} \leq \sqrt{\ex} \, \detam^{2} \) all the results of 
\ifMLE{Theorem~\ref{TFiWititG2} through Theorem~\ref{TQFiWibias}}{Theorem~\ref{TMLEall}} apply with \( \dltwu_{3} \leq \sqrt{\ex} \, \deta \leq \ex^{3/4} \detam \).
\end{proposition}

\begin{proof}
We systematically check the conditions
of \ifMLE{Theorem~\ref{TFiWititG2} through Theorem~\ref{TQFiWibias}}{Theorem~\ref{TMLEall}}.
Conditions \nameref{Eref} and \nameref{LLref} are granted by \eqref{LuYiPsiGLM} with \( \cdens \) convex.
Next, we check \nameref{EU2ref}
for a matrix \( \VP^{2} \) such that \( \VP^{2} \geq 2 \Var(\nabla \zeta) \) and \( \VP^{2} \geq \DPN^{2} \).

\begin{lemma}
\label{LcheckEGLM}
Let \( \VP^{2} \geq 2 \Var(\nabla \zeta) \) and  \( \VP^{2} \geq \DPN^{2} \).
For \( \BBH \in \Matr_{\dimp} \), define 
\begin{EQA}
	\dimA_{\BBH} 
	&=& 
	\tr\bigl( \BBH \bigr) ,
	\qquad
	\lambda_{\BBH} = \| \BBH \| \, ,
	\qquad
	\rr_{\BBH}
	=
	\zq(\BBH,\xx)
	\leq 
	\sqrt{\dimA_{\BBH}} + \sqrt{2\xx \lambda_{\BBH}} \, .
	\qquad
\label{ujf8uejecy7whwejd}
\end{EQA}
With \( \detam \) from \eqref{diuq3rvyewrhjdsfasd}, let
\begin{EQA}
	1 / (2\detam) 
	& > & 
	\sqrt{\dimA_{\BBH}/\lambda_{\BBH}} \, .
\label{7dhdftgrjfte3jfighr4ynvm}
\end{EQA}
Then \nameref{EU2ref} is fulfilled with for this \( \BBH \) and all \( \xx \) such that
\( \sqrt{2\xx} \leq 1 / (2\detam) - \sqrt{\dimA_{\BBH}/\lambda_{\BBH}} \).
\end{lemma}

\begin{proof}
Let \( \gmb	= \log(2) / \detam \) and \( \xxc \) be given by \eqref{kv7367ehjgruwwcewyde}.
By Theorem~\ref{LHDxibound}, 
\begin{EQA}
	\P \bigl( \| \BBH^{1/2} \VP^{-1} \nabla \zeta \| \geq \zq(\BBH,\xx) \bigr) 
	& \leq &
	3 \ex^{-\xx} \, 
\label{wWeightsde}
\end{EQA}
for all \( \xx \leq \xxc \)%
\ifNL{, and by \eqref{if7h3rhgy4676rfhdsjw}, 
\( \sqrt{\xxc} \geq \gmb/2 - \sqrt{\dimA_{\BBH}/(2\lambda_{\BBH})} \).}{.}
This yields the assertion.
\end{proof}

%

\begin{lemma}
\label{PGLMcond0}
Assume \nameref{Psivdref}. 
Let \( \Upsd \) be from \eqref{txxctde43d4s4srgwcrdf} with \( \rrL \) satisfying \( \detam \, \rrL \leq 1/2 \).
Then for any \( \upsv \in \Upsd \)
\begin{EQA}
	\cdens''\bigl( \langle \Psiv_{i}, \upsv \rangle \bigr)
	& \leq &
	\sqrt{\ex} \,\, \cdens''\bigl( \langle \Psiv_{i}, \upsvd \rangle \bigr) ,
	\qquad
	i=1,\ldots,n,
\label{g3d45de5r3ted5tvuye}
	\\
	\frac{1}{\sqrt{\ex}} \, \IF(\upsvd)
	& \leq &
	\IF(\upsv)
	\leq 
	\sqrt{\ex} \, \IF(\upsvd) .
\label{e12f2Pivf2Piiv12}
\end{EQA}
Also, \eqref{k23CPpHkuf2} holds with \( \deta^{2} \leq \sqrt{\ex} \, \detam^{2} \).
\end{lemma}

\begin{proof}
The function \( \cdens(\ups) = \log(1 + \ex^{\ups}) \) satisfies for all \( \ups \in \R \)
\begin{EQA}
	|\cdens^{(k)}(\ups)|
	& \leq &
	\cdens''(\ups) ,
	\qquad
	k=3,4.
\label{phukulppu43}
\end{EQA}
Indeed, it holds
\begin{EQA}
	\cdens'(\ups)
	&=&
	\frac{\ex^{\ups}}{1 + \ex^{\ups}} \, ,
	\\
	\cdens''(\ups)
	&=&
	\frac{\ex^{\ups}}{(1 + \ex^{\ups})^{2}} \, ,
	\\
	\cdens^{(3)}(\ups)
	&=&
	\frac{\ex^{\ups}}{(1 + \ex^{\ups})^{2}} - \frac{2 \ex^{2\ups}}{(1 + \ex^{\ups})^{3}} \, ,
	\\
	\cdens^{(4)}(\ups)
	&=&
	\frac{\ex^{\ups}}{(1 + \ex^{\ups})^{2}} - \frac{6 \ex^{2\ups}}{(1 + \ex^{\ups})^{3}}
	+ \frac{6 \ex^{3\ups}}{(1 + \ex^{\ups})^{4}} \, .
\label{6314412f2e21e}
\end{EQA}
It is straightforward to see that 
\( | \cdens^{(k)}(\ups)| \leq \cdens''(\ups) \) for \( k = 3,4 \) and any \( \ups \).

Next, we check local variability of \( \cdens''(\ups) \).
Fix \( \upsd < 0 \).
As the function \( \cdens''(\ups) \) is monotonously increasing in \( \ups \), it holds 
\begin{EQA}
	\sup_{|\ups - \upsd| \leq \rhou} \frac{\cdens''(\ups)}{\cdens''(\upsd)}
	& = &
	\frac{\cdens''(\upsd + \rhou)}{\cdens''(\upsd)}
	\leq 
	\ex^{\rhou} \, .
\label{exrhfrfppusg}
\end{EQA}
Putting together \eqref{phukulppu43} and \eqref{exrhfrfppusg} leads to a bound on variability of
\( \IF(\upsv) \) for \( \upsv = \upsvd + \uv \) and \( \| \DPN \uv \| \leq \rrL \).
By definition,
\begin{EQA}
	\IF(\upsv)
	&=&
	\sumi \Psiv_{i} \, \Psiv_{i}^{\T} \, \cdens''\bigl( \langle \Psiv_{i}, \upsv \rangle \bigr) .
\label{0bvlgfiu84ehdf6ew3nh}
\end{EQA}
Now \eqref{diuq3rvyewrhjdsfasd} and \( \detam \, \rrL \leq 1/2 \) imply 
\( \bigl| \langle \Psiv_{i}, \uv \rangle \bigr| 
\leq \| \DPN^{-1} \Psiv_{i} \| \, \| \DPN \uv \| \leq \detam \, \rrL \leq 1/2 \) for each \( i \leq \nsize \),
and \eqref{g3d45de5r3ted5tvuye}, \eqref{e12f2Pivf2Piiv12} follow by \eqref{exrhfrfppusg}.

By definition 
\( \| \IF^{1/2}(\upsv) \zv \|^{2} = \sumi \langle \Psiv_{i}, \zv \rangle^{2} \, \weight_{i}(\upsv) \)
and the use of \eqref{diuq3rvyewrhjdsfasd} yields
\begin{EQA}
	\sumi \langle \Psiv_{i}, \zv \rangle^{4} \, \weight_{i}(\upsv) 
	& = &
	\sumi \langle \Psiv_{i}, \zv \rangle^{2} \, 
	\langle \DPN^{-1} \Psiv_{i}, \DPN \zv \rangle^{2} \, \weight_{i}(\upsv) 
	\\
	& \leq &
	\detam^{2} \| \DPN \zv \|^{2} 
	\sumi \langle \Psiv_{i}, \zv \rangle^{2} \, \weight_{i}(\upsv)
	=
	\detam^{2} \| \DPN \zv \|^{2} \, \| \IF^{1/2}(\upsv) \zv \|^{2} .
\label{iud7ejdufuf6ed6th}
\end{EQA}
Therefore, \( \DPN^{2} = \IF(\upsvd) \geq \ex^{-1/2}\IF(\upsv) \) implies \( \deta^{2} \leq \sqrt{\ex} \, \detam^{2} \).
\end{proof}

Further we check \nameref{LLsT3ref}, \nameref{LLsT4ref} at \( \upsvs \)
with \( \rrL = 3 \rr_{\DPN}/2 \), \( \dltwu_{3} = \sqrt{\ex} \, \deta \), and \( \dltwu_{4} = \sqrt{\ex} \, \deta^{2} \).

\begin{lemma}
\label{PGLMcond}
Assume \nameref{Psivdref}. 
Let 
\( \rrL \) satisfy \( \detam \, \rrL \leq 1/2 \).
Then \nameref{LLsT3ref} and \nameref{LLsT4ref} hold with 
\( \dltwu_{3} = \sqrt{\ex} \, \deta \) and \( \dltwu_{4} = \sqrt{\ex} \, \deta^{2} \) for \( \deta \) from \eqref{k23CPpHkuf2}.
\end{lemma}

\begin{proof}
First we evaluate the derivative \( \nabla^{k} f(\upsv) \) for \( f(\upsv) = \E L(\upsv) \)
over \( \upsv \in \Upsd \).
For any \( \zv \in \R^{\dimp} \), it holds
\begin{EQA}
	\bigl\langle \nabla^{k} f(\upsv), \zv^{\otimes k} \bigr\rangle
	& = &
	\sumi \langle \Psiv_{i}, \zv \rangle^{k} \,\, 
	\cdens^{(k)}\bigl( \langle \Psiv_{i}, \upsv \rangle \bigr) .
\label{Cparfxa2IaGusd}
\end{EQA}
With \( \weight_{i}(\upsv) = \cdens''\bigl( \langle \Psiv_{i}, \upsv \rangle \bigr) \),
we derive by \eqref{diuq3rvyewrhjdsfasd}, \eqref{k23CPpHkuf2}, \eqref{phukulppu43},
and \eqref{e12f2Pivf2Piiv12} 
\begin{EQA}
	&& \nquad
	\bigl| \bigl\langle \nabla^{3} f(\upsv), \zv^{\otimes 3} \bigr\rangle \bigr|
	\leq
	\sumi \bigl| \langle \Psiv_{i}, \zv \rangle \bigr|^{3} \, 
	\cdens''\bigl( \langle \Psiv_{i}, \upsv \rangle \bigr)
	\leq 
	\sqrt{\ex} \sumi \bigl| \langle \Psiv_{i}, \zv \rangle \bigr|^{3} \, \weight_{i}(\upsvd)
	\\
	& \leq &
	\sqrt{\ex} \biggl( 
		\sumi \langle \Psiv_{i}, \zv \rangle^{2} \, \weight_{i}(\upsvd)
	\biggr)^{1/2}
	\biggl( 
		\sumi \langle \Psiv_{i}, \zv \rangle^{4} \, \weight_{i}(\upsvd)
	\biggr)^{1/2}
	\leq 
	\sqrt{\ex} \, \deta \, \| \IF^{1/2}(\upsvd) \zv \| \, \| \DPN \zv \|^{2} 
\label{ph2Piu2Pief2}
\end{EQA}
and \nameref{LLsT3ref} follows with \( \dltwu_{3} = \sqrt{\ex} \, \deta \) by \( \IF(\upsvd) = \DPN^{2} \).
Similarly \nameref{LLsT4ref} holds with \( \dltwu_{4} = \sqrt{\ex} \, \deta^{2} \).
\end{proof}

All the required conditions from Section~\ref{SgenBounds} have been checked 
and the results about the behavior of the pMLE \( \tilde{\upsv}_{\GP} \) apply. 
\end{proof}

\ifAnya{}{
\Section{Logistic regression. Random design}
Suppose that \( \Psiv_{1}, \ldots,\Psiv_{n} \) are i.i.d. random vectors. 
This section shows that with a high probability, 
the empirical design \( \Psiv_{1},\ldots,\Psiv_{n} \) fulfills \eqref{maxiPsiu214DG21} through \eqref{k23CPpHkuf2}.
Fix \( \upsv \) and define \( \Zv_{i} = \sqrt{\cdens''(\langle \Psiv_{i} ,\upsv \rangle)} \,\, \Psiv_{i} \), 
\begin{EQA}
	\HLH^{2}
	=
	\HLH^{2}(\upsv)
	& \eqdef &
	\E \, \bigl\{ \cdens''(\langle \Psiv_{1} ,\upsv \rangle) \, \Psiv_{1} \Psiv_{1}^{\T} \bigr\} 
	=
	\E \, (\Zv_{1} \Zv_{1}^{\T}) \, .
\label{7etfy8iiw2qjwo2f7y3r8y}
\end{EQA}
Let \( \HL^{2} \geq \HLH^{2} \).
Later we use \( \HL^{2} \) given by \( n \, \HL^{2} = n \, \HLH^{2} + \GP^{2} \) leading to
the effective dimension \( \dimG = \tr(\HL^{-2} \HLH^{2}) \).
First, we establish a bound on the maximum of the \( \| \HL^{-1} \Zv_{i} \| \)'s.

\begin{lemma}
\label{LRDlogr}
Assume \( \Psiv_{1}, \ldots,\Psiv_{n} \) i.i.d. sub-gaussian random vectors. Then 
\begin{EQA}
	\P\bigl( \max_{i \leq n} \| \HL^{-1} \Zv_{i} \| \geq \sqrt{\tr (\HL^{-2} \HLH^{2})} + 2 \sqrt{\log n} \bigr)
	& \leq &
	1/n \, .
\label{uyhdfijoweiodjfod}
\end{EQA}
\end{lemma}

\begin{proof}
If \( \Psiv_{i} \) is sub-gaussian then \( \Zv_{i} = \sqrt{\cdens''(\langle \Psiv_{i} ,\upsv \rangle)} \,\, \Psiv_{i} \)
is sub-gaussian as well because \( \cdens''(\ups) \leq 1/4 \) for all \( \ups \).
For each \( i \leq n \), it holds by Theorem~\ref{Tdevboundsharp} for any moderate \( \xx \)
\begin{EQA}
	\P\bigl( \| \HL^{-1} \Zv_{i} \| \geq \sqrt{\tr (\HL^{-2} \HLH^{2})} + \sqrt{2\xx} \bigr)
	& \leq &
	\ex^{-\xx} \, .
\label{iu8weuhriw87y68w2}
\end{EQA}
This implies
\begin{EQA}
	&& \nquad
	\P\Bigl( \max_{i \leq n} \| \HL^{-1} \Zv_{i} \| \geq \sqrt{\tr (\HL^{-2} \HLH^{2})} + \sqrt{4\log n} \Bigr)
	\\
	& \leq &
	\sumi \P\bigl( \| \HL^{-1} \Zv_{i} \| \geq \sqrt{\tr (\HL^{-2} \HLH^{2})} + \sqrt{4\log n} \bigr)
	\leq 
	n \cdot n^{-2} \, 
\label{uyhdfijoweiodjfode}
\end{EQA}
and the assertion follows. 
\end{proof}

The next step is to show that \( n^{-1} \sum_{i} \Zv_{i} \Zv_{i}^{\T} \approx \HLH^{2} \).

\begin{lemma}
\label{Lsnlogr}
Let 
\begin{EQA}
	\sigma^{2}
	& \eqdef &
	\Var(\| \HL^{-1} \Zv_{1} \|^{2})
	=
	\E \| \HL^{-1} \Zv_{1} \|^{4} - \tr (\HL^{-1} \HLH^{2} \HL^{-1})^{2} 
	\ll
	n .
\label{ufc7duje476vbyurgters}
\end{EQA}
It holds for any \( \zz > 0 \) with \( \zz^{2} \sigma^{2} \ll n \)
\begin{EQA}
	\P\biggl( 
		\Bigl\| \HL^{-1} \Bigl( \frac{1}{n} \sum_{i} \Zv_{i} \Zv_{i}^{\T} - \HLH^{2} \Bigr) \HL^{-1} \Bigr\| 
		> \frac{\zz \sigma}{\sqrt{n}} 
	\biggr)
	& \leq &
	\frac{1}{\zz^{2}} \, .
\label{76cjhwe9g7y3hdf6dfyt}
\end{EQA}
\end{lemma}

\begin{proof}
Define
\begin{EQA}[c]
	\Xiv_{n} \eqdef \frac{1}{n} \sum_{i} \xiv_{i} \, ,
	\qquad 
	\xiv_{i} \eqdef \HL^{-1} ( \Zv_{i} \Zv_{i}^{\T} - \HLH^{2} ) \HL^{-1} \, .
\label{7jusdyctr5dfye3d6ther4j}
\end{EQA}
We check that  \( \| \Xiv_{n} \|_{\Fr} \) is small in probability.
As \( \xiv_{i} \) are i.i.d. with \( \E \xiv_{i} = 0 \), it holds 
\begin{EQA}[c]
	\P\biggl( \| \Xiv_{n} \|_{\Fr} > \frac{\zz \sigma}{\sqrt{n}} \biggr)
	\leq 
	\frac{n}{\zz^{2} \sigma^{2}} \E\| \Xiv_{n} \|_{\Fr}^{2}
	=
	\frac{1}{\zz^{2} \sigma^{2}} \E\| \xiv_{1} \|_{\Fr}^{2}
\label{76jeckdc8e8ry8ijhfvik}
\end{EQA}
and
\begin{EQA}
	\E\| \xiv_{1} \|_{\Fr}^{2}
	&=&
	\E \tr\bigl\{ \HL^{-1} ( \Zv_{1} \Zv_{1}^{\T} - \HLH^{2} ) \HL^{-1} \bigr\}^{2}
	=
	\E \tr\bigl\{ (\HL^{-1} \Zv_{1} \Zv_{1}^{\T} \HL^{-1})^{2} - (\HL^{-1} \HLH^{2} \HL^{-1})^{2} \bigr\}
	\\
	&=&
	\E \| \HL^{-1} \Zv_{1} \|^{4} - \tr (\HL^{-1} \HLH^{2} \HL^{-1})^{2}
	=
	\sigma^{2} . 
\label{iy3ri98ujmcos62uynbes}
\end{EQA}
This implies the statement.
\end{proof}

Similarly one can bound the derivatives \( \nabla^{3} f(\upsv) \) and \( \nabla^{4} f(\upsv) \); see \eqref{Cparfxa2IaGusd}.
It holds by \( |\cdens^{(3)}(x)| \leq \cdens^{(2)}(x) \) for any \( \uv \in \R^{\dimp} \)
\begin{EQA}
	\bigl| \bigl\langle \nabla^{3} f(\upsv), \uv^{\otimes 3} \bigr\rangle \bigr|
	& \leq &
	\E \sumi \bigl| \langle \Psiv_{i}, \uv \rangle^{3} \bigr| \,\, 
		\cdens^{(2)}\bigl( \langle \Psiv_{i}, \upsv \rangle \bigr) 
	\leq 
	\max_{i} \bigl| \langle \Psiv_{i}, \uv \rangle \bigr|
	\sumi \langle \Psiv_{i}, \uv \rangle^{2} \,\, 
		\cdens^{(2)}\bigl( \langle \Psiv_{i}, \upsv \rangle \bigr)
	\\
	& \leq &
	\| \HL \uv \| \, \max_{i} \| \HL^{-1} \Psiv_{i} \| \, 
	\sumi \langle \HL^{-1} \Zv_{i}, \HL \uv \rangle^{2} 
	\\
	& \leq &
	\| \HL \uv \|^{3} \, \max_{i} \| \HL^{-1} \Psiv_{i} \| \,\, \biggl\| \HL^{-1} \biggl( \sumi \Zv_{i} \Zv_{i}^{\T} \biggr) \HL^{-1} \biggr\| 
	\, .
\label{Cparfxa2IaGusd3}
\end{EQA}
By \eqref{uyhdfijoweiodjfod} and \eqref{76cjhwe9g7y3hdf6dfyt}, 
on a random set \( \Omega_{\zz} \) with \( \P\bigl( \Omega_{\zz} \bigr) \geq 1 - \zz^{-2} - n^{-1} \), this implies
\begin{EQA}
	n^{-1} \bigl| \bigl\langle \nabla^{3} f(\upsv), \uv^{\otimes 3} \bigr\rangle \bigr|
	& \leq &
	\| \HL \uv \|^{3} (1 + \zz \sigma/\sqrt{n}) \, \Bigl( \sqrt{\tr (\HL^{-2} \HLH^{2})} + 2 \sqrt{\log n} \Bigr) \, .
\label{7jfcu7y4ehfuyf65e3gj}
\end{EQA}
Similarly
\begin{EQA}
	n^{-1} \bigl| \bigl\langle \nabla^{4} f(\upsv), \uv^{\otimes 4} \bigr\rangle \bigr|
	& \leq &
	\| \HL \uv \|^{4} (1 + \zz \sigma/\sqrt{n}) \, \Bigl( \sqrt{\tr (\HL^{-2} \HLH^{2})} + 2 \sqrt{\log n} \Bigr)^{2} \, .
\label{7jfcu7y4ehfuyf65e3gj4}
\end{EQA}
}




\newpage
\appendix

\Chapter{Smooth perturbed optimization}
\label{Slocalsmooth}
This section presents conditions and accurate bounds for smooth perturbed optimization.
We consider minimization of a smooth convex function \( \fs \) 
and study the impact of a linear, quadratic, or smooth perturbation.


%
\Section{Gateaux smoothness and self-concordance}
\label{Ssmoothness}
Below we assume 
the function \( \fs(\upsv) \), \( \upsv \in \Ups \subseteq \R^{\dimp} \) to be strongly convex with a positive definite Hessian 
\( \IFN(\upsv) \eqdef \nabla^{2} \fs(\upsv) \).
Also, assume \( \fs(\upsv) \) three or sometimes even four times Gateaux differentiable in \( \upsv \in \Ups \).
Local smoothness of \( \fs \) will be described by the relative error of the Taylor expansion 
of the third or fourth order.
Define
\begin{EQ}[rcl]
	\dltw_{3}(\upsv,\uv) 
	&=& 
	\fs(\upsv + \uv) - \fs(\upsv) - \langle \nabla \fs(\upsv), \uv \rangle 
	- \frac{1}{2} \langle \nabla^{2} \fs(\upsv), \uv^{\otimes 2} \rangle , 
	\\
	\dltwd_{3}(\upsv,\uv) 
	&=&
	\langle \nabla \fs(\upsv + \uv), \uv \rangle - \langle \nabla \fs(\upsv), \uv \rangle 
	- \langle \nabla^{2} \fs(\upsv), \uv^{\otimes 2} \rangle \, ,
\label{dltw3vufuv12f2ga}
\end{EQ}
and
\begin{EQA}
	\dltw_{4}(\upsv,\uv)
	& = &
	\fs(\upsv + \uv) - \fs(\upsv) - \langle \nabla \fs(\upsv), \uv \rangle 
	- \frac{1}{2} \langle \nabla^{2} \fs(\upsv), \uv^{\otimes 2} \rangle
	- \frac{1}{6} \langle \nabla^{3} \fs(\upsv), \uv^{\otimes 3} \rangle \, .
\label{hvcduywgedfuyg2y1y35e3wweg}
\end{EQA}
Now, for each \( \upsv \), suppose to be given a positive symmetric matrix 
\( \DFN(\upsv) \) 
defining a local metric and a local vicinity around \( \upsv \): for some radius \( \rr \),
\begin{EQA}
	\UVz_{\rr}(\upsv)
	&=&
	\bigl\{ \uv \in \R^{\dimp} \colon \| \DFN(\upsv) \uv \| \leq \rr \bigr\}
\label{ed7sycf7wedwgedq2ftwdfgtv}
\end{EQA}
Local smoothness properties of \( \fs \) at \( \upsv \) are given via the quantities
\begin{EQA}[rcccl]
    \dltwb(\upsv)
    & \eqdef &
    \sup_{\uv \colon \| \DFN(\upsv) \uv \| \leq \rr} \,
    \frac{2|\dltw_{3}(\upsv,\uv)|}{\| \DFN(\upsv) \uv \|^{2}} 
    \,\, ,
    \quad
    \dltwbd(\upsv)
    & \eqdef &
    \sup_{\uv \colon \| \DFN(\upsv) \uv \| \leq \rr} \, \frac{|\dltwd_{3}(\upsv,\uv)|}{\| \DFN(\upsv) \uv \|^{2}} \,\, . 
    \qquad
\label{dtb3u1DG2d3GPg}
\end{EQA}
The definition yields for any \( \uv \) with \( \| \DFN(\upsv) \uv \| \leq \rr \)
\begin{EQ}[rcccl]
	\bigl| \dltw_{3}(\upsv,\uv) \bigr|
	& \leq &
	\frac{\dltwb(\upsv)}{2} \| \DFN(\upsv) \uv \|^{2} 
	\, ,
	\qquad
	\bigl| \dltwd_{3}(\upsv,\uv) \bigr|
	& \leq &
	\dltwbd(\upsv) \| \DFN(\upsv) \uv \|^{2}
	\, .
	\qquad
\label{dta3u1DG2d3GPa1g}
\end{EQ}
%
The approximation results can be improved 
under a third-order upper bound on the error of Taylor expansion. 

\begin{description}
    \item[\label{LL3tref} \( \bb{(\mathcal{T}_{3})} \)]
      \textit{For some \( \dltwu_{3} \)}
\begin{EQA}
	\bigl| \dltw_{3}(\upsv,\uv) \bigr|
	& \leq &
	\frac{\dltwu_{3}}{6} \| \DFN(\upsv) \, \uv \|^{3} \, ,
	\quad
	\bigl| \dltwd_{3}(\upsv,\uv) \bigr|
	\leq 
	\frac{\dltwu_{3}}{2} \| \DFN(\upsv) \, \uv \|^{3} \, ,
	\quad
	\uv \in \UVz_{\rr}(\upsv).
\label{bd3xu16f3uo3st}
\end{EQA}
\end{description}
 
\begin{description}
    \item[\label{LL4tref} \( \bb{(\mathcal{T}_{4})} \)]
      \textit{For some \( \dltwu_{4} \)}
\begin{EQA}
	\bigl| \dltw_{4}(\upsv,\uv) \bigr|
	& \leq &
	\frac{\dltwu_{4}}{24} \| \DFN(\upsv) \, \uv \|^{4} \, ,
	\qquad
	\uv \in \UVz_{\rr}(\upsv).
\label{1mffmxum5st}
\end{EQA}
\end{description}

We also present a version of \nameref{LL3tref} resp. \nameref{LL4tref} in terms of the third (resp. fourth) derivative of \( \fs \).
\begin{description}
    \item[\label{LLsT3ref} \( \bb{(\mathcal{T}_{3}^{*})} \)]
    \emph{\( \fs(\upsv) \) is three times differentiable and 
	}
\begin{EQA}
    \sup_{\uv \colon \| \DFN(\upsv) \uv \| \leq \rr} \,\, \sup_{\zv \in \R^{\dimp}} \,\, 
    \frac{\bigl| \langle \nabla^{3} \fs(\upsv + \uv), \zv^{\otimes 3} \rangle \bigr|}
		 {\| \DFN(\upsv) \zv \|^{3}} 
	& \leq &
	\dltwu_{3} \, .
\label{jcxhydtferyu9j3d6vhew}
\end{EQA}

    \item[\label{LLsT4ref} \( \bb{(\mathcal{T}_{4}^{*})} \)]
    \emph{\( \fs(\upsv) \) is four times differentiable and 
	}
\begin{EQA}
    \sup_{\uv \colon \| \DFN(\upsv) \uv \| \leq \rr} \,\, \sup_{\zv \in \R^{\dimp}} \,\, 
    \frac{\bigl| \langle \nabla^{4} \fs(\upsv + \uv), \zv^{\otimes 4} \rangle \bigr|}
		 {\| \DFN(\upsv) \zv \|^{4}} 
	& \leq &
	\dltwu_{4} \, .
\label{jcxhydtferyu9j3d6vhew4}
\end{EQA}

\end{description}

%
\noindent
By Banach's characterization \cite{Banach1938}, \nameref{LLsT3ref} implies
\begin{EQA}
	\bigl| \langle \nabla^{3} \fs(\upsv + \uv), \zv_{1} \otimes \zv_{2} \otimes \zv_{3} \rangle \bigr|
	& \leq &	 
	\dltwu_{3} \| \DFN(\upsv) \zv_{1} \| \, \| \DFN(\upsv) \zv_{2} \| \, \| \DFN(\upsv) \zv_{3} \| \, 
\label{jbuyfg773jgion94euyyfg}
\end{EQA}
for any \( \uv \) with \( \| \DFN(\upsv) \uv \| \leq \rr \) and all \( \zv_{1} , \zv_{2}, \zv_{3} \in \R^{\dimp} \).
Similarly under \nameref{LLsT4ref}
\begin{EQA}
	&& \hspace{-.7cm}
	\bigl| \langle \nabla^{4} \fs(\upsv + \uv), \zv_{1} \otimes \zv_{2} \otimes \zv_{3} \otimes \zv_{4} \rangle \bigr|
	\\
	& \leq &
	\dltwu_{4} \| \DFN(\upsv) \zv_{1} \| \, \| \DFN(\upsv) \zv_{2} \| \, \| \DFN(\upsv) \zv_{3} \| \, \| \DFN(\upsv) \zv_{4} \| \, ,
	\quad 
	\zv_{1} , \zv_{2}, \zv_{3}, \zv_{4} \in \R^{\dimp} \, .
	\qquad
\label{jbuyfg773jgion94euyyfg4}
\end{EQA}

\begin{lemma}
\label{LdltwLa3t}
Under \nameref{LL3tref} or \nameref{LLsT3ref},
it holds for \( \dltwb(\upsv) \) and \( \dltwbd(\upsv) \) from \eqref{dtb3u1DG2d3GPg}
\begin{EQA}[rcccl]
\label{gtcdsftdffrvsewsea}
	\dltwb(\upsv)
	& \leq &
	\frac{\dltwu_{3} \, \rr}{3 } \, ,
	\qquad
	\dltwbd(\upsv)
	& \leq &
	\frac{\dltwu_{3} \, \rr}{2} 
	\, .
\label{gtcdsftdfvtwdsefhfdvfrvsewseG}
\end{EQA}
\end{lemma}

\begin{proof}
For any \( \uv \in \UVz_{\rr}(\upsv) \) with \( \| \DFN(\upsv) \uv \| \leq \rr \)
\begin{EQA}
	\bigl| \dltw_{3}(\upsv,\uv) \bigr|
	& \leq &
	\frac{\dltwu_{3}}{6} \, \| \DFN(\upsv) \uv \|^{3} 
	\leq 
	\frac{\dltwu_{3} \, \rr}{6} \, \| \DFN(\upsv) \uv \|^{2},
\label{jrgeteteer2234587654}
\end{EQA}
and the bound for \( \dltwb(\upsv) \) follows.
The proof for \( \dltwbd(\upsv) \) is similar.
Under \nameref{LLsT3ref}, apply the third order Taylor expansion to the 
univariate function \( \fs(\upsv + t \uv) \) of \( t \) and 
\nameref{LLsT3ref} with \( \zv \equiv \uv \).
\end{proof}

The values \( \dltwu_{3} \) and \( \dltwu_{4} \) are usually very small.
Some quantitative bounds are given later in this section
under the assumption that the function \( \fs(\upsv) \) can be written in the form \( \fs(\upsv) = n \hL(\upsv) \) 
for a fixed smooth function \( h(\upsv) \) with the Hessian \( \nabla^{2} \hL(\upsv) \). 
The factor \( n \) has meaning of the sample size%
\ifapp{; see \Chname \ref{ScritdimMLE} or \Chname \ref{SGBvM}.}{.}

\begin{description}
    \item[\label{LLtS3ref} \( \bb{(\mathcal{S}_{3}^{*})} \)]
      \emph{ \( \fs(\upsv) = n \hL(\upsv) \) for \( \hL(\upsv) \) three times differentiable and for a metric tensor \( \HL(\upsv) \)
\begin{EQA}
	\sup_{\uv \colon \| \HL(\upsv) \uv \| \leq \rr/\sqrt{n}} 
	\frac{\bigl| \langle \nabla^{3} \hL(\upsv + \uv), \uv^{\otimes 3} \rangle \bigr|}{\| \HL(\upsv) \uv \|^{3}}
	& \leq &
	\hmax_{3} \, .
\end{EQA}
}
    \item[\label{LLtS4ref} \( \bb{(\mathcal{S}_{4}^{*})} \)]
      \emph{ the function \( \hL(\cdot) \) satisfies \nameref{LLtS3ref} and  
\begin{EQA}
	\sup_{\uv \colon \| \HL(\upsv) \uv \| \leq \rr/\sqrt{n}}
	\frac{\bigl| \langle \nabla^{4} \hL(\upsv + \uv), \uv^{\otimes 4} \rangle \bigr|}{\| \HL(\upsv) \uv \|^{4}}
	& \leq &
	\hmax_{4} \, .
\end{EQA}
}
\end{description}

\noindent
\nameref{LLtS3ref} and \nameref{LLtS4ref}
are local versions of the so-called self-concordance condition; see \cite{NeNe1994} and \cite{OsBa2021}.
In fact, they require that each univariate function \( \hL(\upsv + t \uv) \) of \( t \in \R \)
is self-concordant with some universal constants \( \hmax_{3} \) and \( \hmax_{4} \).
Under \nameref{LLtS3ref} and \nameref{LLtS4ref}, with \( \DFN^{2}(\upsv) = n \, \HL^{2}(\upsv) \), the values 
\( \dltw_{3}(\upsv,\uv) \), \( \dltw_{4}(\upsv,\uv) \), and \( \dltwb(\upsv) \), \( \dltwbd(\upsv) \) can be well
bounded.

\begin{lemma}
\label{LdltwLaGP}
Suppose \nameref{LLtS3ref}.
Then 
\nameref{LL3tref} follows with \( \dltwu_{3} = \hmax_{3} n^{-1/2} \).
Moreover, for \( \dltwb(\upsv) \) and \( \dltwbd(\upsv) \) from \eqref{dtb3u1DG2d3GPg}, it holds
\begin{EQA}[rcccl]
	\dltwb(\upsv)
	& \leq &
	\frac{\hmax_{3} \, \rr}{3 n^{1/2}} \, ,
	\qquad
	\dltwbd(\upsv)
	& \leq &
	\frac{\hmax_{3} \, \rr}{2 n^{1/2}} \, .
\label{gtcdsftdfvtwdsefhfdvfrvsewseGP}
\end{EQA}
Also \nameref{LL4tref} follows from \nameref{LLtS4ref} with \( \dltwu_{4} = \hmax_{4} n^{-1} \).
\end{lemma}

\begin{proof}
For any \( \uv \in \UVz_{\rr}(\upsv) \), by third order Taylor expansion,
\begin{EQA}
	|\dltw_{3}(\upsv,\uv)|
	& \leq &
	\sup_{t \in [0,1]}
	\frac{1}{6} \bigl| \langle \nabla^{3} \fs(\upsv + t \uv), \uv^{\otimes 3} \rangle \bigr|
	=
	\frac{n}{6} \, 
	\sup_{t \in [0,1]}
	\bigl| \langle \nabla^{3} \hL(\upsv + t \uv), \uv^{\otimes 3} \rangle \bigr|
	\\
	& \leq &
	\frac{n \, \hmax_{3}}{6} \, \| \HL(\upsv) \uv \|^{3} 
	=
	\frac{n^{-1/2} \, \hmax_{3}}{6} \, \| \DFN(\upsv) \uv \|^{3}
	\leq 
	\frac{n^{-1/2} \, \hmax_{3} \, \rr}{6} \, \| \DFN(\upsv) \uv \|^{2} \, .
\label{jrgeteteer2234587654}
\end{EQA}
This implies \nameref{LL3tref} as well as \eqref{gtcdsftdfvtwdsefhfdvfrvsewseGP}; see \eqref{dta3u1DG2d3GPa1g}.
The statement about \nameref{LL4tref} is similar.
\end{proof}

\def\Bv{\operatorname{B}}
\def\AvmGP{\Avm_{\hspace{-1pt}\GP}}

\Section{Linearly perturbed optimization}
\label{Squadnquad}

Let \( \fs(\upsv) \) be a smooth convex function, 
\begin{EQA}
	\upsvs
	&=&
	\argmin_{\upsv} \fs(\upsv),
	\qquad
	\fs(\upsvs)
	=
	\min_{\upsv} \fs(\upsv) ,
	\qquad
	\IFN = \nabla^{2} \fs(\upsvs) \, .
\label{fg5hg3gf98tkj3dciryt}
\end{EQA}
Later we study the question of how the point of minimum and the value of minimum of \( \fs \) change if we add a linear  
component to \( \fs \).
More precisely, let another function \( \fn(\upsv) \) satisfy for some vector \( \Av \)
\begin{EQA}
	\fn(\upsv) - \fn(\upsvs) 
	&=&
	\bigl\langle \upsv - \upsvs, \Av \bigr\rangle + \fs(\upsv) - \fs(\upsvs) .
\label{4hbh8njoelvt6jwgf09}
\end{EQA}
Define
\begin{EQA}
	\upsvn
	& \eqdef &
	\argmin_{\upsv} \fn(\upsv),
	\qquad
	\fn(\upsvn)
	=
	\min_{\upsv} \fn(\upsv) .
\label{6yc63yhudf7fdy6edgehy} 
\end{EQA}
The aim of the analysis is to evaluate the quantities \( \upsvn - \upsvs \) and
\( \fn(\upsvn) - \fn(\upsvs) \).
First, we consider the case of a quadratic function \( \fs \).

\begin{lemma}
\label{Pquadquad}
Let \( \fs(\upsv) \) be quadratic with \( \nabla^{2} \fs(\upsv) \equiv \IFN \) and let 
\( \fn(\upsv) \) be from \eqref{4hbh8njoelvt6jwgf09}.
Then  
\begin{EQA}
	\upsvn - \upsvs
	&=&
	- \IFN^{-1} \Av,
	\qquad
	\fn(\upsvn) - \fn(\upsvs)
	=
	- \frac{1}{2} \| \IFN^{-1/2} \Av \|^{2} .
\label{kjcjhchdgehydgtdtte35}
\end{EQA}
\end{lemma}

\begin{proof}
If \( \fs(\upsv) \) is quadratic, then under \eqref{4hbh8njoelvt6jwgf09}, \( \fn(\upsv) \) is quadratic as well
with \( \nabla^{2} \fn(\upsv) \equiv \IFN \).
This implies
\begin{EQA}
	\nabla \fn(\upsvs) - \nabla \fn(\upsvn)
	&=&
	\IFN (\upsvs - \upsvn) .
\label{dcudydye67e6dy3wujhds7}
\end{EQA}
Further, \eqref{4hbh8njoelvt6jwgf09} and \( \nabla \fs(\upsvs) = 0 \) yield \( \nabla \fn(\upsvs) = \Av \).
Together with \( \nabla \fn(\upsvn) = 0 \), this implies
\( \upsvn - \upsvs = - \IFN^{-1} \Av \).
The Taylor expansion of \( \fn \) at \( \upsvn \) yields by \( \nabla \fn(\upsvn) = 0 \)
\begin{EQA}
	\fn(\upsvs) - \fn(\upsvn)
	&=&
	\frac{1}{2} \| \IFN^{1/2} (\upsvn - \upsvs) \|^{2}
	=
	\frac{1}{2} \| \IFN^{-1/2} \Av \|^{2} 
\label{8chuctc44wckvcuedje}
\end{EQA}
and the assertion follows.
\end{proof}

\Subsection{{Local concentration}}
Let \( \fs \) satisfy \eqref{dtb3u1DG2d3GPg} at \( \upsvs \) with 
\( \DFN(\upsvs) = \DFN \leq \dmax \, \IFN^{1/2} \) for some \( \dmax > 0 \).
The latter means that the matrix \( \IFN - \dmax^{2} \DFN^{2} \) is positive definite.
The next result describes the concentration properties of \( \upsvn \) from \eqref{6yc63yhudf7fdy6edgehy} in a local elliptic set
\begin{EQA}
	\CA(\rr)
	& \eqdef &
	\{ \upsv \colon \| \IFN^{1/2} (\upsv - \upsvs) \| \leq \rr \} ,
\label{0cudc7e3jfuyvct6eyhgwe}
\end{EQA}
where \( \rr \) is slightly larger than \( \| \IFN^{-1/2} \Av \| \).

\begin{theorem}
\label{Pconcgeneric}
Let \( \fs(\upsv) \) be a strongly convex function with \( \fs(\upsvs) = \min_{\upsv} \fs(\upsv) \)  
and \( \IFN = \nabla^{2} \fs(\upsvs) \).
Let, further, \( \fn(\upsv) \) and \( \fs(\upsv) \) be related by \eqref{4hbh8njoelvt6jwgf09} with some vector \( \Av \).
Fix \( \amax < 1 \) and \( \rrn \) such that \( \| \IFN^{-1/2} \Av \| \leq \amax \, \rrn \).
Suppose now that \( \fs(\upsv) \) satisfy \eqref{dtb3u1DG2d3GPg} for \( \upsv = \upsvs \), 
\( \DFN(\upsvs) = \DFN \leq \dmax \, \IFN^{1/2} \) with some \( \dmax > 0 \) 
and \( \dltwbd \) such that 
\begin{EQA}
	1 - \amax - \dltwbd \dmax^{2}
	& > &
	0 .
\label{rrm23r0ut3ua}
\end{EQA}
Then for \( \upsvn \) from \eqref{6yc63yhudf7fdy6edgehy}, it holds 
\begin{EQA}
	\| \IFN^{1/2} (\upsvn - \upsvs) \|  
	& \leq &
	\rrn \, 
	\quad
	\text{ and }
	\quad
	\| \DFN (\upsvn - \upsvs) \|
	\leq 
	\dmax \, \rrn \, . 
\label{rhDGtuGmusGU0a}
\end{EQA}
\end{theorem}

\begin{proof}
Define \( \DFN_{\dmax} = \dmax^{-1} \DFN \).
Then \( \DFN_{\dmax} \leq \IFN^{1/2} \) and the use of \( \DFN_{\dmax} \) in place of \( \DFN \) 
yields \eqref{dtb3u1DG2d3GPg} with \( \dltwbd \dmax^{2} \) in place of \( \dltwbd \).
This allows us to reduce the proof to \( \dmax = 1 \).
The bound \( \| \IFN^{-1/2} \Av \| \leq \amax \, \rrn \) implies for any \( \uv \)
\begin{EQA}
	\bigl| \langle \Av, \uv \rangle \bigr|
	& = &
	\bigl| \langle \IFN^{-1/2} \Av, \IFN^{1/2} \uv \rangle \bigr|
	\leq 
	\amax \, \rrn \| \IFN^{1/2} \uv \| \, .
\label{LLoDGm1nzua}
\end{EQA}
%
Let \( \upsv \) be a point on the boundary of the set \( \CA(\rrn) \) from \eqref{0cudc7e3jfuyvct6eyhgwe}.
We also write \( \uv = \upsv - \upsvs \).
The idea is to show that the derivative  \( \frac{d}{dt} \fn(\upsvs + t \uv) > 0 \) 
is positive for \( t > 1 \).
Then all the extreme points of \( \fn(\upsv) \) are within \( \CA(\rrn) \).
We use the decomposition
\begin{EQA}
	\fn(\upsvs + \rhot \uv) - \fn(\upsvs)
	&=&
	\langle \Av, \uv \rangle \, \rhot 
	+ \fs(\upsvs + \rhot \uv) - \fs(\upsvs) .
\label{LGtsGtuLGtsa}
\end{EQA}
With \( \fGu(t) = \fs(\upsvs + \rhot \uv) - \fs(\upsvs) + \langle \Av, \uv \rangle \, \rhot \), it holds
\begin{EQA}
	\frac{d}{d \rhot} \fs(\upsvs + \rhot \uv)
	&=&
	- \langle \Av, \uv \rangle + \fGu'(\rhot) .
\label{frddtLtGstua}
\end{EQA}
By definition of \( \upsvs \), it also holds \( \fGu'(0) = \langle \Av, \uv \rangle \).
The identity \( \nabla^{2} \fs(\upsvs) = \IFN \) yields \( \fGu''(0) = \| \IFN^{1/2} \uv \|^{2} \).
Bound \eqref{dta3u1DG2d3GPa1g} implies for \( | \rhot | \leq 1 \)
\begin{EQA}
	\bigl| \fGu'(\rhot) - \fGu'(0) - \rhot \fGu''(0) \bigr|
	& \leq &
	\rhot \, \| \DFN \uv \|^{2} \, \dltwbd \, .
\label{fptfp0fpttfpp13a}
\end{EQA}
For \( \rhot = 1 \), we obtain by \eqref{rrm23r0ut3ua} 
\begin{EQA}
	\fGu'(1) 
	& \geq &
	\langle \Av, \uv \rangle + \| \IFN^{1/2} \uv \|^{2} - \| \DFN \uv \|^{2} \, \dltwbd
	\geq 
	\| \IFN^{1/2} \uv \|^{2} (1 -  \dltwbd - \amax)
	> 0 .
\label{fp1fpp13d3rGa}
\end{EQA}
Moreover, convexity of \( \fGu(\rhot) \) implies that \( \fGu'(\rhot) - \fGu'(0) \) increases in 
\( \rhot \) for \( \rhot > 1 \).
Further, summing up the above derivation yields 
\begin{EQA}
	\frac{d}{dt} \fn(\upsvs + \rhot \uv) \Big|_{\rhot=1}
	& \geq &
	\| \IFN^{1/2} \uv \|^{2} (1 - \amax - \dltwbd)
	> 0 .
\label{ddtLGtstu33a}
\end{EQA}
As \( \frac{d}{d \rhot} \fn(\upsvs + \rhot \uv) \) increases with \( \rhot \geq 1 \) together with 
\( \fGu'(\rhot) \) due to \eqref{frddtLtGstua}, the same applies to all such \( \rhot \).
This implies the assertion.
\end{proof}

\Subsection{{Second-order expansions}}
The result of Theorem~\ref{Pconcgeneric} allows to localize the point \( \upsvn = \argmin_{\upsv} \fn(\upsv) \)
in the local vicinity \( \CA(\rrn) \) of \( \upsvs \).
The use of smoothness properties of \( \fn \) or, equivalently, of \( \fs \), in this vicinity helps to obtain
rather sharp expansions for \( \upsvn - \upsvs \) and for \( \fn(\upsvn) - \fn(\upsvs) \).

\begin{theorem}
\label{PFiWigeneric}
Assume the conditions of Theorem~\ref{Pconcgeneric} with \( \dltwb \leq 1/3 \).
Then
\begin{EQ}[rcl]
    - \frac{\dltwb}{1 - \dmax^{2} \hspm \dltwb} \| \DFN \, \IFN^{-1} \Av \|^{2}
    & \leq &
    2 \fn(\upsvn) - 2 \fn(\upsvs) 
    + \| \IFN^{-1/2} \Av \|^{2}
    \\
    & \leq &
    \frac{\dltwb}{1 + \dmax^{2} \hspm \dltwb} \| \DFN \, \IFN^{-1} \Av \|^{2} \, .
    \qquad 
\label{3d3Af12DGttGa}
\end{EQ}
Also
\begin{EQ}[rcl]
    \| \DFN (\upsvn - \upsvs + \IFN^{-1} \Av) \|
    & \leq &
    \frac{2\sqrt{\dltwb}}{1 - \dmax^{2} \hspm \dltwb} \, \| \DFN \, \IFN^{-1} \Av \| \, ,
    \\
    \| \DFN (\upsvn - \upsvs) \|
    & \leq &
    \frac{1 + 2\sqrt{\dltwb}}{1 - \dmax^{2} \hspm \dltwb} \, \| \DFN \, \IFN^{-1} \Av \| \, .
\label{DGttGtsGDGm13rGa}
\end{EQ}
\end{theorem}

\begin{proof}
As in the proof of Theorem~\ref{Pconcgeneric}, replacing \( \DFN_{\dmax} = \dmax^{-1} \DFN \) 
with \( \DFN \) 
reduces the statement to \( \dmax = 1 \) in view of \( \dmax^{2} \dltwb \DFN^{2} = \dltwb \DFN_{\dmax}^{2} \).
By \eqref{dtb3u1DG2d3GPg}, for any \( \upsv \in \CA(\rrn) \)
\begin{EQA}
	\Bigl| 
		\fs(\upsv) - \fs(\upsvs) - \frac{1}{2} \| \IFN^{1/2} (\upsv - \upsvs) \|^{2} 
	\Bigr|
	& \leq &
	\frac{\dltwb}{2} \| \DFN (\upsv - \upsvs) \|^{2} .
\label{d3GrGELGtsG12}
\end{EQA}
Further, 
\begin{EQA}[rcl]
	&& \nquad
	\fn(\upsv) - \fn(\upsvs) + \frac{1}{2} \| \IFN^{-1/2} \Av \|^{2}
	\\
	&=&
	\bigl\langle \upsv - \upsvs, \Av \bigr\rangle
	+ \fs(\upsv) - \fs(\upsvs) + \frac{1}{2} \| \IFN^{-1/2} \Av \|^{2} 
	\\
	&=&
	\frac{1}{2} \bigl\| \IFN^{1/2} (\upsv - \upsvs) + \IFN^{-1/2} \Av \bigr\|^{2}
	+ \fs(\upsv) - \fs(\upsvs) - \frac{1}{2} \| \IFN^{1/2} (\upsv - \upsvs) \|^{2} .
	\qquad 
\label{12ELGuELusG}
\end{EQA}
As \( \upsvn \in \CA(\rrn) \) and it minimizes \( \fn(\upsv) \), we derive by \eqref{d3GrGELGtsG12}
\begin{EQA}
	&& \nquad
	\fn(\upsvn) - \fn(\upsvs) + \frac{1}{2} \| \IFN^{-1/2} \Av \|^{2}
	=
	\min_{\upsv \in \CA(\rrn)} 
	\Bigl\{ 
		\fn(\upsv) - \fn(\upsvs) + \frac{1}{2} \| \IFN^{-1/2} \Av \|^{2} 
	\Bigr\}
	\\
	& \geq &
	\min_{\upsv \in \CA(\rrn)} 
	\Bigl\{ 
		\frac{1}{2} \bigl\| \IFN^{1/2} (\upsv - \upsvs) + \IFN^{-1/2} \Av \bigr\|^{2} 
		- \frac{\dltwb}{2} \| \DFN (\upsv - \upsvs) \|^{2}
	\Bigr\} .
\label{d3G1212222B} 
\end{EQA}
Denote \( \uv = \IFN^{1/2} (\upsv - \upsvs) \), \( \xiv = \IFN^{-1/2} \Av \), and 
\( \BFN = \IFN^{-1/2} \, \DFN^{2} \, \IFN^{-1/2} \).
As \( \DFN^{2} \leq \IFN \), \( \| \xiv \| < \rr \), and \( \dltwb < 1 \), it holds \( \| \BFN \| \leq 1 \) 
for the operator norm of \( \BFN \) and 
\begin{EQA}
	&& \nquad
	\min_{\upsv \in \CA(\rrn)} \bigl\{ \bigl\| \IFN^{1/2} (\upsv - \upsvs) + \IFN^{-1/2} \Av \bigr\|^{2} 
		- \dltwb \| \DFN (\upsv - \upsvs) \|^{2} \bigr\}
	\\
	&=&
	\min_{\| \uv \| \leq \rr} \bigl\{ \| \uv + \xiv \|^{2} - \dltwb \, \uv^{\T} \BFN \uv \bigr\}
	=
	- \xiv^{\T} \bigl\{ (\Id - \dltwb \, \BFN)^{-1} - \Id \bigr\} \xiv
	\geq 
	- \frac{\dltwb}{1 - \dltwb} \xiv^{\T} \BFN \, \xiv
\label{d7eneyf6g53geygywn}
\end{EQA}
with \( \Id \) being the unit matrix in \( \R^{\dimp} \).
This yields
\begin{EQA}
	\fn(\upsvn) - \fn(\upsvs) + \frac{1}{2} \| \IFN^{-1/2} \Av \|^{2}
	& \geq &
	- \frac{\dltwb}{2(1 - \dltwb)} \| \DFN \, \IFN^{-1} \Av \|^{2} . 
\label{fd3G122B2} 
\end{EQA}
Similarly 
\begin{EQA}
	&& \nquad 
	\fn(\upsvn) - \fn(\upsvs) + \frac{1}{2} \| \IFN^{-1/2} \Av \|^{2}
	\\
	& \leq &
	\min_{\upsv \in \CA(\rrn)} 
	\Bigl\{ 
		\frac{1}{2} \bigl\| \IFN^{1/2} (\upsv - \upsvs) + \IFN^{-1/2} \Av \bigr\|^{2} 
		+ \frac{\dltwb}{2} \| \DFN (\upsv - \upsvs) \|^{2}
	\Bigr\}
	\\
	& \leq &
	\frac{1}{2} \xiv^{\T} \bigl\{ \Id - (\Id + \dltwb \, \BFN)^{-1} \bigr\} \xiv
	\leq 
	\frac{\dltwb }{2(1 + \dltwb)} \, \| \DFN \, \IFN^{-1} \Av \|^{2} . 
	\quad 
\label{fd3G122B2m} 
\end{EQA}
These bounds imply 
\eqref{3d3Af12DGttGa}.

Now we derive similarly to \eqref{12ELGuELusG} that for \( \upsv \in \CA(\rrn) \)
\begin{EQA}
	\fn(\upsv) - \fn(\upsvs) 
	& \geq &
	\bigl\langle \upsv - \upsvs, \Av \bigr\rangle
	+ \frac{1}{2} \| \IFN^{1/2} (\upsv - \upsvs) \|^{2}
	- \frac{\dltwb}{2} \| \DFN (\upsv - \upsvs) \|^{2} .
\label{LGvLGvsGf1d3G2}
\end{EQA}
A particular choice \( \upsv = \upsvn \) yields
\begin{EQA}
	\fn(\upsvn) - \fn(\upsvs) 
	& \geq &
	\bigl\langle \upsvn - \upsvs, \Av \bigr\rangle
	+ \frac{1}{2} \| \IFN^{1/2} (\upsvn - \upsvs) \|^{2}
	- \frac{\dltwb}{2} \| \DFN (\upsvn - \upsvs) \|^{2} .
\label{21GsvtvGDG3G2}
\end{EQA}
Combining this inequality with \eqref{fd3G122B2m} allows to bound
\begin{EQA}
	\bigl\langle \upsvn - \upsvs, \Av \bigr\rangle
	+ \frac{1}{2} \| \IFN^{1/2} (\upsvn - \upsvs) \|^{2}
	- \frac{\dltwb}{2} \| \DFN (\upsvn - \upsvs) \|^{2} 
	& \leq &
	- \frac{1}{2} \xiv^{\T} (\Id + \dltwb \, \BFN)^{-1} \xiv .
\label{2m1DGd3G123G}
\end{EQA}
With 
\( \uvd = \IFN^{1/2} (\upsvn - \upsvs) \), this implies
\begin{EQA}
	2 \bigl\langle \uvd, \xiv \bigr\rangle + {\uvd}^{\T}(1 - \dltwb \BFN) \uvd 
	& \leq &
	- \xiv^{\T} (\Id + \dltwb \, \BFN)^{-1} \xiv \, ,
\label{dtxi2fd1d22}
\end{EQA}
and hence,
\begin{EQA}
	&& \nquad
	\bigl\{ \uvd + (\Id - \dltwb \BFN)^{-1} \xiv \bigr\}^{\T} (\Id - \dltwb \BFN) \bigl\{ \uvd + (\Id - \dltwb \BFN)^{-1} \xiv \bigr\}
	\\
	& \leq &
	\xiv^{\T} \bigl\{ (\Id - \dltwb \, \BFN)^{-1} - (\Id + \dltwb \, \BFN)^{-1} \bigr\} \xiv
	\leq 
	\frac{2 \dltwb}{(1 + \dltwb) (1 - \dltwb)} \, \xiv^{\T} \BFN \, \xiv \, .
\label{uv11wxi22w1w}
\end{EQA}
Introduce \( \| \cdot \|_{\afn} \) by \( \| \xv \|_{\afn}^{2} \eqdef \xv^{\T} (\Id - \dltwb \BFN) \xv \),
\( \xv \in \R^{\dimp} \).
Then
\begin{EQA}
	\| \uvd + (\Id - \dltwb \BFN)^{-1} \xiv \|_{\afn}^{2}
	& \leq &
	\frac{2 \dltwb}{1 - \dltwb^{2}} \, \xiv^{\T} \BFN \, \xiv \, .
\label{dcumwf6vhehe6fbwhfr}
\end{EQA}
As 
\begin{EQA}
	\| \xiv - (\Id - \dltwb \BFN)^{-1} \xiv \|_{\afn}^{2}
	& = &
	\dltwb^{2} (\BFN \xiv)^{\T} (\Id - \dltwb \BFN)^{-1} \BFN \xiv
	\leq 
	\frac{\dltwb^{2}}{1 - \dltwb} \, \| \BFN \xiv \|^{2}
	\leq 
	\frac{\dltwb^{2}}{1 - \dltwb} \, \xiv^{\T} \BFN \, \xiv
\label{c8jjkie74he3tftdy3fy}
\end{EQA}
we conclude for \( \dltwb \leq 1/3 \) by the triangle inequality
\begin{EQA}
	\| \uvd + \xiv \|_{\afn}
	& \leq &
	\biggl( \sqrt{\frac{\dltwb^{2}}{1 - \dltwb}} + \sqrt{\frac{2 \dltwb}{1 - \dltwb^{2}}} \biggr)
	\sqrt{\xiv^{\T} \BFN \, \xiv}
	\leq 
	2 \sqrt{\frac{\dltwb}{1 - \dltwb}} \,\, \sqrt{\xiv^{\T} \BFN \, \xiv} \, ,
\label{uxiBws2w1w31w}
\end{EQA}
and \eqref{DGttGtsGDGm13rGa} follows by \( \Id - \dltwb \BFN \geq (1 - \dltwb) \Id \).
\end{proof}

\begin{remark}
\label{Rfsfnlinp}
The roles of the functions \( \fs \) and \( \fn \) are exchangeable.
In particular, the results from \eqref{DGttGtsGDGm13rGa} apply with
\( \IFN = \nabla^{2} \fn(\upsvn) = \nabla^{2} \fs(\upsvn) \) provided that 
\eqref{dtb3u1DG2d3GPg} is fulfilled at \( \upsv = \upsvn \).
\end{remark}

\Subsection{{Expansions under third order smoothness}}
The results of Theorem~\ref{PFiWigeneric} can be refined if
\( \fs \) satisfies condition \nameref{LL3tref}.

\begin{theorem}
\label{Pconcgeneric2}
Let \( \fs(\upsv) \) be a strongly convex function with \( \fs(\upsvs) = \min_{\upsv} \fs(\upsv) \)  
and \( \IFN = \nabla^{2} \fs(\upsvs) \).
Let \( \fn(\upsv) \) fulfill \eqref{4hbh8njoelvt6jwgf09} with some vector \( \Av \).
Suppose that \( \fs(\upsv) \) follows \nameref{LL3tref} at \( \upsvs \) with 
\( \DFN^{2} \), \( \rrn \), and \( \dltwu_{3} \) such that 
\begin{EQA}
	\DFN^{2} \leq \dmax^{2} \hspm \IFN , 
	\quad \rrn \geq \frac{4\dmax}{3} \, \| \IFN^{-1/2} \Av \| ,
	\quad 
	\dmax^{3} \dltwu_{3} \, \| \IFN^{-1/2} \Av \|
	& < &
	\frac{1}{4} \, .
\label{yxdhewndu7jwnjjuu}
\end{EQA}
Then \( \upsvn = \argmin_{\upsv} \fn(\upsv) \) satisfies 
\begin{EQA}
	\| \IFN^{1/2} (\upsvn - \upsvs) \|  
	& \leq &
	\frac{4}{3} \| \IFN^{-1/2} \Av \| \, ,
	\qquad
	 \| \DFN (\upsvn - \upsvs) \|  
	\leq 
	\frac{4\dmax}{3} \, \| \IFN^{-1/2} \Av \| \,. 
\label{rhDGtuGmusGU0a2}
\end{EQA}
Moreover, 
\begin{EQA}
    \Bigl| 2 \fn(\upsvn) - 2 \fn(\upsvs) + \| \IFN^{-1/2} \Av \|^{2} \Bigr|
    & \leq &
    \frac{\dltwu_{3}}{2} \, \| \DFN \, \IFN^{-1} \Av \|^{3} \, .
    \qquad
\label{3d3Af12DGttGa2}
\end{EQA}
\end{theorem}

\begin{proof}
W.l.o.g. assume \( \dmax = 1 \)
and replace \( \rr \) with \( \frac{4}{3} \| \IFN^{-1/2} \Av \| \).
By \eqref{yxdhewndu7jwnjjuu} \( \dltwu_{3} \, \rr \leq 1/3 \),
and \eqref{yxdhewndu7jwnjjuu} implies \eqref{rrm23r0ut3ua}.
Now \nameref{LL3tref} and Lemma~\ref{LdltwLa3t} ensure \eqref{dtb3u1DG2d3GPg} with 
\( \dltwbd(\upsv) = \dltwu_{3} \, \rr/2 \),
and the first statement follows from Theorem~\ref{Pconcgeneric} with \( \amax = 3/4 \).

As \( \nabla \fs(\upsvs) = 0 \), \nameref{LL3tref} implies for any \( \upsv \in \CA(\rrn) \)
\begin{EQA}
	\Bigl| 
		\fs(\upsv) - \fs(\upsvs) - \frac{1}{2} \| \IFN^{1/2} (\upsv - \upsvs) \|^{2} 
	\Bigr|
	& \leq &
	\frac{\dltwu_{3}}{6} \| \DFN (\upsv - \upsvs) \|^{3}
	\, .
	\qquad \quad
\label{d3GrGELGtsG122}
\end{EQA}
Further, 
\begin{EQA}[rcl]
	&& \hspace{-12pt}
	\fn(\upsv) - \fn(\upsvs) + \frac{1}{2} \| \IFN^{-1/2} \Av \|^{2}
	\\
	&=&
	\bigl\langle \upsv - \upsvs, \Av \bigr\rangle
	+ \fs(\upsv) - \fs(\upsvs) + \frac{1}{2} \| \IFN^{-1/2} \Av \|^{2} 
	\\
	&=&
	\frac{1}{2} \bigl\| \IFN^{1/2} (\upsv - \upsvs) + \IFN^{-1/2} \Av \bigr\|^{2}
	+ \fs(\upsv) - \fs(\upsvs) - \frac{1}{2} \| \IFN^{1/2} (\upsv - \upsvs) \|^{2} .
	\qquad 
\label{12ELGuELusG2}
\end{EQA}
As \( \upsvn \in \CA(\rrn) \) and it minimizes \( \fn(\upsv) \), we derive by \eqref{d3GrGELGtsG122} and Lemma~\ref{Llin23}
with \( \AFN = \IFN^{1/2} \DFN^{-1} \) and \( \afv = \DFN \, \IFN^{-1} \Av \)
\begin{EQA}
	&& \hspace{-12pt}
	2 \fn(\upsvn) - 2 \fn(\upsvs) + \| \IFN^{-1/2} \Av \|^{2}
	=
	\min_{\upsv \in \CA(\rrn)} 
	\Bigl\{ 
		2 \fn(\upsv) - 2 \fn(\upsvs) + \| \IFN^{-1/2} \Av \|^{2} 
	\Bigr\}
	\\
	& \geq &
	\min_{\upsv \in \CA(\rrn)} 
	\Bigl\{ 
		\bigl\| \IFN^{1/2} (\upsv - \upsvs) + \IFN^{-1/2} \Av \bigr\|^{2} 
		- \frac{\dltwu_{3}}{3} \| \DFN (\upsv - \upsvs) \|^{3}
	\Bigr\} 
	\geq 
	- \frac{\dltwu_{3}}{2} \, \| \DFN \, \IFN^{-1} \Av \|^{3} \, .
\label{d3G1212222B} 
\end{EQA}
Similarly 
\begin{EQA}
	&&  
	2 \fn(\upsvn) - 2 \fn(\upsvs) + \| \IFN^{-1/2} \Av \|^{2}
	\\
	&& \quad
	\leq 
	\min_{\upsv \in \CA(\rrn)} 
	\Bigl\{ 
		\bigl\| \IFN^{1/2} (\upsv - \upsvs) + \IFN^{-1/2} \Av \bigr\|^{2} 
		+ \frac{\dltwu_{3}}{3} \| \DFN (\upsv - \upsvs) \|^{3}
	\Bigr\}
	\leq 
	\frac{\dltwu_{3} }{2} \, \| \DFN \, \IFN^{-1} \Av \|^{3} \, . 
	\qquad \quad
\label{fd3G122B2m2} 
\end{EQA}
This implies \eqref{3d3Af12DGttGa2}.
\end{proof}

\def\AFN{U}
\begin{lemma}
\label{Llin23}
Let \( \AFN \geq \Id \).
Fix some \( \rr \) and let \( \afv \in \R^{\dimp} \) satisfy \( (3/4) \rr \leq \| \afv \| \leq \rr \).
If \( \dltwu \, \rr \leq 1/3 \), then
\begin{EQA}
\label{hdcf6wyheuv76e34r35eycv}
	\max_{\| \uv \| \leq \rr} \Bigl( \frac{\dltwu}{3} \| \uv \|^{3} - (\uv - \afv)^{\T} \AFN (\uv - \afv) \Bigr)
	& \leq & 
	\frac{\dltwu}{2} \, \| \afv \|^{3} \, ,
	\\
	\min_{\| \uv \| \leq \rr} \Bigl( \frac{\dltwu}{3} \| \uv \|^{3} + (\uv - \afv)^{\T} \AFN (\uv - \afv) \Bigr)
	& \leq & 
	\frac{\dltwu}{2} \, \| \afv \|^{3} \, .
\label{hdcf6wyheuv76e34r35eycvm}
\end{EQA}
\end{lemma}

\begin{proof}
Replacing \( \| \uv \|^{3} \) with \( \rr \| \uv \|^{2} \) reduces 
the problem to quadratic programming. 
It holds with \( \rho \eqdef \dltwu \rr/3 \) and \( \afv_{\rho} \eqdef (\AFN - \rho \Id)^{-1} \AFN \afv \)
\begin{EQA}
	&& \nquad
	\frac{\dltwu}{3} \| \uv \|^{3} - (\uv - \afv)^{\T} \AFN (\uv - \afv)
	\leq 
	\frac{\dltwu \rr}{3} \| \uv \|^{2} - (\uv - \afv)^{\T} \AFN (\uv - \afv)
	\\
	&=&
	- \uv^{\T} \bigl( \AFN - \rho \Id \bigr) \uv + 2 \uv^{\T} \AFN \afv - \afv^{\T} \AFN \afv
	\\
	&=&
	- (\uv - \afv_{\rho})^{\T} (\AFN - \rho \Id) (\uv - \afv_{\rho}) 
	+ \afv_{\rho}^{\T} (\AFN - \rho \Id) \afv_{\rho} - \afv^{\T} \AFN \afv
	\\
	& \leq &
	\afv^{\T} \bigl\{ \AFN (\AFN - \rho \Id)^{-1} \AFN - \AFN \bigr\} \afv 
	=
	\rho \afv^{\T} \AFN (\AFN - \rho \Id)^{-1} \afv .
\label{f9kht446iffrtednftehy}
\end{EQA}
This implies in view of \( \AFN \geq \Id \), \( \rr \leq (4/3) \| \afv \| \), and 
\( \rho = \dltwu \rrn/3 \leq 1/9 \)
\begin{EQA}
	&& \nquad
	\max_{\| \uv \| \leq \rr} \Bigl( \frac{\dltwu}{3} \| \uv \|^{3} - (\uv - \afv)^{\T} \AFN (\uv - \afv) \Bigr)
	\\
	& \leq &
	\frac{\rho}{1-\rho} \| \afv \|^{2}
	\leq 
	\frac{\dltwu \rr}{3(1-\rho)} \| \afv \|^{2}
	\leq 
	\frac{4\dltwu}{9 (1-\rho)} \| \afv \|^{3}
	\leq 
	\frac{\dltwu}{2} \| \afv \|^{3} \, ,
\label{gydw7guywbudvrgte7yruw}
\end{EQA}
and \eqref{hdcf6wyheuv76e34r35eycv} follows.
Further, \( \dltwu \| \uv \|^{3}/3 \leq \dltwu \rr \| \uv \|^{2}/3 = \rho \| \uv \|^{2} \) 
for \( \| \uv \| \leq \rr \), and
\begin{EQA}
	&& \nquad
	\frac{\dltwu}{3} \| \uv \|^{3} + (\uv - \afv)^{\T} \AFN (\uv - \afv) 
	\leq  
	\rho \| \uv \|^{2} + (\uv - \afv)^{\T} \AFN (\uv - \afv) 
	\, .
\label{7jdfvy5433feugywdjw7krfu}
\end{EQA}
The global minimum of the latter function is attained at \( \uv_{\rho} \eqdef (\AFN + \rho \Id)^{-1} \AFN \afv \).
As \( \| \uv_{\rho} \| \leq \| \afv \| \leq \rr \) and \( (3/4) \rr \leq \| \afv \| \), this implies
\begin{EQA}[rcl]
	&& \nquad
	\min_{\| \uv \| \leq \rr} 
	\Bigl( \rho \| \uv \|^{2} + (\uv - \afv)^{\T} \AFN (\uv - \afv) \Bigr)
	= 
	\frac{\dltwu \rr}{3} \| \uv_{\rho} \|^{2} + (\uv_{\rho} - \afv)^{\T} \AFN (\uv_{\rho} - \afv)
	\\
	& \leq &
	\afv^{\T} \bigl\{ \AFN - \AFN (\AFN + \rho \Id)^{-1} \AFN \bigr\} \afv 
	=
	\rho \afv^{\T} \AFN (\AFN + \rho \Id)^{-1} \afv 
	\leq 
	\rho  \| \afv \|^{2}
	\leq 
	\frac{4 \dltwu}{9}  \| \afv \|^{3} \, ,
\end{EQA}
and \eqref{hdcf6wyheuv76e34r35eycvm} follows as well.
\end{proof}


\Subsection{{Advanced approximation under locally uniform smoothness}}
The bounds of Theorem~\ref{Pconcgeneric2} can be made more accurate if \( \fs \) follows
\nameref{LLsT3ref} and \nameref{LLsT4ref} and one can apply Taylor's expansion around any point close to \( \upsvs \).
In particular, the improved results do not involve the value \( \| \IFN^{-1/2} \Av \| \) which can be large or even 
infinite in some situation; see Section~\ref{Slinquadr}.


\begin{theorem}
\label{PFiWigeneric2}
Let \( \fs(\upsv) \) be a strongly convex function with \( \fs(\upsvs) = \min_{\upsv} \fs(\upsv) \)  
and \( \IFN = \nabla^{2} \fs(\upsvs) \).
Assume \nameref{LLsT3ref} at \( \upsvs \) with \( \DFN^{2} \), \( \rrn \), and \( \dltwu_{3} \) such that 
\begin{EQA}[c]
	\DFN^{2} \leq \dmax^{2} \hspm \IFN ,
	\quad
	\rrn \geq \frac{3}{2} \| \DFN \, \IFN^{-1} \Av \| \, ,
	\quad
	\dmax^{2} \hspm \dltwu_{3} \, \| \DFN \, \IFN^{-1} \Av \| < \frac{4}{9} \, .
\label{8difiyfc54wrboer7bjfr}
\end{EQA}
Then 
\begin{EQA}[rcl]
\label{dvue6554d5rdtehes}
    \| \DFN (\upsvn - \upsvs) \|
    & \leq &
    \frac{3}{2} \| \DFN \, \IFN^{-1} \Av \| \, ,
    \\
    \| \DFN^{-1} \IFN (\upsvn - \upsvs + \IFN^{-1} \Av) \|
    & \leq &
    \frac{3\dltwu_{3}}{4} \| \DFN \, \IFN^{-1} \Av \|^{2} 
	\, .
\label{DGttGtsGDGm13rGa2}
\end{EQA}
\end{theorem}

\begin{proof}
W.l.o.g. assume \( \dmax = 1 \).
If the function \( \fs \) is quadratic and convex with the minimum at \( \upsvs \) then the linearly perturbed function
\( \fn \) is also quadratic and convex with the minimum at \( \upsvr = \upsvs - \IFN^{-1} \Av \).
In general, the point \( \upsvr \) is not the minimizer of \( \fn \), however, it is very close to \( \upsvn \).
We use that \( \nabla \fs(\upsvs) = 0 \) 
and \( \nabla^{2} \fs(\upsvs) = \IFN \).
The main step of the proof is given by the next lemma.

\begin{lemma}
\label{Ldltw4s}
Assume \nameref{LLsT3ref} at \( \upsv \)
and let \( \UVz_{\rr} = \{ \uv \colon \| \DFN \uv \| \leq \rr \} \).
Then 
\begin{EQA}[ccl]
	\bigl\| \DFN^{-1} \bigl\{ \nabla \fs(\upsv + \uv) - \nabla \fs(\upsv) 
	- \langle \nabla^{2} \fs(\upsv), \uv \rangle \bigr\} 
	\bigr\|
	& \leq &
	\frac{\dltwu_{3}}{2} \, \| \DFN \uv \|^{2} \, ,
	\quad
	\uv \in \UVz_{\rr} \, .
	\qquad
\label{y6sdjsdy7erwmcuecuid}
\label{y6sdjsdy7erwmcuecuid2}
\end{EQA}
Also for all \( \uv, \uv_{1} \in \UVz_{\rr} \) 
\begin{EQ}[rcl]
	\bigl\| \DFN^{-1} \bigl\{ \nabla^{2} \fs(\upsv + \uv_{1}) - \nabla^{2} \fs(\upsv + \uv) \bigr\} \DFN^{-1}	\bigr\|
	& \leq &
	\dltwu_{3} \, \| \DFN (\uv_{1} - \uv) \|
\label{jhfuy7f7dfyedye663eh}
\end{EQ}
and
\begin{EQA}
	\bigl\| \DFN^{-1} \bigl\{ \nabla \fs(\upsv + \uv_{1}) - \nabla \fs(\upsv + \uv) - \nabla^{2} \fs(\upsv) (\uv_{1} - \uv) \bigr\}
	\bigr\|
	& \leq &
	\frac{3\dltwu_{3}}{2} \, \| \DFN (\uv_{1} - \uv) \|^{2} \, .
	\qquad
\label{6dcfujcu8ed8edsudyf5tre35}
\end{EQA}
Moreover, under \nameref{LLsT4ref}, for any \( \uv \in \UVz_{\rr} \),
\begin{EQA}
	\bigl\| \DFN^{-1} \bigl\{ \nabla \fs(\upsv + \uv) - \nabla \fs(\upsv) 
	- \langle \nabla^{2} \fs(\upsv), \uv \rangle 
	- \frac{1}{2} \langle \nabla^{3} \fs(\upsv), \uv^{\otimes 2} \rangle \bigr\} \bigr\|
	& \leq &
	\frac{\dltwu_{4}}{6} \, \| \DFN \uv \|^{3} \, .
	\qquad
\label{y6sdjsdy7erwmcuecuid4}
\end{EQA}
\end{lemma}

\begin{proof}
Fix \( \uv \in \UVz_{\rr} \) and any vector \( \wv \in \R^{\dimp} \).
For \( t \in [0,1] \), denote
\begin{EQA}
	h(t)
	& \eqdef &
	\bigl\langle \nabla \fs(\upsv + t\uv) - \nabla \fs(\upsv) 
	- t \langle \nabla^{2} \fs(\upsv), \uv \rangle, \, \wv \bigl\rangle
	\, .
\label{7dcjef6chjfer7v54etghf}
\end{EQA}
Then \( h(0) = 0 \), \( h'(0) = 0 \), and 
\nameref{LLsT3ref} and \eqref{jbuyfg773jgion94euyyfg} imply 
\begin{EQA}[rcl]
	|h''(t)|
	&=&
	\bigl| \langle \nabla^{3} \fs(\upsv + t \uv), \uv^{\otimes 2} \otimes \wv \rangle \bigr|
	\leq 
	\dltwu_{3} \| \DFN \uv \|^{2} \, \| \DFN \wv \|
	\, .
\end{EQA}
With \( \av \eqdef \nabla \fs(\upsv + \uv) - \nabla \fs(\upsv) 
	- \langle \nabla^{2} \fs(\upsv), \uv \rangle \), this yields
\begin{EQA}
	|\langle \av,\wv \rangle| 
	&=& 
	|h(1)| \leq \frac{\dltwu_{3}}{2} \| \DFN \uv \|^{2} \, \| \DFN \wv \|
	\, ,
	\\
	\| \DFN^{-1} \av \|
	&=&
	\sup_{\| \wv \| = 1} \bigl| \langle \DFN^{-1} \av,\wv \rangle \bigr|
	=
	\sup_{\| \wv \| = 1} \bigl| \langle \av,\DFN^{-1} \wv \rangle \bigr|
	\leq 
	\frac{\dltwu_{3}}{2} \, \| \DFN \uv \|^{2} 
	\, ,
\label{udfjd6vhfwe36vneo}
\end{EQA}
and the first statement follows. 
For \eqref{y6sdjsdy7erwmcuecuid4}, apply
\begin{EQA}
	\av
	& \eqdef &
	\nabla \fs(\upsv + \uv) - \nabla \fs(\upsv) 
	- \langle \nabla^{2} \fs(\upsv), \uv \rangle 
	- \frac{1}{2} \langle \nabla^{3} \fs(\upsv), \uv^{\otimes 2} \rangle
	\, ,
	\\
	h(t)
	& \eqdef &
	\bigl\langle \nabla \fs(\upsv + t\uv) - \nabla \fs(\upsv) 
	- t \langle \nabla^{2} \fs(\upsv), \uv \rangle 
	- \frac{t^{2}}{2} \langle \nabla^{3} \fs(\upsv), \uv^{\otimes 2} \rangle, \wv \bigl\rangle 
	\, ,
\label{t6dtwsghwesyfyghe322w2w}
\end{EQA} 
and use \nameref{LLsT4ref} and \eqref{jbuyfg773jgion94euyyfg4} instead of \nameref{LLsT3ref} and \eqref{jbuyfg773jgion94euyyfg}.

Further, with \( \Bv_{1} \eqdef \nabla^{2} \fs(\upsv + \uv_{1}) - \nabla^{2} \fs(\upsv + \uv) \) and \( \Delta = \uv_{1} - \uv \), 
by \nameref{LLsT3ref}, for any \( \wv \in \R^{\dimp} \) and some \( t \in [0,1] \),
\begin{EQA}
	&& \nquad
	\bigl| \langle \DFN^{-1} \bigl\{ \nabla^{2} \fs(\upsv + \uv_{1}) - \nabla^{2} \fs(\upsv + \uv) \bigr\} \, \DFN^{-1}, \wv^{\otimes 2} \rangle \bigr|
	=
	\bigl| \langle \Bv_{1}, (\DFN^{-1} \wv)^{\otimes 2} \rangle \bigr|
	\\
	& = &
	\bigl| \bigl\langle \nabla^{3} \fs(\upsv + \uv + t \Delta), \Delta \otimes (\DFN^{-1} \wv)^{\otimes 2} \bigr\rangle \bigr|
	\leq 
	\dltwu_{3} \| \DFN \Delta \| \, \| \wv \|^{2}
	\, .
\label{d76jnwef8j2egtftftyjr5}
\end{EQA}
This proves \eqref{jhfuy7f7dfyedye663eh}.
Similarly, 
for some \( t \in [0,1] \)
\begin{EQA}
	&& \nquad
	\bigl| \bigl\langle 
		\DFN^{-1} \bigl\{ \nabla \fs(\upsv + \uv_{1}) - \nabla \fs(\upsv + \uv) - \nabla^{2} \fs(\upsv + \uv) \Delta \bigr\},\wv  
	\bigr\rangle \bigr|
	\\
	&=&
	\frac{1}{2} \bigl| \bigl\langle \nabla^{3} \fs(\upsv + \uv + t \Delta), \Delta \otimes \Delta \otimes \DFN^{-1} \wv \bigr\rangle
	\bigr|
	\leq 
	\frac{\dltwu_{3}}{2} \| \DFN \Delta \|^{2} \, \| \wv \|
\label{ufhfyt3hbgvbigy4jfiu}
\end{EQA}
and with \( \Bv = \nabla^{2} \fs(\upsv + \uv) - \nabla^{2} \fs(\upsv) \), by \eqref{jhfuy7f7dfyedye663eh},
\begin{EQA}
	\bigl\| \DFN^{-1} \Bv \Delta \bigr\| 
	& \leq &
	\| \DFN^{-1} \Bv \, \DFN^{-1} \| \,\, \| \DFN \Delta \| 
	\leq 
	\dltwu_{3} \| \DFN \Delta \|^{2} \, .
\label{u8ifke0gfjw23gfsd4gy}
\end{EQA}
This completes the proof of \eqref{6dcfujcu8ed8edsudyf5tre35}.
\end{proof}

Now we prove \eqref{DGttGtsGDGm13rGa2}.
W.l.o.g. assume \( \| \DFN \, \IFN^{-1} \Av \| = 2 \rrn/3 \).
Lemma~\ref{Ldltw4s}, \eqref{y6sdjsdy7erwmcuecuid}, applied with \( \upsv = \upsvs \) and \( \uv = \IFN^{-1} \Av \) 
yields for \( \upsvr = \upsvs - \IFN^{-1} \Av \)
\begin{EQA}
	\bigl\| \DFN^{-1} \nabla \fn(\upsvr) \bigr\|
	&=&
	\bigl\| \DFN^{-1} \{ \nabla \fs(\upsvs - \IFN^{-1} \Av) - \nabla \fs(\upsvs) - \Av \} \bigr\| 
	\leq
	\frac{\dltwu_{3}}{2} \| \DFN \, \IFN^{-1} \Av \|^{2} \, .
	\qquad
\label{stedsyteuhyfnmgu3}
\end{EQA}
As \( \| \DFN \, \IFN^{-1} \Av \| = 2 \rrn/3 \), 
condition \nameref{LLsT3ref} can be applied 
in the \( \rrn/3 \)-vicinity of \( \upsvr \).
We show that \( \fn(\upsv) \) attains its minimum within this vicinity.
Fix any \( \upsv \) on its boundary, i.e. with \( \| \DFN (\upsv - \upsvr) \| = \rrn/3 \) and denote 
\( \Delta = \upsv - \upsvr \).
By \eqref{6dcfujcu8ed8edsudyf5tre35} 
\begin{EQA}[c]
	\bigl\| \DFN^{-1} \{ \nabla \fn(\upsv) - \nabla \fn(\upsvr) - \IFN \Delta \} \bigr\|
	=
	\bigl\| \DFN^{-1} \{ \nabla \fs(\upsv) - \nabla \fs(\upsvr) - \IFN \Delta \} \bigr\|
	\leq 
	\frac{3 \dltwu_{3}}{2} \| \DFN \Delta \|^{2} \, .
\label{7djdxhes5ewtwee6e6eed}
\end{EQA}
In particular, this and \eqref{stedsyteuhyfnmgu3} yield
\begin{EQA}
	\bigl\| \DFN^{-1} \{ \nabla \fn(\upsvr + \Delta) - \IFN \Delta \} \bigr\|
	& \leq &
	2 \dltwu_{3} \| \DFN \Delta \|^{2} \, .
\label{7hewjfuv7bh7tur4jwsfycn}
\end{EQA}
For any \( \uv \) with \( \| \uv \| = 1 \), this implies
\begin{EQA}
	\bigl| \bigl\langle \nabla \fn(\upsvr + \Delta) - \IFN \Delta , \DFN^{-1} \uv \bigr\rangle \bigr|
	& \leq &
	2 \dltwu_{3} \| \DFN \Delta \|^{2} \, .
\label{hf9jmw2f7vhehe46fghfnwh}
\end{EQA}
Now consider the function \( h(t) = \fn(\upsvr + t \Delta) \).
Then  
\( h'(t) = \langle \nabla \fn(\upsvr + t \Delta), \Delta \rangle \)
and \eqref{hf9jmw2f7vhehe46fghfnwh} implies with \( \uv = \DFN \Delta/\| \DFN \Delta \| \)
\begin{EQA}
	\bigl| \langle \nabla \fn(\upsvr + \Delta), \Delta \rangle - \| \IFN^{1/2} \Delta \|^{2} \bigr|
	& \leq &
	2 \dltwu_{3} \| \DFN \Delta \|^{3} \, .
\label{ufudhdyf5d53egru7e4u3}
\end{EQA}
As \( \IFN \geq \DFN^{2} \), this yields
\begin{EQA}
	h'(1)
	& \geq &
	\| \DFN \Delta \|^{2} - 2 \dltwu_{3} \| \DFN \Delta \|^{3}  .
\label{67vhyfdhcyeyghevy2hegtie3}
\end{EQA}
Similarly \( - h'(-1) \geq \| \DFN \Delta \|^{2} - 2 \dltwu_{3} \| \DFN \Delta \|^{3} \)
and \eqref{stedsyteuhyfnmgu3} yields by \( \| \DFN \, \IFN^{-1} \Av \| = \frac{2}{3}\rrn \)
\begin{EQA}
	|h'(0)|
	=
	\bigl| \langle \nabla \fn(\upsvr), \Delta \rangle \bigr|
	& \leq &
	\frac{\dltwu_{3}}{2} \| \DFN \, \IFN^{-1} \Av \|^{2} \, \| \DFN \Delta \| 
	\leq 
	\frac{2\dltwu_{3}}{9} \, \rr^{2} \, \| \DFN \Delta \|\, .
\label{yhev7h3bfgreewewweddh2}
\end{EQA}
Due to \eqref{67vhyfdhcyeyghevy2hegtie3}, \eqref{yhev7h3bfgreewewweddh2}, 
\( \| \DFN \Delta \| = \rrn/3 \),  \( \dltwu_{3} \| \DFN \, \IFN^{-1} \Av \| \leq 4/9 \), 
and \( \| \DFN \, \IFN^{-1} \Av \| = 2\rrn/3 \)
\begin{EQA}[rcl]
	h'(1) - |h'(0)|
	& \geq &
	\| \DFN \Delta \|^{2} - \frac{2\dltwu_{3}}{9} \rr^{2} \, \| \DFN \Delta \| - 2 \dltwu_{3} \| \DFN \Delta \|^{3} 
	\\
	&=&
	\| \DFN \Delta \| \, \rrn \Bigl( \frac{1}{3} - \frac{2\dltwu_{3} \, \rrn}{9} - \frac{2 \dltwu_{3} \, \rrn}{9} \Bigr)
	>
	0 
	\, .
\label{thewuvf7ehehctdebenvjyw}
\end{EQA}
Similarly \( - h'(-1) > |h'(0)| \), and 
convexity of \( \fn(\cdot) \) ensures that \( t^{*} = \argmin_{t} h(t) \) satisfies \( |t^{*}| \leq 1 \).
We summarize that \( \upsvn = \argmin_{\upsv} \fn(\upsv) \) satisfies \( \| \DFN \, (\upsvn - \upsvr) \| \leq \rrn/3 \)
while \( \| \DFN (\upsvr - \upsvs) \| = \| \DFN \, \IFN^{-1} \Av \| = 2 \rrn/3 \),
thus yielding \eqref{dvue6554d5rdtehes}.

We can now apply \nameref{LLsT3ref} at \( \upsvn \) for checking \eqref{DGttGtsGDGm13rGa2}.
As \( \nabla \fn(\upsvn) = 0 \), by \eqref{stedsyteuhyfnmgu3}
\begin{EQA}
	\| \DFN^{-1} \{ \nabla \fn(\upsvs - \IFN^{-1} \Av) - \nabla \fn(\upsvn) \} \|
	& \leq &
	\frac{\dltwu_{3}}{2} \| \DFN \, \IFN^{-1} \Av \|^{2} \, .
\label{fhy5345etvty46dgw3}
\end{EQA}
By \eqref{6dcfujcu8ed8edsudyf5tre35} of Lemma~\ref{Ldltw4s}, it holds with \( \Delta = \upsvs - \IFN^{-1} \Av - \upsvn \)
\begin{EQA}
	\bigl\| \DFN^{-1} \{ \nabla \fn(\upsvs - \IFN^{-1} \Av) - \nabla \fn(\upsvn) - \nabla^{2} \fn(\upsvs) \Delta \} \bigr\| 
	& \leq &
	\frac{3 \dltwu_{3}}{2} \| \DFN \Delta \|^{2} \, .
\label{ygtdtdt55636ytv2wg3}
\end{EQA}
Combining with \eqref{fhy5345etvty46dgw3} yields 
\begin{EQA}
	\| \DFN^{-1} \IFN \Delta \|
	& \leq &
	\frac{3 \dltwu_{3}}{2} \| \DFN \Delta \|^{2} + \frac{\dltwu_{3}}{2} \| \DFN \, \IFN^{-1} \Av \|^{2} 
	\leq 
	\frac{3 \dltwu_{3}}{2} \| \DFN^{-1} \IFN \Delta \|^{2} + \frac{\dltwu_{3}}{2} \| \DFN \, \IFN^{-1} \Av \|^{2} \, . 
\label{f7f7fv7f66e6ehb3wyc3}
\end{EQA}
As \( 2x \leq \alpha x^{2} + \beta \) with \( \alpha = 3 \dltwu_{3} \), \( \beta = \dltwu_{3} \| \DFN \, \IFN^{-1} \Av \|^{2} \), and
\( x = \| \DFN^{-1} \IFN \Delta \| \in (0,1/\alpha) \) implies \( x \leq \beta/(2 - \alpha\beta) \),
this yields
\begin{EQA}
	\| \DFN^{-1} \IFN (\upsvn - \upsvs + \IFN^{-1} \Av) \|
	& \leq &
	\frac{\dltwu_{3}}{2 - 3 \dltwu_{3}^{2} \| \DFN \, \IFN^{-1} \Av \|^{2}} \| \DFN \, \IFN^{-1} \Av \|^{2} 
\label{fiufu7df56rgyhvnghbvt3}
\end{EQA}
and \eqref{DGttGtsGDGm13rGa2} follows by 
\( \dltwu_{3} \| \DFN \, \IFN^{-1} \Av \| \leq 4/9 \).
\end{proof}

\begin{remark}
\label{Rbiasgeneric}
As in Remark~\ref{dtb3u1DG2d3GPg}, \( \fs \) and \( \fn \) can be exchanged.
In particular, \eqref{DGttGtsGDGm13rGa2} applies with \( \IFN = \IFN(\upsvn) \) provided that 
\nameref{LLsT3ref} is fulfilled at \( \upsvn \).
\end{remark}

\Subsection{{Fourth-order expansions}}
Under fourth-order condition \nameref{LLsT4ref},  expansion \eqref{DGttGtsGDGm13rGa2} can further be refined.  

\begin{theorem}
\label{Pconcgeneric4}
Let \( \fs(\upsv) \) be a strongly convex function 
satisfying 
\nameref{LLsT3ref} and \nameref{LLsT4ref} at \( \upsvs = \argmin_{\upsv} \fs(\upsv) \) with some 
\( \DFN^{2} \), \( \dltwu_{3} \), \( \dltwu_{4} \), and \( \rrn \) such that 
\begin{EQA}[c]
	\DFN^{2} \leq \dmax^{2} \hspm \IFN \, ,
	\;\;
	\rrn \geq \frac{3}{2} \| \DFN \, \IFN^{-1} \Av \| \, ,
	\;\;
	\dmax^{2} \hspm \dltwu_{3} \| \DFN \, \IFN^{-1} \Av \| < \frac{4}{9} \, ,
	\;\;
	\dmax^{2} \hspm \dltwu_{4} \| \DFN \, \IFN^{-1} \Av \|^{2} < \frac{1}{3} \, .
	\qquad
\label{8difiyfc54wrbosT4}
\end{EQA}
Let \( \fn(\upsv) \) fulfill \eqref{4hbh8njoelvt6jwgf09} with some vector \( \Av \) and \( \fn(\upsvn) = \min_{\upsv} \fn(\upsv) \).
Then \( \| \DFN (\upsvn - \upsvs) \| \leq (3/2) \| \DFN \, \IFN^{-1} \Av \| \). 
Further, define
\begin{EQA}
	\avn 
	&=&
	- \IFN^{-1} \{ \Av + \nabla \Tens(\IFN^{-1} \Av) \} \, ,
\label{8vfjvr43223efryfuweef}
\end{EQA}
where \( \Tens(\uv) = \frac{1}{6} \langle \nabla^{3} \fs(\upsvs), \uv^{\otimes 3} \rangle \) for \( \uv \in \R^{\dimp} \).
Then
\begin{EQ}[rcl]
    \| \DFN^{-1} \IFN (\upsvn - \upsvs - \avn) \|
    & \leq &
    (\dltwu_{4}/2 + \dmax^{2} \hspm \dltwu_{3}^{2}) \, \| \DFN \, \IFN^{-1} \Av \|^{3} \, .
\label{DGttGtsGDGm13rGa4}
\end{EQ}
Also
\begin{EQA}
    && \nquad
    \Bigl| \fn(\upsvn) - \fn(\upsvs) + \frac{1}{2} \| \IFN^{-1/2} \Av \|^{2} + \Tens(\IFN^{-1} \Av) \Bigr|
    \\
    & \leq &
    \frac{\dltwu_{4} + 4 \dmax^{2} \hspm \dltwu_{3}^{2}}{8} \| \DFN \, \IFN^{-1} \Av \|^{4} 
    + \frac{\dmax^{2} \, (\dltwu_{4} + 2 \dmax^{2} \hspm \dltwu_{3}^{2})^{2} }{4} \, \| \DFN \, \IFN^{-1} \Av \|^{6} \, 
    \qquad
\label{3d3Af12DGttGa4}
\end{EQA}
and 
\begin{EQA}
	\bigl| \Tens(\IFN^{-1} \Av) \bigr|
	& \leq &
	\frac{\dltwu_{3}}{6} \| \DFN \, \IFN^{-1} \Av \|^{3} \, .
\label{u7jcc7e45hfiobvioeye6yhy}
\end{EQA}
\end{theorem}

\begin{proof}
W.l.o.g. assume \( \dmax = 1 \) and \( \upsvs = 0 \).
Theorem~\ref{PFiWigeneric2} yields \eqref{DGttGtsGDGm13rGa2}.
Later we use that \( \Tens(- \uv) = - \Tens(\uv) \) while \( \nabla \Tens(- \uv) = \nabla \Tens(\uv) \).
By \nameref{LLsT3ref} and \eqref{jbuyfg773jgion94euyyfg}
\begin{EQA}
	&& \nquad
	\| \DFN^{-1} \, \IFN (\avn + \IFN^{-1} \Av) \|
	=
	\| \DFN^{-1} \, \nabla \Tens(\IFN^{-1} \Av) \|
	\\
	&=&
	\sup_{\| \uv \| = 1} 3 \bigl| \langle \Tens, \IFN^{-1} \Av \otimes \IFN^{-1} \Av \otimes \DFN^{-1} \uv \rangle \bigr|
	\leq 
	\frac{\dltwu_{3}}{2} \| \DFN \, \IFN^{-1} \Av \|^{2} \, .
	\qquad
\label{bhvfwfdsdxexsdwsvwe33a}
\end{EQA}
As \( \DFN^{-1} \, \IFN \geq \IFN^{1/2} \geq \DFN \), this implies by \( \dltwu_{3} \| \DFN \, \IFN^{-1} \Av \| \leq 4/9 \)
\begin{EQA}[rcl]
	\| \DFN \avn \|
	& \leq &
	\| \DFN \, \IFN^{-1} \Av \| + \| \DFN \, \IFN^{-1}  \, \nabla \Tens(\IFN^{-1} \Av) \|
	\\
	& \leq &
	\Bigl( 1 + \frac{\dltwu_{3}}{2} \| \DFN \, \IFN^{-1} \Av \| \Bigr) \| \DFN \, \IFN^{-1} \Av \|
	\leq 
	\frac{11}{9} \, \| \DFN \, \IFN^{-1} \Av \| 
	\qquad
\label{iuvchycvf6e64rygh322}
\end{EQA}
and
\begin{EQA}
	\| \IFN^{1/2} \avn + \IFN^{-1/2} \Av \|
	& \leq &
	\frac{\dltwu_{3}}{2} \| \DFN \, \IFN^{-1} \Av \|^{2} \, .
\label{iuvchycvf6e64ryghF}
\end{EQA}
Next, again by \nameref{LLsT3ref}, for any \( \wv \)
\begin{EQA}
	\| \DFN^{-1} \, \nabla^{2} \Tens(\wv) \, \DFN^{-1} \|
	&=&
	\sup_{\| \uv \| = 1} 6 \bigl| \langle \Tens, \wv \otimes (\DFN^{-1} \uv)^{\otimes 2} \rangle \bigr|
	\leq 
	\dltwu_{3} \| \DFN \wv \| \, .
\label{dunj3df7y76e4hf743j}
\end{EQA}
The tensor \( \nabla^{2} \Tens(\uv) \) is linear in \( \uv \), hence
\( \| \nabla^{2} \Tens(\uv) \| \) is convex in \( \uv \) and
\begin{EQA}
	&& \nquad
	\sup_{t \in [0,1]} \| \DFN^{-1} \, \nabla^{2} \Tens(t \avn - (1-t) \IFN^{-1} \Av) \, \DFN^{-1} \| 
	\\
	&=&
	\max\{ \| \DFN^{-1} \, \nabla^{2} \Tens(- \IFN^{-1} \Av) \, \DFN^{-1} \|, \| \DFN^{-1} \nabla^{2} \Tens(\avn) \DFN^{-1} \| \}
	\\
	& \leq &
	\dltwu_{3} \, \max\{ \| \DFN \, \IFN^{-1} \Av \|, \| \DFN \avn \| \} \, .
\label{huyd76hj3fyt7yfj4e}
\end{EQA}
Based on \eqref{iuvchycvf6e64rygh322}, assume \( \| \DFN \, \IFN^{-1} \Av \| \leq \| \DFN \avn \| \leq (11/9) \| \DFN \, \IFN^{-1} \Av \| \).
Then \eqref{bhvfwfdsdxexsdwsvwe33a} yields by \( \nabla \Tens(\uv) = \nabla \Tens(-\uv) \)
\begin{EQA}
	&& \nquad
	\| \DFN^{-1} \nabla \Tens(\avn) - \DFN^{-1} \nabla \Tens(\IFN^{-1} \Av) \|
	=
	\| \DFN^{-1} \nabla \Tens(\avn) - \DFN^{-1} \nabla \Tens(- \IFN^{-1} \Av) \|
	\\
	& \leq &
	\sup_{t \in [0,1]} \| \DFN^{-1} \, \nabla^{2} \Tens(t \avn - (1-t) \IFN^{-1} \Av) \, \DFN^{-1} \| \,\, 
	\| \DFN \, \IFN^{-1} (\avn + \IFN^{-1} \Av) \|
	\\
	& \leq &
	\frac{\dltwu_{3}^{2}}{2} \, \| \DFN \, \IFN^{-1} \Av \|^{2} \, \| \DFN \avn \|
	\leq 
	\frac{2\dltwu_{3}^{2}}{3} \, \| \DFN \, \IFN^{-1} \Av \|^{3}\, .
\label{ukjikio3278eu7grt64rsa}
\end{EQA}
As \( \nabla^{2} \fs(0) = \IFN \) and
\( \nabla \Tens(\avn) = \frac{1}{2} \langle \nabla^{3} \fs(0),\avn \otimes \avn \rangle \),
by \eqref{y6sdjsdy7erwmcuecuid4} with  \( \upsv = 0 \) and \( \uv = \av \) 
and by \eqref{iuvchycvf6e64rygh322} 
\begin{EQA}
	&& \nquad
	\bigl\| \DFN^{-1} \{ \nabla \fs(\avn) - \IFN \avn - \nabla \Tens(\avn) \} \bigr\| 
	\\
	& \leq &
	\frac{\dltwu_{4}}{6} \| \DFN \avn \|^{3} 
	\leq 
	\frac{(11/9)^{3}\dltwu_{4}}{6} \| \DFN \, \IFN^{-1} \Av \|^{3}
	\leq 
	\frac{\dltwu_{4}}{3} \| \DFN \, \IFN^{-1} \Av \|^{3} \, .
\label{stedsyteuhyfnmgu4}
\end{EQA}
Next we bound \( \bigl\| \DFN^{-1} \{ \nabla \fn(\avn) - \nabla \fn(\upsvn) \} \bigr\| \).
As \( \nabla \fn(\upsvn) = 0 \), \eqref{4hbh8njoelvt6jwgf09} and \eqref{8vfjvr43223efryfuweef} imply 
\begin{EQA}[rcl]
	&& \nquad
	\bigl\| \DFN^{-1} \{ \nabla \fn(\avn) - \nabla \fn(\upsvn) \} \bigr\|
	=
	\bigl\| \DFN^{-1} \nabla \fn(\avn) \bigr\|
	=
	\bigl\| \DFN^{-1} \{ \nabla \fn(\avn) - \IFN \avn - \nabla \Tens(\IFN^{-1} \Av) - \Av \} \bigr\|
	\\
	& \leq &
	\bigl\| \DFN^{-1} \{ \nabla \fs(\avn) - \IFN \avn - \nabla \Tens(\avn) \} \bigr\| 
	+ \| \DFN^{-1} \{ \nabla \Tens(\avn) - \nabla \Tens(\IFN^{-1} \Av) \} \|
	\leq 
	\err_{1} \, ,
\label{fhy5345etvty46dgw35w3a}
\end{EQA}
where 
\begin{EQA}[c]
	\err_{1}
	\eqdef 
	\frac{\dltwu_{4} + 2 \dltwu_{3}^{2}}{3} \| \DFN \, \IFN^{-1} \Av \|^{3} 
\label{fhjvmvgt4judjewtrwgd}
\end{EQA}
and by \eqref{8difiyfc54wrbosT4} 
\begin{EQA}
	3 \dltwu_{3} \, \err_{1}
	&=&
	\dltwu_{3} \| \DFN \, \IFN^{-1} \Av \| 
	\Bigl( 
		\dltwu_{4} \| \DFN \, \IFN^{-1} \Av \|^{2} + 2 \dltwu_{3}^{2} \| \DFN \, \IFN^{-1} \Av \|^{2} 
	\Bigr)
	\leq 
	\frac{4}{9} \Bigl( \frac{1}{3} + \frac{32}{81} \Bigr)
	< 
	\frac{1}{3} \, .
	\qquad
\label{fhjvmvgt4judjewtrwg}
\end{EQA}
Further, \( \nabla^{2} \fn(0) = \nabla^{2} \fs(0) = \IFN \), 
and \eqref{6dcfujcu8ed8edsudyf5tre35} of Lemma~\ref{Ldltw4s} implies 
\begin{EQA}
	&& \nquad
	\bigl\| \DFN^{-1} \{ \nabla \fn(\avn) - \nabla \fn(\upsvn) - \IFN (\avn - \upsvn) \} \bigr\|
	\\
	&=&
	\bigl\| \DFN^{-1} \{ \nabla \fs(\avn) - \nabla \fs(\upsvn) - \IFN (\avn - \upsvn) \} \bigr\|
	\leq 
	\frac{3\dltwu_{3}}{2} \| \DFN (\avn - \upsvn) \|^{2} .
\label{ygtdtdt55636ytv2wgdfa}
\end{EQA}
Combining with \eqref{fhy5345etvty46dgw35w3a} yields in view of \( \DFN^{2} \leq \IFN \)
\begin{EQA}
	\| \DFN^{-1} \IFN (\avn - \upsvn) \|
	& \leq &
	\frac{3\dltwu_{3}}{2} \| \DFN (\avn - \upsvn) \|^{2} + \err_{1}
	\leq 
	\frac{3\dltwu_{3}}{2} \| \DFN^{-1} \IFN (\avn - \upsvn) \|^{2} + \err_{1} \, .
\label{ufjvfchyeghdftdf67dejh}
\end{EQA}
As \( 2x \leq \alpha x^{2} + \beta \) with \( \alpha = 3 \dltwu_{3} \), \( \beta = 2 \err_{1} \), and
\( x \in (0,1/\alpha) \) implies \( x \leq \beta/(2 - \alpha\beta) \), 
we conclude by \eqref{fhjvmvgt4judjewtrwg}
\begin{EQA}
	\| \DFN^{-1} \IFN (\avn - \upsvn) \|
	& \leq &
	\frac{\err_{1}}{1 - 3 \dltwu_{3} \,\err_{1}}  
	\leq 
	\frac{\dltwu_{4} + 2 \dltwu_{3}^{2}}{2} \| \DFN \, \IFN^{-1} \Av \|^{3} \, ,
\label{fiufu7df56rgyhvnghbvtwea}
\end{EQA}
and \eqref{DGttGtsGDGm13rGa4} follows.

Next we bound \( \fn(\upsvn) - \fn(0) \).
By \eqref{iuvchycvf6e64ryghF} and \( \DFN^{2} \leq \IFN \)
\begin{EQA}
	\frac{1}{2} \| \IFN^{-1/2} \Av \|^{2} + \langle \Av, \avn \rangle + \frac{1}{2} \| \IFN^{1/2} \avn \|^{2}
	&=&
	\frac{1}{2} \| \IFN^{1/2} \avn + \IFN^{-1/2} \Av \|^{2}
	\leq 
	\frac{\dltwu_{3}^{2}}{8} \| \DFN \, \IFN^{-1} \Av \|^{4} \, .
\label{7djdfuyv7gh7tur5tuy}
\end{EQA}
This together with \( \nabla \fs(0) = 0 \), \( \nabla^{2} \fs(0) = \IFN \), \nameref{LLsT4ref}, 
and \eqref{iuvchycvf6e64rygh322} implies
\begin{EQA}
	&& \nquad
	\Bigl| \fn(\avn) - \fn(0) + \frac{1}{2} \| \IFN^{-1/2} \Av \|^{2} - \Tens(\avn) \Bigr|
	\\
	&=&
	\Bigl| \fs(\avn) - \fs(0) + \langle \Av, \avn \rangle + \frac{1}{2} \| \IFN^{-1/2} \Av \|^{2} - \Tens(\avn) \Bigr|
	\\
	& \leq &
	\Bigl| \fs(\avn) - \fs(0) - \frac{1}{2} \| \IFN^{1/2} \avn \|^{2} - \Tens(\avn) \Bigr| 
	+ \frac{\dltwu_{3}^{2}}{8} \| \DFN \, \IFN^{-1} \Av \|^{4}
	\\
	& \leq &
	\frac{\dltwu_{4}}{24} \| \DFN \, \avn \|^{4} + \frac{\dltwu_{3}^{2}}{8} \| \DFN \, \IFN^{-1} \Av \|^{4} 
	\leq 
	\Bigl( \frac{\dltwu_{4}}{10} + \frac{\dltwu_{3}^{2}}{8} \Bigr) \| \DFN \, \IFN^{-1} \Av \|^{4} \, .
\label{ghd6w2hehfyttet2hbdgfwa}
\end{EQA}
Further, by \( \nabla \fn(\upsvn) = 0 \) and \( \nabla^{2} \fn(\cdot) \equiv \nabla^{2} \fs(\cdot) \), it holds for some 
\( \upsv \in [\avn,\upsvn] \) 
\begin{EQA}
	2 \bigl| \fn(\avn) - \fn(\upsvn) \bigr|
	& = &
	\bigl| \langle \nabla^{2} \fs(\upsv) , (\avn - \upsvn)^{\otimes 2} \rangle \bigr| \, .
\label{dy6eh3hft636yhg3ffa}
\end{EQA}
As \( \| \DFN \avn \| \leq \rrn = \frac{3}{2} \| \DFN \, \IFN^{-1} \Av \| \) and \( \| \DFN \upsvn \| \leq \rrn \), also 
\( \| \DFN \upsv \| \leq \rrn \).
The use of \( \nabla^{2} \fs(0) = \IFN \geq \DFN^{2} \) and \eqref{jhfuy7f7dfyedye663eh} yields by 
\( \dltwu_{3} \| \DFN \, \IFN^{-1} \Av \| < \frac{4}{9} \) and \eqref{fiufu7df56rgyhvnghbvtwea}
\begin{EQA}
	&& \nquad
	2 \bigl| \fn(\avn) - \fn(\upsvn) \bigr|
	\leq 
	\| \IFN^{1/2} (\avn - \upsvn) \|^{2}
	+ \bigl| \bigl\langle \nabla^{2} \fs(\upsv) - \nabla^{2} \fs(0), (\avn - \upsvn)^{\otimes 2} \bigr\rangle \bigr|	
	\\
	& \leq &
	(1 + \dltwu_{3} \rr) \| \IFN^{1/2} (\avn - \upsvn) \|^{2} 
	\leq 
	\frac{(5/3) (\dltwu_{4} + 2 \dltwu_{3}^{2})^{2}}{4} \, \| \DFN \, \IFN^{-1} \Av \|^{6}
	\, .
\label{f8mn4erf6ffyruyn4e3u8fhw3g}
\end{EQA}
As \( \Tens(\avn) = - \Tens(- \avn) \), it holds with \( \Delta \eqdef \IFN^{-1} \nabla \Tens(\IFN^{-1} \Av) \)
for some \( t \in [0,1] \)
\begin{EQA}
	&& \nquad
	\bigl| \Tens(\avn) + \Tens(\IFN^{-1} \Av) \bigr|
	=
	\bigl| \Tens(\IFN^{-1} \Av + \Delta) - \Tens(\IFN^{-1} \Av) \bigr|
	=
	\bigl| \bigl\langle \nabla \Tens(\IFN^{-1} \Av + t \Delta ), \Delta \bigr\rangle \bigr| \,\, 
	\\
	& \leq &
	\frac{\dltwu_{3}}{2} \| \DFN (\IFN^{-1} \Av + t \Delta) \|^{2} \, \| \DFN \Delta \|
	=
	\frac{\dltwu_{3}}{2} \| \DFN \, \IFN^{-1} \Av + t \, \DFN \Delta \|^{2} \, \| \DFN \Delta \|
	\, .
\label{7jc6hw3f6hw3fuvne8dgwx}
\end{EQA}
Similarly to \eqref{bhvfwfdsdxexsdwsvwe33a}, it holds
\( \| \DFN \Delta \| \leq \| \DFN^{-1} \nabla \Tens(\IFN^{-1} \Av) \| \leq (\dltwu_{3}/2) \| \DFN \, \IFN^{-1} \Av \|^{2} \), and
by \( \dltwu_{3} \| \DFN \, \IFN^{-1} \Av \| \leq 1/2 \)
\begin{EQA}
	&& \nquad
	\bigl| \Tens(\avn) + \Tens(\IFN^{-1} \Av) \bigr|
	\leq 
	\frac{(5/4)^{2} \dltwu_{3}^{2} }{4} \| \DFN \, \IFN^{-1} \Av \|^{4} \, .
\label{dc7hhbejrfweugdf7weuhduw}
\end{EQA}
Summing up the obtained bounds yields \eqref{3d3Af12DGttGa4}.
\eqref{u7jcc7e45hfiobvioeye6yhy} follows from \nameref{LLsT3ref}.
\end{proof}

\Section{Quadratic penalization}
\label{Slinquadr}
Here we discuss the case when \( \fn(\upsv) - \fs(\upsv) \) is quadratic.
The general case can be reduced to the situation with \( \fn(\upsv) = \fs(\upsv) + \| \GP \upsv \|^{2}/2 \).
To make the dependence of \( \GP \) more explicit, denote 
\( \fG(\upsv) \eqdef \fs(\upsv) + \| \GP \upsv \|^{2}/2 \),
\begin{EQA}
	\upsvs 
	&=& 
	\argmin_{\upsv} \fs(\upsv),
	\quad
	\upsvs_{\GP} = \argmin_{\upsv} \fG(\upsv) 
	=
	\argmin_{\upsv} \bigl\{ \fs(\upsv) + \| \GP \upsv \|^{2}/2 \bigr\}.
	\qquad
\label{8cvkfc9fujf6jnmcer4cd}
\end{EQA}
We study the bias \( \upsvs_{\GP} - \upsvs \) induced by this penalization.
To get some intuition, consider first the case of a quadratic function \( \fs(\upsv) \).

\begin{lemma}
\label{Lbiasquadgen}
Let \( \fs(\upsv) \) be quadratic with \( \IFN \equiv \nabla^{2} \fs(\upsv) \).
Denote \( \IFN_{\GP} = \IFN + \GP^{2} \).
Then \( \upsvs_{\GP} \) from \eqref{8cvkfc9fujf6jnmcer4cd} satisfies
\begin{EQA}[rcccl]
	\upsvs_{\GP} - \upsvs
	&=&
	- \IFN_{\GP}^{-1} \GP^{2} \upsvs,
	\qquad
	\fG(\upsvs_{\GP}) - \fG(\upsvs)
	&=&
	- \frac{1}{2} \| \IFN_{\GP}^{-1/2} \GP^{2} \upsvs \|^{2} \, .
\label{kjfydf554dfwertdgwdf}
\end{EQA}
\end{lemma} 

\begin{proof}
By definition, \( \fG(\upsv) \) is quadratic with \( \nabla^{2} \fG(\upsv) \equiv \IFN_{\GP} \) and
\begin{EQA}
	\nabla \fG(\upsvs_{\GP}) - \nabla \fG(\upsvs)
	&=&
	\IFN_{\GP} \, (\upsvs_{\GP} - \upsvs) .
\label{dcudydye67e6dy3wujhdqu}
\end{EQA}
Further, \( \nabla \fs(\upsvs) = 0 \) yielding \( \nabla \fG(\upsvs) = \GP^{2} \upsvs \).
Together with \( \nabla \fG(\upsvs_{\GP}) = 0 \), this implies
\( \upsvs_{\GP} - \upsvs = - \IFN_{\GP}^{-1} \GP^{2} \upsvs \).
The Taylor expansion of \( \fG \) at \( \upsvs_{\GP} \) yields 
\begin{EQA}
	\fG(\upsvs) - \fG(\upsvs_{\GP})
	&=&
	\frac{1}{2} \| \IFN_{\GP}^{1/2} (\upsvs - \upsvs_{\GP}) \|^{2}
	=
	\frac{1}{2} \| \IFN_{\GP}^{-1/2} \GP^{2} \upsvs \|^{2} 
\label{8chuctc44wckvcuedjequ}
\end{EQA}
and the assertion follows.
\end{proof}

Now we turn to the general case with \( \fs \) satisfying \nameref{LLsT3ref}.
Define
\begin{EQA}
	\IFN_{\GP}
	\eqdef
	\nabla^{2} \fG(\upsvs),
	\quad
	\bias_{\GP}
	& \eqdef & 
	\| \DFN \, \IFN_{\GP}^{-1} \GP^{2} \upsvs  \| 
	\, .
\label{fd9dfhy4ye6fuydfrerf}
\end{EQA}

\begin{theorem}
\label{Pbiasgeneric} 
Let \( \fG(\upsv) = \fs(\upsv) + \| \GP \upsv \|^{2}/2 \) be strongly convex
and follow \nameref{LLsT3ref} with some \( \DFN^{2} \), \( \dltwu_{3} \), and \( \rrn \) satisfying for \( \dmax > 0 \)
\begin{EQA}[c]
	\DFN^{2} \leq \dmax^{2} \hspm \IFN_{\GP} \, ,
	\qquad 
	\rrn \geq 3 \bias_{\GP}/2 \, ,
	\qquad
	\dmax^{2} \hspm \dltwu_{3} \, \bias_{\GP} < 4/9 
	\, .
\label{7fjgjgvuvy44erd52f}
\end{EQA}
Then 
\begin{EQA}[c]
	\| \DFN (\upsvs_{\GP} - \upsvs) \| \leq 3 \bias_{\GP}/2 .
\label{odf6fdyr6e4deuewjug}
\end{EQA}
Moreover, 
\begin{EQA}[rcl]
	\bigl\| \DFN^{-1} \IFN_{\GP} (\upsvs_{\GP} - \upsvs + \IFN_{\GP}^{-1} \GP^{2} \upsvs) \bigr\|
	& \leq &
	\frac{3\dltwu_{3}}{4} \, \bias_{\GP}^{2}
	\, ,
\label{11ma3eaelDebbgen}
	\\
	\Bigl| 2 \fG(\upsvs_{\GP}) - 2 \fG(\upsvs) + \| \IFN_{\GP}^{-1/2} \GP^{2} \upsvs \|^{2} \Bigr|
	& \leq &
	\frac{\dltwu_{3}}{2} \, \bias_{\GP}^{3}
	\, .
\label{7ywsjhd7wjhjdbiui84kje}
\end{EQA}
\end{theorem}

\begin{proof}
Define \( \fn_{\GP}(\upsv) \) by
\begin{EQA}
	\fn_{\GP}(\upsv) - \fn_{\GP}(\upsvs_{\GP})
	& = &
	\fG(\upsv) - \fG(\upsvs_{\GP}) - \langle \GP^{2} \upsvs, \upsv - \upsvs_{\GP} \rangle .
\label{f7nhw3dghydfte5w35}
\end{EQA}
The function \( \fG \) is convex, the same holds for \( \fn_{\GP} \) from \eqref{f7nhw3dghydfte5w35}.
Moreover, \( \nabla \fG(\upsvs) = \GP^{2} \upsvs \) yields \( \nabla \fn_{\GP}(\upsvs) = - \GP^{2} \upsvs + \GP^{2} \upsvs = 0 \).
Hence, \( \upsvs = \argmin \fn_{\GP}(\upsv) \) and \( \fG(\upsv) \) is a linear perturbation \eqref{4hbh8njoelvt6jwgf09} 
of \( \fn_{\GP} \) with \( \Av = \GP^{2} \upsvs \).
Now the results follow from Theorems~\ref{PFiWigeneric2} and \eqref{3d3Af12DGttGa2} of Theorem~\ref{Pconcgeneric2} applied 
with \( \fs(\upsv) = \fG(\upsv) - \langle \Av,\upsv \rangle \),
\( \fn(\upsv) = \fG(\upsv) \), \( \Av = \GP^{2} \upsvs \), and \( \nabla^{2} \fs(\upsvs) = \IFN_{\GP} \).
\end{proof}


The bound can be further improved under fourth-order smoothness of \( \fs \) 
using the results of Theorem~\ref{Pconcgeneric4}.

\begin{theorem}
\label{Pbiasgeneric4} 
Let \( \fs(\upsv) \) be strongly convex and \( \upsvs = \argmin_{\upsv} \fs(\upsv) \).
Let also
\( \fs(\upsv) \) follow \nameref{LLsT3ref} and \nameref{LLsT4ref} with \( \IFN_{\GP} = \nabla^{2} \fs(\upsvs) + \GP^{2} \) and some 
\( \DFN^{2} \), \( \dltwu_{3} \), \( \dltwu_{4} \), and \( \rrn \) satisfying 
\begin{EQA}[c]
	\DFN^{2} \leq \dmax^{2} \hspm \IFN_{\GP} \, ,
	\quad
	\rrn = \frac{3}{2} \bias_{\GP} \, ,
	\quad
	\dmax^{2} \hspm \dltwu_{3} \, \bias_{\GP} < \frac{4}{9} \, ,
	\quad
	\dmax^{2} \hspm \dltwu_{4} \, \bias_{\GP}^{2} < \frac{1}{3} \, 
	\qquad
\label{8difiyfc54wrbosT4G}
\end{EQA}
for \( \bias_{\GP} \) from \eqref{fd9dfhy4ye6fuydfrerf}.
Then \eqref{odf6fdyr6e4deuewjug} holds.
Furthermore, 
define
\begin{EQA}
	\bvn_{\GP}
	&=&
	- \IFN_{\GP}^{-1} \{ \GP^{2} \upsvs + \nabla \Tens(\IFN_{\GP}^{-1} \GP^{2} \upsvs) \} \, 
\label{8vfjvr43223efryfuweefG}
\end{EQA}
with \( \Tens(\uv) = \frac{1}{6} \langle \nabla^{3} \fs(\upsvs), \uv^{\otimes 3} \rangle \) and
\( \nabla \Tens = \frac{1}{2} \langle \nabla^{3} \fs(\upsvs), \uv^{\otimes 2} \rangle \).
Then 
\begin{EQA}[c]
	\| \DFN (\bvn_{\GP} - \IFN_{\GP}^{-1} \GP^{2} \upsvs) \|
	\leq 
	\frac{\dltwu_{3}}{2} \bias_{\GP}^{2} 
	\, ,
\label{iuvchycvf6e64rygh322b}
\end{EQA}
and
\begin{EQA}
    \| \DFN^{-1} \IFN_{\GP} (\upsvs_{\GP} - \upsvs - \bvn_{\GP}) \|
    & \leq &
    \frac{\dltwu_{4} + 2 \dmax^{2} \hspm \dltwu_{3}^{2}}{2} \, \bias_{\GP}^{3} \, ,
\label{DGttGtsGDGm13rGa4b}
    \\
    \Bigl| 
    	\fG(\upsvs_{\GP}) - \fG(\upsvs) 
		&+& \frac{1}{2} \| \IFN_{\GP}^{-1/2} \GP^{2} \upsvs \|^{2} 
		+ \Tens(\IFN_{\GP}^{-1} \GP^{2} \upsvs) 
    \Bigr|
    \\
    & \leq &
    \frac{\dltwu_{4} + 4 \dmax^{2} \hspm \dltwu_{3}^{2}}{8} \bias_{\GP}^{4} 
    + \frac{\dmax^{2} \, (\dltwu_{4} + 2 \dmax^{2} \hspm \dltwu_{3}^{2})^{2} }{4} \, \bias_{\GP}^{6} \, .
    \qquad
\label{3d3Af12DGttGa4b}
\end{EQA}
\end{theorem}

\begin{proof}
We apply Theorem~\ref{Pconcgeneric4} and use that for \( \avn \) from \eqref{8vfjvr43223efryfuweef}, it holds
\( \avn = \bvn_{\GP} \).
Also use that
\( \nabla^{3} \fs(\upsvs) = \nabla^{3} \fG(\upsvs) = \nabla^{3} \fn_{\GP}(\upsvs) \).
\end{proof}

\Section{A smooth penalty}
\label{Slinsmooth}
The case of a general smooth penalty \( \pent_{\GP}(\upsv) \) can be studied similarly to the quadratic case.
Denote 
\( \fG(\upsv) \eqdef \fs(\upsv) + \pent_{\GP}(\upsv) \),
\begin{EQA}
	\upsvs 
	&=& 
	\argmin_{\upsv} \fs(\upsv),
	\quad
	\upsvs_{\GP} = \argmin_{\upsv} \fG(\upsv) 
	=
	\argmin_{\upsv} \bigl\{ \fs(\upsv) + \pent_{\GP}(\upsv) \bigr\}.
\label{8cvkfc9fujf6jnmcer4pen}
\end{EQA}
We study the bias \( \upsvs_{\GP} - \upsvs \) induced by this penalization.
The statement of Theorem~\ref{Pbiasgeneric} and its proof can be easily extended to this situation.
Define
\begin{EQA}
	\IFN_{\GP}
	\eqdef
	\nabla^{2} \fG(\upsvs),
	\quad
	\bias_{\GP}
	& \eqdef & 
	\| \DFN \, \IFN_{\GP}^{-1} \AvmGP  \| \, ,
	\quad
	\AvmGP  \eqdef \nabla \pent_{\GP}(\upsvs) \, .
\label{fd9dfhy4ye6fuydfrerfG}
\end{EQA} 


\begin{theorem}
\label{Pbiaspen} 
Let \( \fG(\upsv) = \fs(\upsv) + \pent_{\GP}(\upsv) \) be strongly convex
and follow \nameref{LLsT3ref} at \( \upsvs \)
with some \( \DFN^{2} \), \( \dltwu_{3} \), and \( \rrn \) satisfying for \( \dmax > 0 \)
\begin{EQA}[c]
	\DFN^{2} \leq \dmax^{2} \hspm \IFN_{\GP} \, ,
	\qquad 
	\rrn \geq 3 \bias_{\GP}/2 \, ,
	\qquad
	\dmax^{2} \hspm \dltwu_{3} \, \bias_{\GP} < 4/9 ,
\label{7fjgjgvuvy44erd52fpen}
\end{EQA}
Then 
\begin{EQA}[c]
	\| \DFN (\upsvs_{\GP} - \upsvs) \| \leq 3 \bias_{\GP}/2 .
\label{odf6fdyr6e4deuewjugG}
\end{EQA}
Moreover, 
\begin{EQA}[rcl]
	\bigl\| \DFN^{-1} \IFN_{\GP} (\upsvs_{\GP} - \upsvs + \IFN_{\GP}^{-1} \AvmGP ) \bigr\|
	& \leq &
	\frac{3\dltwu_{3}}{4} \, \bias_{\GP}^{2}
	\, ,
\label{11ma3eaelDebbpen}
	\\
	\Bigl| 2 \fG(\upsvs_{\GP}) - 2 \fG(\upsvs) + \| \IFN_{\GP}^{-1/2} \AvmGP  \|^{2} \Bigr|
	& \leq &
	\frac{\dltwu_{3}}{2} \, \bias_{\GP}^{3}
	\, .
\label{7ywsjhd7wjhjdbiui84pen}
\end{EQA}
If, in addition, \( \fG(\upsv) \) satisfies \nameref{LLsT4ref} and \( \dmax^{2} \hspm \dltwu_{4} \, \bias_{\GP}^{2} < \frac{1}{3} \), 
then with \( \Tens_{\GP}(\uv) = \frac{1}{6} \langle \nabla^{3} \fG(\upsvs), \uv^{\otimes 3} \rangle \),
\( \nabla \Tens_{\GP} = \frac{1}{2} \langle \nabla^{3} \fG(\upsvs), \uv^{\otimes 2} \rangle \), and
\begin{EQA}
	\bvn_{\GP}
	&=&
	- \IFN_{\GP}^{-1} \{ \AvmGP  + \nabla \Tens_{\GP}(\IFN_{\GP}^{-1} \AvmGP ) \} \, ,
\label{8vfjvr43223efryfuweefG}
\end{EQA}
it holds
\begin{EQA}
    \| \DFN^{-1} \IFN_{\GP} (\upsvs_{\GP} - \upsvs - \bvn_{\GP}) \|
    & \leq &
    \frac{\dltwu_{4} + 2 \dmax^{2} \hspm \dltwu_{3}^{2}}{2} \, \bias_{\GP}^{3} \, ,
\label{DGttGtsGDGm13rGa4bpen}
    \\
    \Bigl| \fG(\upsvs_{\GP}) - \fG(\upsvs) 
    &+& \frac{1}{2} \| \IFN_{\GP}^{-1/2} \AvmGP  \|^{2} + \Tens_{\GP}(\IFN_{\GP}^{-1} \AvmGP ) \Bigr|
    \\
    & \leq &
    \frac{\dltwu_{4} + 4 \dmax^{2} \hspm \dltwu_{3}^{2}}{8} \bias_{\GP}^{4} 
    + \frac{\dmax^{2} \, (\dltwu_{4} + 2 \dmax^{2} \hspm \dltwu_{3}^{2})^{2} }{4} \, \bias_{\GP}^{6} \, .
    \qquad
\label{3d3Af12DGttGa4bpen}
\end{EQA}
\end{theorem}

\begin{proof}
Consider \( \fn_{\GP}(\upsv) = \fG(\upsv) - \langle \nabla \pent_{\GP}(\upsvs), \upsv \rangle \).
Then \( \fn_{\GP} \) is strongly convex and \( \nabla \fn_{\GP}(\upsvs) = 0 \) yielding
\( \upsvs = \argmin_{\upsv} \fn_{\GP}(\upsv) \).
Also, \( \fG(\upsv) \) is a linear perturbation of \( \fn_{\GP}(\upsv) \) with 
\( \Av = \AvmGP  = \nabla \pent_{\GP}(\upsvs) \).
Now all the statements of Theorem~\ref{Pbiasgeneric} and Theorem~\ref{Pbiasgeneric4} apply to \( \upsvs_{\GP} \)
with obvious changes.
\end{proof}

\Chapter{Deviation bounds for quadratic forms}
\label{Sdevboundgen}
Here we collect some useful results from probability theory mainly concerning 
Gaussian and non-Gaussian quadratic forms.


\def\gaussB{\gaussv_{\BBH}}

\Section{Moments of a Gaussian quadratic form}
\label{SmomentqfG}
Let \( \gaussv \) be standard normal in \( \R^{\dimp} \) for \( \dimp \leq \infty \).
Given a self-adjoint trace operator \( \BBH \), consider a quadratic form 
\( \bigl\langle \BBH \gaussv, \gaussv \bigr\rangle \).

\begin{lemma}
\label{Gaussmoments}
It holds \( \E \bigl\langle \BBH \gaussv, \gaussv \bigr\rangle = \tr \BBH \).
Moreover, 
\begin{EQA}
	\E \bigl( \bigl\langle \BBH \gaussv, \gaussv \bigr\rangle - \tr \BBH \bigr)^{2}
	&=&
	2 \tr \BBH^{2}  ,
	\\
	\E \bigl( \bigl\langle \BBH \gaussv, \gaussv \bigr\rangle - \tr \BBH \bigr)^{3}
	&=&
	8 \tr \BBH^{3} ,
	\\
	\E \bigl( \bigl\langle \BBH \gaussv, \gaussv \bigr\rangle - \tr \BBH \bigr)^{4}
	&=&
	48 \tr \BBH^{4} + 12 (\tr \BBH^{2})^{2} ,
	\\
	\E \bigl( \bigl\langle \BBH \gaussv, \gaussv \bigr\rangle - \tr \BBH \bigr)^{5}
	&=&
	512 \tr \BBH^{5} + 32 \tr \BBH^{2} \, \tr \BBH^{3} ,
\label{2pG2trD2DGm22m2}
\end{EQA}
and
\begin{EQA}
	\E \bigl\langle \BBH \gaussv, \gaussv \bigr\rangle^{2}
	&=&
	(\tr \BBH)^{2} + 2 \tr \BBH^{2},
	\\
	\E \bigl\langle \BBH \gaussv, \gaussv \bigr\rangle^{3}
	& = &
	(\tr \BBH)^{3} + 6 \tr \BBH \,\, \tr \BBH^{2} + 8 \tr \BBH^{3} ,
	\\
	\E \bigl\langle \BBH \gaussv, \gaussv \bigr\rangle^{4}
	& = &
	(\tr \BBH)^{4} + 12 (\tr \BBH)^{2} \tr \BBH^{2}
	+ 32 (\tr \BBH) \tr \BBH^{3}
	+ 48 \tr \BBH^{4} + 12 (\tr \BBH^{2})^{2} ,
\label{2pG2trD2DGm22m2}
	\\
	\Var \bigl\langle \BBH \gaussv, \gaussv \bigr\rangle^{2}
	& = &
	8 (\tr \BBH)^{2} \tr \BBH^{2}
	+ 32 (\tr \BBH) \tr \BBH^{3}
	+ 48 \tr \BBH^{4} + 8 (\tr \BBH^{2})^{2} .
\label{2pG2trD2DGm22m4}
\end{EQA}
Moreover, if \( \BBH \leq \Id_{\dimp} \) and \( \dimH = \tr \BBH \), then \( \tr \BBH^{m} \leq \dimH \| \BBH \|^{m-1} \) for 
\( m \geq 1 \) and
\begin{EQA}[rcccl]
	\E \bigl\langle \BBH \gaussv, \gaussv \bigr\rangle^{2}
	& \leq &
	\dimH^{2} + 2 \dimH \| \BBH \|
	&\leq &
	(\dimH + \| \BBH \|)^{2},
	\\
	\E \bigl\langle \BBH \gaussv, \gaussv \bigr\rangle^{3}
	& \leq &
	\dimH^{3} + 6 \dimH^{2} \| \BBH \| + 8 \dimH \| \BBH \|^{2}
	&\leq &
	(\dimH + 2 \| \BBH \|)^{3},
	\\
	\E \bigl\langle \BBH \gaussv, \gaussv \bigr\rangle^{4}
	& \leq &
	\dimH^{4} + 12 \dimH^{3} \| \BBH \|
	+ 44 \dimH^{2} \| \BBH \|^{2}
	+ 48 \dimH \| \BBH \|^{3}
	&\leq &
	(\dimH + 3 \| \BBH \|)^{4},
	\\
	\E \bigl\langle \BBH \gaussv, \gaussv \bigr\rangle^{5}
	& \leq &
	\dimH^{5} + 20 \dimH^{4} \| \BBH \| + 140 \dimH^{3} \| \BBH \|^{2} 
	+ 272 \dimH^{2} \| \BBH \|^{3} + 512 \dimH \| \BBH \|^{4}
	& \leq &
	(\dimH + 4 \| \BBH \|)^{5} \, .
\label{2pG2trD2DGm22m2}
	\\
	\Var \bigl\langle \BBH \gaussv, \gaussv \bigr\rangle^{2}
	& \leq &
	8 \dimH^{3} + 40 \dimH^{2} \| \BBH \| + 48 \dimH \| \BBH \|^{2}.
\label{2pG2trD2DGm22m4}
\end{EQA}
Finally,
\begin{EQA}
	\E (\gaussv \gaussv^{\T} - \Id_{\dimp}) \BBH (\gaussv \gaussv^{\T} - \Id_{\dimp}) 
	&=&
	\BBH + \tr (\BBH) \Id_{\dimp}
\label{njt66777888723fdgy}
\end{EQA}
yielding
\begin{EQA}
	\E \| \BBH (\gaussv \gaussv^{\T} - \Id_{\dimp}) \|_{\Fr}^{2}
	&=&
	(\tr \BBH)^{2} + \tr \BBH^{2} .
\label{njt66777888723fdgyf}
\end{EQA}
\end{lemma}

\begin{proof}
Let \( \gaussv \) be standard normal in \( \R^{\dimp} \).
The same holds for \( \Uv \gaussv \) for any orthogonal transform \( \Uv \) in \( \R^{\dimp} \).
The use of the spectral decomposition \( \BBH = \Uv^{\T} \Lambda \Uv \) with \( \Uv \) orthonormal and 
\( \Lambda \) diagonal enables us to represent 
\( \bigl\langle \BBH \gaussv, \gaussv \bigr\rangle = \bigl\langle \Lambda \Uv \gaussv, \Uv \gaussv \bigr\rangle \) and thus,
to reduce the statements to the case when \( \BBH \) is diagonal: \( \BBH = \diag(\supH_{1},\supH_{2},\ldots,\supH_{\dimp}) \).
Then 
\begin{EQA}
	\xi
	\eqdef
	\bigl\langle \BBH \gaussv, \gaussv \bigr\rangle - \tr \BBH
	&=&
	\sum_{j=1}^{\dimp} \supH_{j} (\gauss_{j}^{2} - 1) ,
\label{j1ljgj2m1}
\end{EQA}
where \( \gauss_{j} \) are i.i.d. standard normal. 
This easily yields with \( \dimH_{m} = \tr (\BBH^{m}) \)
\begin{EQA}
	\E \xi^{2}
	&=&
	\sum_{j=1}^{\dimp} \supH_{j}^{2} \E (\gauss_{j}^{2} - 1)^{2}
	=
	\E \chi^{2} \, \tr \BBH^{2} 
	=
	2 \dimH_{2} \, ,
	\\
	\E \xi^{3}
	&=&
	\sum_{j=1}^{\dimp} \supH_{j}^{3} \E (\gauss_{j}^{2} - 1)^{3}
	=
	\E \chi^{3} \, \tr \BBH^{3} 
	=
	8 \dimH_{3} \, ,
	\\
	\E \xi^{4}
	&=&
	\sum_{j=1}^{\dimp} \supH_{j}^{4} (\gauss_{j}^{2} - 1)^{4}
	+ \sum_{i\neq j} \supH_{i}^{2} \supH_{j}^{2} \E (\gauss_{i}^{2} - 1)^{2} \E (\gauss_{j}^{2} - 1)^{2}
	\\
	&=&
	\bigl( \E \chi^{4} - 3 (\E \chi^{2})^{2} \bigr) \tr \BBH^{4} + 3 (\E \chi^{2} \, \tr \BBH^{2})^{2}
	=
	48 \dimH_{4} + 12 \dimH_{2}^{2} ,
	\\
	\E \xi^{5}
	&=&
	\sum_{j=1}^{\dimp} \supH_{j}^{5} (\gauss_{j}^{2} - 1)^{5}
	+ \sum_{i\neq j} \supH_{i}^{2} \supH_{j}^{3} \E (\gauss_{i}^{2} - 1)^{2} \E (\gauss_{j}^{2} - 1)^{3}
	\\
	&=&
	\bigl\{ \E (\gauss^{2} - 1)^{5} - \E (\gauss^{2} - 1)^{2} \E (\gauss^{2} - 1)^{3} \bigr\} \, \tr \BBH^{5} 
	+ \E (\gauss^{2} - 1)^{2} \E (\gauss^{2} - 1)^{3} \, \tr \BBH^{2} \, \tr \BBH^{3}
	\\
	&=&
	512 \dimH_{5} + 32 \dimH_{2} \, \dimH_{3} \, .
\label{2pG2trD2DGm22m2}
\end{EQA}
and
\begin{EQA}
	\E \bigl\langle \BBH \gaussv, \gaussv \bigr\rangle^{2}
	&=&
	\bigl( \E \bigl\langle \BBH \gaussv, \gaussv \bigr\rangle \bigr)^{2} 
	+ \E \xi^{2}
	= 
	\dimH^{2} + 2 \dimH_{2} \, ,
	\\
	\E \bigl\langle \BBH \gaussv, \gaussv \bigr\rangle^{3}
	& = &
	\E ( \xi + \dimH )^{3}
	=
	\dimH^{3} + \E \xi^{3} + 3 \dimH \,\, \E \xi^{2}
	=
	\dimH^{3} + 6 \dimH \,\, \dimH_{2} + 8 \dimH_{3} ,
	\\
	\E \bigl\langle \BBH \gaussv, \gaussv \bigr\rangle^{4}
	& = &
	\E \bigl( \xi + \dimH \bigr)^{4}
	=
	\dimH^{4} + 6 \dimH^{2} \E \xi^{2} + 4 \dimH \, \E \xi^{3} + \E \xi^{4}
	\\
	&=&
	\dimH^{4} + 12 \dimH^{2} \, \dimH_{2}
	+ 32 \dimH \, \dimH_{3}
	+ 48 \dimH_{4} + 12 \dimH_{2}^{2} ,
\label{2pG2trD2DGm22m2}
\end{EQA}
and 
\begin{EQA}
	&& \nquad
	\Var \bigl\langle \BBH \gaussv, \gaussv \bigr\rangle^{2}
	= 
	\E ( \xi + \dimH )^{4}
	- \bigl( \dimH^{2} + 2 \dimH_{2} \bigr)^{2}
	\\
	&=&
	\dimH^{4} + 6 \dimH^{2} \E \xi^{2} + 4 \dimH \, \E \xi^{3} + \E \xi^{4}
	- \bigl( \dimH^{2} + 2 \dimH_{2} \bigr)^{2}
	= 
	8 \dimH^{2} \, \dimH_{2}
	+ 32 \dimH \, \dimH_{3}
	+ 48 \dimH_{4} + 8 \dimH_{2}^{2} \, .
\label{2pG2trD2DGm22m4}
\end{EQA}
Also
\begin{EQA}
	\E \bigl\langle \BBH \gaussv, \gaussv \bigr\rangle^{5}
	& = &
	\E \bigl( \xi + \dimH \bigr)^{5}
	=
	\dimH^{5} + 10 \dimH^{3} \E \xi^{2} + 10 \dimH^{2} \, \E \xi^{3} + 5 \dimH \E \xi^{4}
	+ \E \xi^{5}
	\\
	&=&
	\dimH^{5} + 20 \dimH^{3} \, \dimH_{2} + 80 \dimH^{2} \dimH_{3} 
	+ 5 \dimH (48 \dimH_{4} + 12 \dimH_{2}^{2})
	+ 512 \dimH_{5} + 32 \dimH_{2} \, \dimH_{3} \, .
\label{2pG2trD2DGm22m25}
\end{EQA}
Assume \( \| \BBH \| = 1 \) yielding \( \dimH_{m} \leq \dimH \).
Then
\begin{EQA}
	\E \bigl\langle \BBH \gaussv, \gaussv \bigr\rangle^{2}
	& \leq &
	\dimH^{2} + 2 \dimH
	\leq 
	(\dimH + 1)^{2} \, ,
	\\
	\E \bigl\langle \BBH \gaussv, \gaussv \bigr\rangle^{3} 
	& \leq &
	\dimH^{3} + 6 \dimH^{2} + 8 \dimH
	\leq 
	(\dimH + 2)^{3} ,
	\\
	\E \bigl\langle \BBH \gaussv, \gaussv \bigr\rangle^{4}
	& \leq &
	\dimH^{4} + 12 \dimH^{3} + 44 \dimH^{2} + 48 \dimH 
	\leq 
	(\dimH + 3)^{4} ,
	\\
	\E \bigl\langle \BBH \gaussv, \gaussv \bigr\rangle^{5}
	& \leq &
	\dimH^{5} + 20 \dimH^{4} + 140 \dimH^{3} + 272 \dimH^{2} + 512 \dimH 
	\leq 
	(\dimH + 4)^{5} \, .
\label{2pG2trD2DGm22m25}
\end{EQA}
For the last result of the lemma, observe that with \( \BBH = \diag(\supH_{1},\supH_{2},\ldots,\supH_{\dimp}) \), 
\begin{EQA}
	\E \| \BBH^{1/2} (\gaussv \gaussv^{\T} - \Id_{\dimp}) \BBH^{1/2} \|_{\Fr}^{2}
	&=&
	\sum_{i,j=1}^{\dimp} \supH_{i} \supH_{j} \E (\gauss_{i} \gauss_{j} - \delta_{i,j})^{2}
	=
	\left( \sum_{i=1}^{\dimp} \supH_{i} \right)^{2} + \sum_{i=1}^{\dimp} \supH_{i}^{2} 
\label{ikcduywjwsiv98emdvuw}
\end{EQA}
and assertion \eqref{njt66777888723fdgyf} follows.
\end{proof}

\begin{lemma}
\label{L2tensmoment}
For any vector \( \cv \in \R^{\dimp} \) and \( \gaussB \sim \ND(0,\BBH) \), it holds
with \( \Delta = \BBH \cv \)
\begin{EQA}[c]
	\E \bigl\{ 
			\langle \cv, \gaussB \rangle^{2} \, \gaussB \, \gaussB^{\T} 
	\bigr\}
	=
	2 \Delta \Delta^{\T} + \diag (\Delta_{1}^{2},\ldots,\Delta_{\dimp}^{2})
	\, .
\end{EQA}
\end{lemma}

\begin{proof}
Without loss of generality, assume the matrix \( \BBH \) to be diagonal:
\( \BBH = \diag(\lambda_{1}, \ldots, \lambda_{\dimp}) \).
Then \( \gaussB = (\lambda_{1}^{1/2} \gauss_{1}, \ldots, \lambda_{\dimp}^{1/2} \gauss_{\dimp}) \) for i.i.d. standard Gaussian \( \gauss_{j} \), and
\begin{EQA}[c]
	\langle \cv, \gaussB \rangle^{2}
	=
	\sum_{i,j=1}^{\dimp} c_{i} \, c_{j} \lambda_{i}^{1/2} \lambda_{j}^{1/2} \gauss_{i} \gauss_{j}
	\,
\end{EQA}
and for any pair \( i \neq j \)
\begin{EQA}
	\E \bigl\{ 
			\langle \cv, \gaussB \rangle^{2} \, \lambda_{i}^{1/2} \gauss_{i}
			\,\, \lambda_{j}^{1/2} \gauss_{j} 
	\bigr\}
	&=&
	2 c_{i} \, c_{j} \lambda_{i} \lambda_{j}
	=
	2 \Delta_{i} \, \Delta_{j}
	\, .
\label{uhcu83eknkfiv774fvk}
\end{EQA}
Similarly, for any \( j \)
\begin{EQA}
	\E \bigl\{ 
			\langle \cv, \gaussB \rangle^{2} \, \lambda_{j} \gauss_{j}^{2}
	\bigr\}
	&=&
	c_{j}^{2} \lambda_{j}^{2} \E \gauss_{j}^{4}
	=
	3 c_{j}^{2} \lambda_{j}^{2}
	=
	3 \Delta_{j}^{2}
	\, .
\label{ugcd7y83d8hie4fhr7t7efg}
\end{EQA}
These two identities imply the statement of the lemma.
\end{proof}

Now we compute the exponential moments of centered and non-centered quadratic forms.

\begin{lemma}
\label{Lqfexpmom}
Let \( \| \BBH \| = \supH \) and \( \gaussv \sim \ND(0,\Id_{\dimp}) \).
Then for any \( \mu \in (0,\supH^{-1}) \), 
\begin{EQA}
	\E \exp \Bigl\{ \frac{\mu}{2} \langle \BBH \gaussv, \gaussv \rangle \Bigr\}
	&=&
	\det(\Id_{\dimp} - \mu \BBH)^{-1/2} \, .
\label{m2v241m41m}
\end{EQA}
Moreover, with \( \dimH = \tr \BBH \) and \( \vH^{2} = \tr \BBH^{2} \)
\begin{EQA}
	\log \E \exp \Bigl\{ \frac{\mu}{2} \bigl( \langle \BBH \gaussv, \gaussv \rangle - \dimH \bigr) \Bigr\}
	& \leq &
	\frac{\mu^{2} \vH^{2}}{4 (1 - \supH \mu)} \, .
\label{m2v241m41mb}
\end{EQA}
If \( \BBH \) is positive semidefinite, \( \supH_{j} \geq 0 \), then 
\begin{EQA}
	\log \E \exp \Bigl\{ - \frac{\mu}{2} \bigl( \langle \BBH \gaussv, \gaussv \rangle - \dimH \bigr) \Bigr\}
	& \leq &
	\frac{\mu^{2} \vH^{2}}{4} \, .
\label{m2v241m41mbn}
\end{EQA}
For any complex valued \( \muH \) with \( \supH |\muH| < 1 \),
\begin{EQA}
	\biggl| \log \E \exp \Bigl\{ 
			\frac{\mu}{2} \bigl( \langle \BBH \gaussv, \gaussv \rangle - \dimH \bigr) - \frac{\muH^{2} \tr \BBH^{2}}{4}
		\Bigr\} 
	\biggr|
	& \leq &
	\frac{\supH |\mu|^{3} \vH^{2} }{6 (1 - \supH |\mu|)} \, .
\label{vbu7j3hg8hryhghyidegwdg}
\end{EQA}
\end{lemma}

\begin{proof}
W.l.o.g. assume \( \supH = 1 \).
Let \( \supH_{j} \) be the eigenvalues of \( \BBH \), \( |\supH_{j}| \leq 1 \).
As in Lemma~\ref{Gaussmoments}, one can reduce the statement to the case of a diagonal matrix 
\( \BBH = \diag\bigl( \supH_{j} \bigr) \). 
Then \( \langle \BBH \gaussv, \gaussv \rangle = \sum_{j=1}^{\dimp} \supH_{j} \gauss_{j}^{2} \) and 
by independence of the \( \gauss_{j} \)'s
\begin{EQA}
	&& \nquad
	\E \Bigl\{ \frac{\mu}{2} \langle \BBH \gaussv, \gaussv \rangle  \Bigr\}
	=
	\prod_{j=1}^{\dimp} \E \exp \Bigl( \frac{\mu}{2} \supH_{j} \eps_{j}^{2} \Bigr)
	=
	\prod_{j=1}^{\dimp} \frac{1}{\sqrt{1 - \mu \supH_{j}}} 
	=
	\det \bigl( \Id_{\dimp} - \mu \BBH \bigr)^{-1/2} .
\label{dOImuBm12EB}
\end{EQA}
Below we use the simple bounds: 
\begin{EQ}[rcl]
\label{lo1uusk2iukkp}
	- \log(1 - u) - u
	&=&
	\sum_{k=2}^{\infty} \frac{u^{k}}{k}
	\leq 
	\frac{u^{2}}{2} \sum_{k=0}^{\infty} u^{k} 
	=
	\frac{u^{2}}{2 (1 - u)} \, ,
	\qquad 
	u \in (0,1),
	\\
	- \log(1 - u) + u
	&=&
	\sum_{k=2}^{\infty} \frac{u^{k}}{k}
	\leq 
	\frac{u^{2}}{2} \, ,
	\qquad \qquad
	u \in (-1,0).
\label{lo1uusk2iukk}
\end{EQ}
Now it holds for \( \mu > 0 \)
\begin{EQA}
	&& \nquad
	\log \E \Bigl\{ \frac{\mu}{2} \bigl( \langle \BBH \gaussv, \gaussv \rangle - \dimH \bigr) \Bigr\}
	=
	\log \det(\Id_{\dimp} - \mu \BBH)^{-1/2} - \frac{\mu \, \dimH}{2}
	\\
	&=&
	- \frac{1}{2} \sum_{j=1}^{\dimp} \bigl\{ \log(1 - \mu \supH_{j}) + \mu \supH_{j} \bigr\}
	\leq 
	\sum_{j=1}^{\dimp} \frac{\mu^{2} \supH_{j}^{2}}{4 (1 - \mu \supH_{j})} 
	\leq 
	\frac{\mu^{2} \vH^{2}}{4 (1 - \mu \supH)} \, .
\label{m2v241m4mj1pd}
\end{EQA}
Similarly for any complex \( \mu \) with \( |\mu| \supH < 1 \)
\begin{EQA}
	&& \nquad
	\left| 
		\log \E \Bigl\{ \frac{\mu}{2} \bigl( \langle \BBH \gaussv, \gaussv \rangle - \dimH \bigr) 
		- \frac{\muH^{2} \tr \BBH^{2}}{4} \Bigr\}
	\right|
	=
	\left| \log \det(\Id_{\dimp} - \mu \BBH)^{-1/2} - \frac{\mu \, \dimH}{2} - \frac{\muH^{2} \tr \BBH^{2}}{4} \right|
	\\
	&=&
	\frac{1}{2} \left| 
		\sum_{j=1}^{\dimp} \biggl\{ \log(1 - \mu \supH_{j}) - \mu \supH_{j} - \frac{\mu^{2} \supH_{j}^{2}}{2} \biggr\} 
	\right|
	\leq 
	\sum_{j=1}^{\dimp} \frac{|\mu \supH_{j}|^{3}}{6 (1 - \supH |\mu|)} 
	=
	\frac{| \mu |^{3} \supH \vH^{2}}{6 (1 - \supH |\mu|)} \, .
\label{m2v241m4mj1pd}
\end{EQA}
Statement \eqref{m2v241m41mbn} can be proved similarly.
\end{proof}

Now we consider the case of a non-centered quadratic form
\( \langle \BBH \gaussv,\gaussv \rangle/2 + \langle \Av,\gaussv \rangle \) for a fixed vector \( \Av \).

\begin{lemma}
\label{Lexpmomnoncen}
Let \( \| \BBH \| = \supH < 1 \). 
Then for any \( \Av \)
\begin{EQA}
	\E \exp\Bigl\{ \frac{1}{2}\langle \BBH \gaussv,\gaussv \rangle + \langle \Av,\gaussv \rangle \Bigr\}
	&=&
	\exp\Bigl\{ \frac{\| (\Id_{\dimp} - \BBH)^{-1/2} \Av \|^{2}}{2} \Bigr\} \, \det(\Id_{\dimp} - \BBH)^{-1/2} .
\label{EeBf12BggA}
\end{EQA}
Moreover, for any \( \mu \in (0,1) \)
\begin{EQA}
	&& \nquad
	\log \E \exp\Bigl\{ 
		\frac{\mu}{2} \bigl( \langle \BBH \gaussv,\gaussv \rangle - \dimH \bigr) + \langle \Av,\gaussv \rangle 
	\Bigr\}
	\\
	&=&
	\frac{\| (\Id_{\dimp} - \mu \BBH)^{-1/2} \Av \|^{2}}{2} + \log \det(\Id_{\dimp} - \mu \BBH)^{-1/2} - \mu \, \dimH 
	\\
	& \leq &
	\frac{\| (\Id_{\dimp} - \mu \BBH)^{-1/2} \Av \|^{2}}{2} + \frac{\mu^{2} \vH^{2}}{4 (1 - \supH \mu)} \, .
\label{EeBf12BggAmu}
\end{EQA}
\end{lemma}

\begin{proof}
Denote \( \av = (\Id_{\dimp} - \BBH)^{-1/2} \Av \). 
It holds by change of variables \( (\Id_{\dimp} - \BBH)^{1/2} \xv = \uv \) for \( \CONSTi_{\dimp} = (2\pi)^{-\dimp/2} \)
\begin{EQA}
	&& \nquad
	\E \exp\Bigl\{ \frac{1}{2}\langle \BBH \gaussv,\gaussv \rangle + \langle \Av,\gaussv \rangle \Bigr\}
	=
	\CONSTi_{\dimp}
	\int \exp\Bigl\{ - \frac{1}{2}\langle (\Id_{\dimp} - \BBH) \xv,\xv \rangle + \langle \Av,\xv \rangle \Bigr\} d\xv
	\\
	&=&
	\CONSTi_{\dimp}
	\det(\Id_{\dimp} - \BBH)^{-1/2}
	\int \exp\Bigl\{ - \frac{1}{2} \| \uv \|^{2} + \langle \av,\uv \rangle \Bigr\} d\uv
	=
	\det(\Id_{\dimp} - \BBH)^{-1/2} \, 	\ex^{\| \av \|^{2}/2}  	.
\label{EeBf12BggAp}
\end{EQA}
The last inequality \eqref{EeBf12BggAmu} follows by \eqref{m2v241m41mb}.
\end{proof}

\Section{Deviation bounds for Gaussian quadratic forms}
\label{SdevboundGauss}
The next result explains the concentration effect of \( \| \QP \xiv \|^{2} \)
for a centered Gaussian vector \( \xiv \sim \ND(0,\HVB^{2}) \) and a linear operator \( \QP \colon \R^{\dimp} \to \R^{\dimq} \),
\( \dimp,\dimq \leq \infty \).
We use a version from \cite{laurentmassart2000}.
For completeness, we present a simple proof.

\begin{theorem}
\label{TexpbLGA}
\label{Lxiv2LD}
\label{Cuvepsuv0}
Let \( \xiv \sim \ND(0,\HVB^{2}) \) be a Gaussian element in \( \R^{\dimp} \) and let
\( \QP \colon \R^{\dimp} \to \R^{\dimq} \) be such that \( \BBH = \QP \HVB^{2} \QP^{\T} \) 
is a trace operator in \( \R^{\dimq} \).
Then with \( \dimH = \tr(\BBH) \), \( \vH^{2} = \tr(\BBH^{2}) \), and \( \supH = \| \BBH \| \),
it holds for any \( \xx \geq 0 \)
\begin{EQA}
\label{Pxiv2dimAvp12}
	\P\Bigl( \| \QP \xiv \|^{2} - \dimH > 2 \vH \, \sqrt{\xx} + 2 \supH \xx \Bigr)
	& \leq &
	\ex^{-\xx} ,
	\\
	\P\Bigl( \| \QP \xiv \|^{2} - \dimH \leq - 2 \vH \, \sqrt{\xx} \Bigr)
	& \leq &
	\ex^{-\xx} .
\label{Pxiv2dimAvp12m}
\end{EQA}
It also implies 
\begin{EQA}
	\P\bigl( \bigl| \| \QP \xiv \|^{2} - \dimH \bigr| > \zq_{2}(\BBH,\xx) \bigr)
	& \leq &
	2 \ex^{-\xx} ,
\label{PxivTBBdimA2vp}
\end{EQA}
with
\begin{EQA}
	\zq_{2}(\BBH,\xx)
	& \eqdef &
	2 \vH \, \sqrt{\xx} + 2 \supH \xx \,\, .
\label{zqdefGQF}
\end{EQA}
%
\end{theorem}

\begin{proof}
We use the identity in distribution \( \| \QP \xiv \|^{2} \eqd \langle \BBH \gaussv, \gaussv \rangle \) with
 \( \gaussv \sim \ND(0,\Id_{\dimq}) \).
Markov's inequality yields for any \( \mu > 0 \)
\begin{EQA}
	\P\Bigl( \langle \BBH \gaussv, \gaussv \rangle - \dimH > \zq_{2}(\BBH,\xx) \Bigr)
	& \leq &
	\E \exp \Bigl( \frac{\mu}{2} \bigl( \langle \BBH \gaussv, \gaussv \rangle - \dimH \bigr) - \frac{\mu \, \zq_{2}(\BBH,\xx)}{2} 
	\Bigr) \, .
\label{PBggiz2E2mz2}
\end{EQA}
Given \( \xx > 0 \), fix \( \mu < 1/\supH \) by the equation
\begin{EQA}
	\frac{\mu}{1 - \supH \mu} 
	&=&
	\frac{2 \sqrt{\xx}}{\vH} \, 
	\quad \text{ or } \quad
	\mu^{-1} 
	=
	\supH + \frac{\vH}{2 \sqrt{\xx}} \, .
\label{1v2sxm12m1m}
\end{EQA}
By \eqref{m2v241m41mb}
\begin{EQA}
	&& \nquad
	\log \E \Bigl\{ \frac{\mu}{2} \bigl( \langle \BBH \gaussv, \gaussv \rangle - \dimH \bigr) \Bigr\}
	\leq 
	\frac{\mu^{2} \vH^{2}}{4 (1 - \supH \mu)} \, .
\label{m2v241m4mj1p}
\end{EQA}
For \eqref{Pxiv2dimAvp12}, it remains to check that the choice \( \mu \) by \eqref{1v2sxm12m1m} yields
\begin{EQA}
	\frac{\mu^{2} \vH^{2}}{4 (1 - \supH \mu)} - \frac{\mu \, \zq_{2}(\BBH,\xx)}{2}
	& = &
	\frac{\mu^{2} \vH^{2}}{4 (1 - \supH \mu)} - \mu \bigl( \vH \sqrt{\xx} + \supH \xx \bigr)
	=
	\mu \Bigl( \frac{\vH \sqrt{\xx}}{2} - \vH \sqrt{\xx} - \supH \xx \Bigr)
	=
	- \xx .
\label{m2vA241muz2}
\end{EQA}
The bound \eqref{Pxiv2dimAvp12m} is obtained similarly from Markov's inequality 
applied to \( - \langle \BBH \gaussv, \gaussv \rangle + \dimH \) with \( \mu = 2 \vH^{-1} \sqrt{\xx} \).
The use of \eqref{m2v241m41mbn} yields
\begin{EQA}
	&& \nquad
	\P\Bigl( \langle \BBH \gaussv, \gaussv \rangle - \dimH < - 2 \vH \sqrt{\xx} \Bigr)
	\leq
	\E \exp \Bigl\{ \frac{\mu}{2} \bigl( - \langle \BBH \gaussv, \gaussv \rangle + \dimH \bigr) - \mu \, \vH \sqrt{\xx} 
	\Bigr\}
	\\
	& \leq &
	\exp \Bigl( \frac{\mu^{2} \vH^{2}}{4} - \mu \, \vH \sqrt{\xx} \Bigr) 
	=
	\ex^{-\xx} \, 
\label{PBggiz2E2mz2}
\end{EQA}
as required.
\end{proof}

\begin{corollary}
\label{CTexpbLGAd}
Assume the conditions of Theorem~\ref{TexpbLGA}.
Then for \( \zq > \vH \)
\begin{EQA}
	\P\bigl( \bigl| \| \QP \xiv \|^{2} - \dimH \bigr| \ge \zq \bigr)
	& \leq &
	2 \exp\biggl\{ - \frac{\zq^{2}}{\bigl( \vH + \sqrt{\vH^{2} + 2 \supH \zq} \bigr)^{2}} \biggr\}
	\leq 
	2 \exp\biggl( - \frac{\zq^{2}}{4\vH^{2} + 4 \supH \zq} \biggr) .
	\qquad
	\qquad
\label{3z2spsp2z3z2}
\end{EQA}
\end{corollary}

\begin{proof}
Given \( \zq \), define \( \xx \) by 
\( 2 \vH \sqrt{\xx} + 2 \supH \xx = \zq \) or 
\( 2 \supH \sqrt{\xx} = \sqrt{\vH^{2} + 2 \supH \zq} - \vH \).
Then
\begin{EQA}
	\P\bigl( \| \QP \xiv \|^{2} - \dimH \ge \zq \bigr)
	& \leq &
	\ex^{-\xx} 
	=
	\exp\biggl\{ - \frac{\bigl( \sqrt{\vH^{2} + 2 \supH \zq} - \vH \bigr)^{2}}{4 \supH^{2}} \biggr\}
	=
	\exp\biggl\{ - \frac{\zq^{2}}{\bigl( \vH + \sqrt{\vH^{2} + 2 \supH \zq} \bigr)^{2}} \biggr\}.
\label{3emzmsp22z2c}
\end{EQA}
This yields \eqref{3z2spsp2z3z2} by direct calculus.
\end{proof}

Of course, bound \eqref{3z2spsp2z3z2} is sensible only if \( \zq \gg \vH \).

\ifadap{}{
\begin{corollary}
\label{RsochpHsA}
Assume the conditions of Theorem~\ref{TexpbLGA}.
If also \( \BBH \geq 0 \), then 
\begin{EQA}
\label{Pxiv2dimAxx12}
	\P\Bigl( \| \QP \xiv \|^{2} \geq \zq^{2}(\BBH,\xx) \Bigr)
	& \leq &
	\ex^{-\xx} 
\end{EQA}
with 
\begin{EQA}
	\zq^{2}(\BBH,\xx)
	& \eqdef &
	\dimH + 2 \vH \, \sqrt{\xx} + 2 \supH \xx
	\leq 
	\bigl( \sqrt{\dimH} + \sqrt{2 \supH \xx} \bigr)^{2} \, .
\label{zzxxppdBlroBB}
\end{EQA}
Also
\begin{EQA}
	\P\Bigl( \| \QP \xiv \|^{2} - \dimH < - 2 \vH \, \sqrt{\xx} \Bigr)
	& \leq &
	\ex^{-\xx} .
\label{Pxiv2dimAvp12d}
\end{EQA}
\end{corollary}

\begin{proof}
The definition implies \( \vH^{2} \leq \dimH \supH \)
yielding the statement of the corollary.
\end{proof}

As a special case, we present a bound for the chi-squared distribution 
corresponding to \( \QP = \HVB^{2} = \Id_{\dimp} \), \( \dimp < \infty \).
Then \( \BBH = \Id_{\dimp} \), \( \tr (\BBH) = \dimp \), \( \tr(\BBH^{2}) = \dimp \) and \( \supH(\BBH) = 1 \).

\begin{corollary}
\label{Cchi2p}
Let \( \gaussv \) be a standard normal vector in \( \R^{\dimp} \).
Then for any \( \xx > 0 \)
\begin{EQA}[ccl]
\label{Pxi2pm2px}
	\P\bigl( \| \gaussv \|^{2} \geq \dimp + 2 \sqrt{\dimp \, \xx} + 2 \xx \bigr)
	& \leq &
	\ex^{-\xx},
	\\
	\P\bigl( \| \gaussv \| \,\,  \geq \sqrt{\dimp} + \sqrt{2 \xx} \bigr)
	& \leq &
	\ex^{-\xx} ,
\label{Pxi2pm2px12}
	\\
	\P\bigl( \| \gaussv \|^{2} \leq \dimp - 2 \sqrt{\dimp \, \xx} \bigr)
	& \leq &
	\ex^{-\xx}	.
\label{Pxi2pm2px22}
\end{EQA}
\end{corollary}
}

The bound of Theorem~\ref{TexpbLGA} 
can be represented as a usual deviation bound.

\begin{theorem}
\label{CTexpbLGA}
Assume the conditions of Theorem~\ref{TexpbLGA}.
For \( \yy > 0 \), define
\begin{EQA}
	\xx(\yy)
	& \eqdef &
	\frac{(\sqrt{\yy + \dimH} - \sqrt{\dimH})^{2}}{4 \supH} \, .
\label{iuvfiiow3kboieheuf}
\end{EQA}
Then
\begin{EQA}
	\P\bigl( \| \QP \xiv \|^{2} \ge \dimH + \yy \bigr)
	& \leq &
	\ex^{- \xx(\yy)} ,
\label{3emzmsp22z2}
	\\
	\E \bigl\{ (\| \QP \xiv \|^{2} - \dimH) \Ind\bigl( \| \QP \xiv \|^{2} \ge \dimH + \yy \bigr) \bigr\}
	& \leq &
	2 \Bigl( \frac{\yy + \dimH}{\supH \, \xx(\yy)} \Bigr)^{1/2} \, \, 
	\ex^{- \xx(\yy)} \, .
	\qquad
	\quad
\label{3emzmsp22z2e}
\end{EQA}
Moreover, let \( \muH > 0 \) fulfill \( \rexH =  \muH \supH + \muH \sqrt{\supH \dimH / \xx(\yy)} < 1 \). 
Then 
\begin{EQA}
	\E \bigl\{ \ex^{\muH (\| \QP \xiv \|^{2} - \dimH)/2} \Ind( \| \QP \xiv \|^{2} \ge \dimH + \yy) \bigr\}
	& \leq &
	\frac{1}{1 - \rexH} \, \exp\{ - (1 - \rexH) \xx(\yy) \} \, .
	\qquad
\label{llkknbononjm9hig4e}
\end{EQA}
\end{theorem}

\begin{proof}
Normalizing by \( \supH \) reduces the statements to the case with \( \supH = 1 \).
Define \( \eta = \| \QP \xiv \|^{2} - \dimH \) 
and
\begin{EQA}
	\zq(\xx)
	&=&
	2 \sqrt{\dimH \, \xx} + 2 \xx .
\label{0kmuy765433udgswhhh}
\end{EQA}
Then by \eqref{Pxiv2dimAvp12} \( \P(\eta \geq \zq(\xx)) \leq \ex^{-\xx} \).
Inverting the relation \eqref{0kmuy765433udgswhhh} yields
\begin{EQA}
	\xx(\zq)
	&=&
	\frac{1}{4} \bigl( \sqrt{\zq + \dimH} - \sqrt{\dimH} \bigr)^{2}
\label{jkv78fdjryfgsdfghgj}
\end{EQA}
and \eqref{3emzmsp22z2} follows by applying \( \zq = \yy \).
Further, 
\begin{EQA}
	\E \bigl\{ \eta \Ind(\eta \geq \yy) \bigr\}
	&=&
	\int_{\yy}^{\infty} \P(\eta \geq \zq) \, d\zq
	\leq 
	\int_{\yy}^{\infty} \ex^{ - \xx(\zq) } \, d\zq
	= 
	\int_{\xx(\yy)}^{\infty} \ex^{-\xx} \, \zq'(\xx) \, d\xx \, .
\label{zEe2Iezz2c2H23}
\end{EQA} 
As \( \zq'(\xx) = 2 + \sqrt{\dimH/\xx} \) monotonously decreases with \( \xx \), we derive
\begin{EQA}
	\E \bigl\{ \eta \Ind(\eta \geq \yy) \bigr\}
	& \leq &
	\zq'(\xx(\yy)) \ex^{-\xx(\yy)}
	=
	\frac{1}{\xx'(\yy)} \, \ex^{- \xx(\yy)}
	=
	\frac{4 \sqrt{\yy + \dimH}}{\sqrt{\yy + \dimH} - \sqrt{\dimH}} \, \ex^{- \xx(\yy)}
\label{e7ygv76bgughytuj}
\end{EQA}
and \eqref{3emzmsp22z2e} follows.

In a similar way, define \( \zqe(\xx) \) from the relation
\( \muH^{-1} \log \zqe(\xx) = \sqrt{\dimH \, \xx} + \xx \) yielding
\begin{EQA}
	\zqe(\xx)
	&=&
	\exp \bigl( \muH \sqrt{\dimH \, \xx} + \muH \, \xx \bigr) .
\label{jvcjjuvue37r6gtur4r}
\end{EQA}
The inverse relation reads
\begin{EQA}
	\xxe(\zqe)
	&=&
	\bigl( \sqrt{\muH^{-1} \log \zqe + \dimH/4} - \sqrt{\dimH/4} \bigr)^{2} .
\label{jkv78fdjryfgsdfghgjex}
\end{EQA}
Then with \( \xx(\yy) = \xxe(\ex^{\muH \yy/2}) = \bigl( \sqrt{\yy + \dimH} - \sqrt{\dimH} \bigr)^{2}/4 \)
\begin{EQA}
	\E \bigl\{ \ex^{\muH \eta/2} \Ind(\eta \geq \yy) \bigr\}
	&=&
	\int_{\ex^{\muH \yy/2}}^{\infty} \P(\ex^{\muH \eta/2} \geq \zqe) \, d\zqe
	=
	\int_{\ex^{\muH \yy/2}}^{\infty} \P(\eta \geq 2\muH^{-1} \log \zqe) \, d\zqe
	\\
	& \leq &
	\int_{\ex^{\muH \yy/2}}^{\infty} \ex^{ - \xxe(\zqe) } \, d\zqe
	= 
	\int_{\xx(\yy)}^{\infty} \ex^{-\xx} \, \zqe'(\xx) \, d\xx .
\label{zEe2Iezz2c2H23}
\end{EQA} 
Further, in view of \( \muH + 0.5 \,\muH \sqrt{\dimH/\xx} < \muH + \muH \sqrt{\dimH / \xx(\yy)} = \rexH < 1 \) for 
\( \xx \geq \xx(\yy) \), it holds
\begin{EQA}
	\zqe'(\xx)
	&=&
	\bigl( \muH + 0.5 \, \muH \sqrt{\dimH/\xx} \bigr) \exp \bigl( \muH \sqrt{\dimH \, \xx} + \muH \, \xx \bigr) 
	\leq 
	\exp \bigl( \muH \, \xx \sqrt{\dimH / \xx(\yy)} + \muH \, \xx \bigr)
	=
	\exp (\rexH \, \xx) 
\label{jcuyu3ww3jbkihjitwedk}
\end{EQA}
and  
\begin{EQA}
	\E \bigl\{ \ex^{\muH \eta/2} \Ind(\eta \geq \yy) \bigr\}
	& \leq &
	\int_{\xx(\yy)}^{\infty} \ex^{-(1 - \rexH)\xx} \, d\xx 
	=
	\frac{1}{1 - \rexH} \, \ex^{-(1 - \rexH)\xx(\yy)} \, 
\label{zEe2Iezz2c2H23}
\end{EQA} 
and \eqref{llkknbononjm9hig4e} follows.
\end{proof}

\def\Xvt{\tilde{\Xv}}
\def\muHb{\bar{\muH}}
\def\gaussv{\gammav}
\def\rexH{\omega}
\def\CONSTgmb{\CONSTi_{\hspace{-1pt}X}}
\def\cdensX{\cdens_{\hspace{-1pt} X}}
\def\charfX{\charf_{\hspace{-1pt} X}}

\Section{Deviation bounds for sub-gaussian quadratic forms}
\label{Sprobabquad}
\label{SdevboundnonGauss}
This section collects some probability bounds for sub-gaussian quadratic forms.

\Subsection{A rough upper bound}
\label{SsubGausstail}

Let \( \xiv \) be a random vector in \( \R^{\dimp} \) with \( \E \xiv = 0 \).
We suppose that there exists a positive symmetric operator \( \HVB \) in \( \R^{\dimp} \) such that
\begin{EQA}
	\log \E \exp \bigl( \langle \uv, \HVB^{-1} \xiv \rangle \bigr)
	& \leq &
	\frac{\| \uv \|^{2}}{2} \, ,
	\qquad 
	\uv \in \R^{\dimp} .
\label{devboundinf}
\end{EQA}
In the Gaussian case, one can take \( \HVB^{2} = \Var(\xiv) \).
In general, \( \HVB^{2} \geq \Var(\xiv) \).
We consider a quadratic form \( \| \xiv \|^{2} \), where 
\( \xiv \) satisfies \eqref{devboundinf}.
We show that under \eqref{devboundinf}, the quadratic form \( \| \xiv \|^{2} \)
follows the same upper deviation bound 
\( \P\bigl( \| \xiv \|^{2} \geq \tr(\BBH) + 2 \sqrt{\xx \tr(\BBH^{2})} + 2 \xx \| \BBH \| \bigr) \leq \ex^{-\xx} \) 
with \( \BBH = \HVB^{2} \) as in the Gaussian case.

\begin{theorem}[\cite{HKZ2012}]
\label{Tdevboundinf}
Suppose \eqref{devboundinf}. 
Then for any \( \muH < 1/\| \BBH \| \)
\begin{EQA}
	\E \exp\Bigl( \frac{\muH}{2} \| \xiv \|^{2} \Bigr)
	& \leq &
	\exp \Bigl( \frac{\muH^{2} \tr(\BBH^{2})}{4 (1 - \| \BBH \| \muH)} + \frac{\muH \, \tr(\BBH)}{2} \Bigr)
\label{23iov96oklg7yf532tyfu}
\end{EQA}
and for any \( \xx > 0 \)
\begin{EQA}
	\P\bigl( \| \xiv \|^{2} > \tr(\BBH) + 2 \sqrt{\xx \tr(\BBH^{2})} + 2 \xx \| \BBH \| 
	\bigr)
	& \leq &
	\ex^{-\xx} .
\label{PxivbzzBBroBinf}
\end{EQA}
\end{theorem}

Statement \eqref{PxivbzzBBroBinf} looks identical to the upper bound in \ifapp{\eqref{Pxiv2dimAvp12}}{\eqref{v7hfte535ghewjfyw}}, 
however, there is an essential difference:
\( \tr(\BBH) \) can be much larger than 
\( \E \| \xiv \|^{2} = \tr \Var(\xiv) \) if \( \HVB^{2} \gg \Var(\xiv) \).
The result from \eqref{PxivbzzBBroBinf} is not accurate enough
for supporting the concentration property that \( \| \xiv \|^{2} \) concentrates
around its expectation \( \E \| \xiv \|^{2} \).
The next section presents some sufficient conditions for obtaining sharp Gaussian-like deviation bounds. 

\Subsection{Concentration of the squared norm of a sub-gaussian vector} 
\label{Sdevboundsharp}
Let \( \Xv \) be a centered random vector in \( \R^{\dimp} \).
We study concentration property of the squared norm \( \| \QP \Xv \|^{2} \) for 
a linear mapping \( \QP \colon \R^{\dimp} \to \R^{\dimq} \).
The aim is to establish the results similar to 
\ifapp{\eqref{Pxiv2dimAvp12}}{\eqref{v7hfte535ghewjfyw}} 
with \( \BBH = \QP \Var(\Xv) \QP^{\T} \). 
Later we assume the following condition on the moment-generating function \( \E \ex^{ \langle \uv, \Xv \rangle } \).

\begin{description}
\item[\( \bb{(\cdensX)} \)\label{gmbref}]
\emph{A random vector \( \Xv \in \R^{\dimp} \) satisfies \( \E \Xv = 0 \), \( \Var(\Xv) \leq  \Id_{\dimp} \). 
The function \( \cdensX(\uv) \eqdef \log \E \ex^{ \langle \uv, \Xv \rangle } \) is finite and fulfills 
for some \( \CONSTgmb \)}
\begin{EQA}
	\cdensX(\uv)
	\eqdef
	\log \E \ex^{ \langle \uv, \Xv \rangle }
	& \leq &
	\frac{\CONSTgmb \| \uv \|^{2}}{2} \, ,
	\qquad
	\uv \in \R^{\dimp} \, .
\label{devboundinfgmb}
\end{EQA}
\end{description}

The condition \( \Var(\Xv) \leq  \Id_{\dimp} \) is only for convenience.
One can drop it by rescaling \( \Xv \) and \( \QP \).
The constant \( \CONSTgmb \) can be quite large, it does not show up in the leading term of the obtained bound.
Also, we will only use \eqref{devboundinfgmb} for \( \| \uv \| \geq \gmn \) for some large \( \gmn \).
For \( \| \uv \| \leq \gmn \), we use smoothness properties of \( \cdensX(\uv) \).

The bounds in \ifapp{\eqref{Pxiv2dimAvp12}}{\eqref{v7hfte535ghewjfyw}} 
and in \eqref{PxivbzzBBroBinf} are uniform in the sense that they apply for all \( \xx \) and all \( \BBH \).
The results of this section are limited to a high dimensional situation
with \( \tr(\BBH^{2}) \gg \CONSTgmb \| \QP \QP^{\T} \| \) and apply only for 
\( \xx \ll \tr(\BBH^{2})/(\CONSTgmb \| \QP \QP^{\T} \|) \).
As compensation for this constraint, the bounds are surprisingly sharp.
In fact, they perfectly replicate bounds \ifapp{\eqref{Pxiv2dimAvp12}}{\eqref{v7hfte535ghewjfyw}} 
from the Gaussian case,
the upper and lower quantiles are exactly as in \ifapp{\eqref{Pxiv2dimAvp12}}{\eqref{v7hfte535ghewjfyw}} 
and the deviation probability 
is increased from  \( \ex^{-\xx} \) to \( (1 + \Delta_{\muH}) \ex^{-\xx} \) for a small value \( \Delta_{\muH} \).
For larger \( \xx \), one can still apply rough upper bound \eqref{PxivbzzBBroBinf} involving \( \CONSTgmb \). 

With \( \gaussv \) standard normal in \( \R^{\dimq} \), define
the \emph{effective trace} of \( \QP \) as
\begin{EQ}[rcl]
	\dimQ
	& \eqdef &
	\frac{\E \| \QP^{\T} \gaussv \|^{2}}{\| \QP \QP^{\T} \|}
	=
	\frac{\tr(\QP \QP^{\T})}{\| \QP \QP^{\T} \|} 
	\, .
\label{vrtgnfgih77jrfbegdhyd}
\end{EQ}
%
For \( \wv \in \R^{\dimp} \), define a measure \( \P_{\wv} \) and the corresponding expectation \( \E_{\wv} \) such that for any r.v. \( \eta \)
\begin{EQA}
	\E_{\wv} \, \eta
	& \eqdef &
	\frac{\E (\eta \, \ex^{\langle \wv, \Xv \rangle})}{\E \ex^{\langle \wv, \Xv \rangle}} \, .	
\label{hcxuyjhgjvgui85ww3fg}
\end{EQA}
Also fix some \( \gmn > 0 \) and introduce 
\begin{EQA}[rcl]
	\dltwu_{3}
	& \eqdef &
	\sup_{\| \wv \| \leq \gmn} \, \sup_{\uv \in \R^{\dimp}}
	\frac{1}{\| \uv \|^{3}}
	\bigl| 
		\E_{\wv} \langle \uv, \Xv - \E_{\wv} \Xv \rangle^{3} 
	\bigr| \, ,
\label{7bvmt3g8rf62hjgkhgu3}
 	\\
	\dltwu_{4}
	& \eqdef &
	\sup_{\| \wv \| \leq \gmn} \, \sup_{\uv \in \R^{\dimp}} 
	\frac{1}{\| \uv \|^{4}}
	\bigl| 
		\E_{\wv} \langle \uv, \Xv - \E_{\wv} \Xv \rangle^{4} 
		- 3 \bigl\{ \E_{\wv} \langle \uv, \Xv - \E_{\wv} \Xv \rangle^{2} \bigr\}^{2} 
	\bigr| \, .
	\qquad
\label{7bvmt3g8rf62hjgkhgu}
\end{EQA}
The quantities \( \dltwu_{3} \) and \( \dltwu_{4} \) depend on the distribution of \( \Xv \) and \( \gmn \).
However, they are typically not only finite but also very small.
E.g. for \( \Xv \) Gaussian they just vanish.
If \( \Xv \) is a normalized sum of independent centred random vectors \( \xiv_{1},\ldots,\xiv_{n} \) then 
\( \dltwu_{3} \asymp n^{- 1/2} \) and \( \dltwu_{4} \asymp n^{- 1} \); see Section~\ref{Ssumiiddb}. 

First, we present an upper bound which extends 
\ifapp{\eqref{Pxiv2dimAvp12}}{\eqref{v7hfte535ghewjfyw}} to the non-Gaussian case. 

\begin{theorem}
\label{Tdevboundsharp}
Let \( \Xv \) satisfy \( \E \Xv = 0 \) and condition \nameref{gmbref}.
For any linear mapping \( \QP \colon \R^{\dimp} \to \R^{\dimq} \), define \( \BBH = \QP \Var(\Xv) \QP^{\T} \).
Let \( \gmn \) and \( \dltwu_{3} \) from \eqref{7bvmt3g8rf62hjgkhgu3} fulfill 
\( \gmn^{2} \geq 3 \dimQ \) and \( \gmn \, \dltwu_{3} \leq 2/3 \).
Then for any \( \xx > 0 \) with 
\( \sqrt{4 \xx } \leq \sqrt{\tr (\BBH^{2})}/(3 \CONSTgmb \| \QP \QP^{\T} \|^{2}) \), it holds
\begin{EQA}
	\P\bigl( \| \QP \Xv \|^{2} > \tr (\BBH) + 2 \sqrt{\xx \tr (\BBH^{2})} + 2\xx \| \BBH \| \bigr)
	& \leq &
	(1 + \Delta_{\muH}) \ex^{-\xx} \, ,
	\qquad
\label{PxivbzzBBroBinfsh}
\end{EQA}
where  \( \muH = \muH(\xx) \) is given by 
\( \muH^{-1} = \| \BBH \| + \sqrt{ \tr (\BBH^{2})/(4\xx)} \) and \( \Delta_{\muH} \) 
depends on \( \dltwu_{3} \), \( \dltwu_{4} \), \( \dimQ \) only and
will be given explicitly in the proof.
Moreover, \( \Delta_{\muH} \ll 1 \) under \( \tr (\BBH^{2}) \gg \| \QP \QP^{\T} \|^{2} \) and 
\( (\dltwu_{3}^{2} + \dltwu_{4}) \, \dimQ^{2} \ll 1 \).
\end{theorem}

\begin{remark}
The statement of Theorem~\ref{Tdevboundsharp} looks a bit technical, however, 
the main message is straightforward and useful:
for moderate \( \xx \)-values, the Gaussian upper quantiles \( \tr (\BBH) + 2 \sqrt{\xx \tr (\BBH^{2})} + 2\xx \| \BBH \| \)
ensure the nominal deviation probability \( \ex^{-\xx} \) even if \( \Xv \) is not Gaussian.
\end{remark}

For getting lower deviation bounds, 
in place of condition \nameref{gmbref} on the moment-generating function \( \E \exp\bigl( \langle \uv,\Xv \rangle \bigr) \),
we need a condition on the characteristic function \( \E \exp\bigl( \imi \langle \uv,\Xv \rangle \bigr) \).
Namely, we assume that it does not vanish and its logarithm is bounded on the ball \( \| \uv \| \leq \gmb \).

\begin{description}
\item[\( \bb{(\charfX)} \)\label{gmbiref}]
\emph{For some fixed \( \gmb \) and \( \CONSTchar \), the function 
\( \charfX(\uv) = \log \E \, \ex^{ \imi \langle \uv, \Xv \rangle } \) satisfies}
\begin{EQA}
	|\charfX(\uv)|
	=
	|\log \E \, \ex^{ \imi \langle \uv, \Xv \rangle }|
	& \leq &
	\CONSTchar \, ,
	\qquad
	\| \uv \| \leq \gmb \, .
\label{devboundinfgmbi}
\end{EQA}
\end{description}

Note that this condition can easily be ensured by replacing \( \Xv \) with \( \Xv + \alp \gaussv \) for any positive \( \alp \)
and \( \gaussv \sim \ND(0,\Id_{\dimp}) \).
The constant \( \CONSTchar \) is unimportant, it does not show up in our results.
It, however, enables us to define similarly to \eqref{7bvmt3g8rf62hjgkhgu}
\begin{EQA}[rcl]
	\dltwu_{4}
	& \eqdef &
	\sup_{\| \wv \| \leq \gmn} \,\, \sup_{\uv \in \R^{\dimp}}
	\frac{1}{\| \uv \|^{4}}
	\bigl| 
		\E_{\imi \wv} \langle \imi \uv, \Xv - \E_{\imi \wv} \Xv \rangle^{4} 
		- 3 \bigl\{ \E_{\imi \wv} \langle \imi \uv, \Xv - \E_{\imi \wv} \Xv \rangle^{2} \bigr\}^{2} 
	\bigr| \, .
	\qquad \quad
\label{7bvmt3g8rf62hjgkhgui}
\end{EQA}
The values \( \dltwu_{4} \) in \eqref{7bvmt3g8rf62hjgkhgu} and \eqref{7bvmt3g8rf62hjgkhgui} might be different,
however, we use the same notation without risk of confusion.

\begin{theorem}
\label{Texpquadroi}
Let \( \Xv \) satisfy \( \E \Xv = 0 \), \( \Var(\Xv) \leq  \Id_{\dimp} \).
Let also \( \QP \colon \R^{\dimp} \to \R^{\dimq} \) be a linear mapping,
\( \dimQ = \tr(\QP \QP^{\T}) \), \( \BBH = \QP \Var(\Xv) \QP^{\T} \). 
Assume \nameref{gmbiref} for some \( \gmn \) with \( \gmn^{2} \geq 3 \dimQ^{2} \).
Let also \( \dltwu_{3} \) be given by \eqref{7bvmt3g8rf62hjgkhgu3} and 
\( \gmn \, \dltwu_{3} \leq 2/3 \).
Then for any \( \xx \leq \tr(\BBH^{2})/4 \) 
\begin{EQA}
	\P\bigl( \| \QP \Xv \|^{2} < \tr (\BBH) - 2 \sqrt{\xx \tr(\BBH^{2})} \bigr)
	& \leq & 
	(2 + \err + \rho_{\muH}) \ex^{-\xx} \, ,
\label{PxivbzzBBroBinfshi}
\end{EQA}
where \( \muH \eqdef 2 \sqrt{\xx/\tr(\BBH^{2})} \) and
\begin{EQA}
	\rho_{\muH} 
	\eqdef
	\P\biggl( \| \QP^{\T} \gaussv \|^{2} \geq \frac{4 \muH^{-1} \dimQ^{2}}{\tr(\BBH^{2})} \biggr) 
	& \leq &
	\exp \Bigl( - \frac{4 \dimQ^{2}}{\tr(\BBH^{2})} \Bigr)
	\, .
	\qquad
\label{u8cnyhbnkjmjoyt8re3}
\end{EQA}
The value \( \err \) is given in the proof of Proposition~\ref{Texpquadro} and it is small
under \( \tr (\BBH^{2}) \gg \| \QP \QP^{\T} \| \) and \( (\dltwu_{3}^{2} + \dltwu_{4}) \, \dimQ \ll 1 \).
\end{theorem}

\Subsection{Sum of i.i.d. random vectors}
\label{Ssumiiddb}
Here we specify the obtained results to the case when \( \Xv = n^{-1/2} \sumi \xiv_{i} \) and 
\( \xiv_{i} \) are i.i.d. in \( \R^{\dimp} \) with \( \E \xiv_{i} = 0 \) and \( \Var(\xiv_{i}) = \Sigma \leq \Id_{\dimp} \).
In fact, only independence of the \( \xiv_{i} \)'s is used provided that all the moment conditions later on
are satisfied uniformly over \( i \leq n \). 
However, the formulation is slightly simplified in the i.i.d case.
Let some \( \QP \colon \R^{\dimp} \to \R^{\dimq} \) be fixed.
With \( \BBH = \QP \Sigma \QP^{\T} \), it holds \( \dimH = \E \| \QP \Xv \|^{2} = \tr (\BBH) \).
We study the concentration property for \( \| \QP \Xv \|^{2} \).
The goal is to apply Theorem~\ref{Tdevboundsharp} and Theorem~\ref{Texpquadroi} 
claiming that \( \| \QP \Xv \|^{2} - \dimH \) can be sandwiched between 
\( - 2 \sqrt{\xx \tr(\BBH^{2})} \) and \( 2 \sqrt{\xx \tr(\BBH^{2})} + 2\xx \| \BBH \| \) with probability at least \( 1 - 3 \ex^{-\xx} \).
The major required condition is sub-gaussian behavior of \( \xiv_{1} \).
The conditions are summarized here.

\begin{description}
\item[\( \bb{(\xiv_{1})} \)\label{gmb1ref}]
\emph{A random vector \( \xiv_{1} \in \R^{\dimp} \) satisfies \( \E \xiv_{1} = 0 \), \( \Var(\xiv_{1}) = \Sigma \leq  \Id_{\dimp} \). 
Also}
\begin{enumerate}
	\item \emph{The function \( \cdens_{\xiv}(\uv) \eqdef \log \E \ex^{ \langle \uv, \xiv_{1} \rangle } \) is finite and fulfills 
for some \( \CONSTgmb \)}
\begin{EQA}
	&&
	\cdens_{\xiv}(\uv)
	\eqdef
	\log \E \ex^{ \langle \uv, \xiv_{1} \rangle }
	\leq 
	\frac{\CONSTgmb \| \uv \|^{2}}{2} \, ,
	\qquad
	\uv \in \R^{\dimp} \, .
\label{devboundinfgmb1}
\end{EQA}
	\item
	\emph{For \( \rhogmn > 0 \) and some constants \( \hmax_{3} \) and \( \hmax_{4} \), it holds
	with \( \E_{\wv} \) from \eqref{hcxuyjhgjvgui85ww3fg}}
\begin{EQA}
	&&
	\sup_{\| \wv \| \leq \rhogmn} \,\, \sup_{\uv \in \R^{\dimp}}
	\frac{1}{\| \uv \|^{3}}
	\bigl| 
		\E_{\wv} \langle \uv, \xiv_{1} \rangle^{3} 
	\bigr| 
	\leq 
	\hmax_{3} \, ;
\label{d98k3efu7yvb67r4hfidke}
	\\
	&&
	\sup_{\| \wv \| \leq \rhogmn} \,\, \sup_{\uv \in \R^{\dimp}}
	\frac{1}{\| \uv \|^{4}}
	\bigl| 
		\E_{\wv} \langle \uv, \xiv_{1} - \E_{\wv} \xiv_{1} \rangle^{4} 
		- 3 \bigl\{ \E_{\wv} \langle \uv, \xiv_{1} - \E_{\wv} \xiv_{1} \rangle^{2} \bigr\}^{2} 
	\bigr| 
	\leq 
	\hmax_{4} \, .
\label{7bvmt3g8rf62hjgkhguiid}
\end{EQA}

	\item \emph{The function \( \log \E \, \ex^{ \imi \langle \uv, \xiv_{1} \rangle } \) is well defined and}
\begin{EQ}[c]
	\sup_{\| \wv \| \leq \rhogmn} \,\, \sup_{\uv \in \R^{\dimp}}
	\frac{1}{\| \uv \|^{4}}
	\bigl| 
		\E_{\imi \wv} \langle \imi \uv, \xiv_{1} - \E_{\imi \wv} \xiv_{1} \rangle^{4} 
		- 3 \bigl\{ \E_{\imi \wv} \langle \imi \uv, \xiv_{1} - \E_{\imi \wv} \xiv_{1} \rangle^{2} \bigr\}^{2} 
	\bigr| 
	\leq 
	\hmax_{4} \, .
\label{7bvmt3g8rf62hjgkhguiidi}
\end{EQ}
\end{enumerate}
\end{description}

We are now well prepared to state the result for the i.i.d. case. 

\begin{theorem}
\label{TnormXiid}
Let \( \Xv = n^{-1/2} \sumi \xiv_{i} \), where \( \xiv_{i} \) are i.i.d. in \( \R^{\dimp} \) satisfying
\( \E \xiv_{1} = 0 \) and \( \Var(\xiv_{1}) = \Sigma \leq \Id_{\dimp} \), and condition \nameref{gmb1ref}.
For a fixed \( \QP \), assume \( n \rhogmn^{2} \geq 3 \dimQ \) and \( n \gg \dimQ^{2} \).
Then with \( \BBH = \QP \Sigma \QP^{\T} \), it holds
\begin{EQA}[lcl]
	\P\Bigl( \| \QP \Xv \|^{2} - \tr (\BBH) > 2 \sqrt{\xx \tr (\BBH^{2})} + 2\xx \| \BBH \| \Bigr)
	& \leq &
	(1 + \Delta_{\muH}) \ex^{-\xx} \, ,
	\quad \text{ if } \quad
	\sqrt{4 \xx } \leq \frac{\sqrt{\tr (\BBH^{2})}}{3 \CONSTgmb \| \QP \QP^{\T} \|} \, ,
\label{hvb7ruehft6wjwiqws}
	\\
	\P\Bigl( \| \QP \Xv \|^{2} - \tr (\BBH) < - 2 \sqrt{\xx \tr (\BBH^{2})} \Bigr)
	& \leq &
	(2 + \Delta_{\muH}) \ex^{-\xx} \, ,
	\quad \text{ if } \quad
	\xx \leq \frac{\tr(\BBH^{2})}{4 \| \QP \QP^{\T} \|^{2}} ,
\label{gc82jwjv7w3f5wetjwy}
\end{EQA}
where \( \Delta_{\muH} \lesssim n^{-1} \dimQ^{2} \).
\end{theorem}

\Subsection{Range of applicability, critical dimension}
\label{Srangesqnorm}
This section discusses the range of applicability of the presented results, in particular, 
of the concentration property.
It was already mentioned earlier that concentration of the squared norm \( \| \QP \Xv \|^{2} \) 
is only possible in a high dimensional situation, even for \( \Xv \) Gaussian.
This condition can be written as \( \tr(\BBH^{2})/\| \QP \QP^{\T} \|^{2} \gg 1 \).
In our results, this condition is further detailed.
For instance, bound \eqref{PxivbzzBBroBinfsh} of Theorem~\ref{Tdevboundsharp} is only meaningful 
if \( \tr (\BBH^{2}) \gg \CONSTgmb^{2} \| \QP \QP^{\T} \|^{2} \).
This is the only place where \( \CONSTgmb \) shows up.

Another important quantity is the value \( \Delta_{\muH} \).
It should be small to make the presented results meaningful.
A sufficient condition for this property are
\( (\dltwu_{3}^{2} + \dltwu_{4}) \, \dimQ^{2} \ll 1 \).
For the i.i.d. case, this condition transforms into 
``critical dimension'' condition \( \dimQ^{2} \ll n \).
%
Recent results from \cite{katsevich2023tight} indicate that Laplace approximation could fail if 
\( \dimQ^{2} \ll n \) is not fulfilled even for a simple generalized linear model.
One can guess that a further relaxation of the ``critical dimension'' condition \( \dimQ^{2} \ll n \) is not possible and approximation 
\( \P\bigl( \| \QP \Xv \| > \zq(\BBH,\xx) \bigr) \approx \P\bigl( \| \QP \Xvt \| > \zq(\BBH,\xx) \bigr) \)
with a standard Gaussian vector \( \Xvt \) can fail if \( \dimQ^{2} \gg n \).


\ifadap{
\Section{Proofs of the main results}
\label{SqfnGproofs}
This section collects the proofs of the main results from Section~\ref{SdevboundnonGauss}.
}{}
\Subsection{Proof of Theorem~\ref{Tdevboundinf}}

Let \( \gaussv \) be standard Gaussian in \( \R^{\dimq} \) under \( \Egs \) conditionally on \( \xiv \).
For \( \muH \in (0,1) \), 
\begin{EQA}
	\E \exp\bigl( \muH \| \xiv \|^{2} / 2 \bigr)
	&=&
	\E \, \Egs \, \exp\bigl( \muH^{1/2} \langle \HVB \gammav, \HVB^{-1} \xiv \rangle \bigr) ,
\label{Egexmu12lHm1B12}
\end{EQA} 
Application of Fubini's theorem, \eqref{devboundinf}, and \eqref{m2v241m4mj1p} yields
\begin{EQA}
	\E \exp\Bigl( \frac{\muH}{2} \| \xiv \|^{2} \Bigr)
	& \leq &
	\Egs \, \exp\Bigl( \frac{\muH}{2} \| \HVB \gammav \|^{2} \Bigr)
	\leq 
	\exp \Bigl( \frac{\muH^{2} \tr (\BBH^{2})}{4 (1 - \muH \| \BBH \|)} + \frac{\muH \tr (\BBH)}{2} \Bigr) .
\label{wBmu2vA241mmutrBi}
\end{EQA}
Now the bound follows by Theorem~\ref{TexpbLGA} as in the Gaussian case.

\Subsection{Proof of Theorem~\ref{Tdevboundsharp}}
Normalizing by \( \| \QP \| \) reduces the statement to \( \| \QP \| = 1 \) and \( \dimQ = \tr(\QP \QP^{\T}) \) 
which will be supposed later.
This also implies \( \| \BBH \| = \| \QP \Var(\Xv) \QP^{\T} \| \leq 1 \).
The key step of the proof is the following statement.

\begin{proposition}
\label{Texpquadro}
Assume the conditions of Theorem~\ref{Tdevboundsharp} and \( \| \QP \| = 1 \).
If \( \muH > 0 \) satisfies 
\begin{EQA}[c]
	\CONSTgmb \, \muH \leq 1/3 
	\, ,
\label{uckdfyg6h6t43hfgy2}
\end{EQA}
then it holds 
\begin{EQA}
	\bigl| \E \exp ( \muH \| \QP \Xv \|^{2}/2) - \det (\Id_{\dimq} - \muH \BBH)^{-1/2} \bigr|
	& \leq &
	\Delta_{\muH} \det (\Id_{\dimq} - \muH \BBH)^{-1/2} \, 
	\qquad
\label{jcxu785t83w5ffr4ehjk}
\end{EQA}
for some constant \( \Delta_{\muH} \) such that \( \Delta_{\muH} \ll 1 \) under 
\( \dimQ \gg 1 \), \( (\dltwu_{3}^{2} + \dltwu_{4}) \, \dimQ^{2} \ll 1 \);
see the proof for a closed-form representation.
\end{proposition}

\begin{proof}
We use \eqref{Egexmu12lHm1B12} and Fubini theorem: with \( \Egs = \E_{\gauss \sim \ND(0,\Id_{\dimq})} \)
\begin{EQA}
	\E \exp\bigl( \muH \| \QP \Xv \|^{2} / 2 \bigr)
	&=&
	\E \, \Egs \, \exp\bigl( \muH^{1/2} \langle \QP^{\T} \gaussv, \Xv \rangle \bigr) 
	=
	\Egs \, \exp \cdensX(\muH^{1/2} \QP^{\T} \gaussv).
	\qquad
	\qquad
\label{Egexmu12lHm1B12X}
\end{EQA}
Further, redefine \( \gmn^{2} = 3 \dimQ \) and apply the decomposition
\begin{EQA}
	\Egs \, \exp \cdensX(\muH^{1/2} \QP^{\T} \gaussv)
	&=&
	\Egs \, \exp \cdensX(\muH^{1/2} \QP^{\T} \gaussv) \Ind(\| \muH^{1/2} \QP^{\T} \gaussv \| \leq \gmn)
	\\
	&& 
	+ \, \Egs \, \exp \cdensX(\muH^{1/2} \QP^{\T} \gaussv) \Ind(\| \muH^{1/2} \QP^{\T} \gaussv \| > \gmn) .
	\qquad
\label{pocvuw3jes45q4wtgdhjkj}
\end{EQA}
Each summand here will be bounded separately starting from the second one.
Define
\begin{EQA}[rcccl]
	\zz_{\muH} 
	& \eqdef &
	\frac{1}{4} \Bigl( \sqrt{\CONSTgmb^{-1} \muH^{-1} \gmn^{2}} - \sqrt{\dimQ} \Bigr)^{2} ,
	\qquad
	\rexH_{\muH} 
	& \eqdef &  
	\CONSTgmb \, \muH + \CONSTgmb \, \muH \sqrt{\dimQ / \zz_{\muH}} \, .
\label{jhdfuyewjwkifvgu4rtgy}
\end{EQA}
Then \eqref{uckdfyg6h6t43hfgy2} ensures that \( \zz_{\muH} \geq \bigl( \sqrt{9 \dimQ} - \sqrt{\dimQ} \bigr)^{2}/4 = \dimQ \) 
and \( \rexH_{\muH} \leq 2/3 \).
By \eqref{devboundinfgmb} and \eqref{llkknbononjm9hig4e} of Theorem~\ref{CTexpbLGA}, it holds under 
the condition \( \rexH_{\muH} \leq 2/3 \) 
\begin{EQA}
	&& \nquad
	\Egs \, \exp \cdensX(\muH^{1/2} \QP^{\T} \gaussv) \Ind(\| \muH^{1/2} \QP^{\T} \gaussv \| > \gmn)
	\\
	& \leq &
	\Egs \, \exp\bigl( \CONSTgmb \, \muH \| \QP^{\T} \gaussv \|^{2}/2 \bigr) 
	\Ind(\| \QP^{\T} \gaussv \|^{2} > \muH^{-1} \gmn^{2})
	\\
	& = &
	\exp \bigl( \CONSTgmb \, \muH \, \dimQ/2 \bigr) \,
	\Egs \, \exp\bigl( \CONSTgmb \, \muH (\| \QP^{\T} \gaussv \|^{2} - \dimQ)/2 \bigr) 
	\Ind(\| \QP^{\T} \gaussv \|^{2} > \muH^{-1} \gmn^{2})
	\\
	& \leq &
	\frac{1}{1 - \rexH_{\muH}} \, \exp\{ \CONSTgmb \, \muH \, \dimQ/2 - (1 - \rexH_{\muH}) \zz_{\muH} \} \, .
\label{hbcyedtcdfbghe653ddfh}
\end{EQA}
Note that \( \rexH_{\muH} \leq 2/3 \), \( \zz_{\muH} \geq \dimQ \), and \( \CONSTgmb \, \muH \leq 1/3 \) imply
\begin{EQA}
	\frac{1}{1 - \rexH_{\muH}} \, \exp\{ \CONSTgmb \, \muH \, \dimQ/2 - (1 - \rexH_{\muH}) \zz_{\muH} \}
	& \leq &
	3 \ex^{- \dimQ/6} \, .
\label{y7vcu8eiu36gvhrujeivq}
\end{EQA}
This inequality helps to bound the second term of \eqref{pocvuw3jes45q4wtgdhjkj} corresponding to the event 
\( \{ \| \muH^{1/2} \QP^{\T} \gaussv \| > \gmn \} \).
For the first term, we apply the results on Laplace approximation from Section~\ref{SlocalLaplace}.
First we check that \( \cdensX(\uv) \) satisfies conditions \nameref{l2l3sref} and \nameref{l2l4ref}:
\begin{EQA}
	|\langle \nabla^{3} \cdensX(\xv), \uv^{\otimes 3} \rangle|
	& \leq &
	\dltwu_{3} \| \uv \|^{3} \, , \quad
	\uv \in \R^{\dimp} ,
\label{jcu8dfuerg74jufuiuer3}
\end{EQA}
and
\begin{EQA}
	|\dltw_{4}(\uv)|
	& \eqdef &
	\Bigl| \cdensX(\uv) - \frac{1}{2} \langle \cdensX''(0) ,\uv^{\otimes 2} \rangle 
			- \frac{1}{6} \langle \cdensX^{(3)}(0) ,\uv^{\otimes 3} \rangle 
	\Bigr|
	\leq 
	\frac{\dltwu_{4}}{24} \| \uv \|^{4} \, , \quad
	\| \uv \| \leq \gmn \, .
	\qquad
\label{jcu8dfuerg74jufuiuer4}
\end{EQA}
Let us start with the univariate case.
Let \( X \) satisfy \( \E X = 0 \) and \( \E X^{2} \leq \sigma^{2} \).
Define for any \( t \in [0,\gmn] \)
a measure \( \P_{t} \) such that for any random variable \( \eta \)
\begin{EQA}
	\E_{t} \, \eta
	& \eqdef &
	\frac{\E (\eta \, \ex^{t X})}{\E \ex^{t X}} 
	\, .
\label{cu4e37gurreughierwo}
\end{EQA}
Consider \( \cdensX(t) \eqdef \log \E \ex^{t X} \) as a function of \( t \in [0,\gmn] \).
It is well defined and satisfies \( \cdensX(0) = \cdensX'(0) = 0 \),
\( \cdensX''(0) = \E X^{2} \leq \sigma^{2} \), and
\begin{EQA}
	\cdensX^{(3)}(t) 
	&=& \E_{t} (X - \E_{t} X)^{3} \, ,
	\\
	\cdensX^{(4)}(t)
	&=&
	\E_{t} (X - \E_{t} X)^{4} - 3 \bigl\{ \E_{t} (X - \E_{t} X)^{2} \bigr\}^{2} \, .	
\label{hyhboewhye7y6gvu45u8r}
\end{EQA}
Therefore, conditions \nameref{l2l3sref} and \nameref{l2l4ref} follow from 
\eqref{7bvmt3g8rf62hjgkhgu3} and \eqref{7bvmt3g8rf62hjgkhgu}.
The multivariate case can be reduced to the univariate one by fixing a direction \( \uv \in \R^{\dimp} \)
and considering the function \( \cdensX(t \uv) \) of \( t \).

Next, we apply Proposition~\ref{PTensG42} to evaluate the first term on the right-hand side of \eqref{pocvuw3jes45q4wtgdhjkj}.
Define \( \WV = \{ \wv \in \R^{\dimq} \colon \| \muH^{1/2} \QP^{\T} \wv \| \leq \gmn \} \).
Then with \( \gaussv \sim \ND(0,\Id_{\dimq}) \)
\begin{EQA}
	\Egs \, \exp \cdensX(\muH^{1/2} \QP^{\T} \gaussv) \Ind(\| \muH^{1/2} \QP^{\T} \gaussv \| \leq \gmn)
	& = &
	\CONSTi_{\dimq} \int_{\WV} \ex^{\lgd_{\muH}(\wv)} \, d\wv \, ,
\label{g6hnwyvhwweybgfewwel}
\end{EQA}
where \( \CONSTi_{\dimq} = (2\pi)^{-\dimq/2} \) and for \( \wv \in \R^{\dimq} \)
\begin{EQA}
	\lgd_{\muH}(\wv)
	&=&
	\cdensX(\muH^{1/2} \QP^{\T} \wv) - \| \wv \|^{2}/2
\label{vdsfbcttqbcfccsvcgarfcv}
\end{EQA}
so that \( \lgd_{\muH}(0) = 0 \), \( \nabla \lgd_{\muH}(0) = 0 \).
Also, define \( \DPH^{2} = \muH \QP \QP^{\T} \),
\begin{EQA}
\label{vthwctatqcuhqdgdsdsu3}
	\DPH_{\muH}^{2} 
	& \eqdef &
	- \nabla^{2} \lgd_{\muH}(0) 
	=
	- \muH \QP \Var(\Xv) \QP^{\T} + \Id_{\dimq}
	=
	\Id_{\dimq} - \muH \BBH ,
	\\
	\dimL_{\muH}
	& \eqdef &
	\tr \bigl\{ \DPH_{\muH}^{-1} (\muH \QP \QP^{\T}) \DPH_{\muH}^{-1} \bigr\} 
	=
	\muH \, \tr (\DPH_{\muH}^{-2} \, \QP \QP^{\T}),
	\\
	\normG_{\muH}
	& \eqdef &
	\| \DPH_{\muH}^{-1} (\muH \QP \QP^{\T}) \DPH_{\muH}^{-1} \| 
	=
	\muH \, \| \DPH_{\muH}^{-2} \, \QP \QP^{\T} \| \, .
\label{yyheu7futjthb884fhf6e}
\end{EQA}
By \eqref{uckdfyg6h6t43hfgy2} \( \muH \leq \frac{1}{3\CONSTgmb} \leq \frac{1}{3} \) and \( \| \BBH \| \leq \| \QP \| = 1 \) implies 
\( (1 - \muH) \Id_{\dimp} \leq \DPH_{\muH}^{2} \leq \Id_{\dimp} \) so that
\begin{EQA}
	\dimL_{\muH}
	& \leq &
	\frac{\muH}{1 - \muH} \, \dimQ  
	\leq 
	\frac{\dimQ}{2} \, ,
	\qquad
	\normG_{\muH}
	\leq  
	\frac{\muH}{1 - \muH} 
	\leq 
	\frac{1}{2}
	\, .
\label{ycv8w3ikdfbohyietwqr4cwq}
\end{EQA}
The function \( \lgd_{\muH}(\wv) \) inherits smoothness properties of \( \cdensX(\muH^{1/2} \QP^{\T} \wv) \).
In particular,  for any \( \wv \) with \( \| \muH^{1/2} \QP^{\T} \wv \| \leq \gmn \)
\begin{EQA}
	\bigl| \langle \nabla^{3} \lgd_{\muH}(\wv) , \uv^{\otimes 3} \rangle \bigr|
	& \leq &
	\dltwu_{3} \| \muH^{1/2} \QP^{\T} \uv \|^{3} 
	\, ,
\label{ucie939f7fgr5urfcews}
	\\
	\bigl| \langle \nabla^{4} \lgd_{\muH}(\wv) , \uv^{\otimes 4} \rangle \bigr|
	& \leq &
	\dltwu_{4} \| \muH^{1/2} \QP^{\T} \uv \|^{4} 
	\, .
\label{98cvkdy63ug7rtjhhh}
\end{EQA}
Now Proposition~\ref{PTensG42} applied to \( \lgd_{\muH}(\wv) \) yields
\begin{EQA}[c]
	\biggl| 
		\frac{\int_{\WV} \ex^{\lgd_{\muH}(\wv)} \, d\wv 
				- \int_{\WV} \ex^{- \| \DPH_{\muH} \wv \|^{2}/2} \, d\wv}
			 {\int \ex^{- \| \DPH_{\muH} \wv \|^{2}/2} \, d\wv}
	\biggr|
	\leq 
	\err \, .
\label{ghshyhcn2ynnsdjfgnso}
\end{EQA}
The quantity \( \err \) here is computed as follows.
Let \( \Tens(\uv) = \langle \nabla^{3} \lgd_{\muH}(0), \uv^{\otimes 3} \rangle \), 
\( \gaussv_{\muH} \sim \ND(0,\DPH_{\muH}^{-2}) \).
In view of \( \gmn^{2} = 3 \dimQ \) and \eqref{ycv8w3ikdfbohyietwqr4cwq}, it holds
\begin{EQA}
	\grad_{\muH}
	&=&
	\frac{\dltwu_{3} \, \gmn^{2} \sqrt{\normG_{\muH}}}{2} 
	\leq 
	\dltwu_{3} \, \dimQ \, ,
	\\
	\accu_{\muH}^{2} 
	&=& 
	\E \Tens^{2}(\gaussv_{\muH})
	\leq 
	\sqrt{5/12} \, \dltwu_{3} \, \dimL_{\muH}
	\leq 
	\frac{1}{3} \, \dltwu_{3} \, \dimQ \, ,
	\\
	\dltw_{4,\muH}
	& = &
	\EUV \dltw_{4}^{2}(\gaussv_{\muH})
	\leq 
	\frac{1}{24} \, \dltwu_{4} (\dimL_{\muH} + 3 \normG_{\muH})^{2} 
	\leq 
	\frac{1}{96} \, \dltwu_{4} (\dimQ + 3)^{2} \, .
\label{ydfje348bujtuy5r75t}
\end{EQA}
Then
\begin{EQA}[rcccl]
	\Bigl| \err - \frac{\accu_{\muH}^{2}}{2} \Bigr|
	& \leq &
	\accu_{\muH} \, \dltw_{4,\muH} 
	+ \frac{\dltw_{4,\muH}^{2}}{2} + \frac{5}{3} \grad_{\muH}^{3} \exp(\grad_{\muH}^{2}) \, ,
	\;\;
	\err
	& \leq &
	\frac{1}{2} (\accu_{\muH} + \dltw_{4,\muH})^{2} + \frac{5}{3} \grad_{\muH}^{3} \exp (\grad_{\muH}^{2}) \, .
	\qquad \quad
\label{vuf7fbfygtytewgetejj}
\end{EQA}
Furthermore, 
\begin{EQA}
	\rho_{\muH}
	& \eqdef &
	1 -	\frac{\int_{\WV} \ex^{- \| \DPH_{\muH} \wv \|^{2}/2} \, d\wv}{\int \ex^{- \| \DPH_{\muH} \wv \|^{2}/2} \, d\wv} 
	=
	\P\bigl( \| \sqrt{\muH} \QP^{\T} \DPH_{\muH}^{-1} \gaussv \| > \gmn \bigr)
	\, .
	\qquad
\label{yywnscthwsyd553t23e}
\end{EQA}
The use of \( \muH \leq 1/(3\CONSTgmb) \leq 1/3 \) and \( (1 - \muH) \muH^{-1} \gmn^{2} \geq 6 \CONSTgmb \, \dimQ \geq 6 \dimQ \) 
yields
\begin{EQA}
	\rho_{\muH}
	& \leq &
	\P\bigl( \| \QP^{\T} \gaussv \|^{2} > \frac{(1 - \muH) \gmn^{2}}{\muH} \bigr) 
	\leq 
	\P\bigl( \| \QP^{\T} \gaussv \|^{2} > 6 \dimQ \bigr)
	\leq 
	\ex^{- \dimQ/2} \, ;
\label{8cfkw9vire78vjweiowsck7}
\end{EQA}
see \eqref{3emzmsp22z2} of Theorem~\ref{CTexpbLGA} with \( \yy = 5 \dimQ \).
By \eqref{ghshyhcn2ynnsdjfgnso} and \eqref{yywnscthwsyd553t23e}
\begin{EQA}
	\biggl| \frac{\int_{\WV} \ex^{\lgd_{\muH}(\wv)} \, d\wv}{\int \ex^{- \| \DPH_{\muH} \wv \|^{2}/2} \, d\wv} - 1 \biggr|
	& \leq &
	\err + \rho_{\muH} \, .
\label{obmfurfwubyehjeuy}
\end{EQA}
It remains to be noted that
\begin{EQA}
	\CONSTi_{\dimq} \int \ex^{- \| \DPH_{\muH} \wv \|^{2}/2} \, d\wv
	&=&
	\frac{1}{\det \DPH_{\muH}}
	=
	\det(\Id_{\dimq} - \muH \BBH)^{-1/2} 
	\leq 
	1
\label{ivmeyh8gu467gsbewrde}
\end{EQA}
and \eqref{jcxu785t83w5ffr4ehjk} follows from \eqref{hbcyedtcdfbghe653ddfh} and \eqref{obmfurfwubyehjeuy} with
\begin{EQA}
	\Delta_{\muH}
	& \leq &
	\err + \rho_{\muH} + 
	\frac{1}{1 - \rexH_{\muH}} \, \exp\{ \CONSTgmb \, \muH \, \dimQ/2 - (1 - \rexH_{\muH}) \zz_{\muH} \} \, .
	\qquad
\label{jcxu785t83w5ffr4ehjkde}
\end{EQA}
Moreover,  for \( \dimQ \) large, \( \rho_{\muH} \) is small by \eqref{8cfkw9vire78vjweiowsck7},
the exp-term can be bounded by \eqref{y7vcu8eiu36gvhrujeivq},  
while \( \err \) is small provided that \( (\dltwu_{3}^{2} + \dltwu_{4}) \, \dimQ^{2} \) is small.
\end{proof}

Now we can finalize the proof of Theorem~\ref{Tdevboundsharp}.
Upper deviation bounds for \( \| \QP \Xv \|^{2} \) can be derived as in the Gaussian case 
by applying \eqref{jcxu785t83w5ffr4ehjk} with a proper choice of \( \muH \).
Let \( \xx \) satisfy \( \sqrt{4 \xx } \leq \sqrt{\tr (\BBH^{2})}/(3 \CONSTgmb \| \BBH \|) \).
We check \eqref{uckdfyg6h6t43hfgy2} for \( \muH = \muH(\xx) \).
Indeed, the definition \( \muH^{-1} = \| \BBH \| + \sqrt{ \tr (\BBH^{2})/(4\xx)} \) 
implies \( \muH \leq \sqrt{4 \xx / \tr (\BBH^{2})} \). 
Therefore, \( \sqrt{4 \xx } \leq \sqrt{\tr (\BBH^{2})}/(3 \CONSTgmb) \) yields
\( \muH \leq 1/(3 \CONSTgmb) \) and \eqref{uckdfyg6h6t43hfgy2} is fulfilled for \( \gmn^{2} = 3 \dimQ \).
The bound \eqref{PxivbzzBBroBinfsh} follows from \eqref{jcxu785t83w5ffr4ehjk} as in the Gaussian case of Theorem~\ref{TexpbLGA}.

\Subsection{Proof of Theorem~\ref{Texpquadroi}}

This result is based on an approximation 
\( \E \ex^{ - \muH \| \QP \Xv \|^{2}/2} \approx \det (\Id_{\dimq} + \muH \BBH)^{-1/2} \).
This requires an analog of \eqref{jcxu785t83w5ffr4ehjk} for \( \muH \) negative.
With \( \imi = \sqrt{-1} \), the use of \eqref{Egexmu12lHm1B12X} yields
\begin{EQA}
	\E \, \ex^{ - \muH \| \QP \Xv \|^{2} / 2} 
	&=&
	\E \, \Egs \, \ex^{ \imi \sqrt{\muH} \langle \QP^{\T} \gaussv, \Xv \rangle} 
	=
	\Egs \E \, \ex^{ \imi \sqrt{\muH} \langle \QP^{\T} \gaussv,\Xv \rangle} \, .
\label{Egexmu12lHm1B12Xch}
\end{EQA}

\begin{proposition}
\label{Pexpquadroi}
Assume the conditions of Theorem~\ref{Texpquadroi}.
For any \( \muH \in (0,1) \), it holds with \( \BBH = \QP \Var(\Xv) \QP^{\T} \)
\begin{EQ}[rcl]
	\bigl| \E \ex^{ - \muH \| \QP \Xv \|^{2}/2} - \det (\Id_{\dimq} + \muH \BBH)^{-1/2} \bigr|
	& \leq &
	(\err + \rho_{\muH}) \det (\Id_{\dimq} + \muH \BBH)^{-1/2} + \rho_{\muH} \, ,
	\\
	\rho_{\muH} 
	\leq  
	\Pgs\bigl( \| \QP^{\T} \gaussv \|^{2} \geq 4 \muH^{-1} \dimQ \bigr)
	& \leq &
	\exp \Bigl\{ - \frac{\dimQ}{4} (2\muH^{-1/2} - 1)^{2} \Bigr\} \, .
\label{jcxu785t83w5ffr4ehjki}
\end{EQ}
\end{proposition}

\begin{proof}
We follow the proof of Proposition~\ref{Texpquadro} replacing everywhere \( \cdensX(\uv) \) with \( \charfX(\uv) \).
In particular, we start with representation \eqref{Egexmu12lHm1B12Xch} and apply with \( \gmn^{2} = 3 \dimQ \)
\begin{EQA}
	&& \nquad
	\E \, \ex^{ - \muH \| \QP \Xv \|^{2} / 2}
	=
	\Egs \, \ex^{\charfX(\sqrt{\muH} \, \QP^{\T} \gaussv)}
	\\
	&=&
	\Egs \, \ex^{\charfX(\sqrt{\muH} \, \QP^{\T} \gaussv)} \Ind(\| \sqrt{\muH} \, \QP^{\T} \gaussv \| \leq \gmn)
	+ \Egs \, \ex^{\charfX(\sqrt{\muH} \, \QP^{\T} \gaussv)} \Ind(\| \sqrt{\muH} \, \QP^{\T} \gaussv \| > \gmn) .
\label{pocvuw3jes45q4wtgdhjkji}
\end{EQA}
It holds 
\begin{EQA}
	\charfX(0)
	&=&
	0,
	\quad
	\nabla \charfX(0)
	=
	0,
	\quad
	- \nabla^{2} \charfX(0) = \Var(\Xv) \leq \Id_{\dimp} \, .
\label{ycywhwf5e4g3eerter}
\end{EQA}
Moreover, smoothness conditions \eqref{jcu8dfuerg74jufuiuer3}, \eqref{jcu8dfuerg74jufuiuer4} 
are automatically fulfilled for \( \charfX(\uv) \)
with the same \( \dltwu_{3} \) and \( \dltwu_{4} \) from \eqref{7bvmt3g8rf62hjgkhgui}.
The most important observation for the proof is that the bound \eqref{obmfurfwubyehjeuy} continues to apply for 
\begin{EQA}
	\lgd_{\muH}(\wv)
	&=&
	\charfX(\sqrt{\muH} \, \QP^{\T} \wv) - \| \wv \|^{2}/2,
\label{ujc563hfg8h75urdwyehfcn}
\end{EQA}
with \( \err \) from \eqref{vuf7fbfygtytewgetejj} and 
\begin{EQA}[rcccl]
	\DPH_{\muH}^{2} 
	& \eqdef &
	- \nabla^{2} \lgd_{\muH}(0) 
	&=&
	\muH \QP \Var(\Xv) \QP^{\T} + \Id_{\dimq}
	=
	\Id_{\dimq} + \muH \BBH ,
	\\
	\dimL_{\muH}
	& \eqdef &
	\tr \bigl\{ \DPH_{\muH}^{-2} (\muH \QP \QP^{\T}) \bigr\}
	&\leq &
	\frac{\muH}{1 + \muH} \tr (\QP \QP^{\T}) 
	\leq 
	\muH \, \dimQ \, ,
	\\
	\normG_{\muH}
	& \eqdef &
	\| \DPH_{\muH}^{-1} (\muH \QP \QP^{\T}) \DPH_{\muH}^{-1} \| 
	&\leq &
	\frac{\muH}{1 + \muH} 
	\leq 
	\muH \, ,
\label{vthwctatqcuhqdgdsdsu3i}
\end{EQA}
and \( \rho_{\muH} \leq \P\bigl( \| \QP^{\T} \gaussv \|^{2} \geq \muH^{-1} \gmn^{2} \bigr) \); cf. \eqref{yywnscthwsyd553t23e}.
This yields
\begin{EQA}
	\biggl| \Egs \, \ex^{\charfX(\sqrt{\muH} \, \QP^{\T} \gaussv)} \Ind(\| \sqrt{\muH} \, \QP^{\T} \gaussv \| \leq \gmn) 
		- \frac{1}{\det (\Id_{\dimq} + \muH \BBH)^{1/2}} 
	\biggr|
	& \leq &
	\frac{\err + \rho_{\muH}}{\det (\Id_{\dimq} + \muH \BBH)^{1/2}} \, .
\label{0vkr7vhfh4twhswtfhry4y7}
\end{EQA}
Finally we use \( |\ex^{\charfX(\uv)}| \leq 1 \) and thus,
\begin{EQA}
	\bigl| \Egs \, \ex^{\charfX(\sqrt{\muH} \, \QP^{\T} \gaussv)} \Ind(\| \sqrt{\muH} \, \QP^{\T} \gaussv \| > \gmn) \bigr|
	& \leq &
	\P\bigl( \| \sqrt{\muH} \, \QP^{\T} \gaussv \| > \gmn \bigr) 
	=
	\Pgs\bigl( \| \QP^{\T} \gaussv \|^{2} \geq 4 \muH^{-1} \dimQ \bigr)
\label{hdveyhwedvuwuwejdv1334}
\end{EQA}
and \eqref{jcxu785t83w5ffr4ehjki} follows in view of \eqref{3emzmsp22z2} of Theorem~\ref{CTexpbLGA}.
\end{proof}

The proof of Theorem~\ref{Texpquadroi} is similar to the case of Theorem~\ref{Tdevboundsharp}.
By the exponential Chebyshev inequality and \eqref{jcxu785t83w5ffr4ehjki}, it holds with \( \vH^{2} = \tr(\BBH^{2}) \)
\begin{EQA}
	&& \nquad
	\P\bigl( \tr (\BBH) - \| \QP \Xv \|^{2} > 2 \vH \sqrt{\xx} \bigr)
	\leq
	\exp( - \muH \, \vH \sqrt{\xx}) \E \exp \bigl\{ \muH \tr (\BBH)/2 - \muH \| \QP \Xv \|^{2} /2 \bigr\}
	\\
	& \leq &
	\exp(\muH \tr (\BBH)/2 - \muH \, \vH \sqrt{\xx}) 
	\bigl\{ (1 + \err + \rho_{\muH}) \det (\Id_{\dimq} + \muH \BBH)^{-1/2} + \rho_{\muH} \bigr\} .
\label{cuyfr7ygbyuby6r4yfwete}
\end{EQA}
In view of \( x - \log(1+x) \leq x^{2}/2 \) and \( \muH = 2 \vH^{-1} \sqrt{\xx} \), 
as in the proof of Lemma~\ref{Lqfexpmom}
\begin{EQA}
	- \muH \, \vH \sqrt{\xx} + \muH \tr (\BBH)/2 + \log \det(\Id_{\dimq} + \muH \BBH)^{-1/2}
	& \leq &
	- \muH \, \vH \sqrt{\xx} + \muH^{2} \vH^{2}/4
	=
	- \xx .
\label{d7ywh3dfuyt523fhfyes}
\end{EQA}
Also \( \muH \tr (\BBH)/2 - \muH \, \vH \sqrt{\xx} = \vH^{-1} \tr (\BBH) \, \sqrt{\xx} - 2 \xx \leq \vH^{-1} \dimQ \, \sqrt{\xx} - 2 \xx \).
The bound on \( \rho_{\muH} \) in \eqref{u8cnyhbnkjmjoyt8re3} follows from \eqref{3emzmsp22z2} of Theorem~\ref{CTexpbLGA}
in view of \( \dimQ \geq \vH^{2} \) and hence, \( \dimQ \leq \dimQ^{2}/\vH^{2} \).
Finally, observe that
\begin{EQA}
	\rho_{\muH} \, \exp\bigl( \vH^{-1} \dimQ \, \sqrt{\xx} - 2 \xx \bigr)
	& \leq &
	\exp\Bigl( - \frac{\dimQ^{2}}{4\vH^{2}} + \frac{\dimQ \, \sqrt{\xx}}{\vH} - 2 \xx \Bigr)
	\\
	& \leq &
	\exp\Bigl\{ - \Bigl( \frac{\dimQ}{2 \vH} - \sqrt{\xx} \Bigr)^{2} - \xx \Bigr\}
	\leq 
	\ex^{-\xx} \, 
\label{juf7wjejfc6erdyrjvgbiw3}
\end{EQA}
and \eqref{PxivbzzBBroBinfshi} follows as well.

\Subsection{Proof of Theorem~\ref{TnormXiid}}
The definition  and i.i.d structure of the \( \xiv_{i} \)'s yield 
\( \E \langle \uv,\Xv \rangle^{2} =	\E \langle \uv,\xiv_{1} \rangle^{2} \) and 
\begin{EQA}
	\cdensX(\uv)
	&=&
	\log \E \ex^{\langle \Xv,\uv \rangle}
	=
	n \cdens_{\xiv}(n^{-1/2}\uv) 
\label{ufiu3jgibnhi5keuvytr}
\end{EQA}
for any \( \uv \in \R^{\dimp} \),
where \( \cdens_{\xiv}(\uv) \eqdef \log \E \ex^{\langle \xiv_{1},\uv \rangle} \). 
For the derivatives \( \cdensX^{(k)}(\uv) \) of \( \cdensX(\uv) \), this yields
\begin{EQA}
	\cdensX^{(k)}(\uv)
	&=&
	n^{1 - k/2} \cdens_{\xiv}^{(k)}(n^{-1/2}\uv) .
\label{ufiu3jgibnhi5keuvytr}
\end{EQA}
This enables us to derive \eqref{7bvmt3g8rf62hjgkhgu3} and \eqref{7bvmt3g8rf62hjgkhgu} from \nameref{gmb1ref}
for any \( \gmn \) with \( \gmn/\sqrt{n} \leq \rhogmn \) and
\begin{EQA}
	\dltwu_{3}
	& = &
	n^{-1/2} \hmax_{3} \, ,
	\qquad
	\dltwu_{4}
	= 
	n^{-1} \hmax_{4} \, .
\label{gv8iejkwidwweewfg}
\end{EQA}
Moreover, the quantity \( \err \) from \eqref{vuf7fbfygtytewgetejj} satisfies \( \err \lesssim \dimQ^{2}/n \).
Now the upper bound follows from Theorem~\ref{Tdevboundsharp}.
Similar arguments can be used for checking the lower bound by Theorem~\ref{Texpquadroi}.


\def\zqx{\zq_{\circ}}

\Section{Deviation bounds under light exponential tails}
\label{Sdevboundexp}

Let \( \xiv \) be a zero mean random vector in \( \R^{\dimp} \) with covariance \( \Var(\xiv) \) and 
let \( \QP \colon \R^{\dimp} \to \R^{\dimq} \) be a linear mapping.
This section presents some deviation bounds on the norm \( \| \QP \xiv \| \) for the case of light exponential tails of \( \xiv \).
Namely, 

\begin{description}
\item[\label{gmref} \( \bb{(\gmb)} \)]
	\textit{for some fixed \( \gmb > 0 \) and some self-adjoint operator \( \HVB^{2} \) in \( \R^{\dimp} \) with} 
	\( \HVB^{2} \geq \Var(\xiv) \), 
\begin{EQA}[c]
    \cdens(\uv)
    \eqdef
    \log \E \exp\bigl( \langle \uv, \HVB^{-1} \xiv \rangle \bigr)
    \le
    \frac{\| \uv \|^{2}}{2} \, ,
    \qquad
    \uv \in \R^{\dimp}, \, \| \uv \| \le \gmb ,
\label{expgamgm}
\end{EQA}
\end{description}

In fact, it is sufficient to assume that 
\begin{EQA}
	\sup_{\| \uv \| \leq \gmb} \E \exp\bigl( \langle \uv, \HVB^{-1} \xiv \rangle \bigr)
	& \leq &
	\CONST .
\label{sgagEexlgHx}
\end{EQA}
The quantity \( \CONST \) can be very large but it is not important.
Indeed, the function \( \cdens(\uv) \) is analytic on the disk \( \| \uv \| \leq \gmb \), 
and condition \eqref{sgagEexlgHx} implies an analog of \eqref{expgamgm}: 
\begin{EQA}
	\cdens(\uv) 
	& \leq &
	\frac{\| \uv \|^{2}}{2} +\frac{\dltwu_{3} \| \uv \|^{3}}{6} 
	\leq 
	\frac{\| \uv \|^{2}}{2} \Bigl( 1 + \frac{\dltwu_{3} \gmn}{3} \Bigr)
	\, ,
	\qquad
	\| \uv \| \leq \gmn \, ,
\label{0hy544rtgte7hjuwghbj}
\end{EQA}
for a fixed value \( \dltwu_{3} \).
Moreover, reducing \( \gmb \) allows to take \( \HVB^{2} \) equal or close to \( \Var(\xiv) \) and \( \dltwu_{3} \)
close to zero.
The next section presents our main results under \nameref{gmref}.
The proofs are postponed until the end of the section.

\Subsection{Main results}
Let a random vector \( \xiv \) satisfy \( \E \xiv = 0 \) and \nameref{gmref}. 
The goal is to establish possibly sharp deviation bounds on \( \| \QP \xiv \|^{2} \)
for a given linear mapping \( \QP \colon \R^{\dimp} \to \R^{\dimq} \).
Define 
\begin{EQ}[rcl]
	\BBH 
	& \eqdef &
	\QP \HVB^{2} \QP^{\T}, \quad
	\dimH 
	\eqdef  
	\tr(\BBH) , 
	\quad
	\vH^{2} 
	\eqdef
	\tr(\BBH^{2}) , 
	\quad 
	\supH 
	\eqdef 
	\| \BBH \| ,
	\\
	\zq^{2}(\BBH,\xx) 
	& \eqdef & 
	\tr \BBH + 2 \sqrt{\xx \tr(\BBH^{2})} + 2 \xx \| \BBH \|
	=
	\dimH + 2 \vH \sqrt{\xx} + 2 \xx \supH .
\label{h5ete5f5tegvhvcdvhyendf}
\end{EQ} 
Also fix some \( \rhoH < 1 \), a standard choice is \( \rhoH = 1/2 \).
Our main result applies for all \( \xx \) satisfying the condition
\begin{EQA}
	\zq^{2}(\BBH,\xx)
	& \leq &
	\rhoH \, \biggl(\frac{\gmb \sqrt{\supH}}{\muH(\xx)} - \sqrt{\frac{\dimH}{\muH(\xx)}} \biggr)^{2}
\label{kv7367ehjgruwwcewyd}
\end{EQA}
with \( \zq(\BBH,\xx) \) from \eqref{h5ete5f5tegvhvcdvhyendf} and
\( \muH(\xx) \) defined by \( \muH^{-1}(\xx) = 1 + \frac{\vH}{2 \supH \sqrt{\xx}} \);
see \eqref{1v2sxm12m1m}.
One can see that the left hand-side of \eqref{kv7367ehjgruwwcewyd} increases with \( \xx \) while 
the right hand-side decreases.
Therefore, there exists a unique root \( \xxc \) such that with \( \muHc = \muH(\xxc) \)
\begin{EQA}
	\zq^{2}(\BBH,\xxc)
	& = &
	\rhoH \, \biggl(\frac{\gmb \sqrt{\supH}}{\muHc} - \sqrt{\frac{\dimH}{\muHc}} \biggr)^{2} .
\label{kv7367ehjgruwwcewyde}
\end{EQA}
The value \( \xxc \) is important, it describes the \emph{phase transition} effect:
the upper quantile function of \( \| \QP \xiv \| \) exhibits the Gaussian-like behavior for \( \xx \leq  \xxc \), 
while it grows linearly with \( \xx/\gmb \) for \( \xx > \xxc \) as in a sub-exponential case.


\begin{theorem}
\label{Tdevboundgm}
Assume \nameref{gmref}.
Fix \( \xxc \) by \eqref{kv7367ehjgruwwcewyde} for some 
\( \rhoH \leq 1/2 \).
It holds 
\begin{EQA}
    \P\bigl( \| \QP \xiv \| \ge \zq(\BBH,\xx) \bigr)
    & \le &
    3 \ex^{-\xx} ,
    \qquad
    \xx \leq \xxc \, .
\label{PxivbzzBBroB}
\end{EQA}    
For \( \rhoH = 1/2 \), the value \( \xxc \) from \eqref{kv7367ehjgruwwcewyde} fulfills 
\begin{EQA}
	\frac{1}{4} \biggl( \gmb - \sqrt{\frac{2\dimH}{\supH}} \biggr)_{+}^{2}
	\leq 
	\xxc 
	& \leq &
	\frac{\gmb^{2}}{4} \, .
\label{if7h3rhgy4676rfhdsjw}
\end{EQA}
If \( \gmb > \sqrt{2 \dimH / \supH} \) then \( \zqc = \zq(\BBH,\xxc) \) follows
\begin{EQA}
	\gmb \sqrt{\supH /2} - (1 - 2^{-1/2}) \sqrt{\dimH} 
	\leq 
	\zqc
	& \leq &
	\gmb \sqrt{\supH /2} + \sqrt{\dimH} \, .
\label{f9oi4eogjdgtvgyuj4ek}
\end{EQA}
\end{theorem}


The results of Theorem~\ref{Tdevboundgm} state nearly Gaussian deviation bounds for 
the norm of the vector \( \QP \xiv \) satisfying \nameref{gmref}. 
Namely, the Gaussian deviation bound \( \P\bigl( \| \QP \xiv \| \ge \zq(\BBH,\xx) \bigr) \le \ex^{-\xx} \)
from Theorem~\ref{TexpbLGA} applies with the additional factor 3 for all \( \xx \leq \xxc \).
Condition \( \gmb \gg \sqrt{\dimH / \supH} \) is important.
Otherwise, the value \( \xxc \) is not significantly large and the zone \( \xx \leq \xxc \) 
with Gaussian-like quantiles is too narrow.
It turns out that out of this range, the norm \( \| \QP \xiv \| \) exhibits a sub-exponential behavior.

\begin{theorem}
\label{TQPxivlarge}
Assume \nameref{gmref}.
With \( \xxc \) from \eqref{kv7367ehjgruwwcewyde} and \( \zqc = \zq(\BBH,\xxc) \),
set \( \constg = \frac{\sqrt{\rhoH} \, \gmb}{(2 + \sqrt{\rhoH}) \sqrt{\supH}} \).
It holds 
\begin{EQ}[lll]
	&
	\P\bigl( \| \QP \xiv \| > \zqc + \constg^{-1} (\xx - \xxc) \bigr)
	\leq 
	3 \ex^{-\xx} ,
	&
	\quad
	\xx \geq \xxc \, ,
	\\
	&
	\P\bigl( \| \QP \xiv \| > \zq \bigr)
	\leq
	3 \exp\{ -\xxc - \constg (\zq - \zqc) \} ,
	&
	\quad
	\zq \geq \zqc \, .
\label{9vkjv6njbih9t69t3wfg}
\end{EQ}
\end{theorem}

The obtained deviation bounds of Theorem~\ref{Tdevboundgm} and Theorem~\ref{TQPxivlarge}
can be fused into one.
To be more specific, we fix \( \rhoH = 1/2 \).

\begin{corollary}
\label{CTQPxivlarge}
Assume \nameref{gmref}. 
Let \( \xxc \) be defined by \eqref{kv7367ehjgruwwcewyde} with \( \rhoH = 1/2 \).
For all \( \xx > 0 \)
\begin{EQA}
	\P\bigl( \| \QP \xiv \| > \zqc(\BBH,\xx) \bigr)
	& \leq &
	3 \ex^{-\xx} ,
\label{uyfuyerd7e7uhh8yy689t}
\end{EQA}
where with \( \constg \eqdef \frac{\gmb}{(\sqrt{8} + 1) \sqrt{\supH}} \)
and \( \xx \wedge \xxc \eqdef \min\{ \xx , \xxc \} \)
\begin{EQA}
	\zqc(\BBH,\xx)
	& \eqdef &
	\zq(\BBH,\xx \wedge \xxc)
	+ \constg^{-1} (\xx - \xxc)_{+}
	=
	\begin{cases}
		\zq(\BBH,\xx), & \xx \leq \xxc \, ,
		\\
		\zq(\BBH,\xxc) + \dfrac{\xx - \xxc}{\constg} , & \xx > \xxc \, .
	\end{cases}
	\qquad
	\qquad
\label{dyv6ejf8gjwkerih83}
\end{EQA}
Moreover, \( \xxc \) follows \eqref{if7h3rhgy4676rfhdsjw} and \( \zqc = \zq(\BBH,\xxc) \) 
satisfies \eqref{f9oi4eogjdgtvgyuj4ek} provided \( \gmb \geq \sqrt{2\dimH / \supH} \).
\end{corollary}

If \( \gmb \gg \sqrt{\dimH / \supH} \) then \( \xxc \) is large and 
\( \zqc(\BBH,\xx) = \zq(\BBH,\xx) \leq \sqrt{\dimH} + \sqrt{2 \xx \supH} \) for all reasonable \( \xx \).
For \( \gmb < \sqrt{2 \dimH / \supH} \), the accurate bound \eqref{dyv6ejf8gjwkerih83} can be simplified by a linear majorant
which does not involve \( \xxc \).

\begin{theorem}
\label{Llin3maj}
Assume \nameref{gmref}. 
Fix \( \constg = \frac{\gmb }{(\sqrt{8} + 1) \sqrt{\supH}} \).
Then \eqref{uyfuyerd7e7uhh8yy689t} applies with
\begin{EQA}
	\zqc(\BBH,\xx)
	& \leq &
	\sqrt{\dimH} + \frac{\constg}{\sqrt{2}} + \constg^{-1} \xx \, .
\label{uv7ey4hb743h3hiti4ffd}
\end{EQA}
\end{theorem}

The next result provides some upper bounds on the exponential moments of \( \| \QP \xiv \| \).
We distinguish between zones \( \zq \leq \zqc \) and \( \zq > \zqc \) with \( \zqc = \zq(\BBH,\xxc) \);
see \eqref{kv7367ehjgruwwcewyde}.

\begin{theorem}
\label{TQPlargex}
Assume \nameref{gmref}.
Let \( \xxc \) fulfill \eqref{kv7367ehjgruwwcewyde} and \( \zqc = \zq(\BBH,\xxc) \).
For any \( \zq \in [\sqrt{\dimH}, \zqc] \) and any \( \nuH \leq \frac{\zq - \sqrt{\dimH}}{2\sqrt{\supH}} \), it holds
\begin{EQA}
	\E \ex^{\nuH \| \QP \xiv \|} \Ind(\| \QP \xiv \| \geq \zq)
	& \leq &
	6 \exp\Bigl\{  \nuH \zq - \frac{(\zq - \sqrt{\dimH})^{2}}{2 \supH} \Bigr\} .
\label{ufikwk3ei9vgkj4k4giw3wl}
\end{EQA}
Further, for any \( \nuH < \constg \eqdef \frac{\gmb \sqrt{\rhoH}}{\sqrt{\supH} \, (2 + \sqrt{\rhoH})} \)
\begin{EQA}
	\E \ex^{\nuH \| \QP \xiv \| } \Ind(\| \QP \xiv \| > \zqc)
	& \leq &
	\frac{3 \constg}{\constg - \nuH} \, 
	\exp\biggl\{  \nuH \zqc - \frac{(\zqc - \sqrt{\dimH})^{2}}{2 \supH} \biggr\} \, .
\label{jhf7yehruybyrhe3wevire}
\end{EQA}
Moreover, for \( \zq \geq \zqc \)
\begin{EQA}
	\E \ex^{\nuH \| \QP \xiv \| } \Ind(\| \QP \xiv \| > \zq)
	& \leq &
	\frac{3 \constg}{\constg - \nuH} \, 
	\exp \Bigl\{ \nuH \zqc - \frac{(\zqc - \sqrt{\dimH})^{2}}{2 \supH} - (\constg - \nuH) (\zq - \zqc) \Bigr\} \, .
	\qquad
\label{jhf7yehruybyrhe3wevire2}
\end{EQA}
\end{theorem}

\Subsection{Proof of Theorem~\ref{Tdevboundgm}}
By normalization, one can easily reduce the study to the case \( \| \BBH \| = 1 \).
Moreover, replacing \( \xiv \) with \( \HVB^{-1} \xiv \) and \( \QP \) with \( \QP \HVB \) reduces the proof to the situation 
with \( \HVB = \Id_{\dimp} \).
This will be assumed later on.
%
For \( \muH \in (0,1) \) and \( \zzH(\muH) = \gmb / \muH - \sqrt{\dimH/\muH} > 0 \),
define trimming \( t_{\muH}(\uv) \) of \( \uv \in \R^{\dimp} \) as
\begin{EQA}
	t_{\muH}(\uv)
	& \eqdef &
	\begin{cases}
		\uv, & \text{ if } \| \uv \| \leq \zzH(\muH), 
		\\
		\frac{\zzH(\muH)}{\| \uv \|} \, \uv, & \text{ otherwise}. 
	\end{cases}
\label{ojchcy63eyfvey3gcbfkg}
\end{EQA}
By construction \( \| t_{\muH}(\uv) \| \leq \zzH(\muH) \) for all \( \uv \in \R^{\dimp} \). 

\begin{lemma}
\label{LGDBqfexpB}
Assume \nameref{gmref} and let \( \| \BBH \| = 1 \).
Fix \( \muH \in (0,1) \) s.t. \( \zzH(\muH) = \gmb / \muH - \sqrt{\dimH/\muH} > 0 \).
Then with \( t_{\muH}(\cdot) \) from \eqref{ojchcy63eyfvey3gcbfkg} 
\begin{EQA}
	\E \exp\Bigl\{ \frac{\muH}{2} \, t_{\muH}^{2}(\QP \xiv) \Bigr\}
	& \leq &
	2 \exp\{ \Pmuvp(\muH) \} ,
	\qquad
\label{wBmu2vA241mmutrB}
\end{EQA}
where
\begin{EQA}
	\Pmuvp(\muH)
	& \eqdef &
	\frac{\muH^{2} \vH^{2}}{4 (1 - \muH)} + \frac{\muH \, \dimH}{2} \, .
\label{sdfw6rewqrwdwtffwet66}
\end{EQA}
Furthermore, for any \( \zzH < \zzH(\muH) \)
\begin{EQA}
	\P\bigl( \| \QP \xiv \| > \zzH, \| \QP \xiv \| \leq \zzH(\muH) \bigr)
	& \leq &
	2 \exp\Bigl\{ - \frac{\muH \, \zzH^{2}}{2} + \Pmuvp(\muH) \Bigr\} .
\label{9cdkcf736eryghj7y34wwde}
\end{EQA}
\end{lemma}

\begin{proof}
Let us fix any value of \( \xiv \).
We intend to show that
\begin{EQA}
	\exp\Bigl\{ \frac{\muH}{2} \, \| t_{\muH}(\QP \xiv) \|^{2} \Bigr\}
	& \leq &
	2 \Egs \, \exp\{ \muH^{1/2} \gaussv^{\T} t_{\muH}(\QP \xiv) \} .
\label{0xcujfc5nfiuf76w3hfc}
\end{EQA}
Here \( \Egs \) means conditional expectation w.r.t. \( \gaussv \sim \ND(0,\Id_{\dimp}) \) given \( \xiv \).
Obviously, with \( A = \{ \uv \colon \muH^{1/2} \| \QP^{\T} \uv \| \leq \gmb \} \), 
it suffices to check that
\begin{EQA}
	\II_{\muH}(\xiv)
	& \eqdef &
	\Egs \exp \Bigl\{ \muH^{1/2} \gaussv^{\T} t_{\muH}(\QP \xiv) - \frac{\muH}{2} \| t_{\muH}(\QP \xiv) \|^{2} \Bigr\} 
	\Ind(\gaussv \in A)
	\geq 
	1/2 .
\label{ggfiokhyhiu8y6eeeewe}
\end{EQA}
With \( \CONSTi_{\dimp} = (2\pi)^{-\dimp/2} \), it holds
\begin{EQA}
	\II_{\muH}(\xiv)
	&=&
	\CONSTi_{\dimp} \int_{A} \exp \Bigl( 
		\muH^{1/2} \uv^{\T} t_{\muH}(\QP \xiv) 
		- \frac{\muH}{2} \| t_{\muH}(\QP \xiv) \|^{2} - \frac{1}{2} \| \uv \|^{2} 
	\Bigr) d\uv
	\\
	&=&
	\CONSTi_{\dimp} \int_{A} \exp \Bigl( - \frac{1}{2} \| \uv - \muH^{1/2} t_{\muH}(\QP \xiv) \|^{2} \Bigr) d\uv
	=
	\Pgs(\gaussv - \muH^{1/2} t_{\muH}(\QP \xiv) \in A ) .
\label{kj8vfuff7f7f7ywwew}
\end{EQA}
The definition of \( A \) and the condition \( \| t_{\muH}(\QP \xiv) \| \leq \zzH(\muH) \) imply
in view of \( \| \QP \| \leq 1 \)
\begin{EQA}
	&& \nquad
	\Pgs(\gaussv - \muH^{1/2} t_{\muH}(\QP \xiv) \in A )
	=
	\Pgs\bigl( \| \QP^{\T} (\gaussv - \muH^{1/2} t_{\muH}(\QP \xiv)) \| \leq \gmb/\muH^{1/2} \bigr)
	\\
	& \geq &
	\Pgs\bigl( \| \QP^{\T} \gaussv \| \leq \gmb/\muH^{1/2} - \muH^{1/2} \zzH(\muH) \bigr)
	\geq 
	\Pgs\bigl( \| \QP^{\T} \gaussv \| \leq \sqrt{\dimH} \bigr)
	\geq 
	1/2
\label{mnriuvu8frt7r77r7urfu}
\end{EQA}
and \eqref{ggfiokhyhiu8y6eeeewe} follows.
Taking expectation for both sides of \eqref{0xcujfc5nfiuf76w3hfc} and the use of Fubini's theorem yield
\begin{EQA}
	&& \nquad
	\E \exp\Bigl\{ \frac{\muH}{2} \, \| t_{\muH}(\QP \xiv) \|^{2} \Bigr\} 
	\leq 
	2 \Egs \bigl\{ \E \exp\{ \muH^{1/2} \gaussv^{\T} t_{\muH}(\QP \xiv) \} 
	\, \Ind(\muH^{1/2} \| \QP^{\T} \gaussv \| \leq \gmb) \bigr\} .
\label{icjfcufry7fvc6e32w33wr}
\end{EQA}
Obviously, for any \( \uv \in \R^{\dimp} \)
\begin{EQA}
	\exp\{ \uv^{\T} t_{\muH}(\QP \xiv) \} + \exp\{ - \uv^{\T} t_{\muH}(\QP \xiv) \} 
	& \leq &
	\exp\{ \uv^{\T} \QP \xiv \} + \exp\{ - \uv^{\T} \QP \xiv \}
\label{odfu87e76376rhfyt6wh}
\end{EQA} 
and by \eqref{expgamgm} 
\begin{EQA}
	\E \exp\Bigl\{ \frac{\muH}{2} \, \| t_{\muH}(\QP \xiv) \|^{2} \Bigr\}
	& \leq &
	2 \Egs \Bigl\{ \exp\Bigl( \frac{1}{2} \, \| \muH^{1/2} \gaussv^{\T} \QP \|^{2} \Bigr) \, 
		\Ind(\muH^{1/2} \| \QP^{\T} \gaussv \| \leq \gmb) \Bigr\}
	\\
	& \leq &
	2 \Egs \, \exp\Bigl( \frac{1}{2} \, \| \muH^{1/2} \gaussv^{\T} \QP \|^{2} \Bigr)
	=
	2 \det(\Id_{\dimp} - \muH \QP^{\T} \QP)^{-1/2} .
\label{8ji89656hk0wchdne83}
\end{EQA}
We also use that for any \( \muH > 0 \) by \eqref{m2v241m4mj1p},
\begin{EQA}
	\log \det\bigl( \Id - \muH \BBH \bigr)^{-1/2} 
	& \leq &
	\frac{\muH \tr (\BBH)}{2} + \frac{\muH^{2} \tr (\BBH^{2})}{4 (1 - \muH)} 
	=
	\Pmuvp(\muH) \, ,
\label{mu2v241mmiulogIm12}
\end{EQA}
and the first statement follows.
Moreover, by Markov's inequality 
\begin{EQA}
	\P\bigl( \| \QP \xiv \| > \zzH, \| \QP \xiv \| \leq \zzH(\muH) \bigr)
	& \leq &
	\ex^{- \muH \, \zzH^{2}/2} \E \exp\Bigl\{ \frac{\muH}{2} \, \| t_{\muH}(\QP \xiv) \|^{2} \Bigr\} 
	\leq 
	2 \exp\Bigl\{ - \frac{\muH \, \zzH^{2}}{2} + \Pmuvp(\muH) \Bigr\} ,
\label{9cdkcf736eryfvkler4895ty34wwde}
\end{EQA}
and \eqref{9cdkcf736eryghj7y34wwde} follows as well.
\end{proof}

The use of \( \muH = \muH(\xx) \) from \eqref{1v2sxm12m1m} in \eqref{wBmu2vA241mmutrB} yields 
\begin{EQA}
	- \frac{\muH \zq^{2}(\BBH,\xx)}{2} + \Pmuvp(\muH)
	&=&
	- \xx \, ,
\label{uv73hf8er74e7eereew}
\end{EQA}
and similarly to the proof of Theorem~\ref{TexpbLGA}
\begin{EQA}
	\P\Bigl( \| \QP \xiv \|^{2} > \zq^{2}(\BBH,\xx), \,
		\| \QP \xiv \| \leq \zzH(\muH)
	 \Bigr)
	 & \leq &
	 2 \ex^{-\xx} .
\label{2emxPblrHm1B}
\end{EQA}
It remains to consider the probability of large deviation
\( \P\bigl( \| \QP \xiv \| > \zzH(\muH) \bigr) \).

\begin{lemma}
\label{Ldvbetagmb}
Assume \( \| \BBH \| = 1 \).
Given \( \xx > 0 \), fix \( \muH = \muH(\xx) \) and \( \zzH(\muH) = \gmb/\muH - \sqrt{\dimH/\muH} \).
Assume \eqref{kv7367ehjgruwwcewyd} for some \( \rhoH \leq 1/2 \).
Then 
\begin{EQA}
	\P\bigl( \| \QP \xiv \| > \zzH(\muH) \bigr)
	& \leq &
	\ex^{-\xx} .	
\label{yvy3n3eubvuj2tcy3he}
\end{EQA}
\end{lemma}

\begin{proof}
Denote \( \eta = \| \QP \xiv \| \).
By \eqref{2emxPblrHm1B}
\begin{EQA}
	\P\Bigl( \eta > \zq(\BBH,\xx), \,
		\eta \leq \zzH(\muH)
	 \Bigr)
	 & \leq &
	 2 \ex^{-\xx} ,
\label{me22Iezmemcz2}
\end{EQA}
For \( \muH = \muH(\xx) \), it holds \eqref{uv73hf8er74e7eereew} with \( \Pmuvp(\muH) \) 
given by \eqref{sdfw6rewqrwdwtffwet66}.
Bounding the tails of \( \eta \) in the region \( \eta > \zzH(\muH) \) requires another choice of \( \muH \).
Namely, we apply \eqref{9cdkcf736eryghj7y34wwde} with \( \rhoH \muH \) instead of \( \muH \) yielding
\begin{EQA}
	\P\bigl( \eta > \zzH(\muH), \eta \leq \zz(\rhoH \muH) \bigr)
	& \leq &
	2 \exp\Bigl\{ - \frac{\rhoH \muH \, \zz^{2}(\muH)}{2} + \Pmuvp(\rhoH \muH) \Bigr\} .
\label{uv7e3j2kv88eu3e3536}
\end{EQA}
In a similar way, applying \eqref{me22Iezmemcz2} with \( \rhoH^{2} \muH \) in place of \( \muH \)
and using that 
\begin{EQA}
	\rhoH \, \zz(\rhoH \muH)
	&=&
	\gmb/\muH - \sqrt{\rhoH \, \dimH/\muH}
	\leq 
	\zzH(\muH) 
\label{f9i2mg9gj3g2553}
\end{EQA}
yields
\begin{EQA}
	&& \nquad
	\P\bigl( \eta > \zz(\rhoH \muH), \eta \leq \zz(\rhoH^{2} \muH) \bigr)
	\leq 
	2 \exp\Bigl\{ - \frac{\rhoH^{2} \muH \, \zz^{2}(\rhoH \muH)}{2} + \Pmuvp(\rhoH^{2} \muH) \Bigr\} 
	\\
	& \leq &
	2 \exp\Bigl\{ - \frac{\muH \, \zz^{2}(\muH)}{2} + \Pmuvp(\rhoH^{2} \muH) \Bigr\} .
\label{t6dxhywywybcxcttw2q}
\end{EQA}
This trick can be applied again and again yielding in view of \eqref{f9i2mg9gj3g2553}
\begin{EQA}
	\P\bigl( \eta > \zzH(\muH) \bigr)
	& \leq &
	\sum_{k=0}^{\infty} \P\bigl( \eta > \zz(\rhoH^{k} \muH), \eta \leq \zz(\rhoH^{k+1} \muH) \bigr)
	\\
	& \leq &
	\sum_{k=0}^{\infty} 2 \exp\bigl\{ - \rhoH^{k+1} \muH \, \zz^{2}(\rhoH^{k} \muH)/2 + \Pmuvp(\rhoH^{k+1} \muH) \bigr\} 
	\\
	& \leq &
	\sum_{k=0}^{\infty}  2 \exp\bigl\{ - \rhoH^{-k+1} \muH \, \zz^{2}(\muH)/2 + \Pmuvp(\rhoH^{k+1} \muH) \bigr\} .
\label{uenruvueju33j2223figgi}
\end{EQA}
Condition \( \rhoH \, \zz^{2}(\muH) \geq \zq^{2}(\BBH,\muH)/2 \) and \eqref{uv73hf8er74e7eereew} ensure
for \( \rhoH \leq 1/2 \)
\begin{EQA}
	\P\bigl( \eta > \zzH(\muH) \bigr)
	& \leq &
	\sum_{k=0}^{\infty} 2 \exp\bigl\{ - \rhoH^{-k} \muH \, \zq^{2}(\BBH,\muH)/2 + \Pmuvp(\rhoH^{k+1} \muH) \bigr\}
	\\
	& \leq &	
	2 \sum_{k=0}^{\infty}  \exp\bigl\{ \Pmuvp(\rhoH^{k+1} \muH) - \rhoH^{-k} \Pmuvp(\muH) - \rhoH^{-k} \xx \bigr\} 
	\leq 
	\ex^{-\xx} .	
\label{yvy3n3eubvuj2tyche}
\end{EQA}
This yields \eqref{yvy3n3eubvuj2tcy3he}.
\end{proof}
Putting together \eqref{2emxPblrHm1B} and \eqref{yvy3n3eubvuj2tcy3he} yields \eqref{PxivbzzBBroB}.

Now we check \eqref{if7h3rhgy4676rfhdsjw}.
Normalization by \( \supH \) reduces the proof to the case \( \| \BBH \| = \| \QP \HVB^{2} \QP^{\T} \| = 1 \).
We use the simplified bounds 
\( \zq(\BBH,\xx) \leq \sqrt{\dimH} + \sqrt{2\xx} \) and \( \muH^{-1} = 1 + \sqrt{\dimH/(4\xx)} \). 
Now \eqref{kv7367ehjgruwwcewyd} with \( \rhoH = 1/2 \) can be rewritten as
\begin{EQA}
	\gmb
	& \geq &
	\sqrt{\muH \, \dimH} + \muH \sqrt{2} \bigl( \sqrt{\dimH} + \sqrt{2\xx} \bigr) .
\label{hdf6eh3ye4545dwwe}
\end{EQA}
The use of \( \muH = \sqrt{4 \xx} / (\sqrt{4\xx} + \sqrt{\dimH}) \) yields
\begin{EQA}
	\muH \sqrt{2} \bigl( \sqrt{\dimH} + \sqrt{2\xx} \bigr)
	&=&
	\sqrt{8\xx} \,\, \frac{\sqrt{\dimH} + \sqrt{2\xx}}{\sqrt{\dimH} + \sqrt{4\xx}}
	\geq \sqrt{4\xx} \, ,
\label{uviuedu8e347y4rt87hiuy}
\end{EQA}
and \eqref{hdf6eh3ye4545dwwe} is not possible for \( \xx > \gmb^{2}/4 \).
Further, with \( \yy = \sqrt{4\xx}/\gmb \) and \( \alp = \sqrt{\dimH} / \gmb \)
\begin{EQA}
	\frac{\sqrt{\muH \, \dimH} + \muH \sqrt{2} \bigl( \sqrt{\dimH} + \sqrt{2\xx} \bigr)}{\gmb}
	&=&
	\sqrt{\frac{\yy \alp^{2}}{\alp+\yy}} + \frac{\yy (\sqrt{2} \alp + \yy)}{\alp+\yy}
	\leq 
	\alp + \yy + \frac{\yy(\sqrt{2} - 1) \alp}{\alp+\yy} 
	\leq 
	\yy + \sqrt{2} \alp .
\label{jhfy7637ur8thj986eewe}
\end{EQA}
Together with \eqref{hdf6eh3ye4545dwwe}, this yields \( \yy \geq 1 - \sqrt{2} \alp \) and \eqref{if7h3rhgy4676rfhdsjw} follows.
For \eqref{f9oi4eogjdgtvgyuj4ek} we use \( \zqc \leq \sqrt{\dimH} + \sqrt{2 \supH \xxc} \)
and \( \zqc \geq \sqrt{\dimH/2} + \sqrt{2 \supH \xxc} \).

\Subsection{Proof of Theorem~\ref{TQPxivlarge}}
Assume w.l.o.g. \( \supH = 1 \).
First we present an accurate deviation bound, which, however, does not provide a closed form quantile function 
for \( \| \QP \xiv \| \).
Then we show how it implies a rough linear upper bound on this quantile function.
For \( \xxc \) from \eqref{kv7367ehjgruwwcewyde} and \( \xx > \xxc \), fix \( \muH \) by the relation   
\begin{EQA}
	\frac{\rhoH \, \muH \, \zzH^{2}(\muH)}{2} 
	&=& 
	\xx + \Pmuvp(\muH) 
	=
	\xx + \frac{\muH \, \dimH}{2} + \frac{\muH^{2} \vH^{2}}{4 (1 - \muH)} \, ,
\label{jejfye4ye3hfnw3jfeu}
\end{EQA}
where \( \zzH(\muH) = \gmb/\muH - \sqrt{\dimH/\muH} \)%
; cf. \eqref{uv73hf8er74e7eereew}.
It is easy to see that the solution \( \muH \) exists and unique.
Moreover, if \( \xx = \xxc \) then \( \muH = \muHc \) and \( \zzH^{2}(\muHc) = \zq^{2}(\BBH,\xxc) \);
see \eqref{kv7367ehjgruwwcewyde}.
If \( \xx > \xxc \), then \( \muH < \muHc \) and \( \zzH^{2}(\muH) > \zq^{2}(\BBH,\xx) \).

\begin{lemma}
\label{Ldevbosub}
For \( \xx > \xxc \), define \( \muH \) by \eqref{jejfye4ye3hfnw3jfeu}.
Then with \( \zzH(\muH) = \gmb/ \muH - \sqrt{\dimH/\muH} \)
\begin{EQA}
	\P\bigl( \| \QP \xiv \|^{2} > \rhoH \, \zzH^{2}(\muH) \bigr)
	& \leq &
	3 \ex^{-\xx} \, .
\label{uvw7erf38sjcwtet23w3}
\end{EQA}
\end{lemma}
\begin{proof}
We again apply Lemma~\ref{LGDBqfexpB}, however, the choice \( \muH = \muH(\xx) \) from \eqref{1v2sxm12m1m} 
is not possible anymore in view of \( \zq(\BBH,\xx) > \zzH(\muH) \). 
More precisely, for \( \xx \) large, the value \( \muH(\xx) \) approaches one and this choice of \( \muH \) 
yields the value \( \zzH(\muH) \) smaller than we need.
To cope with this problem, we apply \eqref{9cdkcf736eryghj7y34wwde} of Lemma~\ref{LGDBqfexpB} with a sub-optimal \( \muH \) 
from \eqref{jejfye4ye3hfnw3jfeu} ensuring 
\( \rhoH \muH \, \zzH^{2}(\mu) - \Pmuvp(\muH) = \xx \).
By \eqref{9cdkcf736eryghj7y34wwde} of Lemma~\ref{LGDBqfexpB} 
\begin{EQA}
	\P\bigl( \| \QP \xiv \| > \sqrt{\rhoH} \, \zzH(\muH), \| \QP \xiv \| \leq \zzH(\muH) \bigr)
	& \leq &
	2 \exp\Bigl\{ - \frac{\rhoH \muH \, \zzH^{2}(\muH)}{2} + \Pmuvp(\muH) \bigr\}
	=
	2 \ex^{-\xx} .
\label{mcyer6e46rf62wuvyrddde}
\end{EQA}
Repeating the arguments from the proof of Lemma~\ref{Ldvbetagmb} implies
\begin{EQA}
	&& \nquad
	\P\bigl( \| \QP \xiv \|^{2} > \rhoH \, \zzH^{2}(\muH) \bigr)
	\leq 
	\sum_{k=0}^{\infty} 2 \exp\Bigl\{ - \frac{1}{2} \rhoH^{k+1} \muH \, \zz^{2}(\rhoH^{k} \muH) + \Pmuvp(\rhoH^{k} \muH) \Bigr\} 
	\\
	& \leq &
	\sum_{k=0}^{\infty}  2 \exp\Bigl\{ - \frac{1}{2} \rhoH^{-k+1} \muH \, \zz^{2}(\muH) + \Pmuvp(\rhoH^{k} \muH) \Bigr\} 
	\\
	& \leq &
	2 \ex^{-\xx} 
	+ 2 \ex^{-\xx} \sum_{k=1}^{\infty}  
		\exp\Bigl\{ - \frac{1}{2} (\rhoH^{-k} - 1) \rhoH \muH \, \zz^{2}(\muH) + \Pmuvp(\rhoH^{k} \muH) - \Pmuvp(\muH) \Bigr\} 
	\leq 
	3 \ex^{-\xx} .
\label{uenruvueju3kH223figgi}
\end{EQA}
as stated in \eqref{uvw7erf38sjcwtet23w3}.
\end{proof}

\noindent
It remains to evaluate \( \rhoH \, \zzH^{2}(\muH) \) with \( \muH \) from \eqref{jejfye4ye3hfnw3jfeu} and
\( \zzH(\muH) = \gmb/\muH - \sqrt{\dimH/\muH} \).
For \( \muH \leq \muHc \)
\begin{EQA}
	\frac{\rhoH}{2} \biggl( \frac{\gmb}{\sqrt{\muH}} - \sqrt{\dimH} \biggr)^{2}
	&=&
	\xx + \Pmuvp(\muH)
\label{uedu7euyewhfvhfrt4qww}
\end{EQA}
and 
\begin{EQA}
	\frac{\sqrt{\rhoH} \, \gmb}{\sqrt{\muH}}
	& = &
	\sqrt{2\xx + 2\Pmuvp(\muH)} + \sqrt{\rhoH \, \dimH} 
	\, .
\label{iujk37bujreyvjwtcjeww}
\end{EQA}
This results in
\begin{EQA}
	\sqrt{\rhoH} \, \zzH(\muH)
	& = &
	\frac{\sqrt{\rhoH}}{\sqrt{\muH}} \biggl( \frac{\gmb}{\sqrt{\muH}} - \sqrt{\dimH} \biggr)
	\leq 
	\frac{1}{\sqrt{\rhoH} \, \gmb} \bigl( \sqrt{2\xx + 2\Pmuvp(\muH)} + \sqrt{\rhoH \, \dimH} \bigr) \, \sqrt{2\xx + 2\Pmuvp(\muH)} 
	\\
	& \leq &
	\frac{1}{\sqrt{\rhoH} \, \gmb} 
	\bigl( 2\xx + 2\Pmuvp(\muHc) + \sqrt{\rhoH \, \dimH(2\xx + 2\Pmuvp(\muHc))} \bigr) 
	\eqdef
	\zqs(\xx)
	\, .
\label{hv7jwwsuvjqs89c83jejfy}
\end{EQA}
By \eqref{kv7367ehjgruwwcewyde}, this inequality becomes equality for \( \xx = \xxc \) and \( \muH = \muHc \) with 
\( \sqrt{\rhoH} \, \zzH(\muHc) = \zqs(\xxc) = \zq(\BBH,\xxc) \).
Furthermore, the derivative of \( \zqs(\xx) \) w.r.t. \( \xx \) satisfies
\begin{EQA}
	\frac{d}{d\xx} \, \zqs(\xx)
	&=&
	\frac{1}{\sqrt{\rhoH} \, \gmb} \biggl( 2 + \frac{\sqrt{\rhoH \, \dimH}}{\sqrt{2\xx + 2\Pmuvp(\muHc)}} \biggr)
	\leq 
	\frac{1}{\sqrt{\rhoH} \, \gmb} \biggl( 2 + \frac{\sqrt{\rhoH \, \dimH}}{\sqrt{2\xxc + 2\Pmuvp(\muHc)}} \biggr) \, .
\label{igk373jfgho6952wq2q2}
\end{EQA}
Moreover, \( 2\xxc + 2\Pmuvp(\muHc) = \zq^{2}(\BBH,\xxc) \) and
\begin{EQA}
	\frac{d}{d\xx} \, \zqs(\xx)
	& \leq &
	\frac{1}{\sqrt{\rhoH} \, \gmb} \biggl( 2 + \frac{\sqrt{\rhoH \, \dimH}}{\zq(\BBH,\xxc)} \biggr)
	\leq 
	\frac{2 + \sqrt{\rhoH}}{\sqrt{\rhoH} \, \gmb} 
\label{h8jkew3fbo9rfkw3tehw32h}
\end{EQA}
yielding
\begin{EQA}
	\zqs(\xx)
	& \leq &
	\zqs(\xxc) + \frac{2 + \sqrt{\rhoH}}{\sqrt{\rhoH} \, \gmb} (\xx - \xxc)
	=
	\zq(\BBH,\xxc) + \frac{2 + \sqrt{\rhoH}}{\sqrt{\rhoH} \, \gmb} (\xx - \xxc)
\label{ujw289fjrej3gyhr2wmdx}
\end{EQA}
and hence, 
\begin{EQA}
	\sqrt{\rhoH} \, \zzH(\muH)
	& \leq &
	\zq(\BBH,\xxc) + \frac{2 + \sqrt{\rhoH}}{\sqrt{\rhoH} \, \gmb} (\xx - \xxc) 
	=
	\zqc + \frac{\xx - \xxc}{\constg} \, .
\label{9fj36vhje37gjgr33erfg}
\end{EQA}
This implies \eqref{9vkjv6njbih9t69t3wfg}.
%

\Subsection{Proof of Theorem~\ref{Llin3maj}}
As previously, assume \( \supH = 1 \).
We use \( \zq(\BBH,\xxc) \leq \sqrt{\dimH} + \sqrt{2 \xxc} \).
Further, \( \constg^{-1} \xxc - \sqrt{2\xxc} + \constg/\sqrt{2} \geq 0 \) and thus,
\begin{EQA}
	\sqrt{2\xxc} - \constg^{-1} \xxc
	& \leq &
	\constg/\sqrt{2} \, .
\label{chf63ghfnb74hrmbgogh}
\end{EQA}
Therefore, for \( \xx \geq \xxc \), it holds 
\begin{EQA}
	\zqc(\BBH,\xx)
	& = &
	\zq(\BBH,\xxc) + \frac{\xx - \xxc}{\constg}
	\leq 
	\sqrt{\dimH} + \sqrt{2 \xxc} - \frac{\xxc}{\constg} + \frac{\xx}{\constg}
	\leq 
	\sqrt{\dimH} + \frac{\constg}{\sqrt{2}} + \frac{\xx}{\constg} \, .
\label{hvtwhnfkdhjwifmbye}
\end{EQA}
In the zone \( \xx \leq \xxc \), it holds \( \zqc(\BBH,\xx) = \zq(\BBH,\xx)	\leq \sqrt{\dimH} + \sqrt{2 \xx} \)
and it remains to note that
\( \sqrt{2 \xx} \leq \constg/\sqrt{2} + \constg^{-1} \xx \).

\Subsection{Proof of Theorem~\ref{TQPlargex}}
Assume w.o.l.g. \( \supH = 1 \).
First consider \( \zq \geq \zqc \).
By \eqref{9fj36vhje37gjgr33erfg} of Theorem~\ref{TQPxivlarge}, 
it holds with \( \constg = \gmb \sqrt{\rhoH} /(2 + \sqrt{\rhoH}) \) 
and \( \xxc = (\zqc - \sqrt{\dimH})^{2}/2 \)
\begin{EQA}
	\P\bigl( \| \QP \xiv \| \geq \zq \bigr)
	& = &
	\P\bigl( \| \QP \xiv \| \geq \zqc + \zq - \zqc \bigr)
	\leq 
	3 \ex^{- \xxc - \constg (\zq - \zqc)} \, .
\label{jcu7wj3fibh99y6ky6hbgi8ed}
\end{EQA}
In particular, \( \P( \| \QP \xiv \| \geq \zqc) \leq 3 \ex^{-\xxc} \).
Integration by parts yields for \( \nuH < \constg \)
\begin{EQA}
	&& \nquad
	\E \ex^{ \nuH (\| \QP \xiv \| - \zqc) } \Ind(\| \QP \xiv \| > \zqc)
	= 
	- \int_{\zqc}^{\infty} \ex^{\nuH (\zq - \zqc)} d\P(\| \QP \xiv \| \geq \zq)
	\\
	&=&
	\P(\| \QP \xiv \| \geq \zqc)
	+ \nuH \int_{\zqc}^{\infty} \ex^{\nuH (\zq - \zqc)} \, \P(\| \QP \xiv \| \geq \zq) \, d\zq 
	\\
	& \leq &
	3 \ex^{-\xxc} + \nuH \int_{\zqc}^{\infty} \ex^{\nuH (\zq - \zqc)- \xxc - \constg (\zq - \zqc)} \, d\zq
	=
	\left( 3 + \frac{3 \nuH}{\constg - \nuH} \right) \, \ex^{-\xxc} 
	\qquad
\label{hgdft6gh3wef7u7ruj543}
\end{EQA}
and \eqref{jhf7yehruybyrhe3wevire} follows.
Similarly, for \( \zq \geq \zqc \), we derive  \eqref{jhf7yehruybyrhe3wevire2} as follows
\begin{EQA}
	&& \nquad
	\E \ex^{ \nuH \| \QP \xiv \| } \Ind(\| \QP \xiv \| > \zq)
	= 
	- \int_{\zq}^{\infty} \ex^{\nuH t} d\P(\| \QP \xiv \| \geq t)
	\\
	& \leq &
	3 \ex^{\nuH \zqc - \xxc - \constg (\zq - \zqc)} 
	+ \frac{3 \nuH}{\constg - \nuH} \, \ex^{\nuH \zqc - \xxc - \constg (\zq - \zqc) } 
	=
	\frac{3 \constg}{\constg - \nuH} \, \ex^{\nuH \zqc -\xxc - (\constg - \nuH) (\zq - \zqc)} \, .
\label{hgdft6gmmm7u7ruj543}
\end{EQA}
Now fix \( \zqx \) with \( \zqx - \sqrt{\dimH} \geq 2 \nuH \) but \( \zqx \leq \zqc \).
Then 
\begin{EQA}
	&& \nquad
	\E \ex^{\nuH \| \QP \xiv \|} \Ind( \| \QP \xiv \| > \zqx)
	=
	- \int_{\zqx}^{\infty} \ex^{\nuH \zq} d\P(\| \QP \xiv \| \geq \zq)
	\\
	&=&
	\ex^{\nuH \zqx} \P(\| \QP \xiv \| \geq \zqx) 
	+ \nuH \left( \int_{\zqx}^{\zqc} + \int_{\zqc}^{\infty} \right) \ex^{\nuH \zq } \P(\| \QP \xiv \| \geq \zq) d\zq \, .
\label{jcy603fgy7e34e7yf8reyhu7ewg}
\end{EQA}
By \eqref{PxivbzzBBroB}, for any \( \zq \in [\zqx , \zqc] \),
it holds in view of \( \zq(\BBH,\xx) \leq \sqrt{\dimH} + \sqrt{2 \xx} \)
\begin{EQA}
	\P(\| \QP \xiv \| \geq \zq)
	& \leq &
	3 \ex^{- (\zq - \sqrt{\dimH})^{2}/2} .
\label{fduyr2ytvfiuh23}
\end{EQA}
As \( (\nuH \zq - (\zq - \sqrt{\dimH})^{2}/2)' = \nuH - \zq + \sqrt{\dimH} \leq - \nuH \) 
for \( \zq - \sqrt{\dimH} \geq 2 \nuH \), it holds
\begin{EQA}
	\nuH \int_{\zqx}^{\zqc} \ex^{\nuH \zq - (\zq - \sqrt{\dimH})^{2}/2 } d\zq
	& \leq &
	\ex^{ \nuH \zqx - (\zqx - \sqrt{\dimH})^{2}/2 } \, \nuH \int_{\zqx}^{\zqc} \ex^{- \nuH (\zq - \zqx) } d\zq 
	\leq 
	\ex^{ \nuH \zqx - (\zqx - \sqrt{\dimH})^{2}/2 } \, 
\label{gtw8ibh9t659o2wmljw8}
\end{EQA}
and also \( \nuH \zqx - (\zqx - \sqrt{\dimH})^{2}/2 > \nuH \zqc - (\zqc - \sqrt{\dimH})^{2}/2 \).
Putting this together with the above bound on \( \int_{\zqc}^{\infty} \ex^{\nuH \zq } \P(\| \QP \xiv \| \geq \zq) d\zq  \) 
as in \eqref{hgdft6gh3wef7u7ruj543}
completes the proof of \eqref{ufikwk3ei9vgkj4k4giw3wl}.


\Section{Frobenius norm losses for empirical covariance}
\label{SFrobnorm}

Let \( \Xv_{i} \sim \ND(0,\Sigma) \) be i.i.d. zero mean Gaussian vectors in \( \R^{\dimp} \) 
with a covariance matrix \( \Sigma \in \Matr_{\dimp} \).
By \( \Sigmah \) we denote the empirical covariance
\begin{EQA}
	\Sigmah
	& \eqdef &
	\frac{1}{n} \sumi \Xv_{i} \, \Xv_{i}^{\T} \, .
\label{h7hij476wtdysuwjfjti}
\end{EQA}
Our goal is to establish sharp dimension free deviation bounds on the squared Frobenius norm 
\( \| \Sigmah - \Sigma \|_{\Fr}^{2} \):
\begin{EQA}
	\| \Sigmah - \Sigma \|_{\Fr}^{2}
	&=&
	\tr (\Sigmah - \Sigma)^{2}  .
\label{trqtyc6q67ed737duyq}
\end{EQA}
We demonstrate how the general results of Section~\ref{Sdevboundexp} can be used for obtaining accurate deviation bounds 
for \( \| \Sigmah - \Sigma \|_{\Fr}^{2} \) and for supporting the concentration phenomenon.

\Subsection{Upper bounds}
First we establish a tight upper bound on \( \| \Sigmah - \Sigma \|_{\Fr}^{2} \).
We identify the matrix \( \Sigmah \) with the vector in the linear subspace of \( \R^{\dimp \times \dimp} \) 
composed by symmetric matrices.
Our aim is in showing that the quantiles of \( \| \Sigmah - \Sigma \|_{\Fr}^{2} \) mimic well similar quantiles 
of \( \| \Sigmat - \Sigma \|_{\Fr}^{2} \) for a Gaussian matrix \( \Sigmat \) with the same covariance structure as \( \Sigmah \).
Define 
\begin{EQA}
	\dimH(\Sigma) 
	= (\tr \Sigma)^{2} + \tr \Sigma^{2} ,
	\qquad
	\vH^{2}(\Sigma)
	&=&
	\bigl( \tr \Sigma^{2} \bigr)^{2} + \tr \Sigma^{4} .
\label{9ckd63hggy7r67y3ghft63}
\end{EQA}
Later we show that \( \dimH(\Sigma) = \E \| \Sigmah - \Sigma \|_{\Fr}^{2} = \tr \Var(\Sigmat) \) and 
\( \vH^{2}(\Sigma) = \tr \{ \Var(\Sigmat) \}^{2} \) while \( \supH(\Sigma) = \| \Var(\Sigmat) \| = 2 \| \Sigma \|^{2} \).
In our results we implicitly assume a high dimensional situation with \( \dimH(\Sigma) \) large.
The presented bounds also require that \( n \gg \dimH(\Sigma) \).

\begin{theorem}
\label{TFrobGauss}
Assume \( \| \Sigma \| = 1 \) and \( \dimH(\Sigma) < n/8 \).
Given \( \xx \) with \( 4 \sqrt{\xx} < \sqrt{n/8} - \sqrt{\dimH(\Sigma)} \),
fix \( \alpHm < 1 \) by 
\begin{EQA}
	\alpHm (1 - \alpHm) \sqrt{n/8} 
	&=&
	\sqrt{\dimH(\Sigma)} + 4 \sqrt{\xx} \, .
\label{yf7j3e53thnutjkerkin}
\end{EQA}
Then 
\begin{EQA}
	\P\Bigl( 
		n \| \Sigmah - \Sigma \|_{\Fr}^{2} > \frac{1}{1 - \alpHm} \bigl\{ \dimH(\Sigma) + 2 \vH(\Sigma) \sqrt{\xx} + 4 \xx \bigr\} 
	\Bigr)
	& \leq &
	3 \ex^{-\xx} \, .
\label{uc6d5e4w4ew5dtdbiie}
\end{EQA}
\end{theorem}


\Subsection{Lower bounds}
This section presents a lower bound on the Frobenius norm of \( \Sigmah - \Sigma \).
Later in Section~\ref{SFrobconc} we state the concentration phenomenon for \( \| \Sigmah - \Sigma \|_{\Fr}^{2} \).

\begin{theorem}
\label{TFrobGausslo}
Let \( \| \Sigma \| = 1 \) and  \( \dimH(\Sigma) \) and \( \vH(\Sigma) \) be defined by \eqref{9ckd63hggy7r67y3ghft63}.
For \( \xx > 0 \) with \( 2 \sqrt{\xx} \leq \dimH(\Sigma)/\vH(\Sigma) \), define \( \muH = \muH(\xx) = 2 \sqrt{\xx}/\vH(\Sigma) \) and assume that
there is \( \alpH < 1/2 \) satisfying
\begin{EQA}
	\alpH \sqrt{\frac{1 - 2\alpH}{1 - \alpH}}
	& \geq &
	\sqrt{\frac{\muH(\xx)}{n}} \biggl( \sqrt{2 \dimH(\Sigma)} + \frac{\sqrt{2} \, \dimH(\Sigma)}{\vH(\Sigma)} \biggr) .
\label{ujw3jdfcv7823ujfwqyshf}
\end{EQA}
Then
\begin{EQA}
	\P\Bigl( n \| \Sigmah - \Sigma \|_{\Fr}^{2} 
		< \frac{1 - 2 \alpH}{1 - \alpH} \, \dimH(\Sigma) - 2 \vH(\Sigma) \sqrt{\xx} 
	\Bigr)
	& \leq &
	2 \ex^{-\xx} .
\label{9fhwe3hdfv7ye43kbhnikj}
\end{EQA}
\end{theorem}

\Subsection{Concentration of the Frobenius loss}
\label{SFrobconc}
Putting together Theorem~\ref{TFrobGauss} and Theorem~\ref{TFrobGausslo} yields the following corollary.

\begin{corollary}
\label{CTFrobGauss}
Under conditions of Theorem~\ref{TFrobGauss} and Theorem~\ref{TFrobGausslo}, 
it holds for any \( \xx \) resolving \eqref{yf7j3e53thnutjkerkin} and \eqref{ujw3jdfcv7823ujfwqyshf}
on a random set \( \Omega(\xx) \) with \( \P\bigl( \Omega(\xx) \bigr) \geq 1 - 5 \ex^{-\xx} \)
\begin{EQA}
	\frac{1 - 2\alpH}{1 - \alpH} \, \dimH(\Sigma) - 2 \vH(\Sigma) \sqrt{\xx} 
	& \leq &
	n \| \Sigmah - \Sigma \|_{\Fr}^{2} 
	\leq  
	\frac{1}{1 - \alpHm} \bigl\{ \dimH(\Sigma) + 2 \vH(\Sigma) \sqrt{\xx} + 4 \xx \bigr\} .
	\qquad
\label{uc6d5e4w4ew5dtdbiiec}
\end{EQA}
\end{corollary}

This result mimics similar bound of Theorem~\ref{TexpbLGA} for \( \Sigmah \) Gaussian 
and of Theorem~\ref{Tdevboundgm} for \( \Sigmah \) sub-Gaussian.
However, the empirical covariance \( \Sigmah \) is quadratic in the \( \Xv_{i} \)'s and thus, only sub-exponential.
We pay an additional factor \( (1 - \alpHm)^{-1} \) in the upper quantile function 
and the factor \( \frac{1 - 2\alpH}{1 - \alpH} \) in the lower quantile function for this extension.

Further we discuss the concentration phenomenon for the Frobenius error \( n \| \Sigmah - \Sigma \|_{\Fr}^{2} \) 
around its expectation \( \dimH(\Sigma) \).
Even in the Gaussian case, it meets only in high-dimensional situation with \( \dimH(\Sigma) \) large.
As \( \vH^{2}(\Sigma) \leq \dimH(\Sigma) \supH(\Sigma) = 2 \dimH(\Sigma) \), this also implies
\( \vH(\Sigma) \ll \dimH(\Sigma) \).
Statement \eqref{uc6d5e4w4ew5dtdbiiec} can be rewritten as
\begin{EQA}
	- \frac{\alpH \, \dimH(\Sigma)}{1 - \alpH} - 2 \vH(\Sigma) \sqrt{\xx} 
	& \leq &
	n \| \Sigmah - \Sigma \|_{\Fr}^{2} - \dimH(\Sigma)
	\leq  
	\frac{\alpHm \, \dimH(\Sigma)}{1 - \alpHm} + \frac{2 \vH(\Sigma) \sqrt{\xx} + 4 \xx}{1 - \alpHm} \, .
	\qquad
\label{uc6d5e4w4ew5dtdbiiecm}
\end{EQA}
Therefore, concentration effect of the loss \( n \| \Sigmah - \Sigma \|_{\Fr}^{2} \) requires \( \dimH(\Sigma) \) large and
\( \alpH \) and \( \alpHm \) small.
Then for \( \xx \ll \dimH(\Sigma) \), quantiles of \( n \| \Sigmah - \Sigma \|_{\Fr}^{2} - \dimH(\Sigma) \)
are smaller in order than \( \dimH(\Sigma) \).
Definition \eqref{yf7j3e53thnutjkerkin} of \( \alpHm \) ensures \( \alpHm \asymp \sqrt{\dimH(\Sigma)/n} \),
and hence, ``\( \alpHm \ll 1 \)'' is equivalent to ``\( \dimH(\Sigma) \ll n \)''. 
Condition ensuring \( \alpH \ll 1 \) is similar.
To see this, assume  \( \vH^{2}(\Sigma) \asymp \dimH(\Sigma) \).
Then \( \xx \ll \dimH(\Sigma) \) yields \( \muH(\xx) = 2 \sqrt{\xx}/\vH(\Sigma) \ll 1 \)
and definition \eqref{ujw3jdfcv7823ujfwqyshf} of \( \alpH \) implies 
\begin{EQA}
	\alpH
	& \lesssim &
	\sqrt{\frac{\muH}{n}} \biggl( \sqrt{2 \dimH(\Sigma)} + \frac{\sqrt{2} \, \dimH(\Sigma)}{\vH(\Sigma)} \biggr)
	\lesssim 
	\sqrt{\frac{\dimH(\Sigma)}{n}} \, .
\label{ojju8733erfsy7r745w23e3r}
\end{EQA}

\Subsection{Weighted Frobenius norm}
The result can be easily extended to the case of a weighted Frobenius norm.
Consider for any linear mapping \( \KH \colon \R^{\dimp} \to \R^{\dimq} \) the value 
\( n \| \KH (\Sigmah - \Sigma) \KH^{\T} \|_{\Fr}^{2} \).

\begin{theorem}
\label{TFrobGaussA}
Let \( \| \Sigma \| = 1 \) and \( \KH \colon \R^{\dimp} \to \R^{\dimq} \) be a linear operator with 
\( \| \KH \| = \| \KH^{\T} \KH \| = 1 \).
Define \( \Sigma_{\KH} \eqdef \KH \Sigma \KH^{\T} \),
\begin{EQA}
	\dimH_{\KH} 
	& \eqdef & 
	\dimH(\Sigma_{\KH})
	=
	\tr^{2} (\Sigma_{\KH}) + \tr (\Sigma_{\KH})^{2},
	\quad
	\vH_{\KH}^{2} 
	\eqdef 
	\vH^{2}(\Sigma_{\KH})
	=
	\bigl\{ \tr (\Sigma_{\KH}^{2}) \bigr\}^{2} + \tr (\Sigma_{\KH})^{4},
	\qquad
\label{0j3wgyrfg76hy5tjhtdeert}
\end{EQA}
and assume \( \dimH_{\KH} < n/8 \).
The the statements of Theorem~\ref{TFrobGauss} and Theorem~\ref{TFrobGausslo} apply to 
\( n \| \KH (\Sigmah - \Sigma) \KH^{\T} \|_{\Fr}^{2} \) after replacing 
\( \dimH(\Sigma) \) and \( \vH(\Sigma) \) with \( \dimH_{\KH} \) and \( \vH_{\KH} \).
\end{theorem}

\begin{proof}
We can represented  
\begin{EQA}
	\sqrt{n} \, \KH (\Sigmah - \Sigma) \KH^{\T}
	&=&
	\KH \, \Sigma^{1/2} \errS \, \Sigma^{1/2} \KH^{\T}
\label{on0njm9yu7i8erwhedc}
\end{EQA}
with \( \errS \) from \eqref{jufc6jbjuhhjiutifg}.
This reduces the result to the previous case with \( \Sigma_{\KH} = \KH \Sigma \KH^{\T} \) in place of \( \Sigma \).
\end{proof}

\Subsection{Proof of Theorem~\ref{TFrobGauss}}
Each vector \( \gaussv_{i} = \Sigma^{-1/2} \Xv_{i} \) is standard normal.
Define 
\begin{EQA}
	\errS
	&=&
	\frac{1}{n^{1/2}} \sumi (\gaussv_{i} \gaussv_{i}^{\T} - \Id_{\dimp}) .
\label{jufc6jbjuhhjiutifg}
\end{EQA}
We will use the representation \( \Sigmah - \Sigma = n^{-1/2} \Sigma^{1/2} \, \errS \, \Sigma^{1/2} \) and
\begin{EQA}
	n \| \Sigmah - \Sigma \|_{\Fr}^{2}
	&=&
	\tr ( \Sigma^{1/2} \errS \, \Sigma \, \errS \, \Sigma^{1/2} )
	=
	\| \Sigma^{1/2} \errS \, \Sigma^{1/2} \|_{\Fr}^{2} \, .
\label{hvcftw3ytet655rdffhkfj}
\end{EQA}
The main step is in applying Theorem~\ref{Tdevboundgm} to the quadratic form \( \| \QP \errS \|_{\Fr}^{2} \) 
with \( \QP \errS = \Sigma^{1/2} \, \errS \, \Sigma^{1/2} \). 
First check \eqref{expgamgm} for \( \xiv = \errS \).

\begin{lemma}
\label{LexpmgamG}
For any symmetric \( \Gamma \in \Matr_{\dimp} \) with \( \| \Gamma \|_{\Fr} \leq \gmn < \sqrt{n}/2 \), it holds
\begin{EQA}
\label{vfrw44d4d57y8hjhqtstik2}
	\E \langle \Gamma, \errS \rangle^{2}
	&=&
	2 \| \Gamma \|_{\Fr}^{2} \, ,
	\\
	\log \E \exp \langle \Gamma, \errS \rangle
	& \leq &
	\frac{1}{1 - 2 n^{-1/2} \| \Gamma \|} \, \| \Gamma \|_{\Fr}^{2}
	\leq 
	\frac{1}{1 - 2 n^{-1/2} \gmn} \, \| \Gamma \|_{\Fr}^{2} \, .
\label{vfrw44d4d57y8hjhqtstik}
\end{EQA}
\end{lemma}

\begin{proof}
Let us fix any symmetric \( \Gamma \in \Matr_{\dimp} \) with \( \| \Gamma \|_{\Fr} \leq \gmn \).
For the scalar product \( \langle \Gamma, \errS \rangle \), we use the representation
\begin{EQA}
	\langle \Gamma, \errS \rangle
	&=&
	\tr (\Gamma \errS)
	=
	\frac{1}{n^{1/2}} \sumi \bigl\{ \gaussv_{i}^{\T} \Gamma \gaussv_{i} - \E (\gaussv_{i}^{\T} \Gamma \gaussv_{i}) \bigr\} .
\label{5twy7fuuf37dhjuubjbms}
\end{EQA}
Then by independence of the \( \gaussv_{i} \)'s and Lemma~\ref{Gaussmoments}, it holds 
\begin{EQA}
	\E \langle \Gamma, \errS \rangle^{2}
	&=&
	\frac{1}{n} \sumi \E \bigl\{ \gaussv_{i}^{\T} \Gamma \gaussv_{i} - \E (\gaussv_{i}^{\T} \Gamma \gaussv_{i}) \bigr\}^{2}
	=
	2 \tr \Gamma^{2} .
\label{5twy7fuuf37dhjuubjbms2}
\end{EQA}
Now consider the exponential moment of \( \langle \Gamma, \errS \rangle \).
Again, independence of the \( \gaussv_{i} \)'s yields
\begin{EQA}
	\log \E \exp \langle \Gamma, \errS \rangle
	&=&
	\sumi \log \E \exp \frac{\gaussv_{i}^{\T} \Gamma \gaussv_{i}}{\sqrt{n}} - \sqrt{n} \, \tr \Gamma
	\\
	&=&
	\frac{n}{2} \log \det \bigl( \Id_{\dimp} - \frac{2}{\sqrt{n}} \, \Gamma \bigr) - \sqrt{n} \, \tr \Gamma
	\qquad
\label{hw6jhydhqanchyenqqdbsk} 
\end{EQA}
provided that \( 2\Gamma < \sqrt{n} \Id_{\dimp} \).
Moreover, by Lemma~\ref{Lqfexpmom}
\begin{EQA}
	\left| \frac{n}{2} \log \det (\Id_{\dimp} - 2 n^{-1/2} \Gamma) - \sqrt{n} \, \tr \Gamma \right|
	& \leq &
	\frac{\tr \Gamma^{2}}{1 - 2 n^{-1/2} \| \Gamma \|} 
	=
	\frac{\| \Gamma \|_{\Fr}^{2}}{1 - 2 n^{-1/2} \| \Gamma \|} \, ,
	\qquad
\label{hyakboj84gjqtqgxjvu}
\end{EQA}
and the assertion follows in view of \( \| \Gamma \| \leq \| \Gamma \|_{\Fr} \leq \gmn \).
\end{proof}


We now fix \( \gmn = \alpHm \sqrt{n}/2 \).
Then the random matrix \( \xiv = \errS \) follows
condition \eqref{expgamgm} with 
\( \HVB^{2} = 2 (1 - \alpHm)^{-1} \Id \).
This enables us to apply Theorem~\ref{Tdevboundgm} to the quadratic form \( \| \QP \errS \|_{\Fr}^{2} \) 
for \( \QP \errS = \Sigma^{1/2} \, \errS \, \Sigma^{1/2} \). 
By \eqref{vfrw44d4d57y8hjhqtstik2}, it holds \( \Var(\errS) = 2 \Id \).
Now introduce a Gaussian element \( \errSt \) with the same covariance structure. 
One can use \( \errSt = (\Gaussv + \Gaussv^{\T})/\sqrt{2} \), where \( \Gaussv = (\Gauss_{ij}) \) is a random \( \dimp \)-matrix 
with i.i.d. standard normal entries \( \Gauss_{ij} \).
Indeed, for any symmetric \( \dimp \)-matrix \( \Gamma \),
\begin{EQA}
	\E \langle \errSt, \Gamma \rangle^{2}
	&=&
	2 \E \langle \Gaussv,\Gamma \rangle^{2}
	=
	2 .
\label{8jjjii8i8334eeee4kfoof}
\end{EQA}
Statement \eqref{PxivbzzBBroB} of Theorem~\ref{Tdevboundgm} yields nearly the same deviation bounds 
for \( \| \QP \errS \|_{\Fr}^{2} \) as for \( \| \QP \errSt \|_{\Fr}^{2} \) with \( \errSt \sim \ND(0,\Var(\errS)) \).
Theorem~\ref{TexpbLGA} claims 
\begin{EQA}
	\P\bigl( \| \QP \errSt \|_{\Fr}^{2} > \zq^{2}(\BBHt,\xx) \bigr)
	& \leq &
	\ex^{-\xx} \, ,
\label{jf67wyfuvb87e7eufhgghwxu}
\end{EQA}
where \( \BBHt = \Var(\QP \errSt) \) and 
the quantile \( \zq(\BBH,\xx) \) is defined as
\begin{EQA}
	\zq^{2}(\BBH,\xx)
	&=&
	\tr \BBH + 2 \sqrt{\xx \tr(\BBH^{2})} + 2 \xx \| \BBH \| .
\label{fiu47fkw3ydfywe3hdf7y2h}
\end{EQA}

\begin{lemma}
\label{LBBHt}
Let \( \errSt = (\Gaussv + \Gaussv^{\T})/\sqrt{2} \), where \( \Gaussv = (\Gauss_{ij}) \) is a random \( \dimp \)-matrix 
with i.i.d. standard normal entries \( \Gauss_{ij} \).
Consider \( \QP \errSt = \Sigma^{1/2} \, \errSt \, \Sigma^{1/2} \).
It holds for \( \BBHt = \Var(\QP \errSt) \)
\begin{EQA}
	\tr \BBHt
	&=&
	\dimH(\Sigma) \, ,
	\quad
	\tr \BBHt^{2}
	=
	\vH^{2}(\Sigma) \, ,
	\quad
	\| \BBHt \|
	=
	2. 
\label{yx7hqwjxc78g94mg6e6e4heg}
\end{EQA}
\end{lemma}

\begin{proof}
We may assume \( \Sigma = \diag\{ \lambda_{1},\ldots,\lambda_{\dimp} \} \). 
Then it holds by Lemma~\ref{Gaussmoments}
\begin{EQA}
	\| \QP \errSt \|_{\Fr}^{2}
	=
	\| \Sigma^{1/2} \, \errSt \, \Sigma^{1/2} \|_{\Fr}^{2}
	=
	\frac{1}{2} \sum_{i,j=1}^{\dimp} \lambda_{i} \, \lambda_{j} \, (\Gauss_{ij} + \Gauss_{ji})^{2}
	& \eqd &
	2 \sum_{i \leq j} \lambda_{i} \, \lambda_{j} \, \Gauss_{ij}^{2}
	\qquad
\label{udhf6e3ggun8true4hfvy}
\end{EQA}
and thus
\begin{EQA}
	\tr \BBHt
	=
	\E \| \QP \errSt \|_{\Fr}^{2}
	&=&
	2 \sum_{i \leq j} \lambda_{i} \, \lambda_{j}
	=
	\biggl( \sum_{i=1}^{\dimp} \lambda_{i} \biggr)^{2} + \sum_{i=1}^{\dimp} \lambda_{i}^{2}
	=
	\dimH(\Sigma) \, .
\label{yx7hqwjxc78wjuanct6whw}
\end{EQA}
Further we compute \( \vH^{2}(\Sigma) = \tr \BBHt^{2} \).
Note that \( \Var( \| \QP \errSt \|_{\Fr}^{2} ) \neq \Var( \| \QP \errS \|_{\Fr}^{2} ) \).
Due to Lemma~\ref{Gaussmoments}, it holds
\( \vH^{2}(\Sigma) = \Var( \| \QP \errSt \|_{\Fr}^{2} )/2 \) yielding by \eqref{udhf6e3ggun8true4hfvy}
\begin{EQA}
	\vH^{2}(\Sigma)
	&=&
	2 \sum_{i \leq j} \lambda_{i}^{2} \lambda_{j}^{2} \Var (\Gauss_{ij}^{2})	
	=
	2 \sum_{i \neq j} \lambda_{i}^{2} \lambda_{j}^{2} 
	+ 2 \sum_{i=1}^{\dimp} \lambda_{i}^{4} 
	=
	\bigl( \tr \Sigma^{2} \bigr)^{2} + \tr \Sigma^{4} .
\label{ioer49u8fdj42e3u7fdjgt}
\end{EQA}
Finally, \( \Var(\errS) = 2 \Id \) and \( \| \Sigma \| = 1 \) implies \( \supH(\Sigma) = \| \QP \Var(\errS) \QP^{\T} \| = 2 \).
\end{proof}

Now we apply Theorem~\ref{Tdevboundgm} to 
\( n \| \Sigmah - \Sigma \|_{\Fr}^{2} = \| \QP \errS \|_{\Fr}^{2} \).
Following to Lemma~\ref{LexpmgamG}, define \( \BBH = (1 - \nuH)^{-1} \BBHt \).
Then with \( \zq^{2}(\BBH,\xx) \) from \eqref{fiu47fkw3ydfywe3hdf7y2h}
\begin{EQA}
	\P\Bigl( n \| \Sigmah - \Sigma \|_{\Fr}^{2} > \zq^{2}(\BBH,\xx) \Bigr)
	& = &
	\P\bigl( \| \QP \errS \|_{\Fr}^{2} > \zq^{2}(\BBH,\xx) \bigr)
	\leq 
	3 \ex^{-\xx} ,
	\quad
	\xx \leq \xxc \, ,
\label{ncye6e3y7fiuh8ytst5q}
\end{EQA}
and assertion \eqref{uc6d5e4w4ew5dtdbiie} follows in view of Lemma~\ref{LBBHt} and 
\( \zq^{2}(\BBH,\xx) = (1 - \nuH)^{-1} \zq^{2}(\BBHt,\xx) \).
However, it is still necessary to check that the upper bound \eqref{uc6d5e4w4ew5dtdbiie} 
applies for a given \( \xx \).
\eqref{if7h3rhgy4676rfhdsjw} provides a sufficient condition 
\( \gmn/\supH \geq \sqrt{\dimH/\supH} + \sqrt{8\xx} \) with \( \dimH = \dimH(\Sigma)/(1 - \alpHm) \) and
\( \supH = 2/(1 - \alpHm) \) for \( \gmn = \alpHm \sqrt{n} / 2 \).
%
By \eqref{yf7j3e53thnutjkerkin}
\begin{EQA}
	\frac{\gmn}{\supH} - \sqrt{\frac{\dimH}{\supH}} 
	& = & 
	\frac{\alpHm \sqrt{n}}{2\supH} - \sqrt{\frac{\dimH(\Sigma)}{2}} 
	\geq 
	\frac{\alpHm (1 - \alpHm) \sqrt{n}}{4} - \sqrt{\frac{\dimH(\Sigma)}{2}}
	> 
	\sqrt{8\xx} 
\label{nctc6yw3f76h7u87y98}
\end{EQA}
and the result follows.

\Subsection{Proof of Theorem~\ref{TFrobGausslo}}
As in the proof of the upper bound, we apply Markov's inequality
\begin{EQA}
	\P\bigl( n \| \Sigmah - \Sigma \|_{\Fr}^{2} < \zq \bigr)
	& \leq &
	\ex^{\muH \zq/2} \E \exp\Bigl( - \frac{\muH}{2} n \| \Sigmah - \Sigma \|_{\Fr}^{2} \Bigr) .
\label{9vjnegf66ruhiti4rhsw000}
\end{EQA}
However, now we are free to choose any positive \( \muH \).
Later we evaluate the exponential moments of \( - n \| \Sigmah - \Sigma \|_{\Fr}^{2} \) for all \( \muH > 0 \)
and then, given \( \xx \), fix \( \muH \) and \( \zq \) similarly to the Gaussian case to ensure 
the prescribed deviation probability \( \ex^{-\xx} \). 

Denote by \( \Gaussv = (\Gauss_{ij}) \) a random \( \dimp \times \dimp \) matrix 
with i.i.d. standard Gaussian entries \( \Gauss_{ij} \) and \( \Gaussvb \eqdef (\Gaussv + \Gaussv^{\T})/2 \).
Then for any \( \muH > 0 \)
\begin{EQA}
	\exp \bigl( - \muH n \| \Sigmah - \Sigma \|_{\Fr}^{2}/2 \bigr)
	&=&
	\E_{\Gauss} \exp\bigl\{ \imi \sqrt{\muH n} \, \langle \Sigmah - \Sigma, \Gaussv \rangle \bigr\}
	=
	\E_{\Gauss} \exp\bigl\{ \imi \sqrt{\muH n} \, \langle \Sigmah - \Sigma, \Gaussvb \rangle \bigr\} .
\label{9vje6e36tghu4i3bieijd}
\end{EQA}
Therefore, by independence of the \( \Xv_{i} \)'s
\begin{EQA}
\label{yhs7vjeje8bjekwytvya}
	\E \exp \bigl( - \muH n \| \Sigmah - \Sigma \|_{\Fr}^{2}/2 \bigr)
	&=&
	\E_{\Gauss} \E \exp\bigl( \imi \sqrt{\muH n} \, \langle \Sigmah - \Sigma, \Gaussvb \rangle \bigr)
	\\
	&=&
	\E_{\Gauss} \Bigl\{ \E \exp\bigl( \imi \sqrt{\muH/n} \, \langle \Xv_{1} \Xv_{1}^{\T} - \Sigma, \Gaussvb \rangle \bigr) \Bigr\}^{n}
	\\
	&=&
	\E_{\Gauss} \Bigl\{ \E
	\exp\bigl( 
		\imi \sqrt{\muH/n} \, \langle \gaussv \gaussv^{\T} - \Id_{\dimp} \, , \Sigma^{1/2} \, \Gaussvb \, \Sigma^{1/2} \rangle 
	\bigr) 
	\Bigr\}^{n} .
\end{EQA}
Further, by Lemma~\ref{Lqfexpmom}, with \( \BHG = \Sigma^{1/2} \, \Gaussvb \, \Sigma^{1/2} \)
\begin{EQA}
	&& \nquad
	\Bigl\{ \E
		\exp\bigl( \imi \sqrt{\muH/n} \, \langle \gaussv \gaussv^{\T} - \Id_{\dimp} \, , \BHG \rangle \bigr) 
	\Bigr\}^{n}
	\\
	&=&
	\exp\bigl\{ n \log \det\bigl( \Id_{\dimp} - 2 \imi \sqrt{\muH/n} \, \BHG \bigr)^{-1/2} 
	- \imi \sqrt{\muH n} \tr (\BHG) \bigr\} .
\label{yhs7vjeje8bjekwytgm}
\end{EQA}
Let some \( \xx > 0 \) and some \( \alpH \in (0,1/2) \) be fixed.
Define
\begin{EQA}
	\muH 
	& \eqdef &
	\frac{2 \sqrt{\xx} }{\vH(\Sigma)} \, ,
	\qquad
	\muHa
	\eqdef
	\frac{1 - \alpH}{1 - 2\alpH} \, \muH
	=
	\frac{1 - \alpH}{1 - 2\alpH} \,\, \frac{2 \sqrt{\xx} }{\vH(\Sigma)} \, ,
\label{nsducf7enhw3e7yfgy6ryer4b}
\end{EQA}
and introduce a random set \( \Omega(\alpH) \) with
\begin{EQA}
	\Omega(\alpH)
	& \eqdef &
	\bigl\{ \Gaussv \colon 2 \sqrt{\muHa/n} \, \| \BHG \|
	\leq 
	\alpH \bigr\} ,
	\qquad
	\BHG = \Sigma^{1/2} \, (\Gaussv + \Gaussv^{\T}) \, \Sigma^{1/2}/2.
\label{uiduuvbinjoko4e7gyuuh2w2}
\end{EQA}
It holds on \( \Omega(\alpH) \) 
by \eqref{yhs7vjeje8bjekwytgm} similarly to \eqref{vbu7j3hg8hryhghyidegwdg} of Lemma~\ref{Lqfexpmom} 
\begin{EQA}
	&& \nquad
	\E^{n} \exp\bigl\{ \imi \sqrt{\muHa/n} \, \langle \gaussv \gaussv^{\T} - \Id_{\dimp}, \BHG \rangle \bigr\}
	\leq 
	\exp \biggl( - \muHa \tr(\BHG^{2}) + \frac{\muHa \, \alpH \tr(\BHG^{2})}{1 - \alpH} \biggr)
	\\
	&=&
	\exp\biggl( - \frac{1 - 2 \alpH}{1 - \alpH} \, \muHa \, \tr(\BHG^{2}) \biggr)
	=
	\exp\bigl( - \muH \, \tr(\BHG^{2}) \bigr) .
	\qquad
\label{yhs7vjeje8behfydrjekwytgm}
\end{EQA}
Exponential moments of \( \tr(\BHG^{2}) \) from \eqref{uc8ewjfbu7fue77egkbire} under \( \P_{\Gauss} \) can be easily computed.
We proceed assuming \( \Sigma = \diag\{ \lambda_{j} \} \) and
using that \( \Gauss_{ij} + \Gauss_{ji} \sim \ND(0,2) \) for \( i \neq j \), and all 
\( \Gauss_{ij} + \Gauss_{ji} \) are mutually independent for \( i \leq j \).
This implies
\begin{EQA}
	\tr(\BHG^{2})
	&=&
	\frac{1}{4} \sum_{i,j=1}^{\dimp} \lambda_{i} \, \lambda_{j} \, (\Gauss_{ij} + \Gauss_{ji})^{2}
	\, \eqd \,
	\sum_{i \leq j} \lambda_{i} \, \lambda_{j} \, \Gauss_{ij}^{2}
\label{uc8ewjfbu7fue77egkbire}
\end{EQA}
and 
\begin{EQA}
	\E_{\Gauss} \tr(\BHG^{2})
	&=&
	\sum_{i \leq j} \lambda_{i} \, \lambda_{j}
	=
	\frac{\dimH(\Sigma)}{2} \, ,
\label{87wkfv8ieigy5e5tgvuf}
	\\
	\E_{\Gauss} \exp\{ - \muH \tr(\BHG^{2}) \}
	&=&
	\E_{\Gauss} \exp\biggl( - \muH \sum_{i \leq j} \lambda_{i} \, \lambda_{j} \, \Gauss_{ij}^{2}  \biggr)
	=
	\exp \biggl( - \frac{1}{2} \sum_{i \leq j} \log (1 + 2 \muH \, \lambda_{i} \, \lambda_{j}) \biggr) .
\label{jcvuyfjhi8u97iew22g}
\end{EQA}
The latter expression can be evaluated by using \eqref{m2v241m41mbn} of Lemma~\ref{Lqfexpmom}:
\begin{EQA}
	\E_{\Gauss} \exp\{ - \muH \tr(\BHG^{2}) \}
	& \leq &
	\exp \biggl( 
		- \muH \sum_{i \leq j} \lambda_{i} \, \lambda_{j} + \muH^{2} \sum_{i \leq j} \lambda_{i}^{2} \, \lambda_{j}^{2} 
	\biggr) 
	=
	\exp\biggl( - \frac{\muH \, \dimH(\Sigma)}{2} + \frac{\muH^{2} \vH^{2}(\Sigma)}{4} \biggr).
\label{jcvuyfjhbir64gtughew22g}
\end{EQA}
This and \eqref{yhs7vjeje8behfydrjekwytgm} yield 
\begin{EQA}
	\E \exp \biggl( - \frac{\muHa}{2} n \| \Sigmah - \Sigma \|_{\Fr}^{2} \biggr)
	& \leq &
	\P_{\Gauss}\bigl( \Omega(\alpH)^{c} \bigr) 
	+ \exp \biggl( - \frac{\muH \, \dimH(\Sigma)}{2} + \frac{\muH^{2} \vH^{2}(\Sigma)}{4} \biggr) 
\label{ifgo9elrgpgjn9u9tkfywfgi}
\end{EQA}
and for any \( \zq \) by Markov's inequality \eqref{9vjnegf66ruhiti4rhsw000}
\begin{EQA}
	\P\bigl( n \| \Sigmah - \Sigma \|_{\Fr}^{2} < \zq \bigr)
	& \leq &
	\ex^{\muHa \zq/2} \P_{\Gauss}\bigl( \Omega(\alpH)^{c} \bigr) 
	+ \exp\biggl( \frac{\muHa \zq}{2} - \frac{\muH \, \dimH(\Sigma)}{2} + \frac{\muH^{2} \vH^{2}(\Sigma)}{4} \biggr) .
\label{9vjnegf66ruhiti4rhsw}
\end{EQA}
With \( \muH = 2 \sqrt{\xx} / \vH(\Sigma) \), we  
define \( \zq \) by
\begin{EQA}
	\muHa \,\zq
	&=&
	\muH \{ \dimH(\Sigma) - 2 \vH(\Sigma) \sqrt{\xx} \} 
	=
	\frac{2 \sqrt{\xx} \, \dimH(\Sigma)}{\vH(\Sigma)} - 4 \xx \, 
\label{8ujek94jhkmokjkoiowndy}
\end{EQA}
yielding 
\begin{EQA}
	\frac{\muHa \zq}{2} - \frac{\muH \, \dimH(\Sigma)}{2} + \frac{\muH^{2} \vH^{2}(\Sigma)}{4}
	&=&
	\frac{\muH}{2} \bigl\{ \dimH(\Sigma) - 2 \vH(\Sigma) \sqrt{\xx} \bigr\}
	- \frac{\muH \, \dimH(\Sigma)}{2} + \frac{\muH^{2} \vH^{2}(\Sigma)}{4}
	=
	-\xx
\label{nyhwnvu7u7w2hvuyghiuhyiu}
\end{EQA}
and
\begin{EQA}
	\P( n \| \Sigmah - \Sigma \|_{\Fr}^{2} < \zq)
	& \leq &
	\ex^{-\xx} + \ex^{\muHa \zq/2} \P_{\Gauss}\bigl( \Omega(\alpH)^{c} \bigr)
\label{y7uemb8jrboboierkedfcoaq}
\end{EQA}
where
\begin{EQA}
	\zq
	&=&
	\biggl( 1 - \frac{\alpH}{1 - \alpH} \biggr) \bigl\{ \dimH(\Sigma) - 2 \vH(\Sigma) \sqrt{\xx} \bigr\}
	\geq 
	\dimH(\Sigma)
	- \frac{\alpH}{1 - \alpH} \, \dimH(\Sigma)
	- 2 \vH(\Sigma) \sqrt{\xx} \, .
\label{kfuhyg6er4y6bfhggfy3}
\end{EQA}
For bounding the probability of the set \( \Omega(\alpH)^{c} \) from \eqref{uiduuvbinjoko4e7gyuuh2w2}, 
one can apply the advanced results from the random matrix theory.
To keep the proof self-contained, we use a simple bound \( \| \BHG \|^{2} \leq \| \BHG \|_{\Fr}^{2} = \tr(\BHG^{2}) \).
For any matrix \( \Gamma \), it holds 
\begin{EQA}
	\Var\langle \Gaussvb,\Gamma \rangle 
	&=&
	\frac{1}{4} \E \biggl( \sum_{i,j=1}^{\dimp} \Gamma_{ij} (\Gauss_{ij} + \Gauss_{ji}) \biggr)^{2}
	=
	\| \Gamma \|_{\Fr}^{2} \, 
\label{onhr6g7e8gheh3vujweg}
\end{EQA}
yielding \( \| \Var(\Gaussvb) \| \leq 1 \) and \( \| \Var(\BHG) \| \leq 1 \).
Also by \eqref{87wkfv8ieigy5e5tgvuf} \( \E \| \BHG \|_{\Fr}^{2} = \dimH(\Sigma)/2 \).
Therefore, by Theorem~\ref{TexpbLGA} applied to \( \| \BHG \|_{\Fr}^{2} \), it holds for any \( \xxs \)
\begin{EQA}
	\P_{\Gauss} \bigl( \| \BHG \|_{\Fr} > \sqrt{\dimH(\Sigma)/2} + \sqrt{2\xxs} \bigr)
	& \leq &
	\ex^{-\xxs} .
\label{j8jfkg99gue3g7gyw3bjb}
\end{EQA}
By \eqref{nsducf7enhw3e7yfgy6ryer4b} and \eqref{8ujek94jhkmokjkoiowndy}, it holds 
\begin{EQA}
	\xxs
	\eqdef
	\xx + \frac{\muHa \zq}{2} 
	& \leq &
	\frac{\dimH(\Sigma) \sqrt{\xx}}{\vH(\Sigma)} - \xx
	\leq 
	\frac{\dimH^{2}(\Sigma)}{4 \vH^{2}(\Sigma)}
\label{udf73ujerdf8t754rjfcu8wei}
\end{EQA}
and 
\begin{EQA}
	\P_{\Gauss} \biggl( \| \BHG \|_{\Fr} > \sqrt{\frac{\dimH(\Sigma)}{2}} + \frac{\dimH(\Sigma)}{\sqrt{2} \, \vH(\Sigma)} \biggr)
	& \leq &
	\ex^{-\xx - \muHa \zq/2} .
\label{j8jfkg99gvik4r9hijhew3bjb}
\end{EQA}
Therefore, 
by definition \eqref{uiduuvbinjoko4e7gyuuh2w2} and condition \eqref{ujw3jdfcv7823ujfwqyshf}
\begin{EQA}
	\ex^{\muHa \zq/2} \P_{\Gauss}\bigl( \Omega(\alpH)^{c} \bigr)
	& \leq &
	\ex^{\muHa \zq/2} \P\biggl( \| \BHG \|_{\Fr} > \frac{\alpH \sqrt{n}}{2 \sqrt{\muHa}} \biggr)
	\leq 	
	\ex^{-\xx} 
\label{knegyeg67gh47rskerur}
\end{EQA}
and the result follows.


\def\VBHG{\mathbb{V}}
\def\Tnorm{\mathcal{Z}}
\def\GaussG{\mathcal{G}}
\def\tensn{\delta}
\Section{Concentration for a family of second order tensors}
\label{S2tensFr}
Suppose to be given a family of Gaussian quadratic forms
\begin{EQA}
	\TensG_{i}
	&=&
	\sum_{j,k=1}^{\dimp} \Tens_{i,j,k} \, \gauss_{j} \, \gauss_{k} \, ,\quad
	i=1,\ldots,\dimp, 
\label{bvkie48g78h85tjrefigtoit}
\end{EQA}
with standard Gaussian r.v.'s \( \gauss_{j} \).
Without loss of generality assume that each matrix \( \Tens_{i} = (\Tens_{i,j,k})_{j,k \leq \dimp} \) is symmetric.
The value \( \TensG_{i} \) can be written as 
\begin{EQA}
	\TensG_{i} 
	&=& 
	\gaussv^{\T} \Tens_{i} \gaussv 
	= 
	\langle \Tens_{i} \gaussv, \gaussv \rangle 
	=
	\langle \Tens_{i}, \gauss^{\otimes 2} \rangle.
\label{ycjhdfry64y7t7yu78t64}
\end{EQA}
We study the concentration phenomenon of the vector \( \TensG \) around its expectation
in terms of its covariance matrix \( \VBH^{2} = \Var(\TensG) \).
Note that the use of \( \VBH^{2} = \Var(\TensG) \) is not mandatory.
All the results presented later apply with any matrix \( \VBH^{2} \) satisfying \( \VBH^{2} \geq \Var(\TensG) \).
Denote
\begin{EQA}
	\| \Tens \|_{\Fr}^{2}
	& \eqdef &
	\sum_{i,j,k=1}^{\dimp} \Tens_{i,j,k}^{2} \, .
\label{ufy7w3njfiv8re43uy3w6w2}
\end{EQA}
Given \( \uv \in \R^{\dimp} \), define
\begin{EQA}
	\Tens[\uv]
	& \eqdef &
	\sum_{i=1}^{\dimp} u_{i} \, \Tens_{i} \, .
\label{ufhjdeuvbuftjeee3jrjc}
\end{EQA} 
First, describe the covariance structure of \( \TensG \). 

\begin{lemma}
\label{LcovTens}
Denote
\begin{EQA}
	\langle \Tens_{i},\Tens_{\ic} \rangle
	& \eqdef &
	\sum_{j,k=1}^{\dimp} \Tens_{i,j,k} \, \Tens_{\ic,j,k} \, ,
	\qquad
	i,\ic=1,\ldots,\dimp.
\label{9dr5ew53tfvyrhrcgwhcf}
\end{EQA}
Then 
\begin{EQA}
	\VBH^{2}
	& \eqdef &
	\Var(\TensG)
	=
	\bigl( 2 \langle \Tens_{i},\Tens_{\ic} \rangle \bigr)_{i,\ic=1,\ldots,\dimp} \, ,
\label{yhtwttfyg7v53yf7e4cc}
	\\
	\tr \VBH^{2}
	&=&
	2 \sum_{i=1}^{\dimp} \| \Tens_{i} \|_{\Fr}^{2}
	=
	2 \sum_{i,j,k=1}^{\dimp} \Tens_{i,j,k}^{2}
	=
	2 \| \Tens \|_{\Fr}^{2} \, .
\label{kjv8r54jtfu7e367e76fvh}
\end{EQA}
Moreover, 
\begin{EQA}
	\| \VBH \uv \|^{2}
	&=&
	2 \bigl\| \Tens[\uv] \bigr\|_{\Fr}^{2} \, ,
	\qquad
	\uv \in \R^{\dimp} .
\label{cfy6jh9nmk3e365vyekwp}
\end{EQA}
\end{lemma}

\begin{proof}
For any \( i,\ic \), it holds in view of  \( \E (\gauss_{j} \gauss_{k} - \delta_{j,k})^{2} = 1 + \delta_{j,k} \) 
for all \( j,k \leq \dimp \)
\begin{EQA}
	\E (\TensG_{i} - \E \TensG_{i}) (\TensG_{\ic} - \E \TensG_{\ic})
	&=&
	\E \left( 
		\sum_{j,k=1}^{\dimp} \Tens_{i,j,k} \, (\gauss_{j} \gauss_{k} - \delta_{j,k}) 
		\sum_{\jc,\kc=1}^{\dimp} \Tens_{\ic,\jc,\kc} \, (\gauss_{\jc} \gauss_{\kc} - \delta_{\jc,\kc}) 
	\right)
	\\
	&=&
	2 \sum_{j,k=1}^{\dimp} \Tens_{i,j,k} \, \Tens_{\ic,j,k} 
	=
	2 \langle \Tens_{i},\Tens_{\ic} \rangle .
\label{ucviw38vbj4r78fu2122}
\end{EQA}
This yields \eqref{yhtwttfyg7v53yf7e4cc}.
Further
\begin{EQA}
	\tr \VBH^{2}
	&=&
	2 \sum_{i=1}^{\dimp} \langle \Tens_{i},\Tens_{i} \rangle
	=
	2 \sum_{i=1}^{\dimp} \| \Tens_{i} \|_{\Fr}^{2}
	\eqdef
	2 \| \Tens \|_{\Fr}^{2} \, .
\label{kjv8r54jt7e367e76fvhuffg}
\end{EQA}
Similarly, for any \( \uv = (u_{i}) \in \R^{\dimp} \)
\begin{EQA}
	\| \VBH \uv \|^{2}
	&=&
	\uv^{\T} \VBH^{2} \uv
	=
	2 \sum_{i,\ic=1}^{\dimp} u_{i} \, u_{\ic} \langle \Tens_{i},\Tens_{\ic} \rangle
	=
	2 \biggl\| \sum_{i=1}^{\dimp} u_{i} \Tens_{i} \biggr\|_{\Fr}^{2}
	=
	2 \bigl\| \Tens[\uv] \bigr\|_{\Fr}^{2} \, 
	\qquad
\label{yvgk73dhqe12q5rf64eyu0}
\end{EQA}
completing the proof.
\end{proof}

\noindent
Given \( \VBHG^{2} \geq \VBH^{2} \) we characterize regularity of the family \( (\Tens_{i}) \)
by the value \( \tensn \) such that
\begin{EQA}
	\sup_{\uv \colon \| \VBHG \uv \| \leq 1} 2 \, \bigl\| \Tens [\uv] \bigr\|
	& \leq &
	\tensn \, .
\label{jkvhey6fgyebetw3egbgh}
\end{EQA}

\begin{remark}
With \( \VBHG^{2} = \VBH^{2} \), by \eqref{yvgk73dhqe12q5rf64eyu0}, bound \eqref{jkvhey6fgyebetw3egbgh} follows from the condition
\begin{EQA}
	\sqrt{2} \, \bigl\| \Tens [\uv] \bigr\|
	& \leq &
	\tensn \bigl\| \Tens [\uv] \bigr\|_{\Fr} \, , 
	\qquad 
	\uv \in \R^{\dimp} \, .
\label{jkvhey6fgyebetw3egbghS}
\end{EQA}
Clearly this condition meets for \( \tensn = \sqrt{2} \).
We, however, need this condition to be fulfilled with sufficiently small \( \tensn \).
This can be ensured by choosing another matrix \( \VBHG^{2} \geq \VBH^{2} \). 
For instance, with \( \VBHG^{2} = \CONST^{2} \VBH^{2} \), 
the inequalities
\( \sqrt{2} \, \bigl\| \Tens [\uv] \bigr\| \leq \tensn \bigl\| \Tens [\uv] \bigr\|_{\Fr} \) and \( \| \VBHG \uv \| \leq 1 \) imply
\( 2 \, \bigl\| \Tens [\uv] \bigr\| \leq \tensn/\CONST \).
\end{remark}

\Subsection{An upper bound on \( \| \QP (\TensG - \E \TensG) \| \)}
This section presents an upper bound on the norm of \( \QP \errSv \) 
for \( \errSv = \TensG - \E \TensG \) and a linear mapping \( \QP \).
With \( \VBHG^{2} \geq \VBH^{2} \),
define \( \BBH = \QP \VBHG^{2} \QP^{\T} \) and \( \zq^{2}(\BBH,\xx) = \dimH + 2 \vH \sqrt{\xx} + 2 \supH \xx \)
with
\begin{EQ}[rclcl]
	\dimH
	& \eqdef & 
	\tr \BBH
	&=&
	\tr (\QP \VBHG^{2} \QP^{\T}) \, ,
	\\
	\vH^{2}
	&=&
	\tr \BBH^{2}
	&=&
	\tr (\QP \VBHG^{2} \QP^{\T})^{2} \, ,
	\\
	\supH
	& \eqdef &
	\| \BBH \|
	&=&
	\| \QP \VBHG^{2} \QP^{\T} \| \, .
\label{ytuj22w9ciw33y67rfy66}
\end{EQ}
%
A ``high dimensional'' situation means \( \dimH/\supH \) large.
As \( \dimH \supH \geq \vH^{2} \), this implies \( \dimH \gg \vH \). 

\begin{theorem}
\label{Tdevbount3pm}
Assume \eqref{jkvhey6fgyebetw3egbgh} and let \( \gmb \) fulfill \( \tensn \gmb < 1 \).
Given \( \QP \), consider
\begin{EQA}
	\Tnorm
	&=&
	\sqrt{1 - \tensn \gmb} \,\, \| \QP (\TensG - \E \TensG) \| .
\label{jhcvy6e73hjfgihoi56io43}
\end{EQA}
Then with \( \BBH = \QP \VBHG^{2} \QP^{\T} \), \( \dimH,\vH,\supH \) from \eqref{ytuj22w9ciw33y67rfy66}, 
and  \( \xxc \) from \eqref{kv7367ehjgruwwcewyde}, it holds
\begin{EQA}
	\P\bigl( \Tnorm > \zq(\BBH,\xx) \bigr) 
	& \leq &
	3 \ex^{-\xx} ,
	\qquad
	\xx \leq \xxc \, .
\label{ju7cjw20bxdjhwjyh2rtm}
\end{EQA}
Moreover, with \( \zqc = \zq(\BBH,\xxc) \), \( \constg = \frac{\gmb}{\sqrt{\supH} \, (\sqrt{8} + 1)} \),
it holds 
\begin{EQ}[rll]
\label{8vmdfu65e5tgwysjdigir}
	&
	\P\bigl( \Tnorm \geq \zqc + \constg^{-1} (\xx - \xxc) \bigr)
	\leq 
	3 \ex^{ - \xx } ,
	&
	\quad
	\xx \geq \xxc \, ,
	\\
	&
	\P\bigl( \Tnorm \geq \zq \bigr)
	\leq 
	3 \exp \{ - \xxc - \constg (\zq - \zqc) \} ,
	&
	\quad
	\zq \geq \zqc \, .
\end{EQ}
For any \( \zq \leq \zqc \) and \( \nuH \) with \( 2 \nuH \leq \frac{\zq - \sqrt{\dimH}}{\sqrt{\supH}} \), 
it holds
\begin{EQA}
	\E \, \ex^{\nuH \Tnorm} \Ind(\Tnorm \geq \zq)
	& \leq &
	6 \exp\biggl\{  \nuH \zq - \frac{(\zq - \sqrt{\dimH})^{2}}{2 \supH} \biggr\} ,
\label{ufikwk3ei9vgkj4k4giw3wlE}
\end{EQA}
while the condition \( 2 \nuH < \constg = \frac{\gmb}{\sqrt{\supH} \, (\sqrt{8} + 1)} \) ensures 
\begin{EQA}
	\E \, \ex^{\nuH \Tnorm } \Ind(\Tnorm > \zq)
	& \leq &
	6 \, \exp \biggl\{ 
		\nuH \zqc - \frac{(\zqc - \sqrt{\dimH})^{2}}{2 \supH} - (\constg - \nuH) (\zq - \zqc) 
	\biggr\} ,
	\quad
	\zq \geq \zqc \, .
	\qquad
\label{jhf7yehruybyrhe3wevire2E}
\end{EQA}
\end{theorem}

\begin{proof}
Let \( \xiv = \VBHG^{-1} \errSv \).
For any \( \vv \in \R^{\dimp} \) with \( 2 \bigl\| \Tens [\VBHG^{-1} \vv] \bigr\| < 1 \), we check
\begin{EQA}
	\log \E \exp \bigl( \langle \xiv , \vv \bigr\rangle \bigr)
	& \leq & 
	\frac{\| \vv \|^{2}}{2\bigl( 1 - 2 \bigl\| \Tens [\VBHG^{-1} \vv] \bigr\| \bigr)} \, .
	\qquad
\label{929bm5tug78u998ew}
\end{EQA}
Fix \( \vv \in \R^{\dimp} \) and define \( \wv = 2 \VBHG^{-1} \vv \) and \( \Tens [\wv] \) by \eqref{ufhjdeuvbuftjeee3jrjc}.
By Lemma~\ref{Lqfexpmom}, 
\begin{EQA}
	&& \nquad
	\log \E \exp \langle \xiv,\vv \rangle 
	=
	\log \E \exp \langle \errSv,\VBHG^{-1} \vv \rangle 
	\\
	&=&
	\log \E \exp\Bigl( \frac{1}{2} \langle \Tens [\wv], \gaussv^{\otimes 2} \rangle 
		- \frac{1}{2}  \E \langle \Tens [\wv], \gaussv^{\otimes 2} \rangle
	\Bigr)
	\\
	&=&
	\exp \biggl\{ - \frac{\tr (\Tens [\wv])}{2}
	+ \log \det\bigl( \Id_{\dimp} - \Tens [\wv] \bigr)^{-1/2} \biggr\} 
	\leq 
	\frac{\tr (\Tens [\wv])^{2}}{4 (1 - \| \Tens [\wv] \|)} \, .
\label{87jw3433u889wfjhwehgm}
\end{EQA}
By \eqref{cfy6jh9nmk3e365vyekwp}
\begin{EQA}
	\tr (\Tens [\wv])^{2}
	&=&
	\bigl\| 2 \Tens [\VBHG^{-1} \vv] \bigr\|_{\Fr}^{2}
	=
	2 \vv^{\T} \VBHG^{-1} \VBH^{2} \, \VBHG^{-1} \vv
	\leq 
	2 \| \vv \|^{2} ,
	\qquad
\label{iujweujf3fjuejrf45t7m}
\end{EQA}
and \eqref{929bm5tug78u998ew} follows.
If \( \| \vv \| \leq \gmb \), then by \eqref{jkvhey6fgyebetw3egbgh}
\( 2 \bigl\| \Tens [\VBHG^{-1} \vv] \bigr\| \leq \tensn \gmb \), and condition \eqref{expgamgm} is fulfilled
with \( \HVB^{2} = (1 - \tensn \gmb)^{-1} \Id_{\dimp} \).
Now Theorem~\ref{Tdevboundgm} applied to \( \Tnorm = \sqrt{1 - \tensn \gmb} \,\, \| \QP \VBHG \xiv \| \) implies 
\eqref{ju7cjw20bxdjhwjyh2rtm}.
Furthermore, 
Theorem~\ref{TQPxivlarge} with \( \rhoH = 1/2 \) yields \eqref{8vmdfu65e5tgwysjdigir} 
while Corollary~\ref{CTQPxivlarge} yields \eqref{ufikwk3ei9vgkj4k4giw3wlE} and \eqref{jhf7yehruybyrhe3wevire2E}.
\end{proof}

\Subsection{A lower bound}
We also present a lower bound on the quadratic forms \( \| \QP \errSv \|^{2} \).
Here we assume \( \VBHG^{2} = \VBH^{2} \).

\begin{theorem}
\label{Tdevbout3lm}
Assume \eqref{jkvhey6fgyebetw3egbgh} with \( \VBHG^{2} = \VBH^{2} \).
Fix \( \xx \leq \dimH^{2}/(4\vH^{2}) \), \( \muH = 2 \sqrt{\xx}/\vH \), and
\( \alpH < 1/2 \) s.t. 
\begin{EQA}
	\alpH \sqrt{\frac{1 - \alpH}{1 - 2\alpH}}
	& \geq &
	\tensn \sqrt{\dimH} \, \biggl( 1 + \sqrt{\frac{\dimH \supH}{2 \vH^{2}}} \biggr) \, .
\label{gycu3wjmfgvkjigfu4}
\end{EQA}
Then 
\begin{EQA}
	\P\Bigl( \| \QP \errSv \|^{2} - \dimH 
		< - \frac{\alpH \, \dimH}{1 - \alpH} - 2 \vH \sqrt{\xx} 
	\Bigr)
	& \leq &
	2 \ex^{-\xx} .
\label{9fhwe3hikdney76r43kbhnik}
\end{EQA}
\end{theorem}

\begin{proof}
The main step of the proof is a bound on negative exponential moments of \( \| \QP \errSv \|^{2} \).
Then we apply Markov's inequality with a proper choice of the exponent. 
Namely, given \( \muH = 2 \sqrt{\xx}/\vH \) and \( \alpH < 1/2 \), define \( \muHa \) by 
\begin{EQA}
	\frac{1 - 2\alpH}{1 - \alpH} \, \muHa
	&=&
	\muH . 
\label{9vmkfh3yt4hgkjnur44hg}
\end{EQA}
For \( \Gaussv \sim \ND(0,\Id_{\dimp}) \) independent of \( \GaussG \) and \( \imi = \sqrt{-1} \)
\begin{EQA}
	\E \exp \Bigl( - \frac{\muHa}{2} \| \QP \errSv \|^{2} \Bigr) 
	&=&
	\E \, 
	\E_{\Gaussv} \exp \bigl\{ \imi \muHa^{1/2} \bigl\langle \QP \errSv , \Gaussv \bigr\rangle \bigr\}
	\\
	&=&
	\E_{\Gaussv} \, 
	\E \bigl[ \exp \bigl\{ \imi \muHa^{1/2} \bigl\langle \errSv , \QP^{\T} \Gaussv \bigr\rangle \bigr\} \cond \Gaussv \bigr] ,
	\qquad
\label{juvfc7uyeuje3wf7er453}
\end{EQA}
and similarly to \eqref{87jw3433u889wfjhwehgm}, 
it holds by Lemma~\ref{Lqfexpmom} with \( \wv = 2 \muHa^{1/2} \QP^{\T} \Gaussv \)
\begin{EQA}
	\E \left\{ \exp\Bigl( 
		\imi \muHa^{1/2} \langle \errSv, \QP^{\T} \Gaussv \rangle 
	\Bigr) \Cond \Gaussv \right\}
	&=&
	\exp \biggl\{ - \frac{\imi \tr (\Tens [\wv])}{2}
	+ \log \det\bigl( \Id_{\dimp} - \imi \Tens [\wv] \bigr)^{-1/2} \biggr\} .
\label{87jw3433u889wfjhwehgimi}
\end{EQA}
Now introduce a random set
\begin{EQA}
	\Omega(\alpH)
	& \eqdef &
	\bigl\{ 2 \muHa^{1/2} \bigl\| \Tens [\QP^{\T} \Gaussv] \bigr\| \leq \alpH \bigr\} .
\label{hdyt6w6w26fhy3w3wkg}
\end{EQA}
Then \( \| \Tens [\wv] \| \leq \alpH \) on this set and by \eqref{vbu7j3hg8hryhghyidegwdg} of Lemma~\ref{Lqfexpmom}
\begin{EQA}
	\biggl| 
		\log \det\bigl( \Id_{\dimp} - \imi \Tens [\wv] \bigr)^{-1/2} - \frac{\imi \tr (\Tens [\wv])}{2} 
		+ \frac{\tr (\Tens [\wv])^{2}}{4} 
	\biggr|
	& \leq &
	\frac{\alpH \tr (\Tens [\wv])^{2}}{6(1 - \alpH)} \, .
\label{0iuuy54rgeyuduyew3fdr}
\end{EQA}
This implies on \( \Omega(\alpH) \) in view of \eqref{iujweujf3fjuejrf45t7m}
\begin{EQA}
	&& \nquad
	\left| \E \left\{ \exp\Bigl( \imi \muHa^{1/2} \langle \QP \errSv, \Gaussv \rangle \Bigr) \Cond \Gaussv \right\} \right|
	\leq 
	\exp\Bigl( - \frac{(1 - 2 \alpH) \tr (\Tens [\wv])^{2}}{4(1 - \alpH)} \Bigr)
	\\
	&=&
	\exp\Bigl( 
		- \frac{(1 - 2 \alpH) \muHa \tr (2 \Tens [\QP^{\T} \Gaussv])^{2}}{4(1 - \alpH)} 
	\Bigr)
	\\
	&=&
	\exp\Bigl\{ - \frac{\muHa (1 - 2 \alpH)}{1 - \alpH} \, \frac{\| \VBH \QP^{\T} \Gaussv \|^{2}}{2} \Bigr\}
	=
	\exp \Bigl( - \frac{\muH}{2} \| \VBH \QP^{\T} \Gaussv \|^{2} \Bigr) .
\label{ghy6dhwtrrrt26dghe7m}
\end{EQA}
Now, by \eqref{juvfc7uyeuje3wf7er453} and by \eqref{m2v241m41mbn} of Lemma~\ref{Lqfexpmom}
\begin{EQA}
	&& \nquad
	\E \exp \Bigl\{ - \frac{\muHa}{2} \| \QP \errSv \|^{2} \Bigr\} 
	\leq 
	\E_{\Gaussv} \exp\Bigl\{ - \frac{\muH}{2} \| \VBH \QP^{\T} \Gaussv \|^{2} \Bigr\}
	+ \P\bigl( \Omega(\alpH) \bigr)
	\\
	&=&
	\det\bigl( \Id_{\dimp} + \muH \BBH \bigr)^{-1/2} + \P\bigl( \Omega(\alpH) \bigr)
	\leq 
	\exp\Bigl\{ - \frac{\muH \tr (\BBH)}{2} + \frac{\muH^{2} \tr (\BBH^{2})}{4} \Bigr\} 
	+ \P\bigl( \Omega(\alpH) \bigr) .
\label{neu78bjtm6chbweede}
\end{EQA}
For any fixed \( \zq \), by Markov's inequality
\begin{EQA}
	&& \nquad
	\P\bigl( \| \QP \errSv \|^{2} < \zq \bigr)
	\leq
	\exp\Bigl( \frac{\muHa \zq}{2} \Bigr) \E \exp \Bigl( - \frac{\muHa}{2} \| \QP \errSv \|^{2} \Bigr) 
	\\
	& \leq &
	\exp\Bigl\{ \frac{\muHa \zq}{2} - \frac{\muH \tr (\BBH)}{2} + \frac{\muH^{2} \tr (\BBH^{2})}{4} \Bigr\} 
	+ \exp\Bigl( \frac{\muHa \zq}{2} \Bigr) \P\bigl( \Omega(\alpH) \bigr) .
\label{jyhw3jfuy4rhgirfhy3332m}
\end{EQA}
With \( \dimH = \tr \BBH \), \( \vH^{2} = \tr \BBH^{2} \),
and \( \muH = 2 \sqrt{\xx} / \vH \),
define \( \zq \) by
\begin{EQA}
	\frac{\muHa \, \zq}{2}
	&=&
	\frac{\muH}{2} \, \bigl( \dimH - 2 \vH \sqrt{\xx} \bigr) 
	=
	\frac{\dimH \sqrt{\xx}}{\vH} - 2 \xx \, 
\label{8ujekvikjr74rokjkoiowndy}
\end{EQA}
yielding  
\begin{EQA}
	\frac{\muHa \zq}{2} - \frac{\muH \, \dimH}{2} + \frac{\muH^{2} \vH^{2}}{4}
	&=&
	\frac{\muH}{2} \bigl( \dimH - 2 \vH \sqrt{\xx} \bigr)
	- \frac{\muH \, \dimH}{2} + \frac{\muH^{2} \vH^{2}}{4}
	=
	-\xx
\label{nyhwnlfdge634yghiuhyi}
\end{EQA}
while
\begin{EQA}
	\zq 
	&=&
	\frac{1 - 2\alpH}{1 - \alpH} \, (\dimH - 2 \vH \sqrt{\xx})
	\geq 
	\dimH - \frac{\alpH}{1 - \alpH} \, \dimH - 2 \vH \sqrt{\xx} \, .
\label{jfvyuf7y6764gjhwetef}
\end{EQA}
Now we check that
\( \ex^{\muHa\zq /2} \P\bigl( \Omega(\alpH)^{c} \bigr) \leq \ex^{-\xx} \).
By \eqref{jkvhey6fgyebetw3egbgh} 
\( 2 \bigl\| \Tens [\QP^{\T} \Gaussv] \bigr\| \leq \tensn \| \VBH \QP^{\T} \Gaussv \| \), 
and 
\begin{EQA}
	\P\bigl( \Omega(\alpH)^{c} \bigr)
	& \leq &
	\P\bigl( 2 \sqrt{\muHa} \bigl\| \Tens [\QP^{\T} \Gaussv]  \bigr\| > \alpH \bigr)
	\leq 
	\P\bigl( \tensn \sqrt{\muHa} \| \VBH \QP^{\T} \Gaussv \| > \alpH \bigr) .
\label{dyvhb3t6fnruigtkh3yfh}
\end{EQA}
Gaussian deviation bound \eqref{Pxiv2dimAvp12} yields for any \( \xxs > 0 \) by 
\( \| \VBH \QP^{\T} \Gaussv \|^{2} = \Gaussv^{\T} \BBH \Gaussv \)
\begin{EQA}
	\P\bigl( \| \VBH \QP^{\T} \Gaussv \| > \sqrt{\dimH} + \sqrt{2 \xxs \supH} \bigr)
	& \leq &
	\ex^{-\xxs} .
\label{hwt6fyhg4e43euthujr453}
\end{EQA}
By construction, 
\begin{EQA}
	\xx + \frac{\muHa \, \zq}{2}
	&=&
	\xx + \frac{\muH}{2} (\dimH - 2 \vH \sqrt{\xx})
	=
	\frac{\dimH \sqrt{\xx}}{\vH} - \xx
	\leq 
	\frac{\dimH^{2}}{4 \vH^{2}} \, ,
\label{jdy77dy643yry76fwedf}
\end{EQA}
and the use of \( \xxs = \dimH^{2}/(4 \vH^{2}) \) ensures under \eqref{gycu3wjmfgvkjigfu4} 
\begin{EQA}
	\ex^{\muHa \zq/2} \P\bigl( \Omega(\alpH)^{c} \bigr)
	&=&
	\ex^{\muHa \zq/2} \P\biggl( \| \VBH \QP^{\T} \Gaussv \| > \frac{\alpH}{\tensn \sqrt{\muHa}} \biggr)
	\\
	& \leq &
	\ex^{\muHa \zq/2} 
	\P\biggl( \| \VBH \QP^{\T} \Gaussv \| > \sqrt{\dimH} + \frac{\dimH \, \sqrt{\supH}}{\sqrt{2} \, \vH} \biggr)
	\leq 
	\ex^{-\xx} .
\label{bycjq2jfgbyehenwh}
\end{EQA}
Putting this together with \eqref{jyhw3jfuy4rhgirfhy3332m} and \eqref{nyhwnlfdge634yghiuhyi} 
yields \eqref{9fhwe3hikdney76r43kbhnik}.
\end{proof}


\def\TGD{\mathbbmsl{J}}
\def\TG{\mathbb{D}}
\def\GaussD{\gaussv_{\IFL}}

\Section{Some bounds for a third order Gaussian tensor}
\label{SDB3tens}
Let \( \Tens = \bigl( \Tens_{i,j,k} \bigr) \) be a third order symmetric tensor, that is,
\( \Tens_{i,j,k} = \Tens_{\perm(i,j,k)} \) for any permutation \( \perm \) of the triple \( (i,j,k) \).
%
This section present a deviation bound for a Gaussian tensor sum 
\( \Tens(\GaussD) \eqdef \langle \Tens, \GaussD^{\otimes 3} \rangle \) 
for a Gaussian zero mean vector \( \GaussD \sim \ND(0,\IFL^{-1}) \) in \( \R^{\dimp} \).
Much more general results for higher order tensors are available in the literature, see e.g. \cite{GSS2021} and \cite{AW2013} and references therein.
We, however, present an independent self-contained study which delivers finite sample and sharp results.
Later we use notations
\begin{EQA}[ccl]
	\| \Tens \|
	&=&
	\sup_{\| \uv_{1} \| = \| \uv_{2} \| = \| \uv_{3} \| = 1} 
		\bigl| \langle \Tens, \uv_{1} \otimes \uv_{2} \otimes \uv_{3} \rangle \bigr| \, .
\label{jcu7ejufd76ehyew236h}
\end{EQA}
Banach's characterization \cite{Banach1938,nie2017} yields 
\begin{EQA}
	\| \Tens \|
	&=&
	\sup_{\| \uv \| = 1} \bigl| \langle \Tens, \uv^{\otimes 3} \rangle \bigr| \, . 
\label{iv7dcu3921kfvhw7fjew}
\end{EQA}
Define
\begin{EQA}
	\Tens(\uv)
	&=&
	\langle \Tens, \uv^{\otimes 3} \rangle
	=
	\sum_{i,j,k=1}^{\dimp} \Tens_{i,j,k} \, u_{i} \, u_{j} \, u_{k} \, ,
	\qquad
	\uv = (u_{i}) \in \R^{\dimp} \, .
\label{vhjr8vh8r84423rghyrf}
\end{EQA}
Clearly \( \Tens(\uv) \) is a third order polynomial function on \( \R^{\dimp} \).
Define also its gradient \( \nabla \Tens(\uv) \in \R^{\dimp} \).
Each entry of the gradient vector \( \nabla \Tens(\uv) \) is a second order polynomial of \( \uv \).
Symmetricity of \( \Tens \) implies for any \( \uv \in \R^{\dimp} \)
\begin{EQ}[rcccl]
	\nabla \Tens(\uv)
	&=&
	\biggl( 3 \sum_{j,k=1}^{\dimp} \Tens_{i,j,k} \, u_{j} \, u_{k}
	\biggr)_{i=1,\ldots,\dimp} 
	&=&
	3 \bigl( \langle \Tens_{i} \, , \uv \otimes \uv \rangle \bigr)_{i=1,\ldots,\dimp} \,\, ,
	\\
	\nabla^{2} \Tens(\uv)
	&=&
	\biggl( 6 \sum_{i=1}^{\dimp} \Tens_{i,j,k} \, u_{i} 
	\biggr)_{j,k=1,\ldots,\dimp} 
	&=&
	6 \Tens[\uv] \, ,
\label{ho5r9hkrh4ygj32fjur}
\end{EQ}
where \( \Tens_{i} \) is the sub-tensor of order 2 with \( (\Tens_{i})_{j,k} = \Tens_{i,j,k} \) and
\begin{EQA}
	\Tens[\uv]
	& \eqdef &
	\sum_{i=1}^{\dimp} u_{i} \, \Tens_{i} \, .
\label{ufhjdeuvbuftjeee3jrjc3}
\end{EQA} 
Also
\begin{EQA}
	\Tens(\uv)
	&=&
	\frac{1}{3} \langle \nabla \Tens(\uv),\uv \rangle
	=
	\frac{1}{6} \langle \nabla^{2} \Tens(\uv),\uv^{\otimes 2} \rangle \, .
\label{tv8f8e35etfghr67ghj}
\end{EQA}
%
For the norm of the vector \( \nabla \Tens(\uv) \) and of the matrix \( \nabla^{2}\Tens(\uv) \), it holds by \eqref{jcu7ejufd76ehyew236h}
\begin{EQA}
	\| \nabla \Tens(\uv) \|
	&=&
	\sup_{\omegav \in \R^{\dimp} \colon \| \omegav \| = 1} \langle \nabla \Tens(\uv),\omegav \rangle
	=
	\sup_{\omegav \in \R^{\dimp} \colon \| \omegav \| = 1} 3 \langle \Tens,\uv \otimes \uv \otimes \omegav \rangle
	=
	3 \| \Tens \| \, \| \uv \|^{2},
	\\
	\| \nabla^{2} \Tens(\uv) \|
	&=&
	\sup_{\omegav \in \R^{\dimp} \colon \| \omegav \| = 1} 
	\bigl| \bigl\langle \nabla^{2} \Tens(\uv) \, , \omegav \otimes \omegav \bigr\rangle \bigr|
	=
	\sup_{\omegav \in \R^{\dimp} \colon \| \omegav \| = 1} 6 \bigl| \langle \Tens,\uv \otimes \omegav \otimes \omegav \rangle \bigr|
	=
	6 \| \Tens \| \, \| \uv \| .
\label{ucvkje8ij4efui2w65hvgjh}
\end{EQA}

\Subsection{Moments of a Gaussian 3-tensor}
Consider a Gaussian 3-tensor \( \Tens(\gaussv) = \langle \Tens,\gaussv^{\otimes 3} \rangle \).
Define
\begin{EQA}
	\trT_{i}
	&=&
	\sum_{j=1}^{\dimp} \Tens_{i,j,j} 
	=
	\tr \Tens_{i} \, ,
	\qquad
	i = 1,\ldots,n \, .
\label{ydy7w38uw3eiovkeredrgb}
\end{EQA}

\begin{lemma}
\label{LtensFr}
Let \( \Tens = (\Tens_{i,j,k}) \) be a 3-dimensional symmetric tensor in \( \R^{\dimp} \) and 
\( \Tens(\gaussv) = \langle \Tens, \gaussv^{\otimes 3} \rangle \) for \( \gaussv \sim \ND(0,\Id_{\dimp}) \).
With \( \trTv = (\trT_{i}) \in \R^{\dimp} \) and \( \| \Tens \|_{\Fr}^{2} = \sum_{i,j,k = 1}^{\dimp} \Tens_{i,j,k}^{2} \), 
it holds
\begin{EQA}[rcl]
	\E \bigl( \Tens(\gaussv) - 3 \langle \trTv,\gaussv \rangle \bigr)^{2}
	& = &
	6 \| \Tens \|_{\Fr}^{2} \, ,
	\\
	\E \, \Tens^{2}(\gaussv)
	&=&
	6 \| \Tens \|_{\Fr}^{2} + 9 \| \trTv \|^{2} \, .
\label{jhdyuy3e7fghy3r542tdght}
\end{EQA}
\end{lemma}

\begin{proof}
By definition
\begin{EQA}
	&& 
	\Tens(\gaussv) - 3 \langle \trTv,\gaussv \rangle 
	=
	\sum_{i,j,k=1}^{\dimp} \Tens_{i,j,k} \, \gauss_{i} \gauss_{j} \gauss_{k} 
		- 3 \sum_{i=1}^{\dimp} \gauss_{i} \, \sum_{j,k=1}^{\dimp} \Tens_{i,j,k} \, \delta_{j,k} \, .
\label{hfd7j23dmc6hwfgbo}
\end{EQA}
It is easy to see that for each \( i \) by symmetricity of \( \Tens \)
\begin{EQA}
	\E \biggl( \gauss_{i} \sum_{i,j,k=1}^{\dimp} \Tens_{i,j,k} \, \gauss_{i} \gauss_{j} \gauss_{k} \biggr)
	&=&
	3 \sum_{j \in \IIm_{i}} \Tens_{i,j,j} \E (\gauss_{i}^{2} \gauss_{j}^{2}) + \sum_{i=1}^{\dimp} \Tens_{i,i,i} \E \gauss_{i}^{4}
	=
	3 \sum_{j=1}^{\dimp} \Tens_{i,j,j} 
	=
	3 \trT_{i} \, ,
\label{vyehr6bgy4hrtdb34gybn}
\end{EQA}
where the index set \( \IIm_{i} = \{ 1,\ldots,i-1,i+1,\ldots,\dimp \} \) is obtained
by removing the index \( i \) from \( 1,\ldots,\dimp \).
This implies orthogonality 
\begin{EQA}
	\E \bigl\{ 
		\bigl( \Tens(\gaussv) - 3 \langle \trTv,\gaussv \rangle \bigr) \langle \trTv,\gaussv \rangle 
	\bigr\}
	&=&
	0 .
\label{bf6ejh7ehdsnvfetg}
\end{EQA}
Introduce the index set \( \II = \{ (i,j,k) \colon i \neq j \neq k \} \):
\begin{EQA}
	\II
	& \eqdef &
	\{ (i,j,k) \colon \Ind(i=j) + \Ind(i=k) + \Ind(j=k) = 0 \} \, .
\label{iviu38rtfgugh7rdhue}
\end{EQA}
Represent \eqref{hfd7j23dmc6hwfgbo} as
\begin{EQA}
	\Tens(\gaussv) - 3 \langle \trTv,\gaussv \rangle 
	=
	\sum_{\II} \Tens_{i,j,k} \, \gauss_{i} \, \gauss_{j} \,\gauss_{k} 
	+ 3 \sum_{i=1}^{\dimp} \sum_{j \in \IIm_{i}} \Tens_{i,j,j} \, \gauss_{i} (\gauss_{j}^{2} - 1) 
	+ \sum_{i=1}^{\dimp} \Tens_{i,i,i} \, (\gauss_{i}^{3} - 3 \gauss_{i}) \, .
\label{ij3bhj47nedtyvyhwed}
\end{EQA}
All terms in the right hand-side are orthogonal to each other allowing to compute 
\( \E \bigl( \Tens(\gaussv) - 3 \langle \trTv,\gaussv \rangle \bigr)^{2} \):
\begin{EQA}
	&& \nquad\nquad
	\E \bigl( \Tens(\gaussv) - 3 \langle \trTv,\gaussv \rangle \bigr)^{2}
	=
	\E \biggl( \sum_{\II} \Tens_{i,j,k} \, \gauss_{i} \, \gauss_{j} \,\gauss_{k} \biggr)^{2}
	\\
	&&
	+ \, \E \biggl( 
		3 \sum_{i=1}^{\dimp} \sum_{j \in \IIm_{i}} \Tens_{i,j,j} \, \gauss_{i} (\gauss_{j}^{2} - 1) 
	\biggr)^{2}
	+ \E \biggl( \sum_{i=1}^{\dimp} \Tens_{i,i,i} \, (\gauss_{i}^{3} - 3 \gauss_{i}) \biggr)^{2} .
\label{ij3bhjvu37eg47ned2wed}
\end{EQA}
Further, by symmetricity of \( \Tens \)
\begin{EQA}
	&& \nquad
	\E \biggl( \sum_{\II} \Tens_{i,j,k} \, \gauss_{i} \, \gauss_{j} \,\gauss_{k} \biggr)^{2}
	=
	\E \biggl( \sum_{\II} \Tens_{i,j,k} \, \gauss_{i} \, \gauss_{j} \,\gauss_{k} 
	\sum_{\II} \Tens_{\ic,\jc,\kc} \, \gauss_{\ic} \, \gauss_{\jc} \, \gauss_{\kc} \biggr) 
	\\
	&=&
	\E \biggl( \sum_{\II} \Tens_{i,j,k} \, \gauss_{i} \, \gauss_{j} \,\gauss_{k} 
	\sum_{(\ic,\jc,\kc) = \perm(i,j,k)} \Tens_{\ic,\jc,\kc} \, \gauss_{\ic} \, \gauss_{\jc} \, \gauss_{\kc} \biggr)
	=
	6 \sum_{\II} \Tens_{i,j,k}^{2} \, .
\label{hvhfy63trgy6rhdfybhe}
\end{EQA}
Similarly
\begin{EQA}
	\E \biggl( 3 \sum_{i=1}^{\dimp} \sum_{j \in \IIm_{i}} \Tens_{i,j,j} \, \gauss_{i} (\gauss_{j}^{2} - 1) \biggr)^{2}
	&=&
	9 \sum_{i=1}^{\dimp} \sum_{j \in \IIm_{i}} \Tens_{i,j,j}^{2} \, \E \bigl\{ \gauss_{i}^{2} (\gauss_{j}^{2} - 1)^{2} \bigr\}
	=
	18 \sum_{i=1}^{\dimp} \sum_{j \in \IIm_{i}} \Tens_{i,j,j}^{2} \, ,
	\\
	\E \biggl( \sum_{i=1}^{\dimp} \Tens_{i,i,i} \, (\gauss_{i}^{3} - 3 \gauss_{i}) \biggr)^{2}
	&=&
	\sum_{i=1}^{\dimp} \Tens_{i,i,i}^{2} \, \E (\gauss_{i}^{3} - 3 \gauss_{i})^{2}
	=
	6 \sum_{i=1}^{\dimp} \Tens_{i,i,i}^{2}
\label{h6heiubhrhdybferyfn}
\end{EQA}
yielding again by symmetricity of \( \Tens \)
\begin{EQA}
	\E \bigl( \Tens(\gaussv) - 3 \langle \trTv,\gaussv \rangle \bigr)^{2}
	&=&
	6 \sum_{\II} \Tens_{i,j,k}^{2}  
	+ 18 \sum_{i=1}^{\dimp} \sum_{j \in \IIm_{i}} \Tens_{i,j,j}^{2} 
	+ 6 \sum_{i=1}^{\dimp} \Tens_{i,i,i}^{2}
	=
	6 \| \Tens \|_{\Fr}^{2}
\label{jche6e6ghtur5ewgby32}
\end{EQA}
and assertion \eqref{jhdyuy3e7fghy3r542tdght} follows in view of orthogonality \eqref{bf6ejh7ehdsnvfetg}.
\end{proof}

Similarly we study the moments of the scaled gradient vector 
\begin{EQA}
	\TensG
	&=&
	\frac{1}{3} \nabla \Tens(\gaussv) .
\label{gfe98bm3y67vnjnh90t}
\end{EQA}
The entries \( \TensG_{i} \) of \( \TensG \) can be written as 
\( \TensG_{i} = \gaussv^{\T} \Tens_{i} \, \gaussv \); see \eqref{ho5r9hkrh4ygj32fjur}.

\begin{lemma}
\label{LmomGtens}
It holds \( \E \TensG = \trTv \),
\begin{EQA}
	\Var(\TensG)
	&=&
	\VBH^{2}
	=
	\bigl( 2 \langle \Tens_{i},\Tens_{\ic} \rangle \bigr)_{i,\ic=1,\ldots,\dimp} \, ,
\label{u7fyv65eyte56vt6rjikoy}
	\\
	\tr \VBH^{2}
	&=&
	2 \sum_{i=1}^{\dimp} \| \Tens_{i} \|_{\Fr}^{2}
	=
	2 \sum_{i,j,k=1}^{\dimp} \Tens_{i,j,k}^{2}
	=
	2 \| \Tens \|_{\Fr}^{2} \, ,
\label{dyuwe7cvjw24r5cvtjw2b21tc}
	\\
	\E \| \TensG \|^{2}
	&=&
	\| \trTv \|^{2} + 2 \| \Tens \|_{\Fr}^{2} 
	\leq 
	\frac{1}{3} \E \Tens^{2}(\gaussv) \, .
\label{udc6663yv7ryefywjwe7}
\end{EQA}
Moreover, for any \( \uv \in \R^{\dimp} \)
\begin{EQA}
	\| \VBH \uv \|^{2}
	&=&
	2 \bigl\| \Tens[\uv] \bigr\|_{\Fr}^{2} \, .
	\qquad
\label{yvgk73dhqe12q5rf64eyu}
\end{EQA}
\end{lemma}

\begin{proof}
The first statement follows directly from 
\( \E \TensG_{i} = \E \gaussv^{\T} \Tens_{i} \, \gaussv = \tr \Tens_{i} \).
For any \( i,\ic \), it holds in view of  \( \E (\gauss_{j} \gauss_{k} - \delta_{j,k})^{2} = 1 + \delta_{j,k} \) 
for all \( j,k \leq \dimp \)
\begin{EQA}
	\E (\TensG_{i} - \E \TensG_{i}) (\TensG_{\ic} - \E \TensG_{\ic})
	&=&
	\E \left( 
		\sum_{j,k=1}^{\dimp} \Tens_{i,j,k} \, (\gauss_{j} \gauss_{k} - \delta_{j,k}) 
		\sum_{\jc,\kc=1}^{\dimp} \Tens_{\ic,\jc,\kc} \, (\gauss_{\jc} \gauss_{\kc} - \delta_{\jc,\kc}) 
	\right)
	\\
	&=&
	2 \sum_{j,k=1}^{\dimp} \Tens_{i,j,k} \, \Tens_{\ic,j,k} 
	=
	2 \langle \Tens_{i},\Tens_{\ic} \rangle .
\label{ucviw38vbj4r78fu2122}
\end{EQA}
This yields \eqref{u7fyv65eyte56vt6rjikoy}.
Further
\begin{EQA}
	\tr \VBH^{2}
	&=&
	2 \sum_{i=1}^{\dimp} \langle \Tens_{i},\Tens_{i} \rangle
	=
	2 \sum_{i=1}^{\dimp} \| \Tens_{i} \|_{\Fr}^{2}
	\eqdef
	2 \| \Tens \|_{\Fr}^{2} \, ,
\label{kjv8r54jt7e367e76fvhuffg}
\end{EQA}
which proves \eqref{dyuwe7cvjw24r5cvtjw2b21tc}.
Similarly, for any \( \uv = (u_{i}) \in \R^{\dimp} \)
\begin{EQA}
	\| \VBH \uv \|^{2}
	&=&
	\uv^{\T} \VBH^{2} \uv
	=
	2 \sum_{i,\ic=1}^{\dimp} u_{i} \, u_{\ic} \langle \Tens_{i},\Tens_{\ic} \rangle
	=
	2 \biggl\| \sum_{i=1}^{\dimp} u_{i} \Tens_{i} \biggr\|_{\Fr}^{2}
	=
	2 \bigl\| \Tens[\uv] \bigr\|_{\Fr}^{2} \, 
	\qquad
\label{yvgk73dhqe12q5rf64eyupr}
\end{EQA}
completing the proof.
\end{proof}

\Subsection{\( \ell_{3}-\ell_{2} \) condition}
This section introduces a special \( \ell_{3}-\ell_{2} \) condition for a symmetric 3-tensor \( \Tens \).

\begin{description}
    \item[\label{l2l3Tref} \( \bb{(\TG)} \)]
      \textit{For some symmetric \( \dimp \)-matrix \( \TG \) and \( \tensco > 0 \),
      \( \Tens(\uv) = \langle \Tens, \uv^{\otimes 3} \rangle \) fulfills}
\begin{EQA}
	|\Tens(\uv)|
	& \leq &
	\tensco \, \| \TG \uv \|^{3} ,
	\qquad
	\uv \in \R^{\dimp} \, .
\label{7cmvvc7e3hghjj856uiedT}
\end{EQA}
\end{description}


\begin{lemma}
\label{LTensTGm}
Suppose that the tensor \( \Tens \) satisfies \nameref{l2l3Tref}.
Then
\begin{EQA}
	|\langle \Tens,\uv_{1} \otimes \uv_{2} \otimes \uv_{3} \rangle|
	& \leq &
	\tensco \, \| \TG \uv_{1} \| \, \| \TG \uv_{2} \| \, \| \TG \uv_{3} \| \, ,
	\quad
	\uv_{1} \, , \uv_{2} \, , \uv_{3} \in \R^{\dimp} \, ,
	\qquad
\label{h8wxuew7fje4uyruej}
\end{EQA}
and it holds for 
any \( \uv \in \R^{\dimp} \)
\begin{EQA}
\label{8cvmd3r6gfy33gj2wTG2}
	\| \nabla \Tens(\uv) \|
	& \leq &
	3 \tensco \, \| \TG \uv \|^{2} \, \| \TG \| \, ,
	\\
	\Tens[\uv]
	& \leq &
	\tensco \, \| \TG \uv \| \, \TG^{2} \, ,
\label{8cvmd3r6gfy33gj2wfr}
\end{EQA}
yielding
\begin{EQ}[rcl]
	\| \Tens[\uv] \|_{\Fr}^{2}
	& \leq &
	\tensco^{2} \, \| \TG \uv \|^{2} \, \tr (\TG^{4}) \, ,
	\qquad
	\uv \in \R^{\dimp} \, ,
	\\
	\| \Tens \|_{\Fr}^{2}
	& \leq &
	\tensco^{2} \, \tr (\TG^{2}) \, \tr (\TG^{4}) \, .
\label{bhdctsweghfnjggye2qv}
\end{EQ}
Further, for \( \trTv = (\trT_{i}) \in \R^{\dimp} \) with \( \trT_{i} = \tr \Tens_{i} \), it holds
\begin{EQA}
	\| \trTv \|
	& \leq &
	\tensco \, \| \TG \| \, \tr (\TG^{2}) \, ,
\label{hsfd87ewyiewy7er84eewe}
\end{EQA}
The matrix \( \VBH^{2} \) from \eqref{u7fyv65eyte56vt6rjikoy} fulfills
\begin{EQA}
	\VBH^{2}
	& \leq &
	2 \tensco^{2} \tr (\TG^{4}) \, \TG^{2} .
\label{yv8deub78546thi594}
\end{EQA}
It holds for the Gaussian tensor \( \Tens(\gaussv) \) 
\begin{EQA}
	\E \, \Tens^{2}(\gaussv)
	& \leq &
	6 \tensco^{2} \, \tr (\TG^{2}) \, \tr (\TG^{4}) + 9 \tensco^{2} \, \| \TG \|^{2} \, \tr^{2} (\TG^{2})
	\leq 
	15 \tensco^{2} \, \| \TG \|^{2} \, \tr^{2} (\TG^{2}) .
	\qquad
\label{uvjkv62hvub883jdt2jbi}
\end{EQA}
\end{lemma}

\begin{proof}
Define 3-tensor \( \Tens_{\TG} \) by \( \Tens_{\TG}(\uv) = \Tens(\TG^{-1} \uv) \). 
Then condition \eqref{7cmvvc7e3hghjj856uiedT} reads \( |\Tens_{\TG}(\uv)| \leq \tensco \) for all \( \| \uv \| \leq 1 \) while
\eqref{h8wxuew7fje4uyruej} can be written as
\begin{EQA}
	|\langle \Tens_{\TG},\uv_{1} \otimes \uv_{2} \otimes \uv_{3} \rangle|
	& \leq &
	\tensco \, ,
	\quad
	\forall \| \uv_{j} \| \leq 1, \, \, \,
	j=1,2,3.
\label{67che4v5egv5yrr9gjsfyte}
\end{EQA}
The latter holds by Banach's characterization as in \eqref{iv7dcu3921kfvhw7fjew}. 
Further, 
\begin{EQA}
	\| \nabla \Tens(\uv) \|
	&=&
	\sup_{\| \uv_{1} \| = 1} \bigl| \langle \nabla \Tens(\uv),\uv_{1} \rangle \bigr| 
	= 
	\sup_{\| \uv_{1} \| = 1} 3 \bigl| \langle \Tens,\uv \otimes \uv \otimes \uv_{1} \rangle \bigr|
	\\
	& \leq &
	3 \tensco \, \| \TG \uv \|^{2} \, \sup_{\| \uv_{1} \| = 1} \| \TG \uv_{1} \|
	\leq 
	3 \tensco \, \| \TG \uv \|^{2} \, \| \TG \| \, ,
	\\
	\| \Tens[\uv] \|
	&=&
	\sup_{\| \uv_{1} \| = 1} \langle \Tens[\uv],\uv_{1}^{\otimes 2} \rangle 
	= 
	\sup_{\| \uv_{1} \| = 1} \langle \Tens,\uv \otimes \uv_{1} \otimes \uv_{1} \rangle
	\\
	& \leq &
	\tensco \, \| \TG \uv \| \, \sup_{\| \uv_{1} \| = 1} \| \TG \uv_{1} \|^{2}
	\leq 
	\tensco \, \| \TG \uv \| \, \| \TG^{2} \| \, ,
\label{v87buib73jhf7yvbguytgu8}
\end{EQA} 
yielding \eqref{8cvmd3r6gfy33gj2wfr}.
%
Further, \( \langle \trTv, \uv \rangle = \tr \Tens[\uv] \) and by \eqref{8cvmd3r6gfy33gj2wfr}
\begin{EQA}
	\| \trTv \|
	&=&
	\sup_{\| \uv \| = 1} \bigl| \langle \trTv,\uv \rangle \bigr|
	=
	\sup_{\| \uv \| = 1} \bigl| \tr \Tens[\uv] \bigr| 
	\leq 
	\tensco \, \| \TG \| \, \tr (\TG^{2}) \, .
\label{hsfd87ewyiewy7er84eewep}
\end{EQA}
Similarly for \( \uv \in \R^{\dimp} \)
\begin{EQA}
	\bigl\| \Tens [\uv] \bigr\|_{\Fr}^{2}
	&=&	
	\tr (\Tens[\uv]^{2})
	\leq 
	\tensco^{2} \, \| \TG \uv \|^{2} \, \tr (\TG^{4}) \, .
\label{bhdctsweghfnjggye2qvp}
\end{EQA}
Finally, the use of \( \Tens_{i} = \Tens[\ev_{i}] \) for the canonic basis vectors \( \ev_{i} \in \R^{\dimp} \) yields
\begin{EQA}
	\| \Tens \|_{\Fr}^{2}
	&=&
	\sum_{i=1}^{\dimp} \tr (\Tens[\ev_{i}]^{2})
	\leq 
	\tensco^{2} \, \sum_{i=1}^{\dimp} \| \TG \ev_{i} \|^{2} \, \tr (\TG^{4})
	=
	\tensco^{2} \, \tr (\TG^{2}) \, \tr (\TG^{4}) \, ,
\label{bhdctsweghfnjggye2qvp}
\end{EQA}
and \eqref{bhdctsweghfnjggye2qv} follows.
By \eqref{yvgk73dhqe12q5rf64eyu} and \eqref{8cvmd3r6gfy33gj2wfr}, 
it holds for any \( \uv \in \R^{\dimp} \)
\begin{EQA}
	\| \VBH \uv \|^{2}
	&=&
	2 \bigl\| \Tens[\uv] \bigr\|_{\Fr}^{2}
	\leq 
	2 \tensco^{2} \tr (\TG^{4}) \| \TG \uv \|^{2} \, .
\label{yfjf7b6rrtgubuhtweee2w}
\end{EQA}
This yields \eqref{yv8deub78546thi594}.
The obtained bounds lead to \eqref{uvjkv62hvub883jdt2jbi} in view of \eqref{jhdyuy3e7fghy3r542tdght}.
\end{proof}

\Subsection{Colored case}
\label{Stenscolor}
This section extends the established upper bound to the case when the standard Gaussian vector \( \gaussv \)
is replaced by a general zero mean Gaussian vector \( \GaussD \sim \ND(0,\IFL^{-1}) \) for 
a symmetric covariance matrix \( \IFL \).
Then \( \GaussD = \IFL^{-1/2} \gaussv \) with \( \gaussv \) standard normal and 
\( \Tens(\GaussD) = \Tens(\IFL^{-1/2} \gaussv) = \Tenst(\gaussv) \) with \( \Tenst(\uv) = \Tens(\IFL^{-1/2} \uv) \).
If \( \Tens \) satisfies \nameref{l2l3Tref} then \( \Tenst \) does as well
but \( \TG^{2} \) has to be replaced by \( \TGD^{2} = \IFL^{-1/2} \, \TG^{2} \, \IFL^{-1/2} \).

\begin{lemma}
\label{LTensTGmG}
Let \( \Tens(\uv) \) satisfies \nameref{l2l3Tref} with some \( \TG \) and \( \tensco \).
Then \( \Tenst(\uv) = \Tens(\IFL^{-1/2} \uv) \) satisfies \nameref{l2l3Tref} with 
\( \TGD^{2} = \IFL^{-1/2} \, \TG^{2} \, \IFL^{-1/2} \) in place of \( \TG^{2} \) and the same \( \tensco \).
In particular, with \( \Tenst = (\Tenst_{i}) \), \( \trTvt = (\tr \Tenst_{i}) \), and 
\( \VBHt^{2} \eqdef \bigl( 2 \langle \Tenst_{i},\Tenst_{\ic} \rangle \bigr)_{i,\ic=1,\ldots,\dimp} \), it holds
\begin{EQA}
\label{bhdctsweghfnjggye2qvt}
	\| \Tenst \|_{\Fr}^{2}
	& \leq &
	\tensco^{2} \, \tr (\TGD^{2}) \, \tr (\TGD^{4}),
	\\
	\| \trTvt \|
	& \leq &
	\tensco \, \| \TGD \| \, \tr (\TGD^{2}) \, ,
\label{7vjnvc7e3yfgb6ywvbueom}
	\\
	\VBHt^{2}
	& \leq &
	2 \tensco^{2} \, \tr (\TGD^{4}) \, \TGD^{2} ,
\label{7vjnvc7e3yfgb6ywvbueuj}
\end{EQA}
Moreover, for any \( \uv \in \R^{\dimp} \) 
\begin{EQ}[rcl]
	\| \nabla \Tenst(\uv) \|
	& \leq &
	3 \tensco \, \| \TGD \, \uv \|^{2} \, \| \TGD \| \, .
	\\
	\bigl\| \Tenst[\uv] \bigr\|_{\Fr}^{2}
	& \leq &
	\tensco^{2} \, \| \TGD \, \uv \|^{2} \, \tr \TGD^{4} \, ,
\label{vdf8we3hjg7gt4f7yvbg}
\end{EQ}
and
\begin{EQA}
	\E \, \Tens^{2}(\GaussD)
	& \leq &
	15 \tensco^{2} \, \| \TGD \|^{2} \, \tr^{2} (\TGD^{2}) .
\label{uvjkb7nh3f6eyeyudt2jbi}
\end{EQA}
\end{lemma}

\begin{proof}
By definition, for any \( \uv \in \R^{\dimp} \)
\begin{EQA}
	\Tenst(\uv)
	&=&
	\Tens(\IFL^{-1/2} \uv)
	\leq 
	\tensco \| \TG \, \IFL^{-1/2} \uv \|^{3} 
	=
	\tensco \| \TGD \, \uv \|^{3} 
\label{d7vuvuufu7e4e347nsduyc}
\end{EQA}
yielding \nameref{l2l3Tref} for \( \Tenst \).
Now Lemma~\ref{LTensTGmG} enables us to apply the results of 
Lemma~\ref{LTensTGm} with \( \TGD \) in place of \( \TG \).
Finally, for any \( \uv \) with \( \| \TGD \, \uv \| \leq \rrTG \), it holds by \eqref{8cvmd3r6gfy33gj2wTG2}
\begin{EQA}
	\| \nabla \Tenst(\uv) \|
	& \leq &
	3 \tensco \, \| \TGD \, \uv \|^{2} \, \| \TGD \|
	\leq 
	\tensco \, \rrTG^{2} \, \| \TGD \| \, .
\label{f8uvuvuv7e7edgw2hgc7ww}
\end{EQA}
Lemma~\ref{LtensFr} applied to \( \Tens(\GaussD) = \Tenst(\gaussv) \) 
and \eqref{bhdctsweghfnjggye2qvt} and \eqref{7vjnvc7e3yfgb6ywvbueom}
imply \eqref{uvjkb7nh3f6eyeyudt2jbi}.
\end{proof}

\Subsection{Log-Sobolev inequality and Herbst's arguments}
\label{SlogSobolev}
Let \( \dlt(\uv) \) be a smooth function on \( \R^{\dimp} \).
A typical example we have in mind is \( \dlt(\uv) = \Tens(\uv) \), where
\( \Tens \) is a symmetric 3-tensor satisfying \nameref{l2l3Tref} with some \( \TG^{2} \) and \( \tensco \).
Let also \( \GaussD \sim \ND(0,\IFL^{-1}) \).
Our aim is a possibly accurate exponential bound for \( \dlt(\GaussD) \), in particular, for the Gaussian tensor  
\( \Tens(\GaussD) = \langle \Tens,\GaussD^{\otimes 3} \rangle \).
We use \( \dlt(\GaussD) = \dlt(\IFL^{-1/2} \gaussv) = \dltt(\gaussv) \) for \( \dltt(\uv) = \dlt(\IFL^{-1/2} \uv) \)
and \( \gaussv \) standard normal.
The results use a bound on the norm \( \| \nabla \dltt(\uv) \| \) which is hard to verify on the whole domain \( \R^{\dimp} \).
Therefore, we limit the domain of \( \dlt(\uv) \) to a subset \( \UV \) on which the Gaussian measure 
\( \ND(0,\IFL^{-1}) \) well concentrates.
Clearly, for any \( \uv \in \R^{\dimp} \) 
\begin{EQA}
	\nabla \dltt(\uv) 
	&=&
	\IFL^{-1/2} \nabla \dlt(\IFL^{-1/2} \uv) .
\label{uv7yrjrf846yr9df}
\end{EQA}
For the rest of this section, we assume 
that \( \| \nabla \dltt(\uv) \| \) is uniformly bounded over the set of \( \uv \) with \( \IFL^{-1/2} \uv \in \UV \),
where 
\begin{EQA}
	\UV
	& \eqdef &
	\{ \uv \colon \| \TG \uv \| \leq \rr \} .
\label{fc7vhfe367hf78r78ewhc}
\end{EQA}
This applies to \( \dlt(\uv) = \Tens(\uv) \) for a tensor \( \Tens \) satisfying \nameref{l2l3Tref}.
We consider local behavior of \( \dlt(\GaussD) \) for \( \GaussD \sim \ND(0,\IFL^{-1}) \).
With \( \UV \) fixed, introduce the notation 
\begin{EQA}[c]
	\EUV \, \xi \eqdef \E \{ \xi \Ind(\GaussD \in \UV) \} .
\label{yvjby6wehybnbuejwfct}
\end{EQA}
Remind the definition \( \TGD^{2} =	\IFL^{-1/2} \, \TG^{2} \, \IFL^{-1/2} \).
The next lemma explains the choice of the radius \( \rr \) to ensure a concentration effect of \( \GaussD \) on \( \UV \).

\begin{lemma}
\label{PTensG}
For a fixed \( \xx \), set \( \rr = \rr(\xx) = \zq(\TGD^{2},\xx) \) with   
\begin{EQA}
	\zq^{2}(\TGD^{2},\xx) 
	&=& 
	\tr(\TGD^{2}) + 2 \sqrt{\xx \tr(\TGD^{4})} + 2 \xx \| \TGD^{2} \| .
\label{fubnv6e67ehfvyterbue3u}
\end{EQA}
For the set \( \UV \) from \eqref{fc7vhfe367hf78r78ewhc}, suppose
\begin{EQA}
	\sup_{\vv \colon \IFL^{-1/2} \, \vv \, \in \, \UV}
	\| \nabla \dltt(\vv) \|
	=
	\sup_{\uv \, \in \, \UV}
	\| \IFL^{-1/2} \nabla \dlt(\uv) \|
	& \leq &
	\grad \, .
\label{yfuejfivi8eujnbwerwree}
\end{EQA}
Then it holds for \( X = \dlt(\GaussD) - \EUV \, \dlt(\GaussD) \), 
with any \( \muH \) and any integer \( k \)
\begin{EQA}[ccll]
	\EUV \, \ex^{\muH X} 
	& \leq &
	\exp\bigl( \muH^{2} \grad^{2}/2 \bigr) \, ,
	&
\label{v7jdciuwmedfib5i84k}
	\\
	\EUV |X|^{2k} 
	& \leq &
	\CONSTi_{k}^{2} \grad^{2k} \, ,
	&
	\CONSTi_{k}^{2} = 2^{k+1} k! \, .
\label{yjwvue4ufvnwyvne}
\end{EQA}
Also
\begin{EQA}
	\P\Bigl( X > \grad \sqrt{2\xx} \Bigr)
	& \leq &
	2 \ex^{-\xx} .
\label{76djhv83iurfgjhyhP}
\end{EQA}
\end{lemma}

\begin{proof}
With \( \gaussv \sim \ND(0,\Id_{\dimp}) \) and \( \GaussD \sim \ND(0,\IFL^{-1}) \), it holds
\begin{EQA}
	\P(\GaussD \not\in \UV)
	&=&
	\P\bigl( \| \TG \, \IFL^{-1/2} \gaussv \| > \rr \bigr) .
\label{tduc7wquyjhiw6wyhdj}
\end{EQA}
For \( \rr = \zq(\TGD^{2},\xx) \), 
Gaussian concentration bound yields
\begin{EQA}
	\P(\GaussD \not\in \UV)
	=
	\P\bigl( \| \TG \, \IFL^{-1/2} \gaussv \| > \zq(\TGD^{2},\xx) \bigr)
	& \leq &
	\ex^{-\xx}.
\label{xc7y6wqjhc9iwndreae}
\end{EQA}
Further, \( \dlt(\GaussD) = \dlt(\IFL^{-1/2} \gaussv) = \dltt(\gaussv) \) for \( \gaussv \) standard normal
and by \eqref{yfuejfivi8eujnbwerwree}, the norm of the gradient \( \nabla \dltt(\vv) \) is bounded by \( \grad \) 
for all \( \vv \) with \( \IFL^{-1/2} \vv \in \UV \).
The use of log-Sobolev inequality and Herbst's arguments yields \eqref{v7jdciuwmedfib5i84k} for 
\( X = \dlt(\IFL^{-1/2} \gaussv) - \EUV \dlt(\IFL^{-1/2} \gaussv) \);
see Theorem 5.5 in \cite{boucheron2013concentration} or Proposition 5.4.1 in \cite{bakry2013analysis}.
Result \eqref{v7jdciuwmedfib5i84k} also implies the probability bound 
\begin{EQA}
	\P\Bigl( X > \sqrt{2\xx} \, \grad \Bigr)
	& \leq &
	\P(\GaussD \not\in \UV) + \P\Bigl( X > \grad \sqrt{2\xx} \, , \, \GaussD \in \UV \Bigr)
	\leq 
	2 \ex^{-\xx} ;
\label{76djhv83iurfgjhyhy4r}
\end{EQA}
see (5.4.2) in \cite{bakry2013analysis}.
Now Lemma~\ref{Lmomentsexp} and \eqref{v7jdciuwmedfib5i84k} imply \eqref{yjwvue4ufvnwyvne}.
\end{proof}

\ifAnya{}{
\begin{lemma}
\label{Lexmoments}
Let a random variable \( X \) satisfy \( \E \exp (\muH |X|) \leq \CONST \) for some fixed constant \( \CONST \).
Then for any integer \( k \)
\begin{EQA}
	\E |X|^{k}
	& \leq &
	k! \, \muH^{-k} \, \E \exp (\muH |X|) 
	\, .
\label{6cyvhwyvhejhewyex}
\end{EQA}
\end{lemma}

\begin{proof}
Markov inequality imply for any \( u > 0 \) 
\begin{EQA}
	\P\bigl( \muH |X| > u \bigr)
	& \leq &
	\ex^{- u} \E \exp (\muH |X|)
	\leq 
	\CONST \ex^{- u} \, .
\label{gsnvfu3ruhmuhefwex}
\end{EQA}
Hence, 
\begin{EQA}
	\E |\muH X|^{k} 
	&=&
	\int_{0}^{\infty} \P\bigl( |\muH X|^{k} > x \bigr) \, dx
	=
	k \int_{0}^{\infty} x^{k-1} \P\bigl( |\muH X| > x \bigr) \, dx
	\\
	& \leq &
	\CONST k \int_{0}^{\infty} x^{k-1} \ex^{-x} \, dx
	= 
	\CONST k!  
\label{jvt6wehbf7y67e34nex}
\end{EQA}
and the result follows.
\end{proof}
}

\begin{lemma}[\cite{boucheron2013concentration}, Theorem 2.1]
\label{Lmomentsexp}
Let a r.v. \( X \) satisfy 
\( \E \exp (\muH X) \leq \exp(\muH^{2} \grad^{2}/2) \) for all \( \muH \) with some \( \grad^{2} > 0 \).
Then for any integer \( k \)
\begin{EQA}
	\E |X|^{2k}
	& \leq &
	\CONSTi_{k}^{2} \, \grad^{2k} \, ,
	\qquad
	\CONSTi_{k}^{2} = 2^{k+1} k! 
	\, .
\label{6cyvhwyvhejhewycnwsbj}
\end{EQA}
In particular, \( \CONSTi_{1} = 2 \), \( \CONSTi_{2} = 4 \), \( \CONSTi_{3} = \sqrt{96} \leq 10 \),
\( \CONSTi_{4} = 16 \sqrt{3} \leq 28 \).
\end{lemma}

\begin{proof}
Conditions of the lemma and Markov inequality imply for any \( u > 0 \) with \( \mu = u \)
\begin{EQA}
	\P\bigl( X/ \grad > u \bigr)
	& \leq &
	\ex^{- \muH u} \E \exp (\muH X/\grad)
	\leq 
	\exp(- u^{2} /2)
\label{gsnvfu3ruhadsqsgvdfw}
\end{EQA}
and similarly for \( \P\bigl( - X/ \grad > u \bigr) \) hence, 
\begin{EQA}
	\E |X/\grad|^{2k} 
	&=&
	\int_{0}^{\infty} \P\bigl( |X/\grad|^{2k} > x \bigr) \, dx
	=
	2k \int_{0}^{\infty} x^{2k-1} \P\bigl( |X/\grad| > x \bigr) \, dx
	\\
	& \leq &
	4k \int_{0}^{\infty} x^{2k-1} \ex^{-x^{2}/2} \, dx
	= 
	4k \int_{0}^{\infty} (2t)^{k-1} \ex^{-t} \, dt
	=
	2^{k+1} k! 
\label{jvt6wehbf7y67e34nhfdd}
\end{EQA}
as claimed. 
\end{proof}

Let \( X \) satisfy \( \E \exp (\muH X) \leq \exp(\muH^{2} \grad^{2}/2) \) for all \( \muH \) with some \( \grad \) small.
One can expect that \( \ex^{X} \) can be well approximated for \( k \geq 2 \) by
\begin{EQA}
	\EX_{k}(X)
	& \eqdef &
	1 + X + \ldots + \frac{X^{k-1}}{(k-1)!} \, .
\label{tyfcyhwwuv7wj2apc7ew3j}
\end{EQA}

\begin{lemma}
\label{PTensG3}
Let a random variable \( X \) satisfy 
\( \E \exp (\muH X) \leq \exp(\muH^{2} \grad^{2}/2) \) for all \( \muH \) with some \( \grad^{2} > 0 \).
Then for a random variable \( \xi \) such that \( |\xi| \leq 1 \) and any integer \( k \) with \( \CONSTi_{k} \) from \eqref{yjwvue4ufvnwyvne} and \( \EX_{k}(X) \) from \eqref{tyfcyhwwuv7wj2apc7ew3j}
\begin{EQA}
	\bigl| 
		\E \bigl( \ex^{X} - \EX_{k}(X) \bigr) \xi
	\bigr|
	& \leq &
	\frac{\CONSTi_{k}}{k!} \, \grad^{k} \ex^{\grad^{2}} \, .
\label{7vnf53g6fy363ghfughd}
\end{EQA}
In particular, with \( \CONSTi_{2} = 4 \) and \( \CONSTi_{3}^{2} = 96 \)
\begin{EQ}[rcl]
	\bigl| 
		\E \bigl( \ex^{X} - 1 - X \bigr) \xi 
	\bigr|
	& \leq &
	2 \grad^{2} \, \ex^{\grad^{2}} \, ,
	\\
	\bigl| 
		\E \bigl( \ex^{X} - 1 - X - \frac{X^{2}}{2} \bigr) \xi 
	\bigr|
	& \leq &
	\frac{5}{3} \grad^{3} \, \ex^{\grad^{2}} \, .
\label{7vnf53g6fy363ghfughd2}
\end{EQ}
If \( \xi \) is not bounded but \( \E \xi^{2k+2} < \infty \), then with \( \rho = k/(k+1) \) 
\begin{EQA}
	\bigl| 
		\E \bigl( \ex^{X} - \EX_{k}(X) \bigr) \xi 
	\bigr|
	& \leq &
	\frac{\CONSTi_{k+1}^{\rho}}{k!} \, \grad^{k} \, \ex^{\grad^{2}} \, \bigl( \E \xi^{2k+2} \bigr)^{\frac{1}{2k+2}} \, .
\label{7vnf53g6fy363ghfugh62d}
\end{EQA}
In particular, with \( \CONSTi_{3}^{2/3} = 96^{1/3} \leq 4.6 \) and \( \CONSTi_{4}^{3/4} \leq 12 \)
\begin{EQA}
	\bigl| 
		\E \bigl( \ex^{X} - 1 - X \bigr) \xi 
	\bigr|
	& \leq &
	2.3 \, \grad^{2} \, \ex^{\grad^{2}} \, \bigl( \E \xi^{6} \bigr)^{1/6} \, ,
\label{7vnf53g6fy363ghfugh2}
	\\
	\bigl| 
		\E \bigl( \ex^{X} - 1 - X - \frac{X^{2}}{2} \bigr) \xi 
	\bigr|
	& \leq &
	2 \, \grad^{3} \, \ex^{\grad^{2}} \, \bigl( \E \xi^{8} \bigr)^{1/8} \, .
\label{7vnf53g6fy363ghfugh4}
\end{EQA}
\end{lemma}

\begin{proof}
Define 
\begin{EQA}
	\riskt(t)
	& \eqdef &
	\E \bigl\{ \bigl( \ex^{t X} - \EX_{k}(tX) \bigr) \xi \bigr\}  \, .
\label{cuvh46fghb83bbuwbdu}
\end{EQA}
Obviously \( \riskt(0) = \riskt'(0) = \ldots = \riskt^{(k-1)}(0) = 0 \).
The Taylor expansion of order \( k \) yields
\begin{EQA}
	\bigl| \riskt(1) \bigr|
	& \leq &
	\frac{1}{k!} \sup_{t \in [0,1]} |\riskt^{(k)}(t)| .
\label{hjvg7ejfi9i4nretcrwty}
\end{EQA}
Further, 
\begin{EQA}
	\riskt^{(k)}(t)
	& = &
	\E (X^{k} \,  \xi \, \ex^{t X}) \, .
\label{hfywhjfvueybwfsgh}
\end{EQA}
Consider first the case \( |\xi| \leq 1 \) a.s.
By the Cauchy-Schwarz inequality, \eqref{v7jdciuwmedfib5i84k} of Lemma~\ref{PTensG}, and \eqref{6cyvhwyvhejhewycnwsbj},
it holds  for any \( t \in [0,1] \) 
\begin{EQA}
	|\riskt^{(k)}(t)|^{2}
	& \leq &
	\E |X|^{2k} \,\, \E \, \ex^{2 t X} 
	\leq 
	\CONSTi_{k}^{2} \grad^{2k} \ex^{2 \grad^{2}}
	\, .
\label{hgdrbketrdefvbdtctgewb}
\end{EQA}
For a general \( \xi \), in a similar way, it holds with \( \rho = k/(k+1) \)
\begin{EQA}
	|\riskt^{(k)}(t)|^{2}
	& \leq &
	\E (|X|^{2k} \xi^{2}) \,\, \E \, \ex^{2 t X} 
	\leq 
	\bigl( \E |X|^{2k+2} \bigr)^{\rho} \bigl( \E \xi^{2k+2} \bigr)^{1 - \rho} \,\, \E \, \ex^{2 \grad^{2}} 
	\\
	& \leq &
	\CONSTi_{k+1}^{2\rho} \grad^{2k} \ex^{2 \grad^{2}} \, \bigl( \E \xi^{2k+2} \bigr)^{1 - \rho} \, ,
\label{hv6gdffv5egd5wghb87ruwe}
\end{EQA}
and \eqref{7vnf53g6fy363ghfugh62d} follows.
\end{proof}

This result with \( \xi = 1 \) yields an approximation \( \E \ex^{X} \approx 1 + \E X + \E X^{2}/2 \)
and with \( \xi = X \) an approximation \( \E (X \ex^{X}) \approx \E X + \E X^{2} \).

\begin{lemma}
\label{PTensGEX}
Let a random variable \( X \) satisfy 
\( \E \exp (\muH X) \leq \exp(\muH^{2} \grad^{2}/2) \) for all \( \muH \) with some \( \grad^{2} > 0 \).
Then
\begin{EQA}
	\bigl| \E \ex^{X} - 1 - \E X - \E X^{2}/2) \bigr|
	& \leq &
	2 \grad^{3} \ex^{\grad^{2}} \, ,
	\\
	\bigl| \E (X \ex^{X}) - \E X - \E X^{2}) \bigr|
	& \leq &
	5 \grad^{3} \ex^{\grad^{2}} \, .
\label{iucfjhc76ehfyrujgbge}
\end{EQA}
\end{lemma}

\begin{proof}
The first bound follows from \eqref{7vnf53g6fy363ghfugh4} with \( \xi \equiv 1 \).
Further, \eqref{6cyvhwyvhejhewycnwsbj} for \( k=3 \) implies \( \E X^{6} \leq 96 \grad^{6} \) and
\eqref{7vnf53g6fy363ghfugh2} with \( \xi = X \) yields the second bound.
\end{proof}

Now we specify the obtained bounds for two scenarios. 
Let \( \lgd \) be a function on \( \R^{\dimp} \).
First we consider a symmetric 3-tensor \( \Tens \) which can be viewed as third order derivative of \( \lgd \) 
at some point \( \xv \) and define \( X = \Tens(\GaussD) \) for \( \GaussD \sim \ND(0,\DVL^{-2}) \).

\begin{lemma}
\label{PtensGTens}
Let \( \Tens \) be a symmetric 3-tensor \( \Tens \) satisfying \nameref{l2l3Tref} and \( \GaussD \sim \ND(0,\IFL^{-1}) \).
Consider the set \( \UV \) from Lemma~\ref{PTensG}.
Then all the statements of Lemma~\ref{PTensG} and Lemma~\ref{PTensG3} continue to apply with \( X = \Tens(\GaussD) \)
and \( \grad \eqdef 3 \, \tensco \, \rr^{2} \| \TGD \| \).
In particular, it holds for any \( \muH \) and any integer \( k \)
\begin{EQA}[ccl]
	\EUV \, \ex^{\muH \Tens(\GaussD)} 
	& \leq &
	\exp\bigl( \muH^{2} \grad^{2}/2 \bigr) ,
\label{v7jdciuwmedfib5i84kT}
	\\
	\EUV |\Tens(\GaussD)|^{2k} 
	& \leq &
	\CONSTi_{k}^{2} \grad^{2k} \, , \qquad \qquad \CONSTi_{k}^{2} = 2^{k+1} k! \, .
\label{yjwvue4ufvnwyvneT}
\end{EQA}
\end{lemma}

\begin{proof}
Condition \nameref{l2l3Tref} ensures that 
the gradient of \( \Tens(\DVL^{-1} \uv) \) is uniformly bounded by
\( \grad = 3 \, \tensco \, \rr^{2} \| \TGD \| \) on the local set \( \UV \); see Lemma~\ref{LTensTGmG}.
This enables the statements of Lemma~\ref{PTensG} and Lemma~\ref{PTensG3}.
\end{proof}

For the second scenario, let \( \dlt(\uv) \) be  the third order remainder in the Taylor expansion of 
\( \lgd(\xv + \uv) \) at a fixed point \( \xv \):
\begin{EQA}
	\dlt(\uv)
	& \eqdef &
	\lgd(\xv + \uv) - \lgd(\xv) - \langle \nabla \lgd(\xv), \uv \rangle 
	- \frac{1}{2} \langle \nabla^{2} \lgd(\xv), \uv^{\otimes 2} \rangle 
\label{ytvcjf63wehg7ehfyjjf}
\end{EQA}
and consider \( \dlt(\GaussD) \).
In this case we assume that \( \lgd \) satisfies the following condition.
 
\begin{description}
    \item[\label{l2l3gref} \( \bb{(\TG_{3})} \)]
    \textit{For some \( \TG \) and \( \dltwu_{3} > 0 \), it holds 
     for any \( \uv \in \UV = \{ \uv \colon \| \TG \uv \| \leq \rr \} \),
      }
\begin{EQA}
	|\langle \nabla^{3} \lgd(\xv + \uv), \uv_{1}^{\otimes 3} \rangle|
	& \leq &
	\dltwu_{3} \, \| \TG \uv_{1} \|^{3} ,
	\qquad
	\uv_{1} \in \R^{\dimp} .
\label{vjb7hyer5ewre5fty4fghg}
\end{EQA}
\end{description}

Banach's characterization \cite{Banach1938} yields for any \( \uv_{1} , \uv_{2}, \uv_{3} \in \R^{\dimp} \)
\begin{EQA}
	\bigl| \langle \nabla^{3} \lgd(\xv + \uv), \uv_{1} \otimes \uv_{2} \otimes \uv_{3} \rangle \bigr|
	& \leq &	 
	\dltwu_{3} \, \| \TG \uv_{1} \| \, \| \TG \uv_{2} \| \, \| \TG \uv_{3} \| \, ;
\label{jbufiv784ejf76e3n94}
\end{EQA}
see Lemma~\ref{LTensTGm}.

\begin{lemma}
\label{PtensGdlt}
Let a function \( \lgd \) satisfy \nameref{l2l3gref} and \( \GaussD \sim \ND(0,\IFL^{-1}) \).
Define \( \TGD \) by \( \TGD^{2} = \IFL^{-1/2} \, \TG^{2} \, \IFL^{-1/2} \).
Then all the statements of Lemma~\ref{PTensG} and Lemma~\ref{PTensG3} continue to apply with \( X = \dlt(\GaussD) \)
for \( \dlt(\uv) \) from \eqref{ytvcjf63wehg7ehfyjjf}
and \( \grad \eqdef \dltwu_{3} \, \rr^{2} \| \TGD \|/2 \).
\end{lemma}

\begin{proof}
Define \( \dltt(\uv) = \dlt(\IFL^{-1/2} \uv) \). 
Note that \( \IFL^{-1/2} \uv \in \UV \) means \( \| \TG \, \IFL^{-1/2} \uv \| = \| \TGD \uv \| \leq \rr \).
We only have to check that condition \nameref{l2l3gref} implies with \( \grad = \dltwu_{3} \, \rr^{2} \| \TGD \|/2 \)
\begin{EQA}
	\sup_{\uv \colon \| \TGD \uv \| \leq  \rr} \| \nabla \dltt(\uv) \|
	=
	\sup_{\uv \colon \| \TGD \uv \| \leq  \rr} \| \IFL^{-1/2} \nabla \dlt(\uv) \|
	& \leq &
	\grad .
\label{i8foi9eoeg8j4dfgifert9i}
\end{EQA}
Indeed, the Taylor expansion at \( \uv = 0 \) yields by \( \nabla \dlt(0) = 0 \) and \( \nabla^{2} \dlt(0) = 0 \)
\begin{EQA}
	\| \IFL^{-1/2} \nabla \dlt(\uv) \|
	& \leq &
	\sup_{\| \ww \| = 1} \langle \nabla \dlt(\uv), \IFL^{-1/2} \wv \rangle 
	=
	\sup_{\| \ww \| = 1} \langle \nabla \dlt(\uv) - \nabla \dlt(0) - \nabla^{2} \dlt(0) \uv, \IFL^{-1/2} \wv \rangle
	\\
	&=&
	\frac{1}{2} \sup_{\| \ww \| = 1} \langle \nabla^{3} \dlt(t \uv) , \uv \otimes \uv \otimes \IFL^{-1/2} \wv \rangle
\label{bn9m3gf8rjdgtgi8evig}
\end{EQA}
By \eqref{jbufiv784ejf76e3n94}
\begin{EQA}
	\| \IFL^{-1/2} \nabla \dlt(\uv) \|
	& \leq &
	\frac{1}{2} \sup_{\| \ww \| = 1} \langle \nabla^{3} \dlt(t \uv) , \uv \otimes \uv \otimes \IFL^{-1/2} \wv \rangle
	\leq 
	\frac{\dltwu_{3}}{2} \sup_{\| \ww \| = 1} \| \TG \uv \| \, \| \TG \uv \| \, \| \TG \, \IFL^{-1/2} \wv \|
	\\
	& \leq &
	\frac{\dltwu_{3}}{2} \, \rr^{2} \, \sup_{\| \ww \| = 1} \, \| \TGD \wv \|
	=
	\frac{\dltwu_{3}}{2} \, \rr^{2} \, \| \TGD \| \, 
\label{hgf7ehf7ru3nmghkhucdjh}
\end{EQA}
as required.
\end{proof}


\def\gaussF{\gauss_{\IFL}}
\Section{Local Laplace approximation}
\label{SlocalLaplace}

This section presents some bounds on the error of local Laplace approximation.
Let \( \lgd(\xv) \) be a function in a high-dimensional Euclidean space \( \R^{\dimp} \) such that
\( \int \ex^{\lgd(\xv)} \, d\xv = \CONST < \infty \),
where the integral sign \( \int \) without limits means the integral over the whole space \( \R^{\dimp} \).
Then \( \lgd \) determines a distribution \( \PfL \) with the density
\( \CONST^{-1} \ex^{\lgd(\xv)} \).
Let \( \xvs \) be a point of maximum:
\begin{EQA}
	\lgd(\xvs)
	&=&
	\sup_{\uv \in \R^{\dimp}} \lgd(\xvs + \uv) .
\label{scdygw7ytd7wqqsquuqydtdtd}
\end{EQA}
We also assume that \( \lgd(\cdot) \) is at least three time differentiable. 
Introduce the negative Hessian \( \IFL = - \nabla^{2} \lgd(\xvs) \).
Later we assume that the negative Hessian \( \IFL = - \nabla^{2} \lgd(\xvs) \) is sufficiently large
in the sense that the Gaussian measure \( \ND(0,\IFL^{-1}) \) concentrates on a small local set \( \UVL \).
This allows to use a local Taylor expansion for 
\( \lgd(\xvs;\uv) \approx - \| \IFL^{1/2} \uv \|^{2}/2 \) in \( \uv \) on \( \UVL \).
For this local set \( \UVL \), we evaluate the quantity
\begin{EQA}
	\err
	& \eqdef &
	\biggl| \frac{\int_{\UVL} \ex^{\lgd(\xvs + \uv) - \lgd(\xvs)} \, d\uv - \int_{\UVL} \ex^{- \| \IFL^{1/2} \uv \|^{2}/2} \, d\uv} 
		 {\int \ex^{- \| \IFL^{1/2} \uv \|^{2}/2} d\uv} 
	\biggr| \, .
\label{fio9vkmfg763eu7g8gyhffr}
\end{EQA}
As \( \xvs = \argmax_{\xv} \lgd(\xv) \), it holds \( \nabla \lgd(\xvs) = 0 \) and 
\begin{EQA}
	\err
	& = &
	\biggl| \frac{\int_{\UVL} \ex^{\lgd(\xvs;\uv)} \, d\uv - \int_{\UVL} \ex^{- \| \IFL^{1/2} \uv \|^{2}/2} \, d\uv} 
		 {\int \ex^{- \| \IFL^{1/2} \uv \|^{2}/2} d\uv} 
	\biggr| ,
\label{IIgfifxutudufxutdu}
\end{EQA}
where \( \lgd(\xv;\uv) \) is the Bregman divergence 
\begin{EQA}
	\lgd(\xv;\uv)
	&=&
	\lgd(\xv + \uv) - \lgd(\xv) - \bigl\langle \nabla \lgd(\xv), \uv \bigr\rangle .
\label{fxufxpufxfpxu}
\end{EQA}
Assume that \( \lgd(\cdot) \) be a four times continuously differentiable function on \( \R^{\dimp} \).
Consider the remainder of the second and third-order Taylor approximation 
\begin{EQ}[rcl]
	\dltw_{3}(\uv)
	&=&
	\lgd(\xvs;\uv) - 
	\bigl\langle \nabla^{2} \lgd(\xvs) , \uv^{\otimes 2} \bigr\rangle/2 ,
	\\
	\dltw_{4}(\uv)
	&=&
	\lgd(\xvs;\uv) - 
	\bigl\langle \nabla^{2} \lgd(\xvs) , \uv^{\otimes 2} \bigr\rangle/2 
	- 
	\bigl\langle \nabla^{3} \lgd(\xvs), \uv^{\otimes 3} \bigr\rangle / 6 ,
\label{d4fuv1216303}
\end{EQ}
where \( \lgd(\xv;\uv) \) is given by \eqref{fxufxpufxfpxu}.
We will use the decomposition
\begin{EQA}
	\lgd(\xvs;\uv)
	&=&
	- \frac{1}{2} \| \IFL^{-1/2} \uv \|^{2} + \dltw_{3}(\uv)
	=
	- \frac{1}{2} \| \IFL^{-1/2} \uv \|^{2} + \Tens(\uv) + \dltw_{4}(\uv) ,
\label{yvjd7e3jgir63hfgkdd}
\end{EQA} 
where \( \Tens(\uv) = \langle \nabla^{3} \lgd(\xvs), \uv^{\otimes 3} \rangle/6 \) is the third order tensor
corresponding to the third derivative in the fourth order Taylor expansion for \( \lgd(\xvs;\uv) \).
For ease of notation, we skip dependence of \( \Tens \), \( \dltw_{3} \), and \( \dltw_{4} \) on \( \xvs \).

The local set \( \UVL \) will be described using a metric tensor \( \DVL \).
For for some \( \rrL > 0 \),
\begin{EQA}
	\UVL 
	&=& 
	\bigl\{ \uv \colon \| \DVL \uv \| \leq \rrL \bigr\} .
\label{UvTDunm12spT}
\end{EQA}
We assume that \( \lgd \) is sufficiently smooth within \( \UVL \) and satisfies the following conditions.

\begin{description}
    \item[\label{l2l3sref} \( \bb{(\DVL_{3}^{*})} \)]
      \textit{\( \DVL^{2} \leq \IFL \) and for some \( \dltwu_{3} > 0 \), }
\begin{EQA}
	\sup_{\xv \colon \xv - \xvs \in \UVL} |\langle \nabla^{3} \lgd(\xv),\uv^{\otimes 3} \rangle|
	& \leq &
	\dltwu_{3} \, \| \DVL \uv \|^{3} ,
	\qquad
	\uv \in \R^{\dimp} \, .
\label{7cmvvc7e3hghjj856uied}
\end{EQA}
    \item[\label{l2l4ref} \( \bb{(\DVL_{4})} \)]
    \textit{For some \( \dltwu_{4} > 0 \), 
    }
\begin{EQA}
	|\dltw_{4}(\uv)|
	& \leq &
	\frac{\dltwu_{4}}{24} \| \DVL \uv \|^{4} ,
	\qquad
	\uv \in \UVL .
\label{vjb7hyer5ewre5fty4fgh}
\end{EQA}
\end{description}

Expansion \eqref{yvjd7e3jgir63hfgkdd} allows to represent 
the error \( \err \) from \eqref{IIgfifxutudufxutdu} as
\begin{EQA}
	\err
	=
	\frac{\int_{\UVL} \ex^{\lgd(\xvs;\uv)} \, d\uv - \int_{\UVL} \ex^{- \| \IFL^{1/2} \uv \|^{2}/2} \, d\uv} 
		 {\int \ex^{- \| \IFL^{1/2} \uv \|^{2}/2} d\uv} 
	&=&
	\EUV \Bigl[ \bigl\{ \exp \dltw_{3}(\gaussF) - 1 \bigr\} \Bigr] \, ,
\label{ed3le2d3pGd}
\end{EQA}
where \( \gaussF \sim \ND(0,\IFL^{-1}) \) and \( \EUV \xi \) means \( \E \{ \xi \Ind(\gaussF \in \UVL) \} \).

\begin{proposition}
\label{PTensG42}
Assume \nameref{l2l3sref}, \nameref{l2l4ref}.
Define
\( \TGD^{2} = \IFL^{-1/2} \DFL^{2} \, \IFL^{-1/2} \),
\begin{EQA}[c]
	\grad = \dltwu_{3} \, \rr^{2} \| \TGD \|/2 ,
	\qquad
	\accu_{\IFL}^{2} = \E \Tens^{2}(\gaussF) ,
	\qquad 
	\dltw_{4,\IFL} = \EUV \dltw_{4}^{2}(\gaussF) .
\label{•}
\end{EQA}
Then 
\begin{EQA}[rcccl]
	\Bigl| \err - \frac{\accu_{\IFL}^{2}}{2} \Bigr|
	& \leq &
	\accu_{\IFL} \, \dltw_{4,\IFL} 
	+ \frac{\dltw_{4,\IFL}^{2}}{2} + \frac{5}{3} \grad^{3} \ex^{\grad^{2}} \, ,
	\quad
	\err
	& \leq &
	\frac{1}{2} (\accu_{\IFL} + \dltw_{4,\IFL})^{2} + \frac{5}{3} \grad^{3} \ex^{\grad^{2}} \, .
	\qquad
\label{u7yfjeoi8g8tj5udftrugh}
\end{EQA}
Moreover, 
\begin{EQA}[rcl]
\label{gdt6djiu97ie4425g7}
	\dltw_{4,\IFL}
	& \leq &
	\frac{1}{24} \, \dltwu_{4} \, \bigl\{ \tr (\TGD^{2}) + 3 \| \TGD^{2} \| \bigr\}^{2} \, ,
	\\
	\accu_{\IFL}
	& \leq &
	\sqrt{\frac{5}{12}} \,\, \dltwu_{3} \, \| \TGD \| \, \tr (\TGD^{2}) \, .
\label{gdt6djiu97ie4425g77}
\end{EQA}
\end{proposition}

\begin{proof}
%
%
Condition \nameref{l2l3sref} enables us to apply
Lemma~\ref{PtensGdlt} with \( X = \dltwhat_{3}(\gaussF) \) and \( k=3 \).
This yields 
\begin{EQA}
	\Bigl| 
		\EUV \Bigl\{ \Bigl(	\ex^{\dltwhat_{3}(\gaussF)} - 1 - \dltwhat_{3}(\gaussF) 
			- \frac{\dltwhat_{3}^{2}(\gaussF)}{2} \Bigr) g(\gaussF) 
		\Bigr\}
	\Bigr| 
	& \leq &
	\frac{5}{3} \grad^{3} \ex^{\grad^{2}} \, ,
\label{7nhig763hf5bv9edfuejn}
\end{EQA}
for any \( g(\cdot) \) with \( \sup_{\uv \in \UV} |g(\uv)| \leq 1 \).
Further, by \nameref{l2l4ref} and Lemma~\ref{Gaussmoments}
\begin{EQA}
	\EUV \dltw_{4}^{2}(\gaussF)
	& \leq &
	\frac{\dltwu_{4}^{2}}{24^{2}} \, \E \| \DFL \, \IFL^{-1/2} \gaussv \|^{8}
	\leq 
	\frac{\dltwu_{4}^{2}}{24^{2}} \, \bigl\{ \tr(\TGD^{2}) + 3 \| \TGD^{2} \| \bigr\}^{4}
\label{yhdyvcyhwjg8it5j5uwfg}
\end{EQA}
and \eqref{gdt6djiu97ie4425g7} follows.
As \( \dltwhat_{3}(\gaussF) = \Tens(\gaussF) + \dltwhat_{4}(\gaussF) \), it holds
\begin{EQA}[rcl]
	\EUV |\dltwhat_{3}(\gaussF) - \Tens(\gaussF)| 
	&=& 
	\EUV |\dltwhat_{4}(\gaussF)|
	\leq 
	\sqrt{\EUV \dltwhat_{4}^{2}(\gaussF)}
	\leq 
	\sqrt{\EUV \dltw_{4}^{2}(\gaussF)} \, ,
\label{ycsy6ty25cgc53c7vyeetch}
\end{EQA}
and by \( \dltwhat_{4}(\gaussF) = \dltw_{4}(\gaussF) - \EUV \dltw_{4}(\gaussF) \)
\begin{EQA}[rcl]
	&& \nquad
	\EUV \bigl| \dltwhat_{3}^{2}(\gaussF) - \Tens^{2}(\gaussF) \bigr|
	\leq 
	2 \EUV \bigl| \dltwhat_{4}(\gaussF) \, \Tens(\gaussF) \bigr|
	+ \EUV \dltwhat_{4}^{2}(\gaussF)
	\\
	& \leq &
	2 \sqrt{\E \Tens^{2}(\gaussF) \, \EUV \dltwhat_{4}^{2}(\gaussF)}
	+ \EUV \dltwhat_{4}^{2}(\gaussF)
	\leq 
	2 \accu_{\IFL} \, \dltw_{4,\IFL} + \dltw_{4,\IFL}^{2} \, .
\label{ycsy6ty25cgc53c7vyeetch2}
\end{EQA}
%
%
As \( \EUV \dltwhat_{3}(\gaussF) = 0 \), \eqref{7nhig763hf5bv9edfuejn} with \( g(\cdot) \equiv 1 \) 
and \eqref{ycsy6ty25cgc53c7vyeetch2} imply \eqref{u7yfjeoi8g8tj5udftrugh}.
The use of \eqref{uvjkb7nh3f6eyeyudt2jbi} from Lemma~\ref{LTensTGmG} with \( \tensco = \dltwu_{3}/6 \) yields
\begin{EQA}
	\EUV |\Tens(\gaussF)|
	& \leq &
	\sqrt{\EUV \Tens^{2}(\gaussF) }
	\leq 
	\frac{1}{6} \sqrt{15 \dltwu_{3}^{2} \, \| \TGD \|^{2} \, \tr^{2} (\TGD^{2})} 
\label{iujhdsufvhgfreu8478}
\end{EQA}
and \eqref{gdt6djiu97ie4425g77} follows as well. 
\end{proof}


\Section{Deviation bounds for Bernoulli vector sums}
\label{SdevBern}
Let \( Y_{i} \) be independent \( \Bernoulli(\thetas_{i}) \), \( i = 1,\ldots,\nsize \).
We denote \( \Yv = (Y_{i}) \in \R^{\nsize} \).
Weighted sums of the \( Y_{i} \) naturally appear in various statistical tasks 
including classification, binary response models, logistic regression etc. 
Recent applications include e.g. stochastic block modeling; see e.g. \cite{GaZh2017}, \cite{Abbe2018} and references therein,
or ranking from pairwise comparison \cite{ChGa2022} among many others.
We show how the general bounds of Section~\ref{Sdevboundexp} can be used for vector sums of Bernoulli r.v.s.
For a linear mapping \( \Psiv \colon \R^{\nsize} \to \R^{\dimp} \), define
\( \xiv = \Psiv (\Yv - \E \Yv) \).
Below we state some deviation bounds on the squared norm \( \| \xiv \|^{2} \) starting from the univariate case. 

\Subsection{Weighted sums of Bernoulli r.v.'s: univariate case}
Given a collections of weights \( (\weight_{i}) \), define 
\begin{EQA}
	S
	&=&
	\sum_{i=1}^{\nsize} Y_{i} \weight_{i} \, ,
	\\
	\VP^{2}
	&=&
	\Var(S)
	=
	\sum_{i=1}^{\nsize} \thetas_{i}(1-\thetas_{i}) \weight_{i}^{2} \, ,
	\\
	\weights 
	&=& 
	\max_{i} |\weight_{i}| .
\label{Het2Fetetp}
\end{EQA}
First, we state a deviation bound for a centered sum \( S - \E S \).

\begin{proposition}
\label{LBeBvM}
Let \( Y_{i} \) be independent \( \Bernoulli(\thetas_{i}) \)
and \( \weight_{i} \in \R \), \( i = 1,\ldots,\nsize \). 
Then \( S = \sum_{i=1}^{\nsize} Y_{i} \weight_{i} \) satisfies
\begin{EQA}
	\log \E \exp\Bigl\{ \frac{\lambda (S - \E S)}{\VP} \Bigr\}
	& \leq &
	\lambda^{2}	,
	\qquad
	\lambda \leq \frac{\log(2) \VP}{\weights} \, .
\label{llEelY1tsel}
\end{EQA}
Furthermore, suppose that given \( \xx \geq 0 \), 
\begin{EQA}
	\VP 
	& \geq &
	\frac{3}{2} \, \weights \sqrt{\xx} \, .
\label{Het2Fetetpzz}
\end{EQA}
Then
\begin{EQA}
	\P\bigl( \VP^{-1} |S - \E S| \geq 2 \sqrt{\xx} \bigr)
	& \leq &
	2 \ex^{-\xx} .
\label{PHetm1SESs2x}
\end{EQA}
Without \eqref{Het2Fetetpzz}, the bound \eqref{PHetm1SESs2x} applies
with \( \VP \) replaced by \( \VP_{\xx} = \VP \vee (3 \, \weights \sqrt{\xx} \, /2) \).
\end{proposition}

\begin{proof}
Without loss of generality assume \( \weights = 1 \), otherwise just rescale all the weights by the factor \( 1/\weights \).
We use that 
\begin{EQA}
	f(u)
	\eqdef
	\log \E \exp\Bigl\{ u (S - \E S) \Bigr\}
	&=&
	\sum_{i=1}^{\nbin} \Bigl[ \log\bigl( \thetas_{i} \ex^{u \weight_{i}} + 1 - \thetas_{i} \bigr) 
	- u \weight_{i} \thetas_{i} \Bigr] \, .
\label{llEelY1tsela}
\end{EQA}
This is an analytic function of \( u \) for \( |u| \leq \log 2 \) satisfying 
\( f(0) = 0 \), \( f'(0) = 0 \), and, with \( \upss_{i} = \log \thetas_{i} - \log(1-\thetas_{i}) \),
\begin{EQA}
	f''(u)
	&=&
	\sum_{i=1}^{\nbin} 
	\frac{\weight_{i}^{2} \, \thetas_{i}(1-\thetas_{i}) \, \ex^{u \weight_{i}}}
		 {( \thetas_{i} \ex^{u \weight_{i}} + 1 - \thetas_{i})^{2}}
	=
	\sum_{i=1}^{\nbin} \frac{\weight_{i}^{2} \, \ex^{\upss_{i}+u \weight_{i}}}{(\ex^{\upss_{i}+u \weight_{i}} + 1)^{2}} 
	=
	\sum_{i=1}^{\nbin} \theta_{i}(u) \bigl\{ 1-\theta_{i}(u) \bigr\} \weight_{i}^{2}
\label{nu2He2fpput}
\end{EQA}
for \( \theta_{i}(u) = \ex^{\upss_{i}+u \weight_{i}}/(\ex^{\upss_{i}+u \weight_{i}} + 1) \).
Clearly \( \theta_{i}(u) \) and thus, \( \theta_{i}(u) \bigl\{ 1-\theta_{i}(u) \bigr\} \) monotonously increases with \( u \) and 
it holds for \( \thetas_{i} = \theta_{i}(0) \) 
\begin{EQA}
	\theta_{i}(u) \bigl\{ 1-\theta_{i}(u) \bigr\}
	& \leq &
	\ex^{|u|} \, \thetas_{i} \, (1 - \thetas_{i})
	\leq 
	2 \, \thetas_{i} \, (1 - \thetas_{i}) ,
	\qquad
	|u| \leq \log 2.
\label{2ts1tseu1tb}
\end{EQA}
This yields 
\begin{EQA}
	f(u)
	& \leq & 
	\VP^{2} \, u^{2}
	\qquad
	|u| \leq \log 2.
\label{flaHetnu2la22}
\end{EQA}
As \( \xx \leq 4\VP^{2}/9 \), the value \( \lambda = \sqrt{\xx} \) fulfills 
\( \lambda/\VP = \sqrt{\xx}/\VP \leq  \log 2 \leq 2^{-1/2} \).
Now by the exponential Chebyshev inequality 
\begin{EQA}
	\P\Bigl( \VP^{-1}(S - \E S) \geq 2 \sqrt{\xx} \Bigr)
	& \leq &
	\exp\bigl\{ - 2 \lambda \sqrt{\xx} + f(\lambda/\VP) \bigr\}
	\\
	& \leq &
	\exp\bigl( - 2 \lambda \sqrt{\xx} + \lambda^{2} \bigr)
	=
	\ex^{-\xx} .
\label{exmxxsq2xnu2la22}
\end{EQA}
Similarly one can bound \( \E S - S \).
\end{proof}

\Subsection{Deviation bounds for Bernoulli vector sums}

Now we present an upper bound on the norm of a vector \( \xiv = \Psiv (\Yv - \E \Yv) \),
where \( \Psiv \) is a linear mapping \( \Psiv \colon \R^{\nsize} \to \R^{\dimp} \).
It holds 
\begin{EQA}
	\Var(\xiv)
	&=&
	\Var(\Psiv \Yv)
	=
	\Psiv \Var(\Yv) \Psiv^{\T} .
\label{HP2VxiWWYT}
\end{EQA}
We aim at bounding the squared norm \( \| \QP \xiv \|^{2} \) for another linear mapping \( \QP \colon \R^{\dimp} \to \R^{\dimq} \).

\begin{theorem}
\label{LHDxibound}
Let \( Y_{i} \sim \Bernoulli(\thetas_{i}) \), \( i = 1,\ldots, n \).
Consider \( \xiv = \Psiv (\Yv - \E \Yv) \), and let \( \HVB^{2} \geq 2 \Var(\xiv) \).
Define 
\begin{EQA}
	\weights
	&=&
	\max_{i \leq \nsize} \| \HVB^{-1} \Psiv_{i} \| \, ,\qquad
	\gmb
	=
	\log(2) / \weights \, .
\label{WHm1uinfmiwi}
\end{EQA} 
Then with \( \BBH = \QP \HVB^{2} \QP^{\T} \) and \( \zqc(\BBH,\xx) \) from \eqref{dyv6ejf8gjwkerih83}, it holds
\begin{EQA}
    \P\bigl( \| \QP \xiv \| \ge \zqc(\BBH,\xx) \bigr)
    & \le &
    3 \ex^{-\xx} .
\label{PxivbzzBBroB3m}
\end{EQA}    
\end{theorem}

\begin{proof}
We apply the general result of Corollary~\ref{CTQPxivlarge} under conditions \eqref{expgamgm}.
For any vector \( \uv \), consider the scalar product 
\( \langle \HVB^{-1} \xiv,\uv \rangle = \langle \HVB^{-1} \Psiv (\Yv - \E \Yv), \uv \rangle \).
It is obviously a weighted centered sum of the Bernoulli r.v.'s \( Y_{i} - \thetas_{i} \) with
\begin{EQA}
	\Var \langle \HVB^{-1} \xiv,\uv \rangle
	& \leq &
	\| \uv \|^{2}/2 .
\label{VlHm1xiuu22}
\end{EQA}
One can write with \( \eps_{i} = Y_{i} - \thetas_{i} \) and \( \epsv = (\eps_{i}) \)
\begin{EQA}
	\langle \HVB^{-1} \xiv,\uv \rangle
	&=&
	\bigl\langle \epsv, \Psiv^{\T} \HVB^{-1} \uv \bigr\rangle .
\label{lHm1xiuepWH1u}
\end{EQA}
By the Cauchy-Schwarz inequality, it holds 
\begin{EQA}
	\| \Psiv^{\T} \HVB^{-1} \uv \|_{\infty}
	=
	\max_{i} \bigl| (\HVB^{-1} \Psiv_{i})^{\T} \uv \bigr|
	& \leq &
	\weights \| \uv \| .
\label{WHm1uinfmiwi}
\end{EQA} 
Bound \eqref{llEelY1tsel} of Proposition~\ref{LBeBvM} on the exponential moments of \( \langle \HVB^{-1} \xiv,\uv \rangle \) implies 
\begin{EQA}
	\log \E \exp\bigl\{ \langle \HVB^{-1} \xiv,\uv \rangle \bigr\}
	& \leq &
	\| \uv \|^{2}/2	,
	\qquad
	\| \uv \| \leq \log(2) / \weights \, .
\label{llEelY1tsel2}
\end{EQA}
Therefore, \eqref{expgamgm} is fulfilled with 
\( \gmb = \log(2) / \weights \). 
The deviation bound \eqref{uyfuyerd7e7uhh8yy689t} of Corollary~\ref{CTQPxivlarge} yields the assertion.
\end{proof}

\bibliography{exp_ts,listpubm-with-url}

\end{document}